\documentclass[11pt]{article}
\usepackage[margin=1in]{geometry} 
\usepackage{amsmath}        
\usepackage{graphicx}
\usepackage{amssymb}
\usepackage{amsthm}
\usepackage{epstopdf}
\usepackage{amsfonts}
\usepackage[scr]{rsfso}
\usepackage{hyperref}
\usepackage{booktabs}
\usepackage{placeins}
\makeatletter
\def\amsbb{\use@mathgroup \M@U \symAMSb}
\makeatother
\usepackage{bbold}
\usepackage{natbib}
\usepackage{epigraph}
\usepackage{authblk}
\usepackage{esint}
\usepackage{subfig}

\def\d{\delta}
\def\c{\circ}

\def\E{\amsbb{E}}

\def\C{\amsbb{C}ov}
\def\P{\amsbb{P}}
\def\R{\amsbb{R}}
\def\Q{\amsbb{Q}}
\def\N{\amsbb{N}}

\def\H{\amsbb{H}}
\def\G{\amsbb{G}}

\def\D{\mathrm{d}}
\def\Var{\amsbb{V}ar}
\def\C{\amsbb{C}ov}
\def\Cor{\amsbb{C}orr}

\DeclareMathOperator*{\for}{\quad\text{for}\quad}
\DeclareMathOperator*{\Leb}{Leb}

\newcommand\ind[1]{\amsbb{I}_{#1}}

\usepackage{fancyhdr}
\title{A field equation for induction-transduction of activation-deactivation probability on measurable space}

\author[1]{Caleb Deen Bastian\thanks{cbastian@princeton.edu}}
\author[1,2]{Herschel Rabitz\thanks{hrabitz@princeton.edu}}
\affil[1]{Program in Applied and Computational Mathematics, Princeton University, Princeton, NJ. 08544}
\affil[2]{Department of Chemistry, Princeton University, Princeton, NJ. 08544}

\date{\today}  

\begin{document}

\vspace{-6cm}                    
\maketitle

\begin{abstract} Induction-transduction of activating-deactivating points are fundamental mechanisms of action that underlie innumerable systems and phenomena, mathematical, natural, and anthropogenic, and can exhibit complex behaviors such as self-excitation, phase transitions, hysteresis, polarization, periodicity, chaos, wave behavior, geometry, and energy transfer. We describe a class of primitives for induction-transduction based on dynamics on images of marked random counting measures under graphical random transformations. We derive a field equation for the law of the activation-deactivation (Bernoulli) process on an arbitrary measurable space and describe some mechanisms of action on the unit interval.

\end{abstract}

\section{Introduction} In this note we derive, from the theory of random measures and graphs (for review of general and particular concepts, please see \S6 of \cite{cinlar} and \S2-3 of \cite{bastianrg}), a field equation for the transduction of probability law of a spatiotemporal Bernoulli process of expansion and contraction, induced by a graphical random transformation of a random counting measure. The idea is that the system initially starts in a state of zero entropy, with all states deactivation (or all activation). This is at time $0_-$ and is pre-induction. The induction activation event (a `perturbation' or `kick') occurs at time zero ($0$) and transitions the system from a state of zero entropy to a state of non-zero entropy. The resultant activation flow in positive time is the transduced system response, which may increase activation until complete activation, decay until complete deactivation, or anything in between. The response is encoded by the system configuration and the nature of the kick. Note the appearance of the Bernoulli process and its significance in view of the Ornstein isomorphism theorem \citep{ornstein}, and the appearance of the random counting measure, which presides over the process and regulates system cardinality, its projection into space-time, and point correlation. The image of the random measure under graphical random transformation conveys a degree function, which we refer to as a random density field, over a measurable space (here ``space''), where each location of the density field is identified to a random finite count (`density'). If space is a continuum, then the density field is uncountably infinite. The density field is said to be entangled with-by-through the points. The field equation for activation law admits dual interpretations, where it may be solved in continuous or discrete time, and continuous or discrete space, and any combination of the two. The field equation can be provisioned as a wave equation. As a general matter, the graph encodes point density through its density kernel, which encodes the dynamics of the induction-transduction process of activation-deactivation, whose law is given by the field equation. Moreover, the induced stochastic geometry of space-time, which we refer to as stochastic space time, is identified by the maximum entropy isoline, which we refer to as the activation frontier or activation event horizon. 

The field equation fundamentally requires computing spatiotemporal count probabilities of the random counting measure, which are defined in terms of the Laplace functional. Typically these probabilities requires repeated differentiation (unstable) or complex integration (expensive) to attain, but a recent result establishes that the Poisson-type (PT) family of random counting measures are the only members in the canonical non-negative power series (CNNPS) family that are subset invariant (self-similar) \citep{Bastian:2020tb}, e.g., counts in all subspaces are rescaled PT random variables, established through discovery of a generalized additive Cauchy (functional) equation. Such strongly invariant counting probability measures are said to be `bones.' The PT random measures thus exhaust the possibilities in a broad class and are fundamental to practically solving the field equation. Interestingly, their functional equation (and its dual) has recently been shown to lead to the unification of the Cauchy and Golab–Schinzel-Type equations in mathematical analysis \citep{jacek}. 

The field equation admits wave and energy equations. For a certain monotone continuous density kernel on the unit square, which we call central induction-transduction, the corresponding mean density field, which turns out to be related to the Laplace transform of the uniform distribution on the unit interval, is discontinuous at the source, and for this case, the field equation transports the discontinuity in space and time, bending, folding, and making discontinuous the stochastic space-time, and generates energetic waves and negative energy. Note that source discontinuity is unnecessary for central, where central mechanisms having continuous sources are easy to find.

Other density kernels include subcentral, decentral, and local, each a distinct mechanism for induction-transduction and induced stochastic space time. Subcentral is similar to central, in that it bends the stochastic space-time given measure (`mass'), but it does not seem capable of folding nor making it discontinuous, and its stochastic space-time is locally linear. Note that subcentral cannot retrieve central by taking limits. Decentral is homogeneous in space, i.e. there are no spatial dynamics in decentral. Decentral, in discrete-time with activation annihilation and superfluity, transduces chaos. Decentral also exhibits profound hysteresis of the equilibria in the activation-deactivation transduction thresholds, enabling efficient cycling. Local, having a discontinuous density kernel function, also bends, folds, and makes discontinuous the stochastic space-time, but in a different manner than central. 

The field equation can be mapped to physical processes. It is encoded by the counting law for system cardinality (for Poisson, the mean number of points), the spatial domain, the density kernel on the product spatial domain, and the induction-transduction activation-deactivation thresholds. Additional frills for the model include forcings and overloadings on the probability fluxes and making the system open. Suitable wave and energy equations may be provisioned in view of the medium, thinking of space as a kind of elastic membrane, with appropriately dimensionalized constants. Interestingly, the continuous wave formulation of the field equation admits a reduction to the master equation underlying the Dirac equation.

This note is organized as follows. In \S\ref{sec:induction} and \S\ref{sec:transduction} we discuss induction and transduction. Using this material, we introduce the induction-transduction activation-deactivation (ITAD) field equation (\S\ref{sec:itad}). We end with a brief discussion (\S\ref{sec:discuss}). We include an appendix containing exhibits of some solutions of the ITAD field equation for the four considered mechanisms (central, subcentral, decentral, local) and some of their properties, where all continuum solutions are attained using the function `NDSolve' in Mathematica 11.3 on Mac OS X x86. The solutions are stable with increasing precision, until global loss of precision, which occurs around $10^{-15}$ in these examples; the examples are generated relative to precision $10^{-13}$.

\section{Induction}\label{sec:induction} In this section we consider the problem of modeling the distribution of the number of activations in a system of correlated points. The fundamental idea is to identify a point in activation to an active vertex of a random graph. In particular we construct a model of point interactions of the system as a (density) random graph, where edges correspond to density interactions. The degree (density) of a point is the number of density interactions it has with itself and the other points and corresponds to integration. If the degree (integrated interactions) of a point exceeds some threshold, then the point is in activation and the point vertex is said to be active in the graph. Both the mean activation interaction (edge) weight function, here a density kernel, and the system point counting distribution, encode the point activation covariance structure. 

A random counting measure $N$ on $(E,\mathscr{E})$ is identified to a pair of independent deterministic probability measures $(\kappa,\nu)$---$\kappa$ on $(\N_{\ge0},2^{\N_{\ge0}})$ and $\nu$ on $(E,\mathscr{E})$---through the \emph{stone throwing construction} \citep{cinlar,Bastian:2020tb} (STC) \begin{equation}\label{eq:stc}N(A)=N\ind{A}=\int_{E}N(\D x)\ind{A}(x)\equiv\sum_{i}^{K}\ind{A}(X_i)\quad\text{for}\quad A\in\mathscr{E}\end{equation}  where counting variable $K\sim\kappa$ has mean $c>0$ and variance $\delta^2\ge0$ and $\ind{A}$ is a set function. It is denoted the random measure $N=(\kappa,\nu)$ on $(E,\mathscr{E})$ and is said to be formed by the collection of iid points (independency) $\mathbf{X}$. It is also known as the mixed binomial process. Let $\mathscr{E}_{\ge0}$ denote the collection of non-negative $\mathscr{E}$-measurable functions. For functions $f$ in $\mathscr{E}_{\ge0}$, we have \[Nf = \int_EN(\D x)f(x)=\sum_i^Kf(X_i)\] Now we consider the basic STC random graph set-up \citep{bastianrg}: let $M=N\times N$ be the product STC random measure on $(E\times E,\mathscr{E}\otimes\mathscr{E})$ and let $\phi$ be a symmetric Bernoulli transformation from $E\times E$ into $F=\{0,1\}$ defined as $\phi(x,y)\sim\text{Bernoulli}(f(x,y))$ for symmetric $f$ in $(\mathscr{E}\otimes\mathscr{E})_{[0,1]}$. 
 
\subsection{Point activation degree} Let $\phi_x(y)\equiv\phi(x,y)$ and $I$ be the identity function. The point activation degree function (random density field) is given by the image \[d(x) = N\phi(x,\cdot)= (N\circ\phi_x^{-1})I=\sum_i^K\phi(x,X_i)\for x\in E\] Put $f_x(\cdot)\equiv f(x,\cdot)$. The degree function has probability generating function (pgf), mean, variance, covariance, and correlation \begin{align*}\psi_x(t)&\equiv \psi(1-\nu(f_x)+\nu(f_x)t)\\\E d(x) &= c\nu(f_x)\\\Var d(x) &=c\nu(f_x) + (\delta^2-c)(\nu(f_x))^2\\\C(d(x),d(y))&=\begin{cases}c\nu(f_xf_y) + (\delta^2-c)\nu(f_x)\nu(f_y)&\text{if } x\ne y\\c\nu(f_x) + (\delta^2-c)(\nu(f_x))^2&\text{otherwise}\end{cases}\\\Cor(d(x),d(y))&=\begin{cases}\frac{c\nu(f_xf_y) + (\delta^2-c)\nu(f_x)\nu(f_y)}{\sqrt{(c\nu(f_x) + (\delta^2-c)(\nu(f_x))^2)(c\nu(f_y) + (\delta^2-c)(\nu(f_y))^2)}}&\text{if } x\ne y\\1&\text{otherwise}\end{cases}\end{align*}

For $a>0$, both of these choices for $f$ result in decreasing degree functions $d$. In the first case, the normalized degree function  is one at zero, $\nu(f_0)=1$, so these points have maximum mean activation degree. For the second, the normalized degree function is less than one at zero, $\nu(f_0)=\varphi_U(a)<1$, so these points have mean degree less than the maximum.    

An interesting special case is constant density kernel $f(x,y)=p$ and fixed number $c$ of points, $\kappa=\text{Dirac}(c)$. Then $\Cor(d(x),d(y))=\ind{}(x=y)$. Hence for the constant density kernel $f=p$ and a fixed number of points $c$, the activation degrees (densities) are decorrelated $\text{Binomial}(c,p)$ random variables. For $\kappa=\text{Poisson}(c)$, then $\Cor(d(x),d(y))=\ind{}(x=y)+p\ind{}(x\ne y)$. Hence the Poisson random measure ($c=\delta^2$) correlates the density field. More generally, the Poisson-type (PT) random measures correlate the degree function, with binomial ($\delta^2<c$) and negative binomial ($\delta^2>c$) respectively having reduced and increased correlation relative to Poisson, and binomial recovering Dirac as $\delta^2\rightarrow0$. We illustrate the degree function correlations of the constant density kernel for the PT random measures below in Figure~\ref{fig:ptcorr}. The constant density kernel demonstrates that the correlation structure of the point activation density is determined by both the density kernel $f$ and the point counting law $\kappa$. 

\begin{figure}[h!]
\centering
\includegraphics[width=5in]{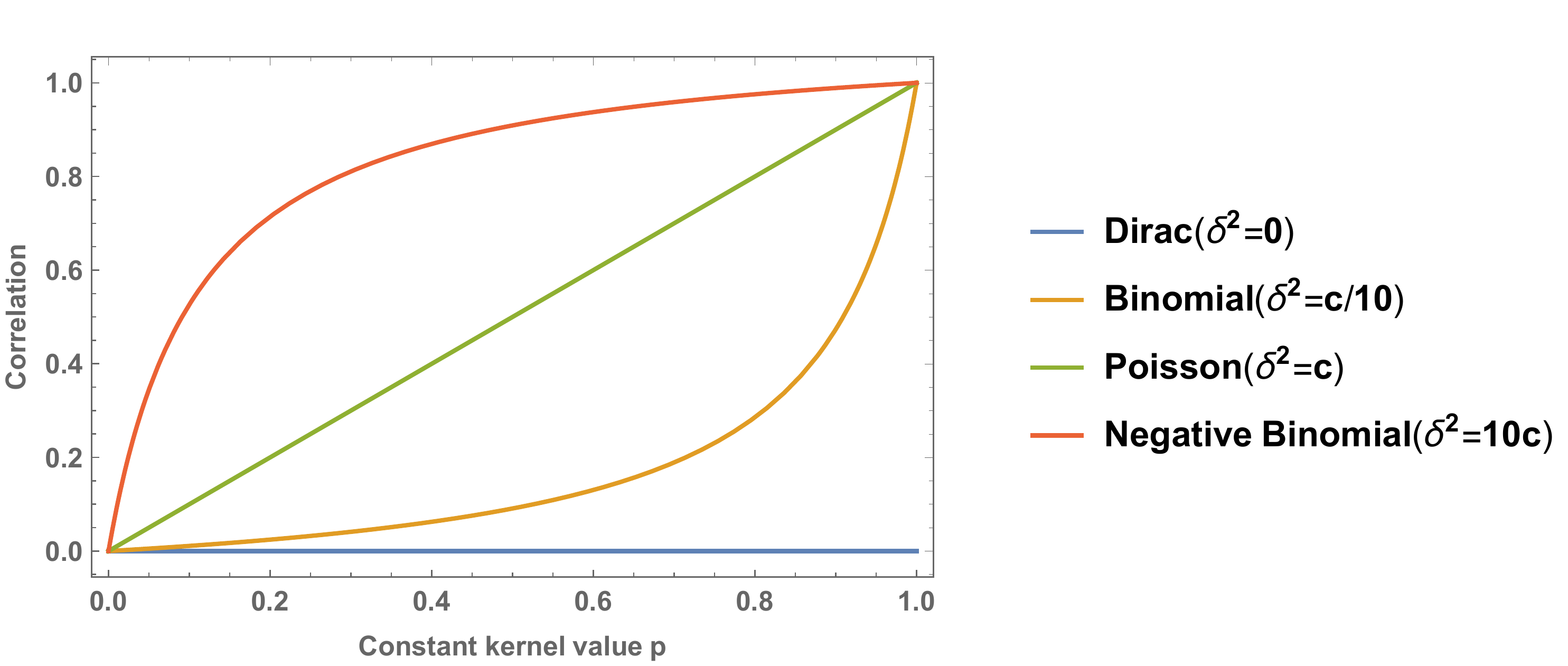}
\caption{Constant density kernel $f(x,y)=p$ degree function (density field) correlations $\Cor(d(x),d(y))$ in $p$ of Poisson-type random measures for $x\ne y$}\label{fig:ptcorr}
\end{figure}

\FloatBarrier

\subsection{Number of points in activation} We say a point $X$ is in activation when $d(X)\ge k$ for some $k\ge1$, and its vertex is said to be $k$-active. The number of $k$-active vertices \[V_k(G) \equiv(N\circ d^{-1})\ind{\ge k}= \sum_i^K\ind{\ge k}(d(X_i))\] corresponds to the number of points in activation and has mean\begin{align*}\E V_k(G) &= c\int_E\nu(\D x)\P(d(x)\ge k) + c\int_E\nu(\D x)\P(d(x)=k-1)f(x,x) \\&= c\nu g_k\end{align*} where \[g_k(x)\equiv\P(d(x)\ge k)+ \P(d(x)=k-1)f(x,x)\for x\in E\] 

\subsection{Point system example} 
\begin{figure}[h!]
\centering
\begingroup
\captionsetup[subfigure]{width=5in,font=normalsize}
\subfloat[$\C(d(x),d(y))$ for $f_2$\label{fig:g1}]{\includegraphics[width=2.5in]{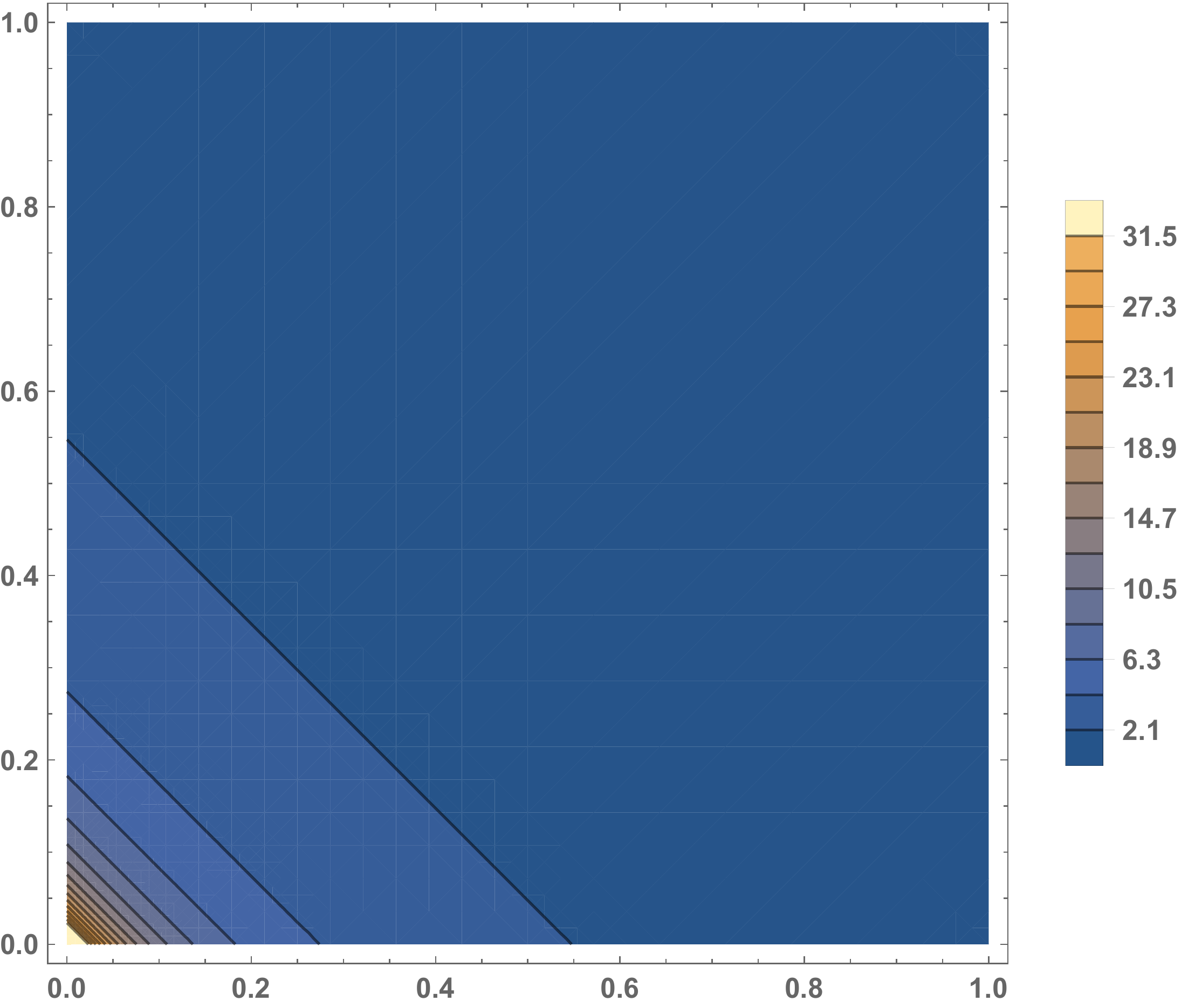}}
\subfloat[$\Cor(d(x),d(y))$ for $f_2$\label{fig:g2}]{\includegraphics[width=2.5in]{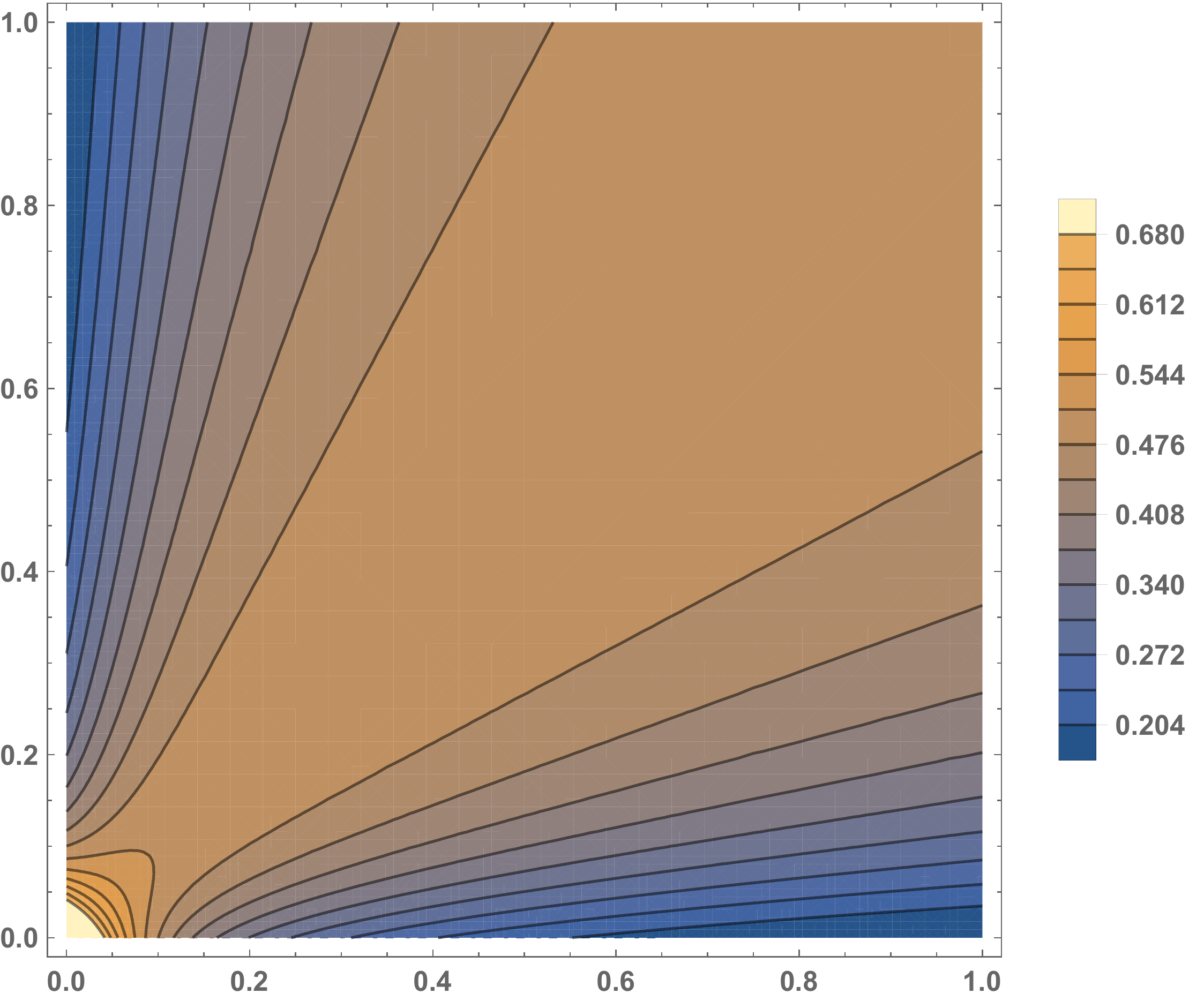}}\\
\subfloat[$\C(d(x),d(y))$ for $f_3$\label{fig:g3}]{\includegraphics[width=2.5in]{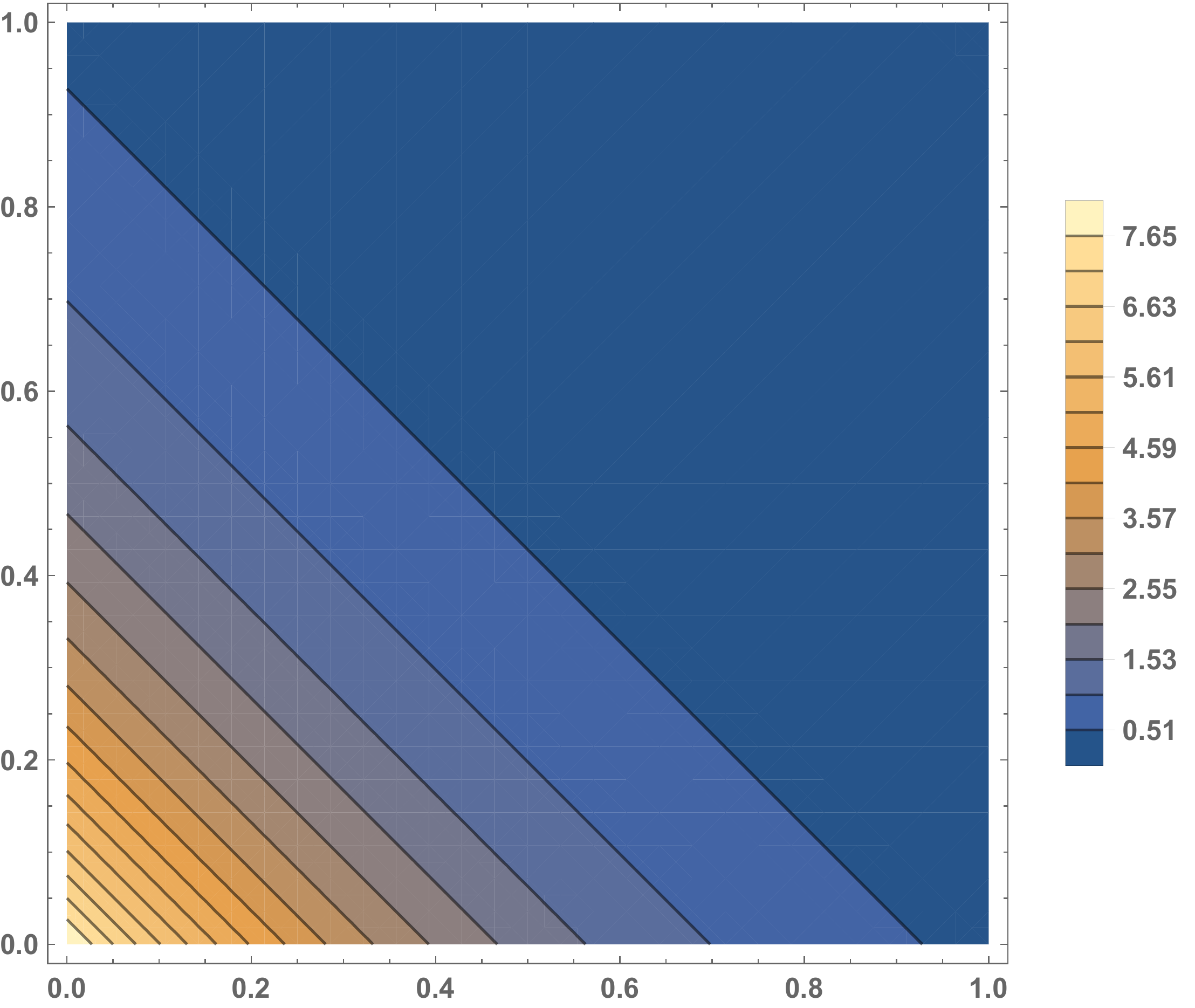}}
\subfloat[$\Cor(d(x),d(y))$ for $f_3$\label{fig:g3}]{\includegraphics[width=2.5in]{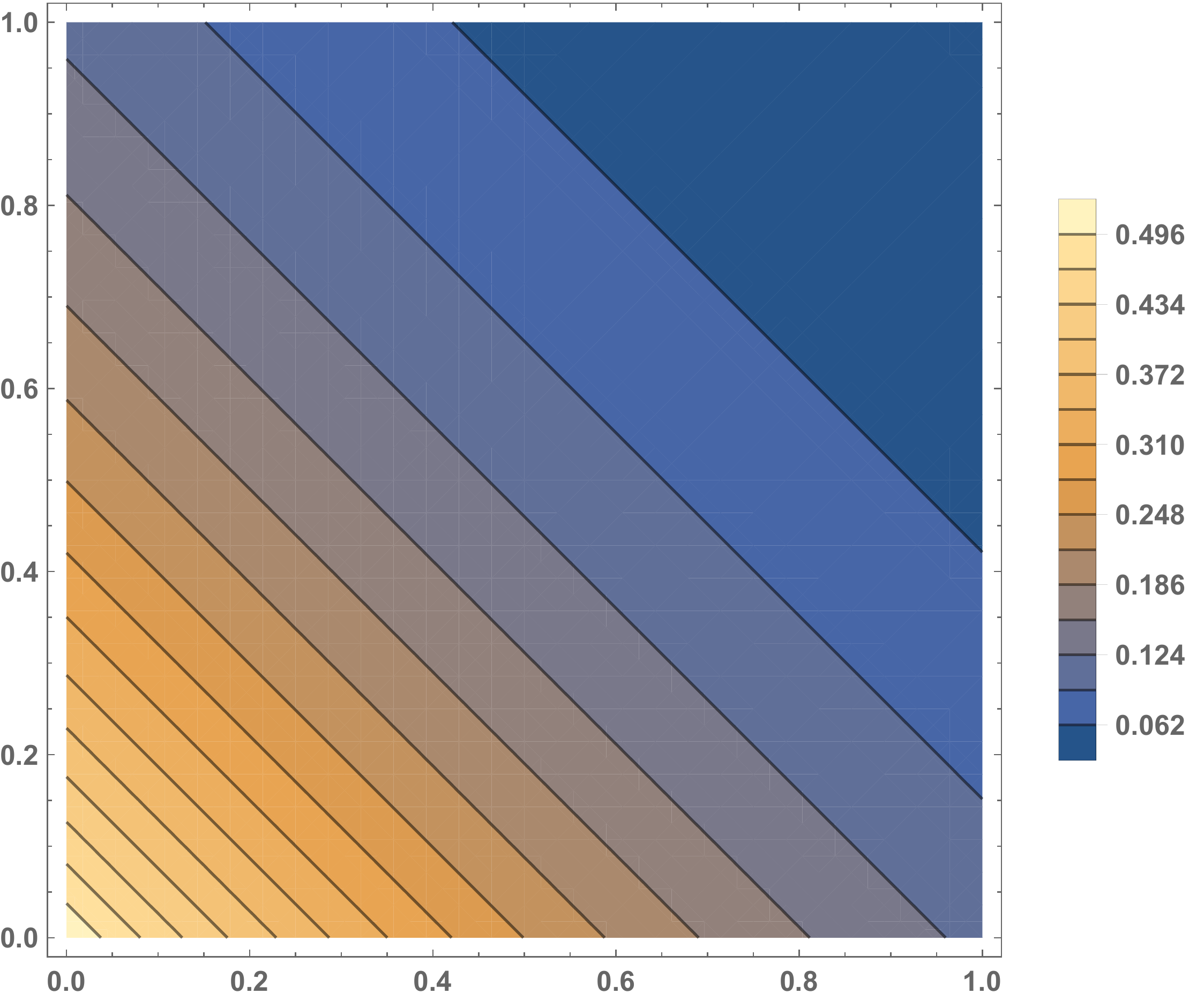}}
\endgroup
\caption{Covariance $\C(d(x),d(y))$ and correlation $\Cor(d(x),d(y))$ of point activation degree $d$ for density kernels $f_2(x,y)=\exp_-a_2xy$ and $f_3(x,y)=\exp_-a_3(x+y)$ having equal total mass $a_1=1/10$ and $\kappa=\text{Poisson}(c=50)$ (shown for $x\ne y$)}\label{fig:covcor}
\end{figure}

\begin{figure}[h!]
\centering
\begingroup
\captionsetup[subfigure]{width=5in,font=normalsize}
\subfloat[$G_1$\label{fig:g1}]{\includegraphics[width=2in]{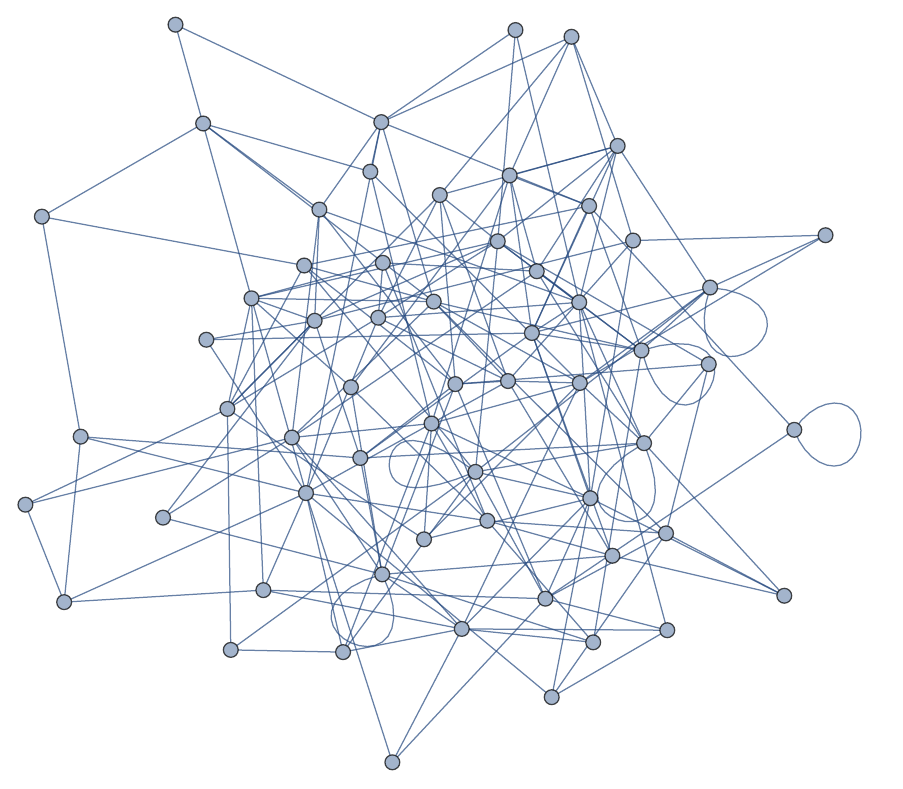}}
\subfloat[$G_2$\label{fig:g2}]{\includegraphics[width=2in]{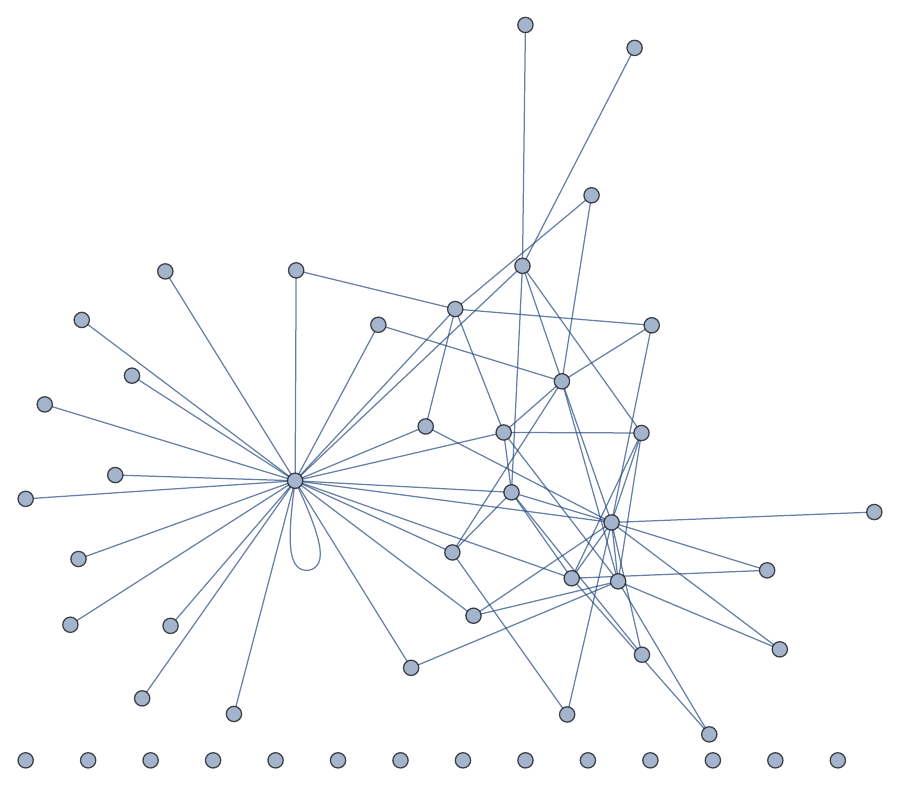}}
\subfloat[$G_3$\label{fig:g3}]{\includegraphics[width=2in]{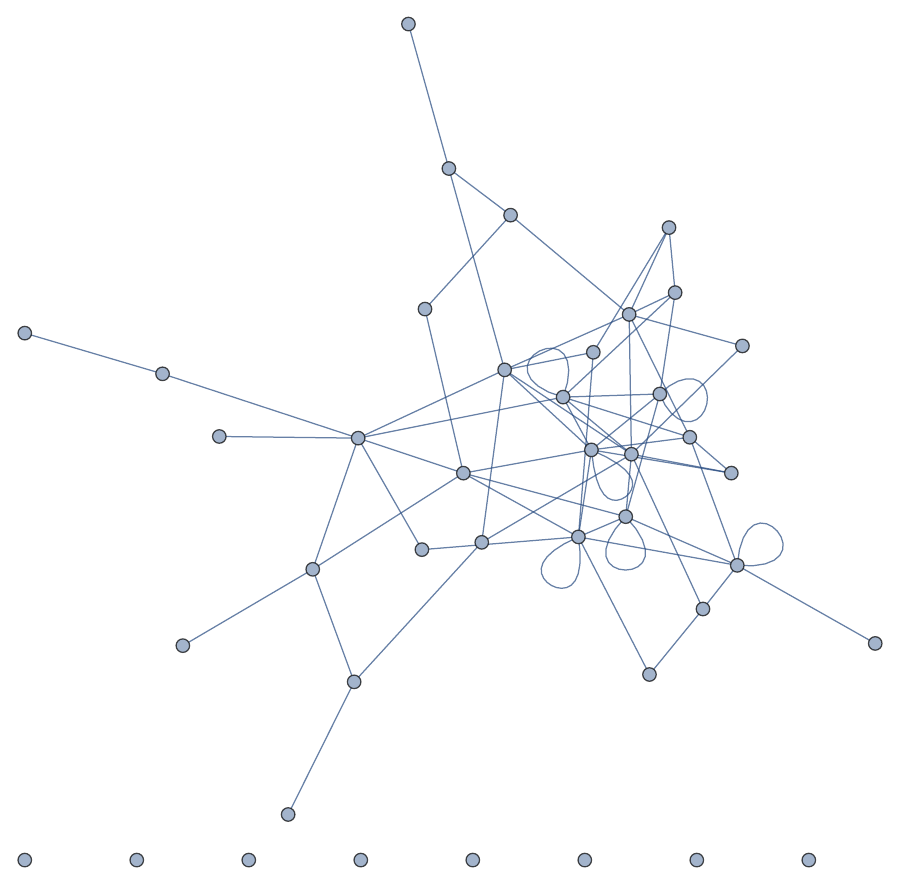}}
\endgroup
\caption{Random realizations of point activation Bernoulli graph $G$ for density kernels $f_1,f_2,f_3$ having equal total mass $a_1=1/10$ and $\kappa=\text{Poisson}(c=50)$}\label{fig:gs}
\end{figure}

Consider a point system of $K\sim\kappa=\text{Poisson}(c)$ points with $c=50$. Consider the Bernoulli random graph $G(\phi)$ with density kernel $f$. We choose $f$ as $f_1(x,y)=a_1$ for $a_1\in(0,1]$, $f_2(x,y)=\exp_-a_2xy$, and $f_3(x,y)=\exp_-a_3(x+y)$, where $a_2$ and $a_3$ are chosen such that $a_1=(\nu\times\nu)f_2=(\nu\times\nu)f_3$, i.e. $f_1, f_2, f_3$ have the same normalized edge density. We take $a_1=1/10$, giving $a_2\simeq 43.4997$ and $a_3\simeq 3.0057$. We simulate $V_k(G)$ 10,000 times for $k=8$ for each kernel, plot the empirical distributions below in Figure~\ref{fig:loss}, and report their sample statistics in Table~\ref{tab:stats}. Despite the fixed mean total number of activation interactions, the resultant random variables $V_k(G)$ have markedly different distributions. Firstly, the constant kernel has the smallest mean and by far the largest probability of zero activations. The kernels $f_2$ and $f_3$ give small probability to zero activations, with the mean and variance of $f_2$ dominated by those of $f_3$, and $f_2$ has the smallest variance of the three considered kernels. We plot degree covariance and correlation for $f_2$ and $f_3$ in Figure~\ref{fig:covcor}, conveying distinct structures. We show random graphs $G_1,G_2,G_3$ corresponding to the three kernels in Figure~\ref{fig:gs}. These illustrate the qualitative differences between the homogeneous structure of $f_1$ and the non-homogenous structures of $f_2$ and $f_3$. For instance, $G_2$ contains one highly connected point and a number of disconnected points, whereas $G_1$ is a connected homogeneous Poisson graph. 

\begin{figure}[h!]
\centering
\includegraphics[width=5in]{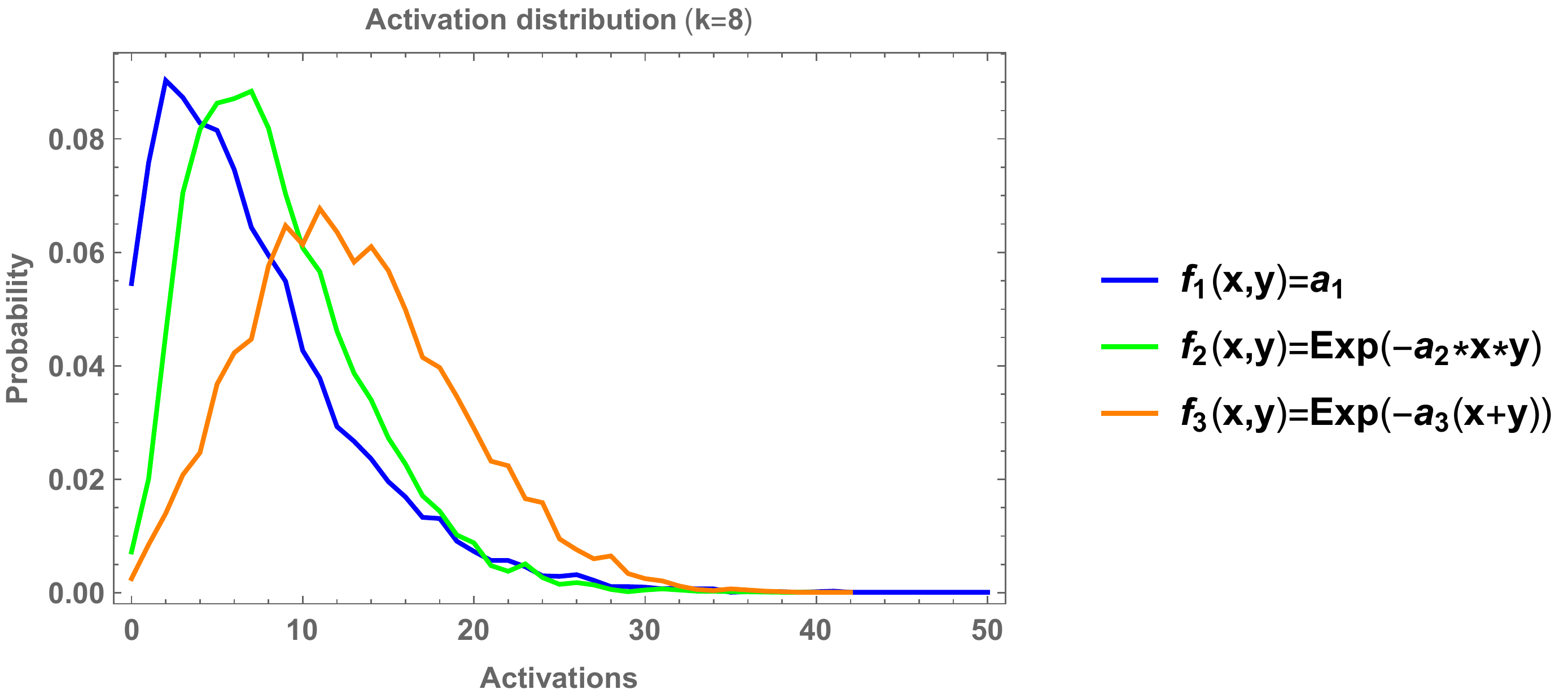}
\caption{Empirical distributions of point activation count $V_8(G)$ for density kernels $f_1,f_2,f_3$ having equal total mass $a_1=1/10$ }\label{fig:loss}
\end{figure}

\FloatBarrier

\begin{table}[h!]
\begin{center}
\begin{tabular}{cll}
\toprule
Kernel & $\E_n V_8(G)$ & $\Var_n V_8(G)$ \\\midrule
$f_1$ & 7.2026 & 35.4671\\
$f_2$ & 8.4552 & 25.2061 \\
$f_3$ & 12.9328 & 38.2605\\
\bottomrule
\end{tabular}
\caption{Sample mean $\E_n$ and variance $\Var_n$ of point activation count $V_8(G)$ for density kernels $f_1,f_2,f_3$ having equal total mass $a_1=1/10$ and $n=10^4$}\label{tab:stats}
\end{center}
\end{table}

\FloatBarrier

\FloatBarrier

\FloatBarrier

\subsection{Activation processes} Consider $N=(\kappa,\nu)$ on $(E,\mathscr{E})=([0,1],\mathscr{B}_{[0,1]})$ with $\nu=\Leb$ formed by independency $\mathbf{X}$ and its product random measure $M=N\times N$ with symmetric Bernoulli transformation $\phi$ and degree function $d(x)=N\phi(x,\cdot)$. Suppose we have mark-space time-set $[0,T]$. The idea is that each point is identified to location in $E=[0,1]$ with distribution $\nu=\Leb$ and a time-of-arrival (TOA) of the location in $[0,T]$ with distribution $\pi=\text{Uniform}[0,T]$. This is a simple marked random measure design where TOA and locations are independent. Let independency $(\mathbf{X},\mathbf{T})$ form the random measure $N^*=(\kappa,\nu\times\pi)$ on $(E\times[0,T],\mathscr{E}\otimes\mathscr{B}_{[0,T]})$. Define the time-dependent point activation degree function \[d_t(x,s) \equiv \sum_i^K\phi(x,X_i)\ind{[0,t]^2}(s,T_i)\for (x,s)\in E\times[0,T],\quad t\in[0,T]\]Then we can define the point activation process as \[V_k^t(G) \equiv \sum_i^K\ind{\ge k}(d_t(X_i,T_i))\for t\in[0,T]\] The mean carries over in a manner similar to those of the system activation distribution based on marked random measures. In this setting, $(d_t(x,s))_{t\in[0,T]}$ is a non-decreasing, pure-jump, counting process for every $(x,s)\in E\times[0,T]$, with initial condition $d_0(x,s)=0$ and terminal condition $d_T(x,s)=d(x)$. Similarly, $(V_k^t(G))_{t\in[0,T]}$ is a non-decreasing, pure-jump, counting process. They are both non-decreasing because they both may be the zero process. An immediate generalization is taking the TOA distribution $\pi$ to be non-uniform. Note that these processes are neither L\'{e}vy nor Markov, even for Poisson $N^*$. A further generalization is letting the activation value and TOA be dependent. 

\subsection{Activation motifs / higher-order activation interactions} We have defined points in activation in terms of their degree functions, formed from edge counts. Instead, we can define points in activation in terms of higher-order edge relations, such as relative to the number of triangles the points participate in or other motifs. 

\subsection{Generalizations} A number of generalizations are possible. The transformation need not be symmetric, yielding directed graphs. Then the points have in- and out-degrees. Also the transformation can be a choice different from Bernoulli, such as Poisson (giving activation multigraphs). Others are randomizing the threshold $k$ in the threshold function $\ind{\ge k}$ or considering a general non-negative (possibly random) function. The mark space law can be described through a transition kernel, encoding conditional dependence on location. The location space can be multi-dimensional, incorporating point attributes, etc. 

\section{Transduction}\label{sec:transduction} 

Let $F=\{D,A\}$ be a mark-space of point activation labels of deactivation $D$ and activation $A$, with distribution $\xi$. Consider the random measure $N=(\kappa,\eta\equiv\nu\times\xi)$ on $(E\times F,\mathscr{E}\otimes\mathscr{F})$ formed by $(\mathbf{X},\mathbf{Y})$ with $\nu=\Leb$. As before, $\mathbf{X}$ are the activation locations. We take $\xi\{D\}=1-\nu(g_k)$ such that the fraction of activated points is $\xi\{A\}=\nu(g_k)$, where $\nu(g_k)$ is the mean normalized activation degree on $(E,\mathscr{E})$. That is, we initialize the fraction of activations to the activation model's normalized number of points in activation, a probability. 

Let $B\equiv\{(D,A),(A,D)\}\subset F\times F$. Consider the symmetric Bernoulli transformation $\Phi_B=\phi\ind{B}$ with symmetric mean function $g$ in $(\mathscr{E}\otimes\mathscr{F}\otimes\mathscr{E}\otimes\mathscr{F})_{[0,1]}$ defined as $g((x,y),(x',y'))=f(x,x')\ind{B}(y,y')$, such that \[\Phi_B((x,y),(x',y'))\sim\text{Bernoulli}(g((x,y),(x',y'))=\begin{cases}\text{Bernoulli}(f(x,x'))&\text{if }(y,y')\in B\\\text{Bernoulli}(0)\equiv\text{Dirac}(0)&\text{otherwise}\end{cases}\] where symmetric $f$ in $(\mathscr{E}\otimes\mathscr{E})_{[0,1]}$ is the activation interaction kernel of the symmetric Bernoulli transformation $\phi$. We have that $\Phi_B((x,y),(x',y'))=\phi(x,x')\ind{B}(y,y')$. 

\subsection{Activation degree} Let $g_{xy}(x',y')\equiv g((x,y),(x',y'))$. Define the activation degree function as \[d_B(x,y)=N\Phi_B((x,y),\cdot)=\sum_i^K\phi(x,X_i)\ind{B}(y,Y_i)\for(x,y)\in E\times F\] The degree has mean \begin{align*}\E d_B(x,y) &= c\eta(g_{xy}) \\&= c\int_E\nu(\D x')f(x,x')\int_F\xi(\D y')\ind{B}(y,y')\\&=c\,\nu(f_x)(\ind{}(y=D)\xi\{A\}+\ind{}(y=A)\xi\{D\})\end{align*} and pgf \[\psi_{xy}^B(t) = \psi(1-\eta(g_{xy})+\eta (g_{xy})t)\] 

\subsection{Number of activations} We say that the deactivation point $(X_i,Y_i=D)$ is in activation if $d_B(X_i,Y_i)\ge l$ for $l\ge 1$. The number of activations is \[V_l^B(G) = \sum_i^K\ind{\ge l}(d_B(X_i,Y_i))\ind{}(Y_i=D)\] with mean \[\E V_l^B(G) = c\nu(q_l^B)\xi\{D\}\] where \[q_l^B(x)\equiv\P(d_B(x,D)\ge l) + \P(d_B(x,D)=l-1)f(x,x)\] and $d_B(x,D)$ has pgf \[\psi_{xD}^B(t)=\psi(1-\nu(f_x)\xi\{A\} + \nu(f_x)\xi\{A\}t)\]

\section{Induction-transduction activation-deactivation field equation}\label{sec:itad}

Consider probability measure $\kappa$ on $(\N_{\ge0},2^{\N_{\ge0}})$ with pgf $\psi$ and probability measure $\nu$ on $(E,\mathscr{E})$. Let $f\in(\mathscr{E}\otimes\mathscr{E})_{[0,1]}$ be an $[0,1]$-valued $\mathscr{E}\otimes\mathscr{E}$-measurable function and put $f_x(y)=f(x,y)$. Here we assume here that $f$ is symmetric. 

\subsection{Spatial-temporal field equation} The spatial-temporal induction-transduction activation-deactivation (ITAD) field equation on $(F,\mathscr{F})=([0,1],\mathscr{B}_{[0,1]})^{\R_{\ge0}\times E}$ is given by the following spatio-temporal ordinary differential equation \begin{align*}(t,x,u)&\in\R_{\ge0}\times E\times[0,1]\\\\d(x)&\text{ has pgf } \psi_x(z)=\psi(1-\nu(f_x)+\nu(f_x)z)\\p_k(x)&=\P(d(x)\ge k)+\P(d(x)=k-1)f(x,x)\\\\d_D(x,u)&\text{ has pgf }\psi_{xu}^D(z)=\psi(1-\nu(f_x)u+\nu(f_x)uz)\\ q_l(x,u)&=\P(d_D(x,u)\ge l)+\P(d_D(x,u)=l-1)f(x,x)\\\\d_A(x,u)&\text{ has pgf }\psi_{xu}^A(t)=\psi(1-\nu(f_x)(1-u)+\nu(f_x)(1-u)t)\\r_m(x,u)&=\P(d_A(x,u)\ge m)+\P(d_A(x,u)=m-1)f(x,x)\\\\ \frac{\partial}{\partial t}\P(t,x)&=q_l(x,\P(t,x))(1-\P(t,x)) - r_m(x,\P(t,x))\P(t,x)\\\P(0,x)&=p_k(x)\end{align*} The induction-transduction activation field equation is retrieved by taking $m=\infty$ such that \begin{align*}\frac{\partial}{\partial t}\P(t,x)&=q_l(x,\P(t,x))(1-\P(t,x))\\\P(0,x)&=p_k(x)\end{align*} 

 For the case of Poisson and unit thresholds $k=l=m=1$, the ITAD field equation is given by \begin{align*}\frac{\partial}{\partial t}\P(t,x)&=\left(1- e^{-c\nu(f_x)\P(t,x)}+e^{-c \nu(f_x)\P(t,x)}f(x,x)\right)(1-\P(t,x))\\&\hspace{1cm}-\left(1- e^{-c\nu(f_x)(1-\P(t,x))}+e^{-c \nu(f_x)(1-\P(t,x))}f(x,x)\right)\P(t,x)\\\P(0,x)&= 1-e^{-c\nu(f_x)}+e^{-c\nu(f_x)}f(x,x)\end{align*} 

The field equation can be solved in discrete-time as \begin{align*}t&\in\N_{\ge0}\\\P(t+1,x)&=\P(t,x)+q_l(x,\P(t,x))(1-\P(t,x)) - r_m(x,\P(t,x))\P(t,x)\\\P(0,x)&=p_k(x)\end{align*}

For example the density function has distribution 
\[\P(d(x)=j) = \frac{1}{j!}\psi_x^{(j)}(0) = \frac{1}{2\pi i}\oint_{C}\frac{\psi_x(t)}{t^{j+1}}\D t\for j\in\N_{\ge0}\] where the contour integral is over the unit circle $C$. These calculations can be unstable or expensive. Towards efficiency, notice that the density pgfs $\psi_x$, $\psi_{xu}^D$, and $\psi_{xu}^A$ are each binomially thinned. If the counting law is closed under binomial thinning, then the density distributions are rescaled versions of themselves, i.e. they are self-similar, and there is no need for differentiation nor complex integration. Alas, it turns out on diffuse space that there are only three self-similar members of the canonical non-negative power series family of distributions to possess closure under binomial thinning, the Poisson-type (PT): Poisson, negative binomial, and binomial. The Dirac measure can be said to be a fourth, as the limit of the binomial distribution. This means that systems having fixed or PT cardinality possess self-similar configurations. Note that the PT random measures are also closed under thinning for atomic measures, but not uniquely so. For completeness the PT family of distributions are given by \[K\sim\begin{cases}\text{Binomial}(n,p)\in\{0,\dotsc,n\}&\text{for}\quad n\in\N_{>0}, p\in(0,1]\\\text{Poisson}(c)\in\N_{\ge0}&\text{for}\quad c\in(0,\infty)\\\text{NegativeBinomial}(r,p)\in\N_{\ge0}&\text{for}\quad r\in\N_{>0}, p\in(0,1]\end{cases}\] and the rescalings to subspace $A\subseteq E$ with mass $a$ are given by \[K_A\sim\begin{cases}\text{Binomial}(n,ap) &\text{if }K\sim\text{Binomial}(n,p)\\\text{Poisson}(ac)&\text{if }K\sim\text{Poisson}(c)\\\text{NegativeBinomial}(r,ap/(1-(1-a)p))&\text{if }K\sim\text{NegativeBinomial}(r,p)\end{cases}\] The Poisson-type distributions can be approximated by Gaussian when the degrees are large. Putting $a_x=\nu(f_x)$, we have (with some abuse of notation) \[d(x)\simeq\begin{cases}\text{Gaussian}(na_xp,na_xp(1-a_xp))&\text{if }K\sim\text{Binomial}(n,p)\\\text{Gaussian}(a_xc,a_xc)&\text{if }K\sim\text{Poisson}(c)\\\text{Gaussian}(\frac{ra_xp}{1-p},\frac{ra_xp(1-(1-a_x)p)}{(1-p)^2})&\text{if }K\sim\text{NegativeBinomial}(r,p).\end{cases}\] 

\subsection{Temporal field equation} 
The (marginal) temporal field equation is attained as the following ordinary differential equation \begin{align*}(t,u)&\in\R_{\ge0}\times[0,1]\\\\p_k&=\int_E\nu(\D x)\left(\P(d(x)\ge k)+\P(d(x)=k-1)f(x,x)\right)\\ q_l(u)&=\int_E\nu(\D x)\left(\P(d_D(x,u)\ge l)+\P(d_D(x,u)=l-1)f(x,x)\right)\\r_m(u)&=\int_E\nu(\D x)\left(\P(d_A(x,u)\ge m)+\P(d_A(x,u)=m-1)f(x,x)\right)\\\\ \frac{\partial}{\partial t}\P(t)&=q_l(\P(t))(1-\P(t)) - r_m(\P(t))\P(t)\\\P(0)&=p_k\end{align*} The equation can be solved in discrete-time as \begin{align*}t&\in\N_{\ge0}\\\P(t+1)&=\P(t)+q_l(\P(t))(1-\P(t)) - r_m(\P(t))\P(t)\\\P(0)&=p_k\end{align*}

\subsection{Induced wave equation} Assuming $\nu$ is diffuse, the associated (first-order) wave equation on $(W,\mathscr{W})=(\R,\mathscr{B}_{\R})^{\R_{\ge0}\times E}$ is given by \begin{align*}\frac{\partial}{\partial t}\Q(t,x)&=C_1\frac{\partial}{\partial x}\P(t,x)\\\Q(0,x)&=0\end{align*} where $C_1>0$ is an appropriately dimensionalized constant, e.g., transport velocity.

\subsection{Energy equation} Consider $E=[0,1]$. The energy equation on $(G,\mathscr{G})=(\R,\mathscr{B}_{\R})^{\R_{\ge0}\times E}$ is given by \begin{align*}\frac{\partial}{\partial t}\G(t,x)&=\frac{1}{2}C_2\left(\frac{\partial}{\partial t}\Q(t,x)\right)^2+\frac{1}{2}C_3\left(\frac{\partial}{\partial x}\Q(t,x)\right)^2\\\G(0,x)&=0\end{align*} where $C_2,C_3>0$ are appropriately dimensionalized constants, e.g, $C_2$ and $C_3$ are the mass density and tension of a unit length vibrating string.

\subsection{Field quantities} The entropy field is \[\H(t,x)=\P(t,x)\log(1/\P(t,x))\] The maximum entropy activation frontier (or activation event horizon) is given by the level set \[\widehat{\H}=\{(t,x)\in\R_{\ge0}\times E: \P(t,x)=1/2\}\] 

\subsection{Forcing and overloading} The field equation can be extended to include forcing and overloading through non-negative forcing parameters $(a,b)$ and overloading parameters $(\alpha,\beta)$, i.e., \[\frac{\partial}{\partial t}\P(t,x)=(\alpha q_l(x,\P(t,x)+a)(1-\P(t,x)) - (\beta r_m(x,\P(t,x))+b)\P(t,x)\] The forcing-overloading pairs for activation $(a,\alpha)$ and deactivation $(b,\beta)$ are chosen such that $\sup_{(x,u)\in E\times [0,1]}(\alpha q_l(x,u)+a)\le 1$ and $\sup_{(x,u)\in E\times [0,1]}(\beta r_m(x,u)+b)\le 1$ respectively hold. 

\subsection{ITAD wave equation} For continuous marginal density functions (`sources'), i.e. $\nu(f_x)$ is continuous on diffuse $(E,\mathscr{E})$, the ITAD field equation can be provisioned as a continuous wave equation. The ITAD wave equation, with forcings $a$ and $b$, is given by \begin{align*} \frac{\partial}{\partial t}\P(t,x) +C_1\frac{\partial}{\partial x}\P(t,x)&=(q_l(x,\P(t,x))+a)(1-\P(t,x)) - (r_m(x,\P(t,x))+b)\P(t,x)\\\P(0,x)&=p_k(x)\end{align*} where $C_1\ge0$ is activation transport velocity. If $C_1=0$, then the standard ITAD field equation is retrieved. The additional term processes the field fluctuations relative to finite activation transport velocity. Note that, for given activation-deactivation transduction thresholds $l$ and $m$, the forcings $a$ and $b$ must be chosen such that $\sup_{(x,u)\in E\times [0,1]}(q_l(x,u)+a)\le 1$ and $\sup_{(x,u)\in E\times [0,1]}(r_m(x,u)+b)\le 1$. Solution stability exhibits a bifurcation in $C_1$, where increasing $C_1$ beyond a certain level causes the ITAD wave equation to blow-up in finite-time. The ITAD wave equation can also be solved in discrete time and/or space, employing  suitable finite differences.

\paragraph{Reduction to the Dirac equation} Note that for infinite activation-deactivation transduction thresholds ($l=m=\infty$), equal forcings ($a=b\in(0,1]$), and positive finite transport velocity ($C_1=v>0$), the continuous ITAD wave equation reduces to the master equation underlying the Dirac equation, c.f. the positive component of equation (2) in \cite{prl}, where $\P_+(t,x)\equiv\P(t,x)$ and $\P_-(t,x)\equiv1-\P(t,x)$, and by duality, its negative component. In particular the reduction conveys the telegrapher's equation \[\frac{\partial}{\partial t}\P_{\pm}(t,x) = -a(\P_{\pm}(t,x)-\P_{\mp}(t,x)) \mp v\frac{\partial}{\partial x}\P_{\pm}(t,x)\] which implies the formulation \[\frac{\partial^2}{\partial t^2}\P_{\pm}(t,x) - v^2\frac{\partial^2}{\partial x^2}\P_{\pm}(t,x) = -2a\frac{\partial}{\partial t}\P_{\pm}(t,x)\] shown to underlie the one-dimensional Dirac equation. This reduction suggests that the ITAD wave equation is a generalization of the telegrapher's equation and that there is some kind of relationship between ITAD theory and physics; however, the physical interpretation of ITAD theory at present eludes us. 

\subsection{Density kernels} We think of $f$ in $(\mathscr{E}\otimes\mathscr{E})_{[0,1]}$ as a density function or kernel. We discuss four archetypical ITAD mechanisms for $E=[0,1]$ and $\nu=\Leb$: central, subcentral, decentral, and local. We denote $\varphi_U(\alpha)=(1-e^{-\alpha})/\alpha$ the Laplace transform of the uniform distribution on $[0,1]$.

\paragraph{Central} The marginal densities of central mechanisms, i.e., $\nu(f_x)$, evaluate to unity at the boundary and decay rapidly. Let $g(x)=\exp_-x$. The canonical central density kernel is $f(x,y) = \exp_-axy=g(axy)$ for $a\ge0$. Recall $f_x(y)\equiv f(x,y)$. Then \begin{align*}\nu(f_xf_y)&=\int_E\nu(\D z)g(axz)g(azy)\\&=\varphi_U(a(x+y))\\\nu(f_x)&=\int_E\nu(\D z)g(axz)\\&=\varphi_U(ax)\\\lim_{x\downarrow0}\nu(f_x)&=1=g(0)\end{align*} $\nu(f_x)$ is discontinuous at $x=0$.  A family of non-canonical central density kernels is defined by $f(x,y)=1/(1+axy)^d$ for $a\ge0, d\ge1$. For $d=1$, we have $\nu(f_x)=\log(1+ax)/(ax)$ with $\lim_{x\downarrow0}\nu(f_x)=1$, which is discontinuous at $x=0$. For $d=2$, this is $f(x,y)=1/(1+axy)^2$ with $\nu(f_x)=1/(1+ax)$ and $\nu(f_0)=1$, which is continuous at $x=0$. For $d\ge3$, the behavior is the same as $d=2$.

\paragraph{Subcentral} The marginal subcentral density functions evaluate to less than unity at the boundary and decay rapidly. Consider subcentral density function $\widetilde{f}(x,y)=\exp_-a(x+y)=\widetilde{f}(x,0)\widetilde{f}(y,0)=f(x,1)f(y,1)=g(ax)g(ay)$ (for $a\ge0$). Then \begin{align*}\nu(\widetilde{f}_x\widetilde{f}_y)&=\int_E\nu(\D z)g(ax)g^2(az)g(ay)\\&=g(ax)g(ay)\int_E\nu(\D z)g^2(az)\\&=e^{-a(x+y)}\varphi_U(2a)\\&=\widetilde{f}(x,y)\varphi_U(2a)\\&=\widetilde{f}(x,0)\widetilde{f}(y,0)\varphi_U(2a)\\&=f(x,1)f(y,1)\varphi_U(2a)\\&=g(ax)g(ay)\varphi_U(2a)\\\nu(\widetilde{f}_x)&=g(ax)\int_E\nu(\D z)g(az)\\&=e^{-ax}\varphi_U(a)\\&=\widetilde{f}(x,0)\varphi_U(a)\\&=f(x,1)\varphi_U(a)\\&=g(ax)\varphi_U(a)\\\nu(\widetilde{f}_0)&=\varphi_U(a)\\\lim_{a\downarrow0}\nu(\widetilde{f}_0)&=1\end{align*} In this case, $\nu(\widetilde{f}_x)=e^{-ax}\varphi_U(a)\simeq (1-ax)\varphi_U(a)$ for small $ax$, recovering linear activation density. Another subcentral density kernel is $f(x,y)=(1+ax)^{-2}(1+ay)^{-2}$ for $a\ge0$. Then $\nu(f_x)=(1+a)^{-1}(1+ax)^{-2}$. Note that $\nu(f_0)=1/(1+a)$ and that $\nu(f_x)$ is continuous at $x=0$.

\paragraph{Decentral} The decentral mechanism is based on the constant density kernel $f(x,y)=p$ for $p\in[0,1]$. 

\paragraph{Local} The local mechanism is based on the local density kernel $f(x,y)=q\ind{}(|x-y|\le r)$ for $(q,r)\in[0,1]\times[0,1]$.

\section{Discussion}\label{sec:discuss} The induction-transduction activation-deactivation (ITAD) field equation represents the spatiotemporal probability of activation-deactivation. The underlying process is Bernoulli. By the Ornstein isomorphism theorem, the Bernoulli process is universal, unique, and the `most random' possible process. The ITAD field equation can be interpreted as the law of a special kind of random Boolean dynamical system where the graphical random transformation induces a density random field, whose exceedance probabilities constitute fluxes of probability. Moreover, this system can be identified to the dynamics of a neural network, where each point is regarded as a neuron that activates (`fires') when its activation exceeds some level, here a perceptron, and deactivates (`stops' firing) similarly. This model is based on graphical random transformations of 2-product random measures. Higher order graphical random transformations of higher-order product random measures may be used to define higher-order ITAD field equations, conveying activation-deactivation through for example triangulations of points, and so on and so forth. The measurable space is general. We have considered the (diffuse) unit interval and square. Atomic measure spaces can be considered, such as in lattices of integers or alphabets of symbols and their density kernels, i.e. induction-transduction through co-prime integers. Some future work include studying the existence, uniqueness, and stability of solutions to the ITAD equations; exploring the physical significance of the given mechanisms and others; developing estimators for the parameters given measurements; and so on. In regard to identification, ideas from stochastic epidemiology may be helpful, e.g., dynamic survival analysis \citep{dsa}. 

\section*{Acknowledgements} The authors thank Andr\'{e} Gontijo Campos at the Max-Planck Institute for Nuclear Physics (Heidelberg, Baden-W\"{u}rttemberg, Germany) for suggesting reduction to the Dirac equation. 

 \bibliographystyle{apalike}

\appendix

\newpage
\section{Induction-transduction}
\begin{figure}[h!]
\centering
\includegraphics[width=4.5in]{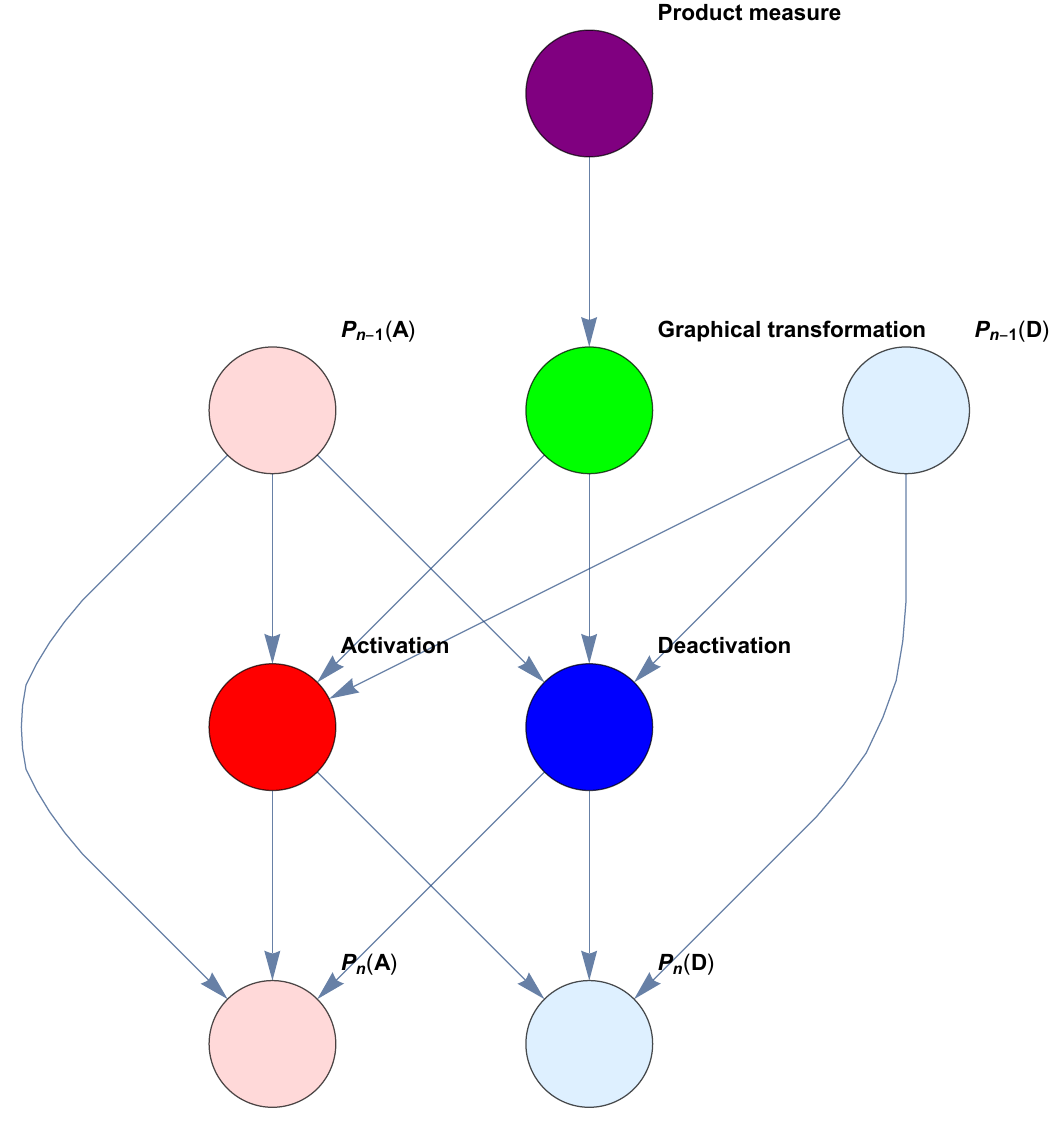}
\caption{Mean-field of primitive in discrete time for induction-transduction}\label{fig:cm}
\end{figure}
\FloatBarrier

\begin{figure}[h!]
\centering
\begingroup
\captionsetup[subfigure]{width=5in,font=normalsize}
\subfloat[Induction activations for constant kernel $p\in\{0.05,0.1,0.15\}$ in threshold $k$\label{fig:r0}]{\includegraphics[width=4.5in]{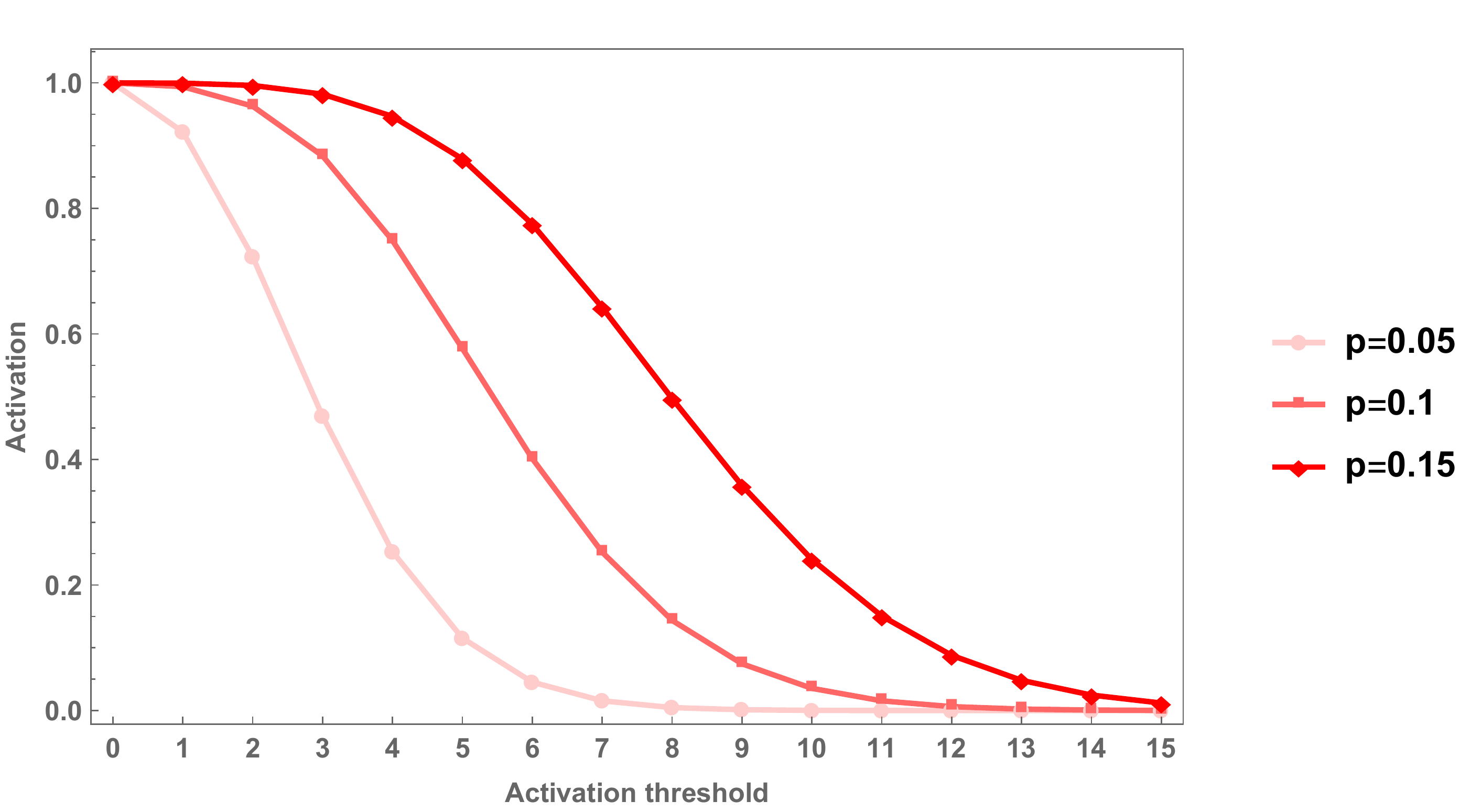}}\\
\subfloat[Induction activations for constant kernel $p$ and functional kernels $q\ind{}(|x-y|\le r)$, $\exp(-a(x+y))$, $\exp(-bxy)$ in threshold $k$\label{fig:r0}]{\includegraphics[width=4.5in]{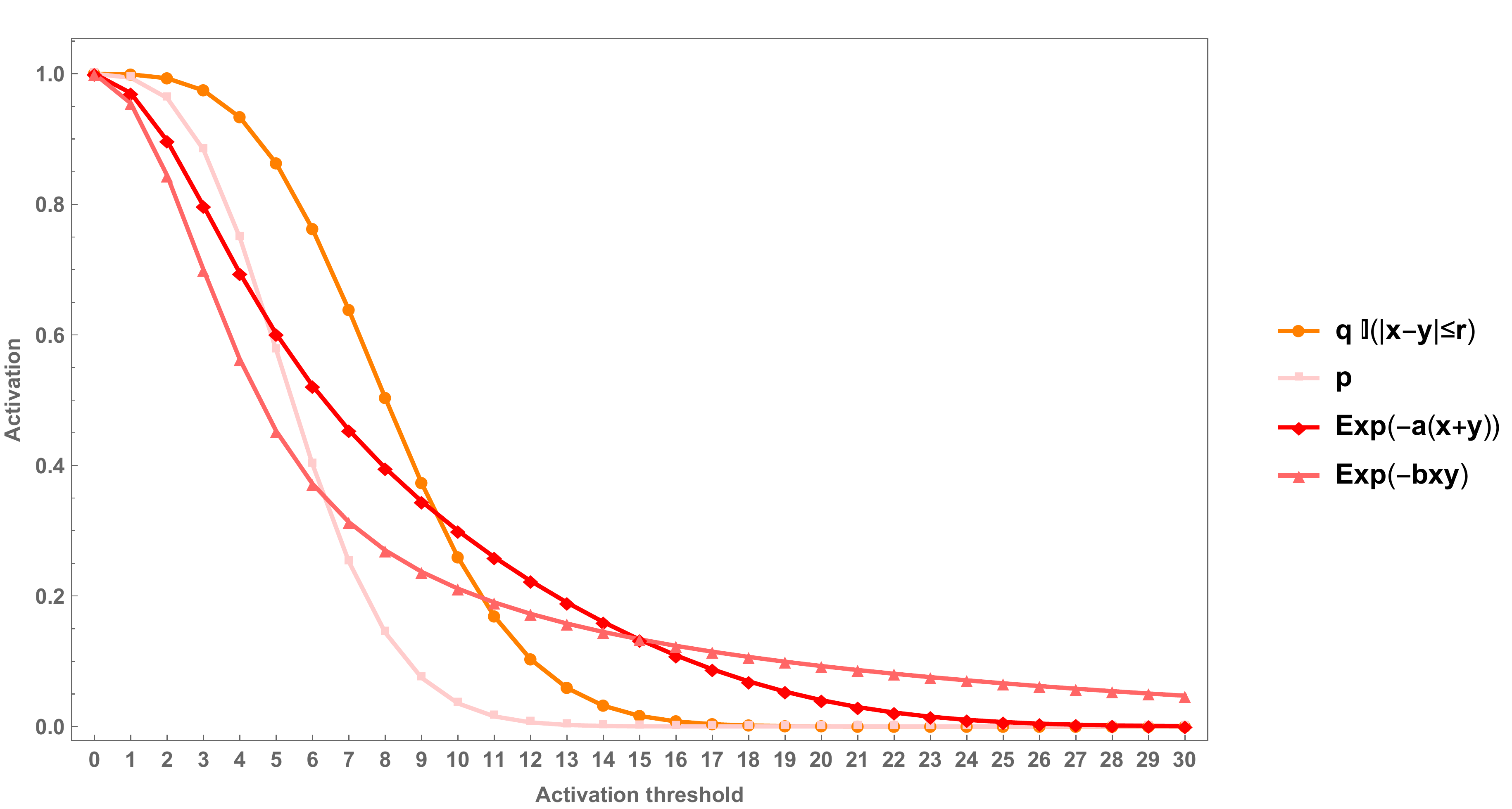}}\\
\subfloat[Induction activations in latent induction for constant kernel $p$ and functional kernels $\exp(-a(x+y))$ and $\exp(-bxy)$ and induction threshold $k=8$\label{fig:r0}]{\includegraphics[width=4.5in]{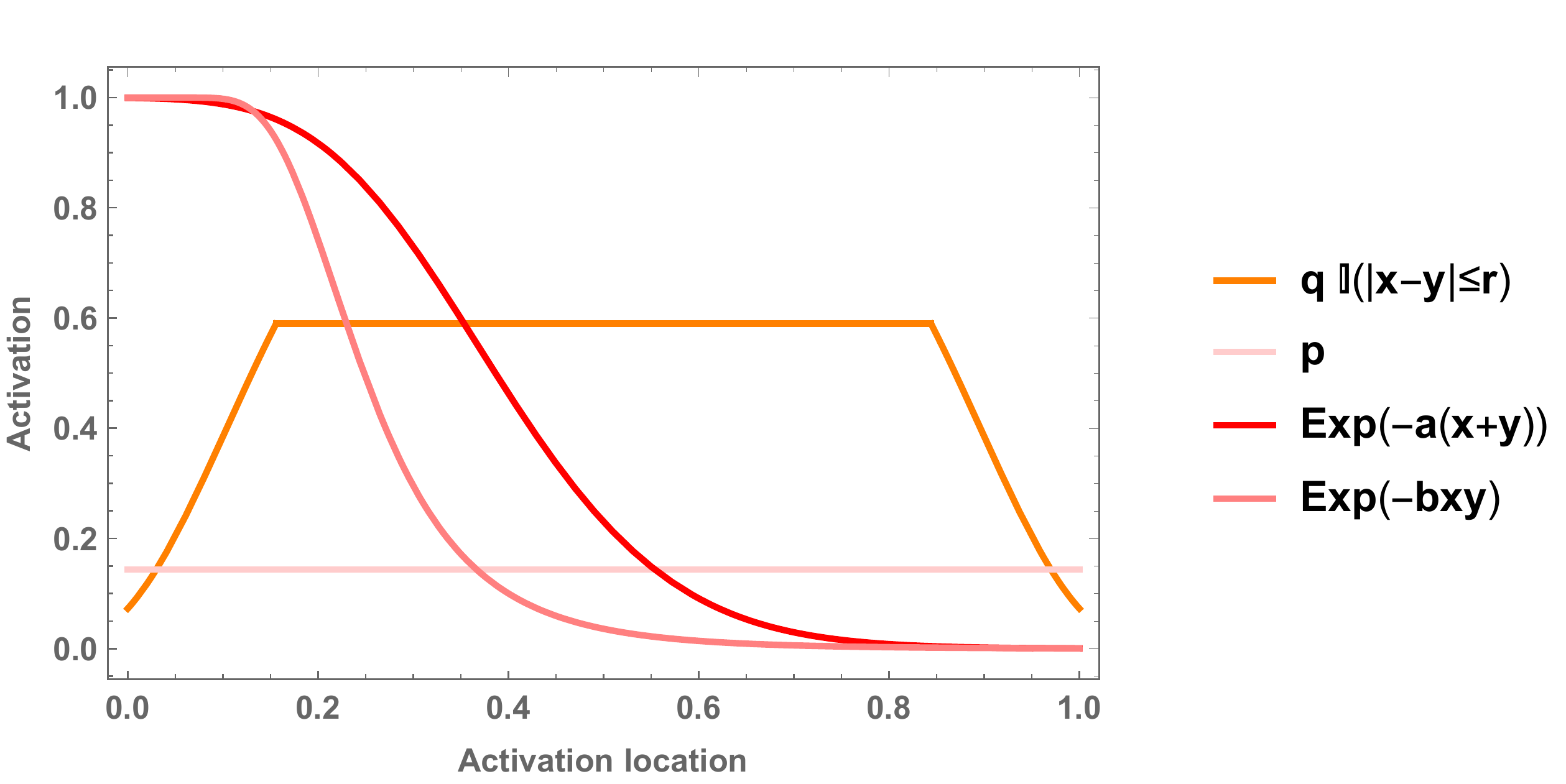}}\\
\endgroup
\caption{Induction activations for $K\sim\text{Poisson}(50)$ and varying density kernels having common kernel mass $p=1/10$ (local: $q=0.5$ and $r$ set) and varying transduction activation thresholds $l\in\{1,2,3\}$}\label{fig:orbitsdeact}
\end{figure}

\FloatBarrier

\begin{figure}[h!]
\centering
\includegraphics[width=6in]{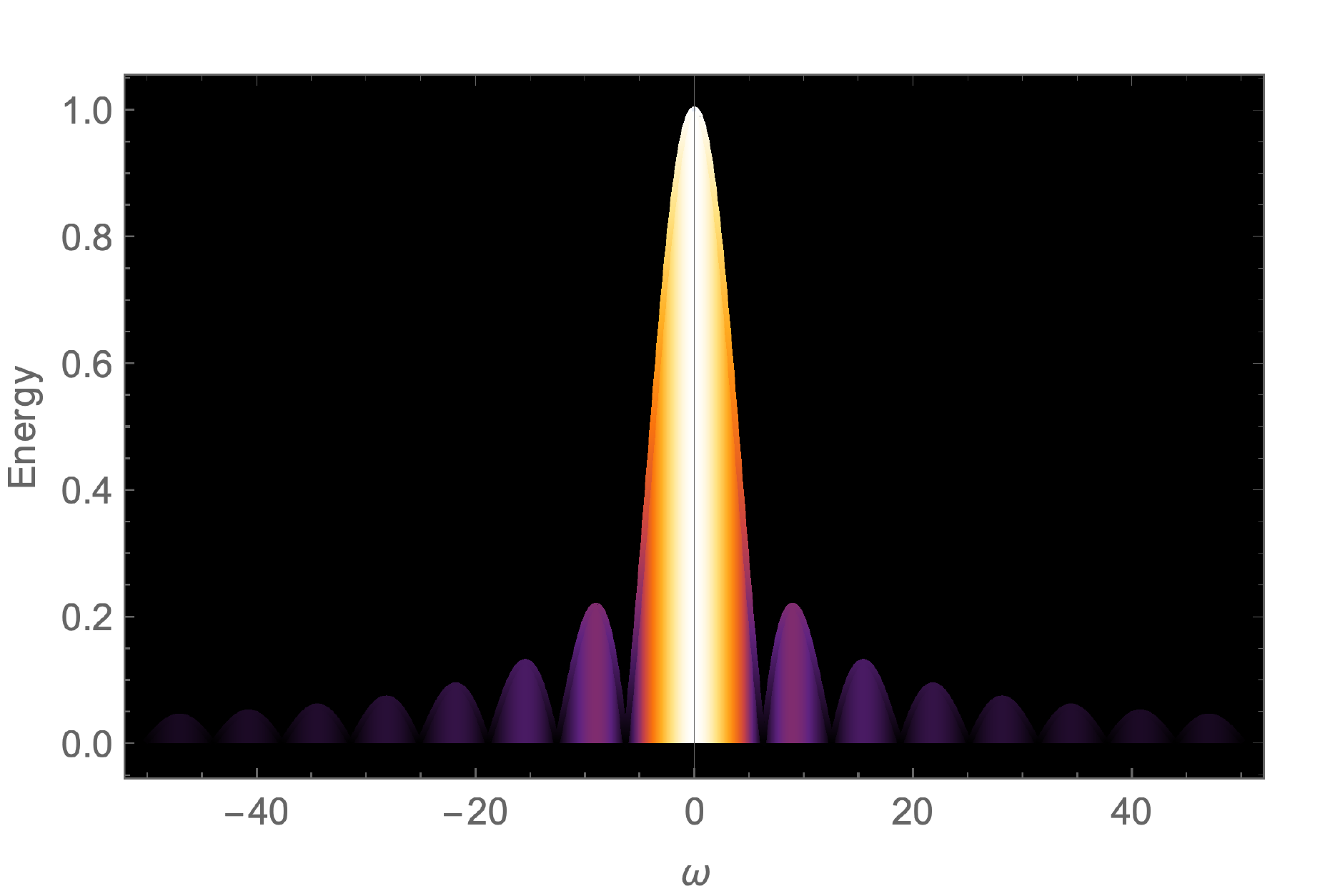}
\caption{Activation location spectra of $\text{Uniform}[0,1]$}\label{fig:cm}
\end{figure}
\FloatBarrier

\begin{figure}[h!]
\centering
\includegraphics[width=7in]{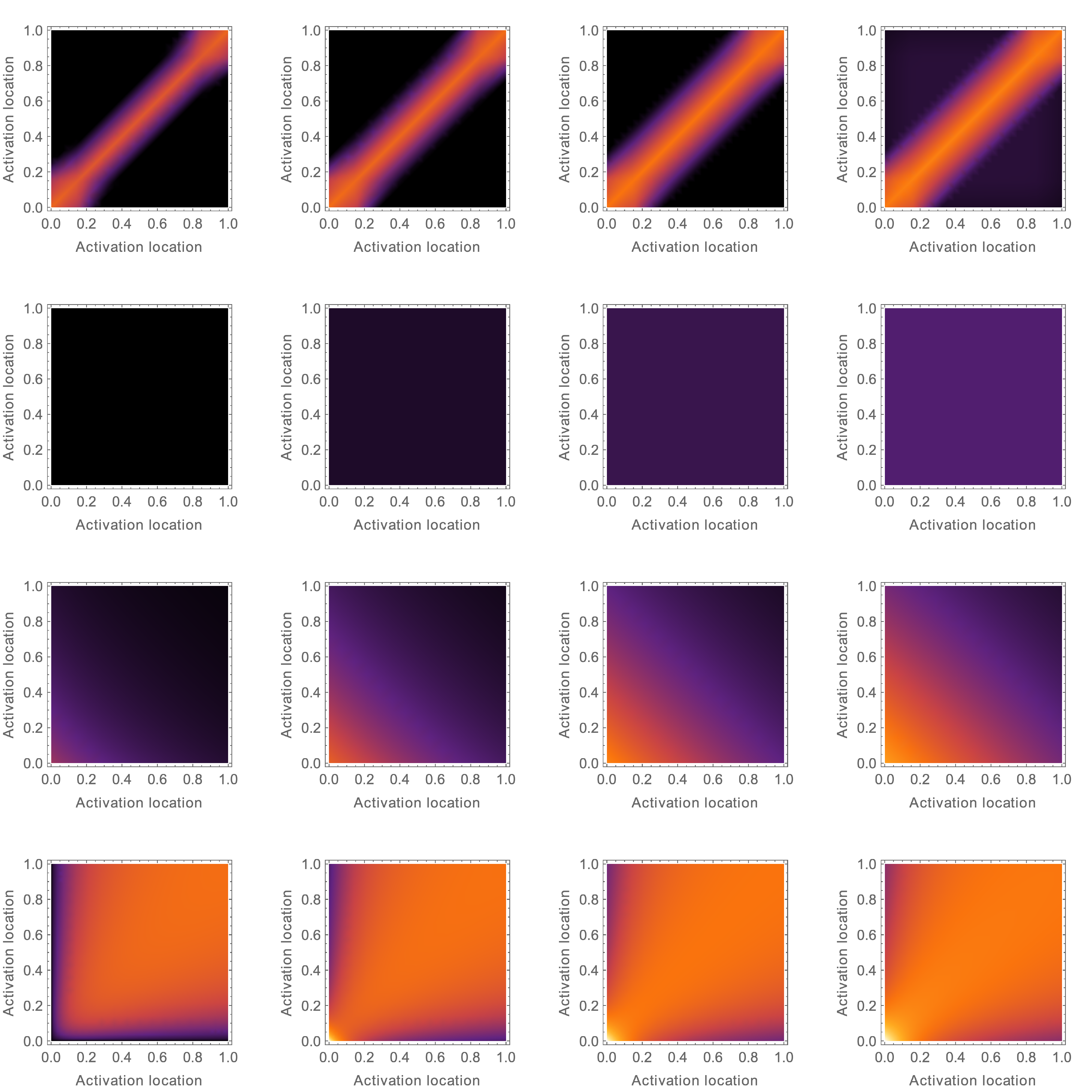}
\caption{Spatial point correlations on $[0,1]\times[0,1]$ for counting mean $c=50$, counting variance $\delta^2\in\{0,25,50,75\}$ (left to right) and varying density kernels (local, decentral, subcentral, central; top to bottom) having common kernel mass $p=1/10$}\label{fig:cm}
\end{figure}
\FloatBarrier

\begin{figure}[h!]
\centering
\includegraphics[width=6in]{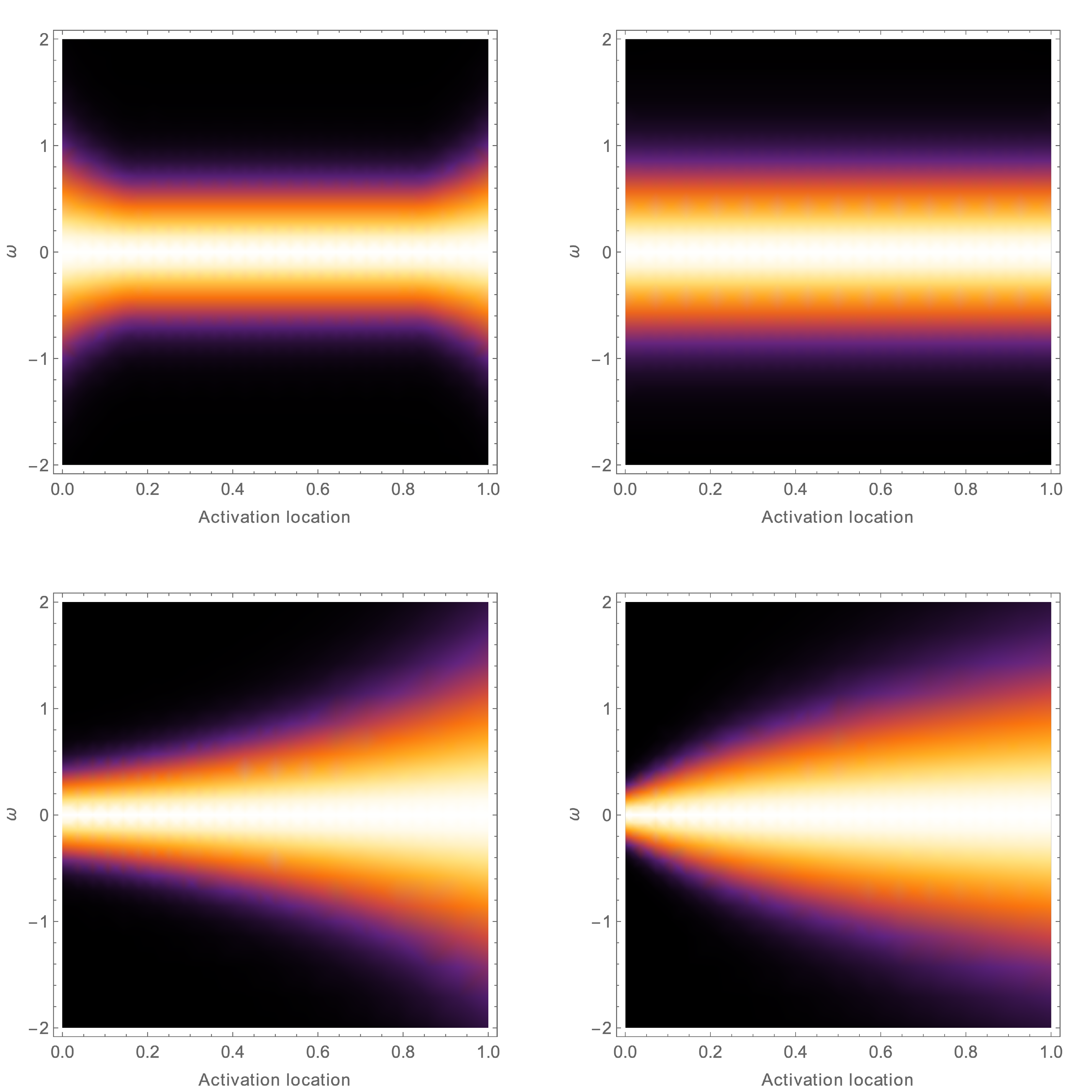}
\caption{Density field spectra (top row: local (L), decentral (R); bottom row: subcentral (L), central (R)) for $K\sim\text{Poisson}(c=50)$ having common kernel mass $p=1/10$}\label{fig:cm}
\end{figure}
\FloatBarrier

\begin{figure}[h!]
\centering
\includegraphics[width=7in]{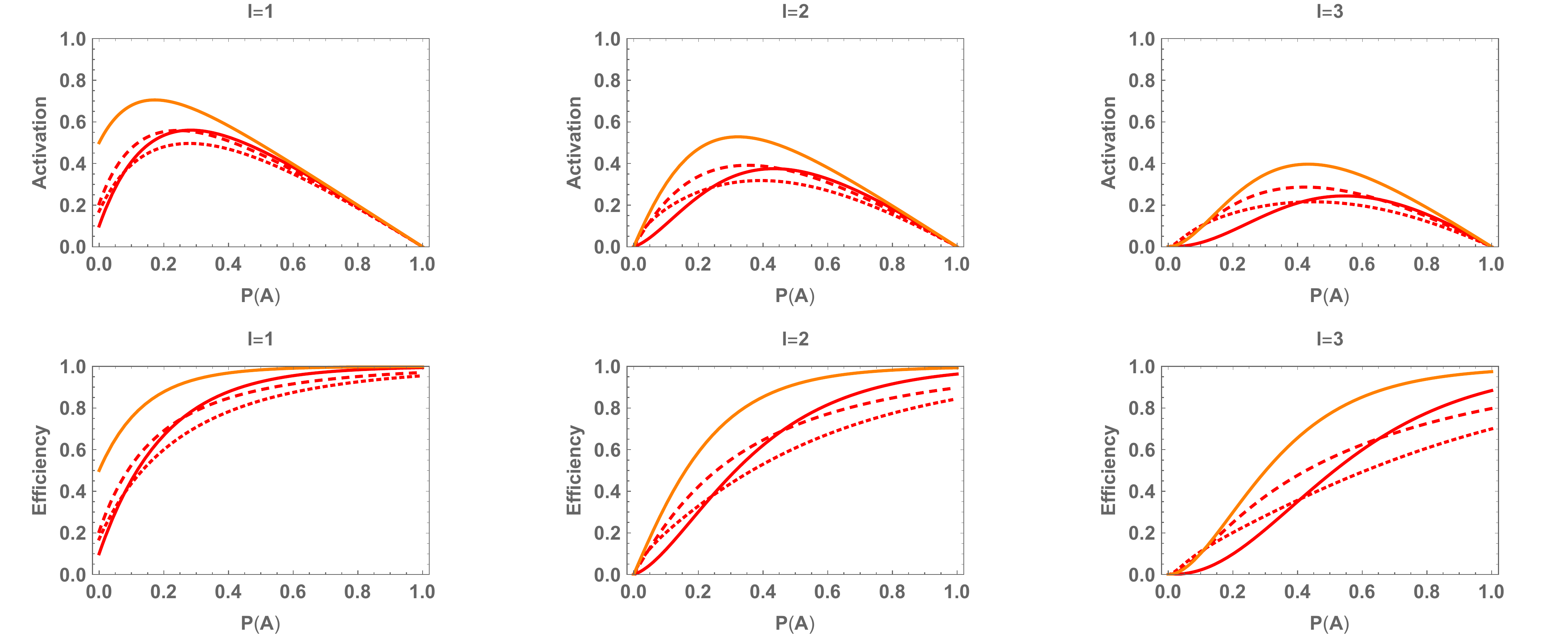}
\caption{Transduction activations and efficiencies for $K\sim\text{Poisson}(50)$ and varying density kernels (orange is local; red solid is decentral; red small dash is central; red large dash is subcentral) having common kernel mass $p=1/10$ and varying transduction activation thresholds $l\in\{1,2,3\}$ (left to right)}\label{fig:cm}
\end{figure}
\FloatBarrier

\begin{figure}[h!]
\centering
\includegraphics[width=7in]{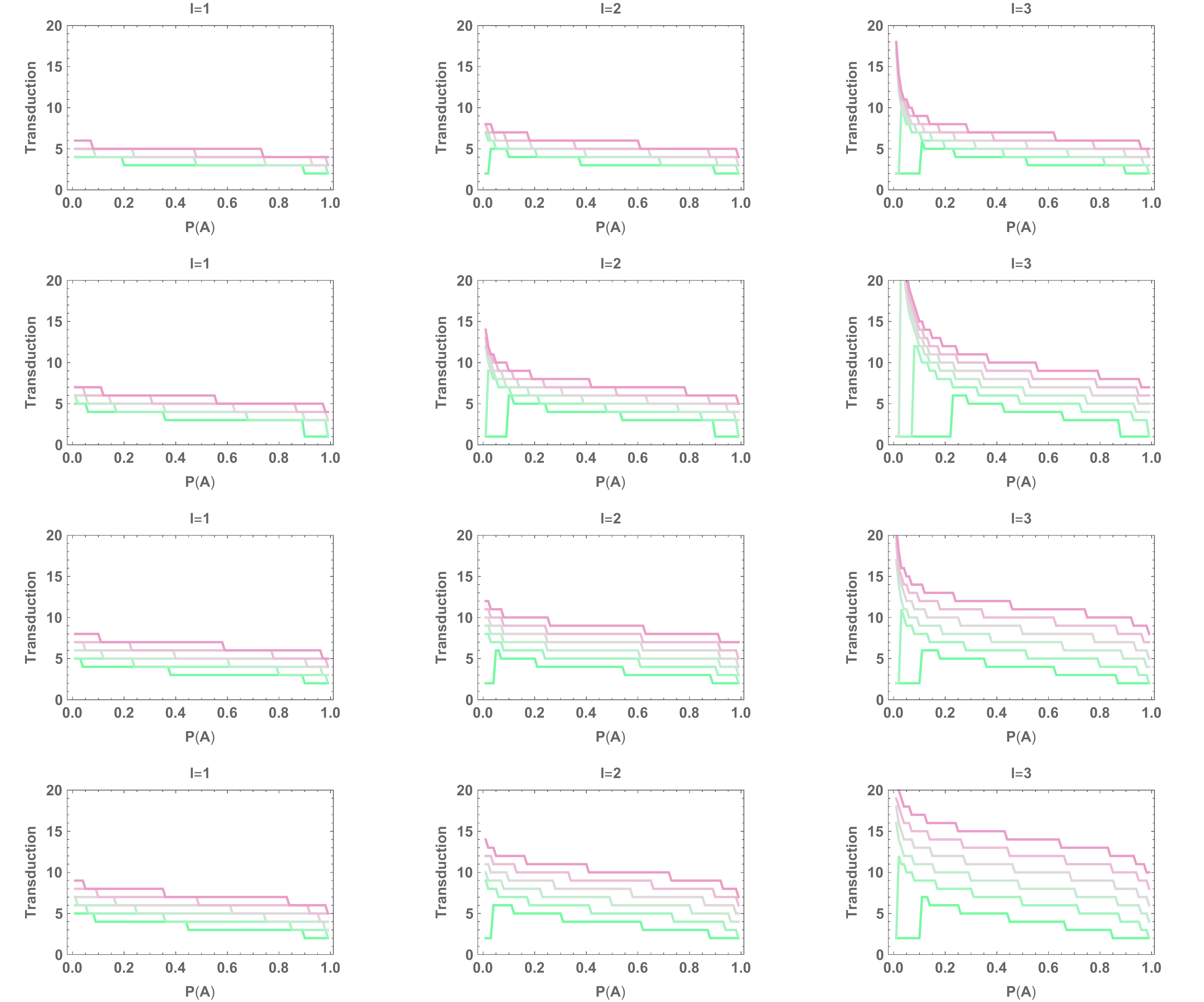}
\caption{Necessary transductions for activation to tolerance $\varepsilon\in\{10^{-1},10^{-2},10^{-3},10^{-5},10^{-6}\}$ (green to pink) on grid $\P(A)=0.01,0.02,\dotsb,0.99$ for $K\sim\text{Poisson}(50)$ and varying density kernels (local, decentral, subcentral, central; top to bottom) having common kernel mass $p=1/10$ and varying transduction activation thresholds $l\in\{1,2,3\}$ (left to right)}\label{fig:cm}
\end{figure}
\FloatBarrier

\begin{figure}[h!]
\centering
\includegraphics[width=7in]{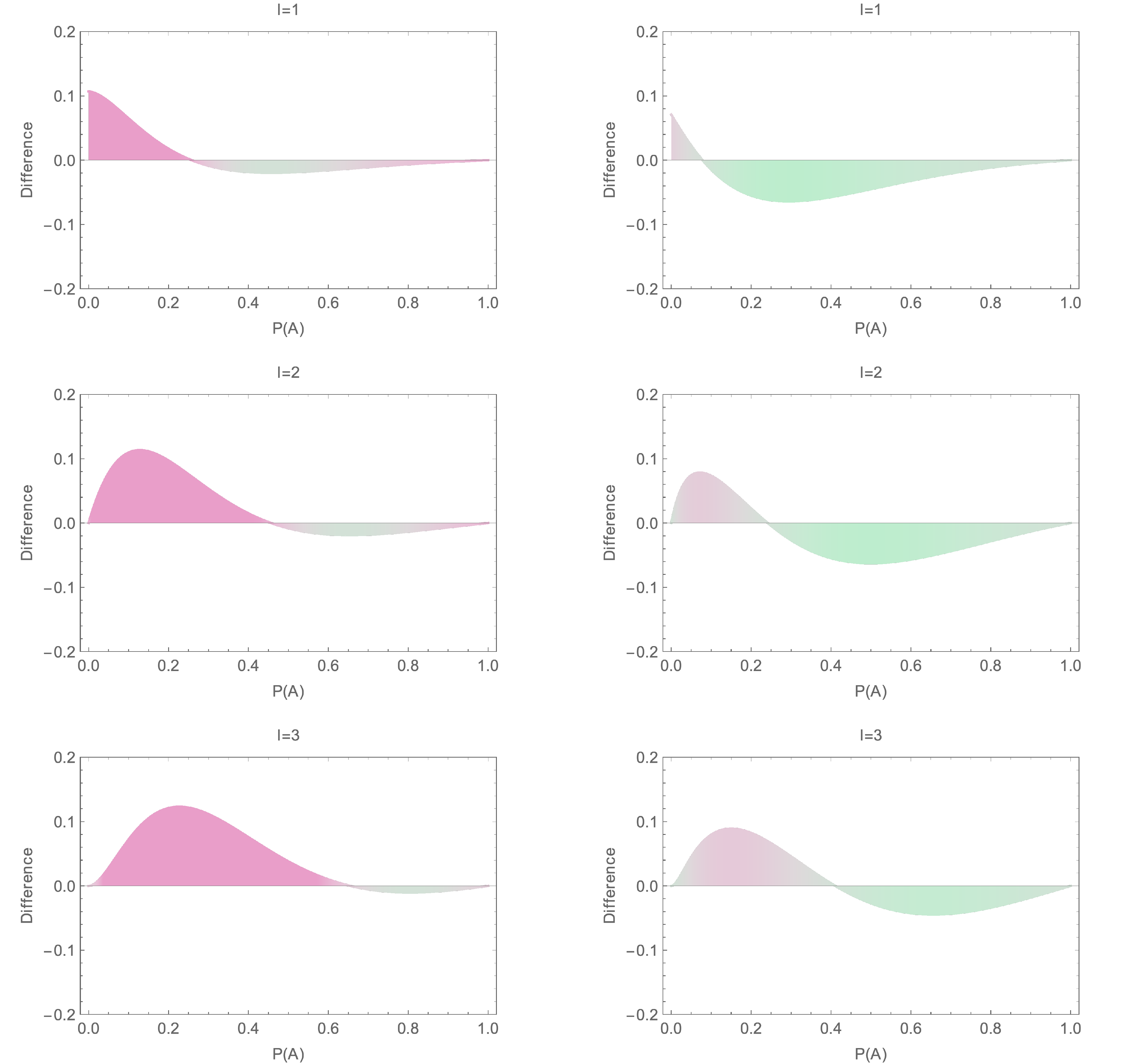}
\caption{Differences relative to decentral kernel (left subcentral, right central; $l=1,2,3$ top to bottom)}\label{fig:cm}
\end{figure}
\FloatBarrier

\begin{table}[h!]
\begin{center}
\begin{tabular}{ccccccccc}
\toprule
&Local & Local &Decentral & Decentral & Subcentral & Subcentral & Central & Central\\
Threshold & Capacity & Gain & Capacity & Gain & Capacity & Gain & Capacity & Gain \\\midrule
$l=1$ & 0.4386 & 4.386 & 0.3557 & 3.557 &0.3609 & 3.609& 0.3277  &3.277 \\
$l=2$ & 0.3272 & 3.272 & 0.2300 & 2.300 & 0.2527  & 2.527& 0.2132  & 2.132  \\
$l=3$ & 0.2363 & 2.363 & 0.1396 & 1.396 & 0.1802  & 1.802 & 0.1447 & 1.447\\
\bottomrule
\end{tabular}
\caption{Transduction activation capacities and efficiencies in activation threshold $l\in\{1,2,3\}$ relative to $K\sim\text{Poisson}(c=50)$ and density kernels having common kernel mass $p=1/10$}\label{tab:stats}
\end{center}
\end{table}

\FloatBarrier

\begin{figure}[h!]
\centering
\includegraphics[width=6in]{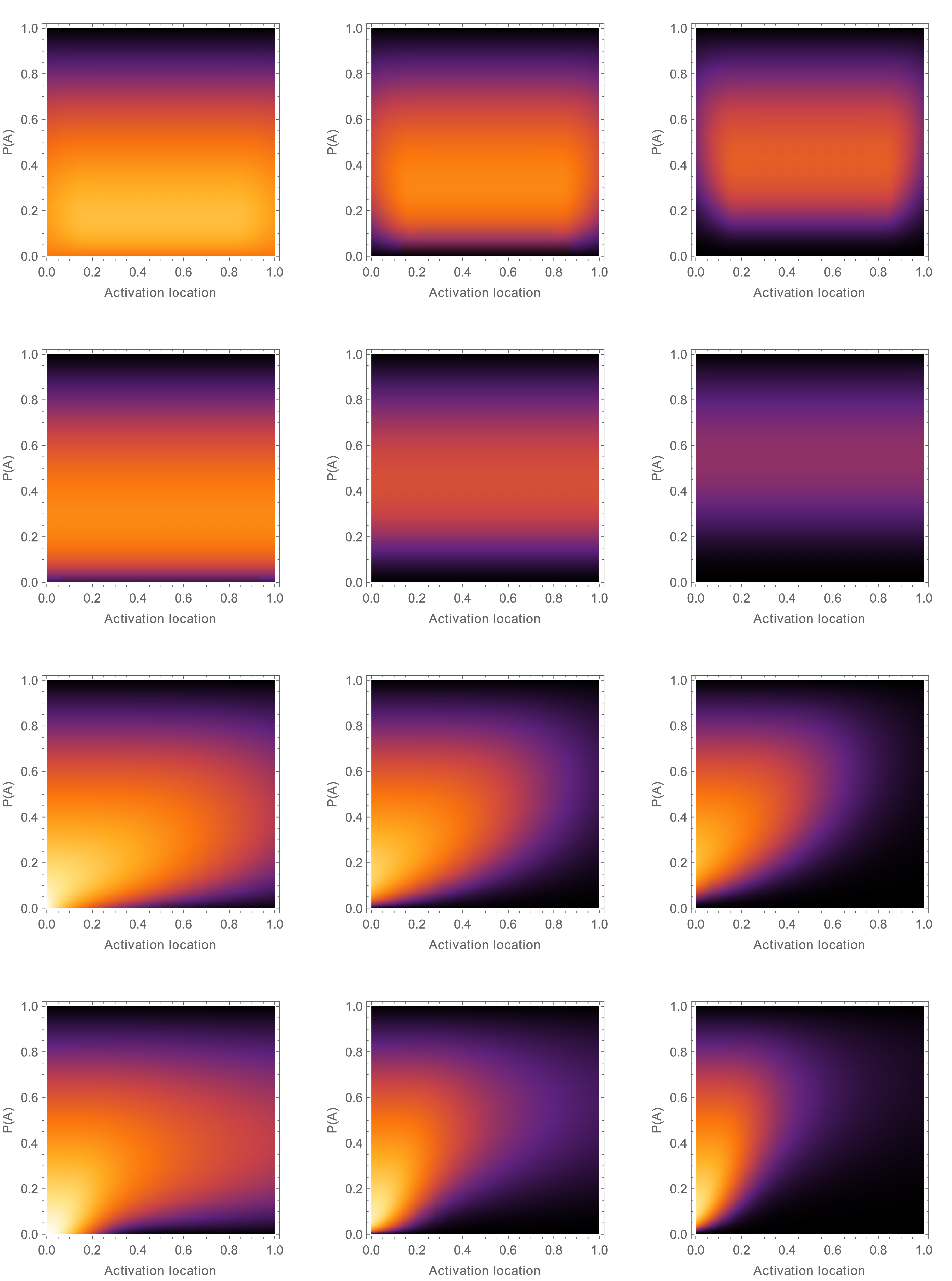}
\caption{Transduction activations for $K\sim\text{Poisson}(50)$ and varying density kernels (local, decentral, subcentral, central; top to bottom) having common kernel mass $p=1/10$ and varying transduction activation thresholds $l\in\{1,2,3\}$ (left to right)}\label{fig:cm}
\end{figure}
\FloatBarrier

\begin{figure}[h!]
\centering
\includegraphics[width=6in]{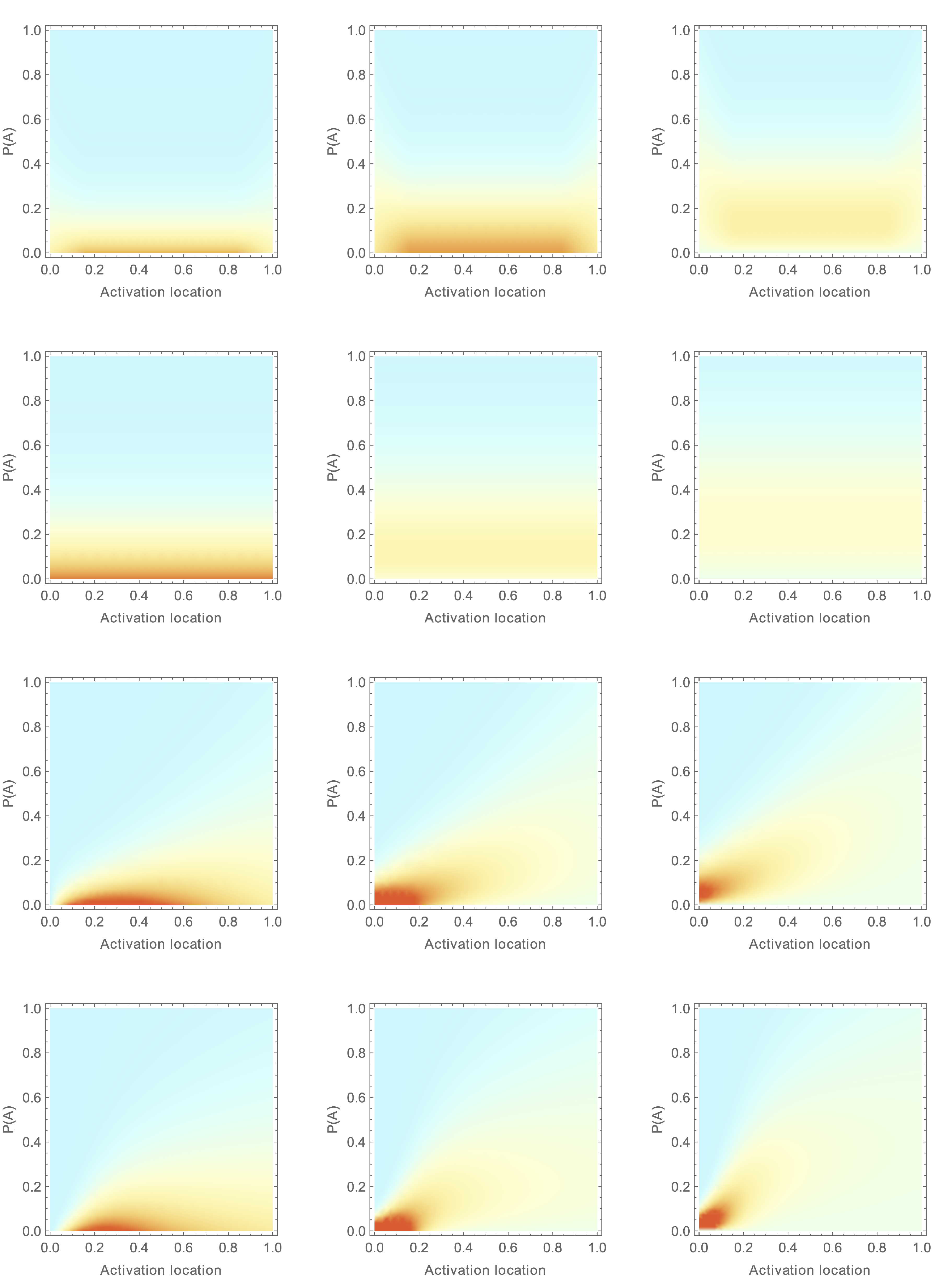}
\caption{Transduction activation flux relative to activation for $K\sim\text{Poisson}(50)$ and varying density kernels (local, decentral, subcentral, central; top to bottom) having common kernel mass $p=1/10$ and varying transduction activation thresholds $l\in\{1,2,3\}$ (left to right)}\label{fig:cm}
\end{figure}
\FloatBarrier

\begin{figure}[h!]
\centering
\includegraphics[width=6in]{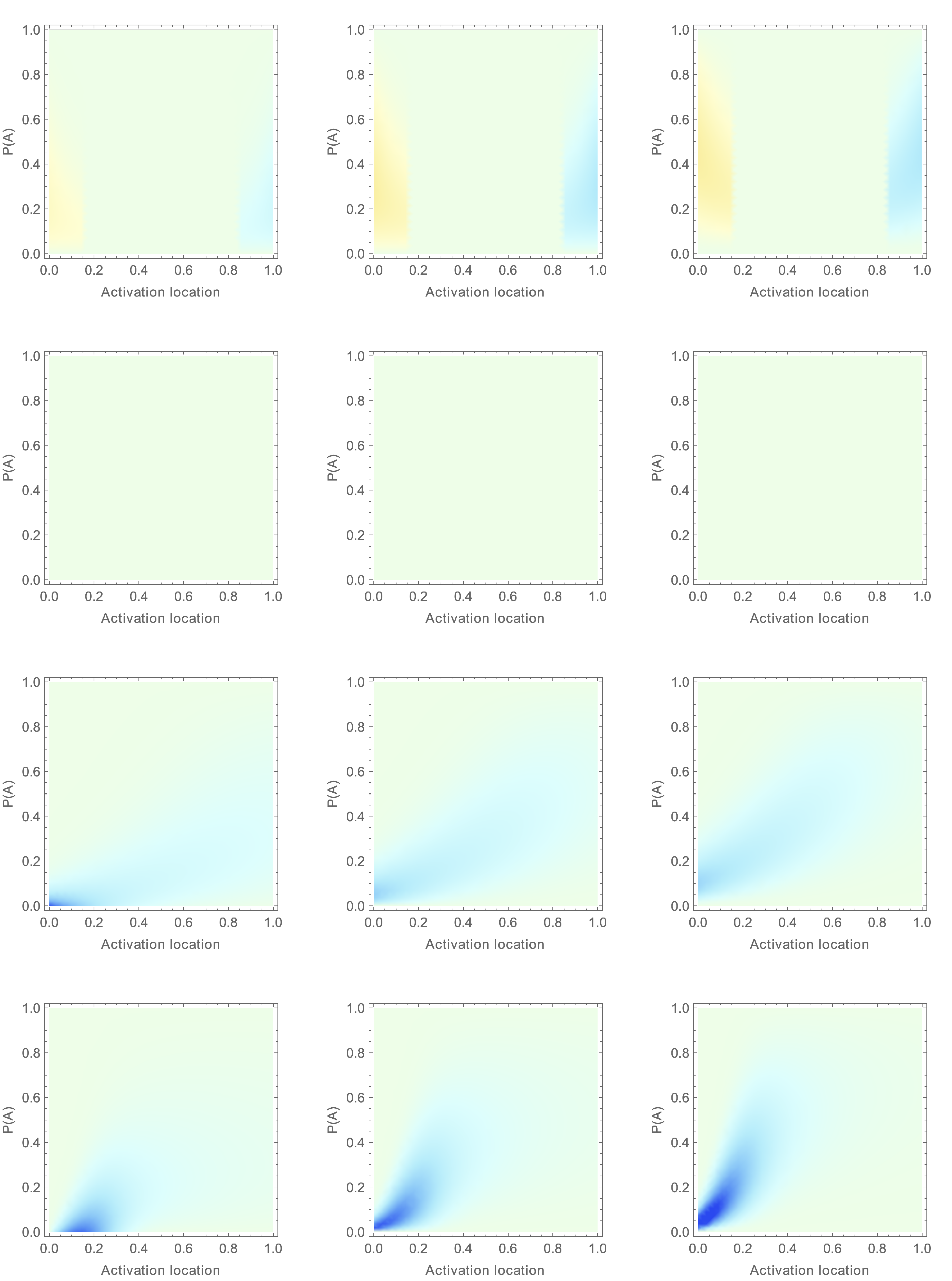}
\caption{Transduction activation flux relative to location for $K\sim\text{Poisson}(50)$ and varying density kernels (local, decentral, subcentral, central; top to bottom) having common kernel mass $p=1/10$ and varying transduction activation thresholds $l\in\{1,2,3\}$ (left to right)}\label{fig:cm}
\end{figure}
\FloatBarrier

\begin{figure}[h!]
\centering
\includegraphics[width=6in]{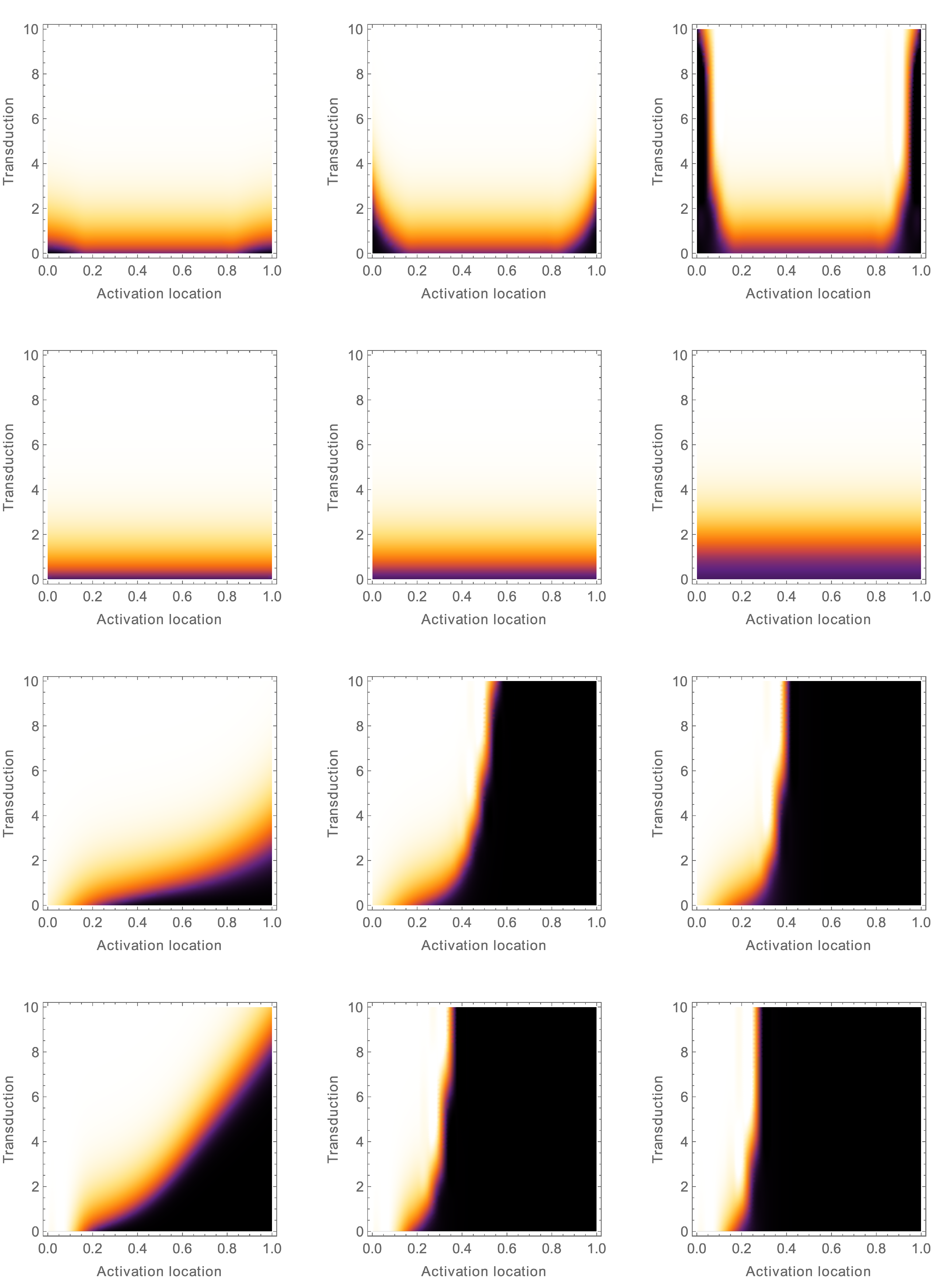}
\caption{Transduction activation in location and time for $K\sim\text{Poisson}(50)$ and varying density kernels (local, decentral, subcentral, central; top to bottom) having common kernel mass $p=1/10$, initializing induction activation thresholds $(11,8,14,13)$, and varying transduction activation thresholds $l\in\{1,2,3\}$}\label{fig:cm}
\end{figure}
\FloatBarrier

\begin{figure}[h!]
\centering
\includegraphics[width=6in]{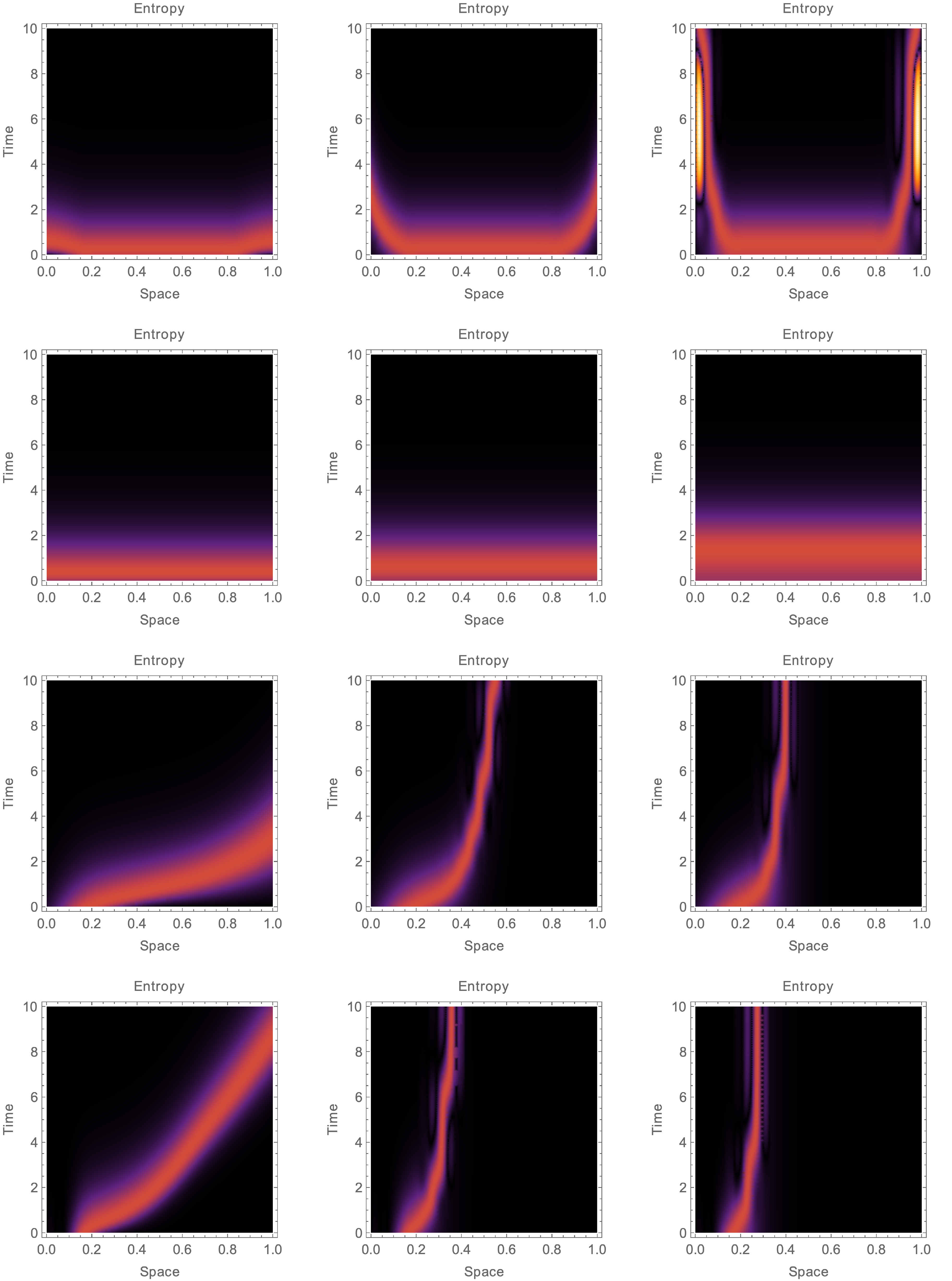}
\caption{Transduction entropy in location and time for $K\sim\text{Poisson}(50)$ and varying density kernels (local, decentral, subcentral, central; top to bottom) having common kernel mass $p=1/10$, initializing induction activation thresholds $(11,14,13)$, and varying transduction activation thresholds $l\in\{1,2,3\}$}\label{fig:cm}
\end{figure}
\FloatBarrier

\newpage
\section{Central induction-transduction activation-deactivation} For these figures, we take $E=[0,1]$, $\nu=\Leb$, and the density kernel as $f(x,y)=\exp_-axy$ on $E\times E$, where $a$ is chosen such that \[(\nu\times\nu)f = \int_{[0,1]^2}\D x\D y e^{-axy}= p\] and use the standard ITAD field equation, i.e. $C_1=0$. 

\begin{figure}[h!]
\centering
\begingroup
\captionsetup[subfigure]{width=5in,font=normalsize}
\subfloat[Activation with threshold $l\in\{1,2,3\}$\label{fig:r0}]{\includegraphics[width=3.5in]{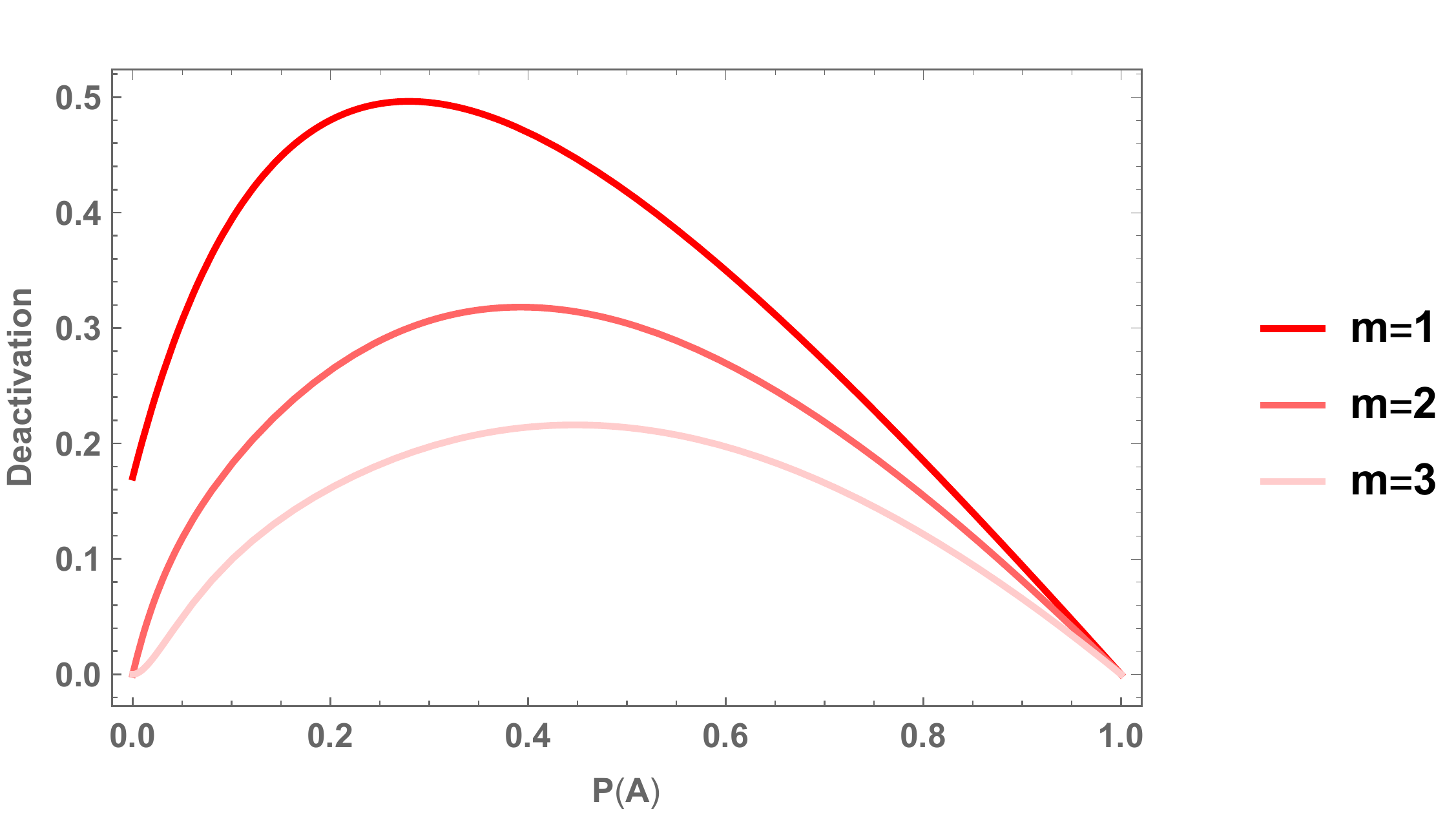}}
\subfloat[Deactivation with threshold $m\in\{1,2,3\}$\label{fig:r1}]{\includegraphics[width=3.5in]{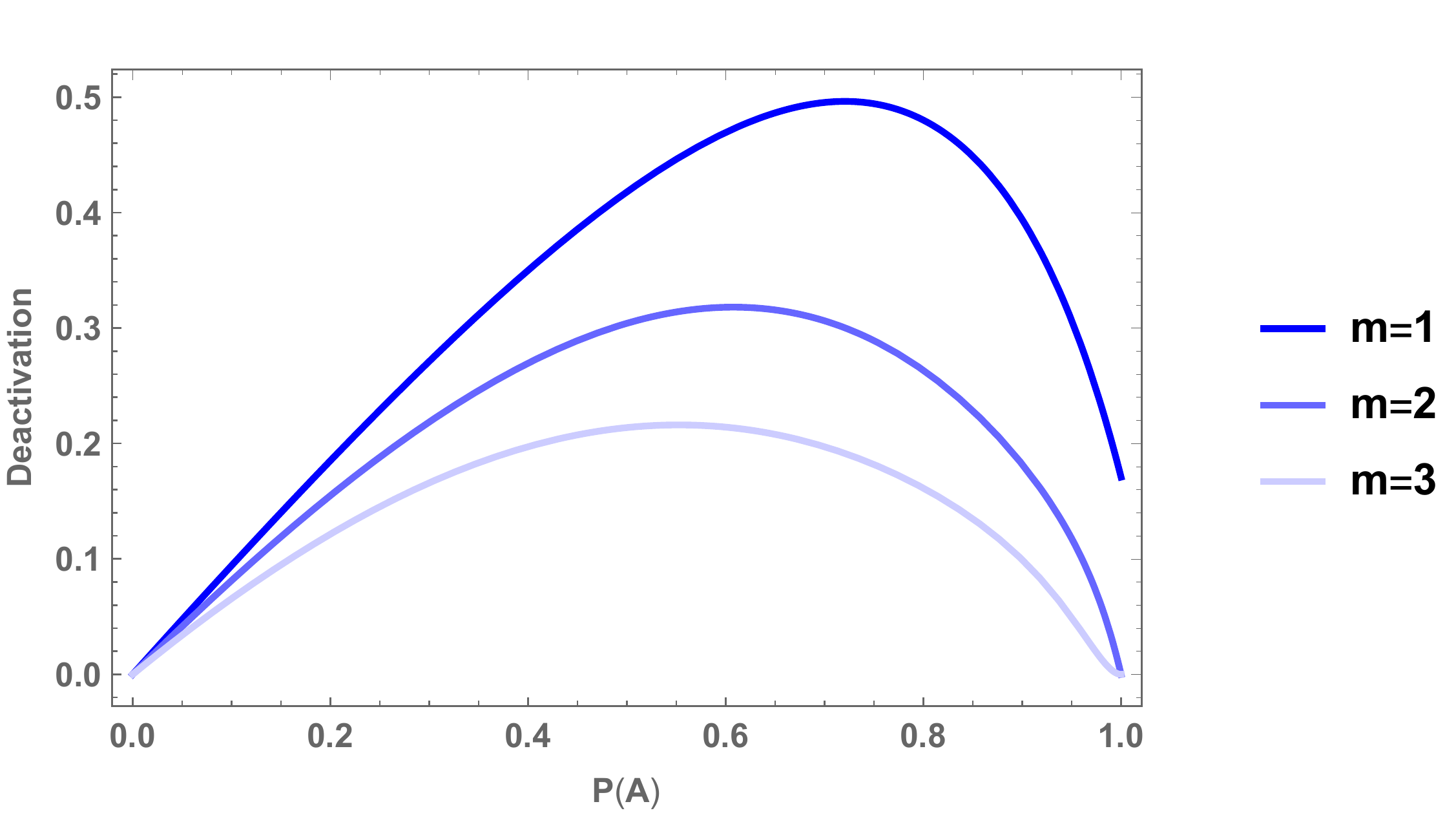}}\\
\subfloat[Net activation for $(l,m)\in\{1,2,3\}\times\{1,2,3\}$\label{fig:r2}]{\includegraphics[width=6.5in]{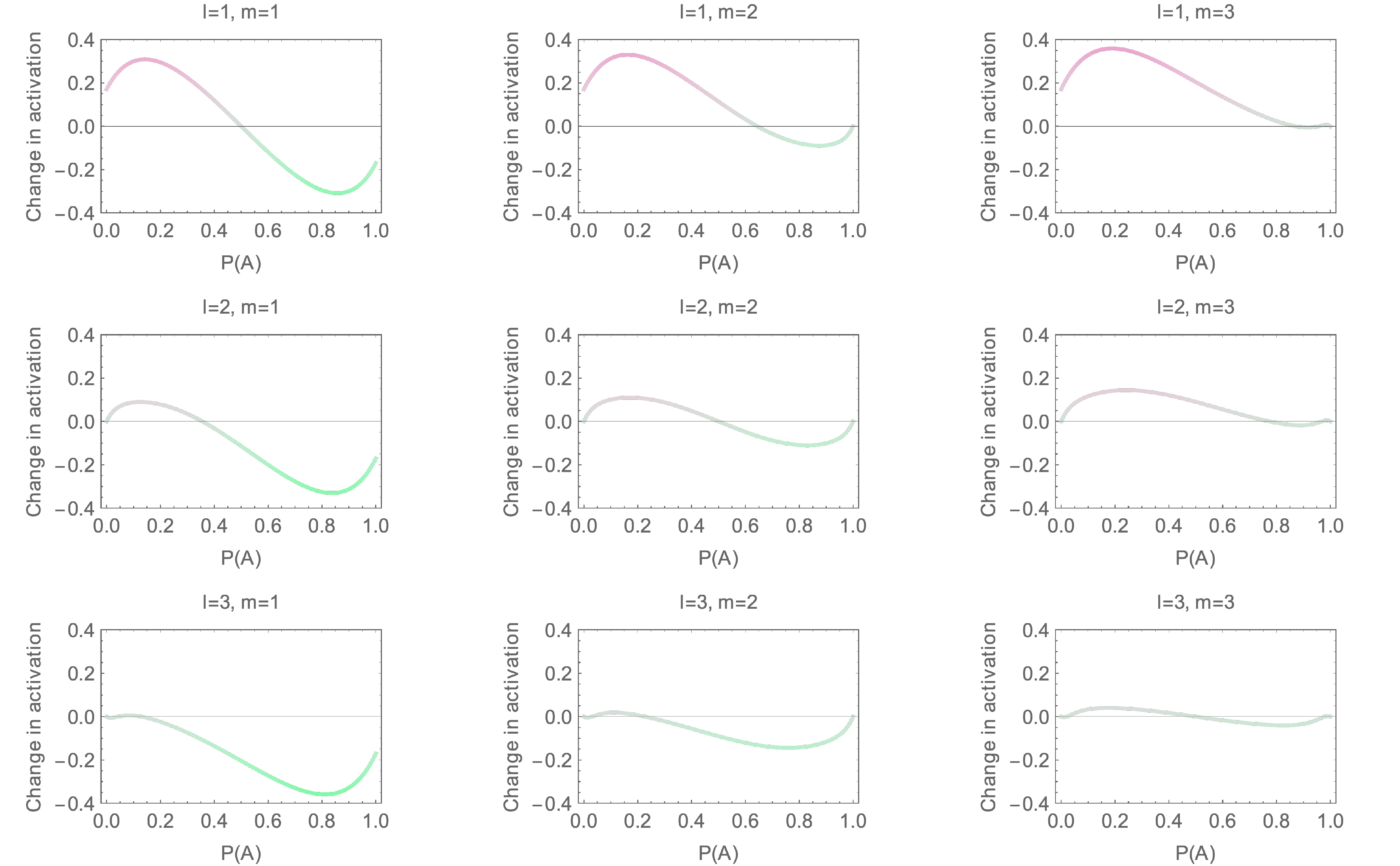}}
\endgroup
\caption{Central transduction in activation, deactivation, and net activation, in activation threshold $l$ and deactivation probability $r$ (blue is zero, red is one), for $K\sim\text{Poisson}(c=50)$ with the density kernel mass $p=1/10$}\label{fig:orbitsdeact}
\end{figure}

\begin{figure}[h!]
\centering
\includegraphics[width=6.5in]{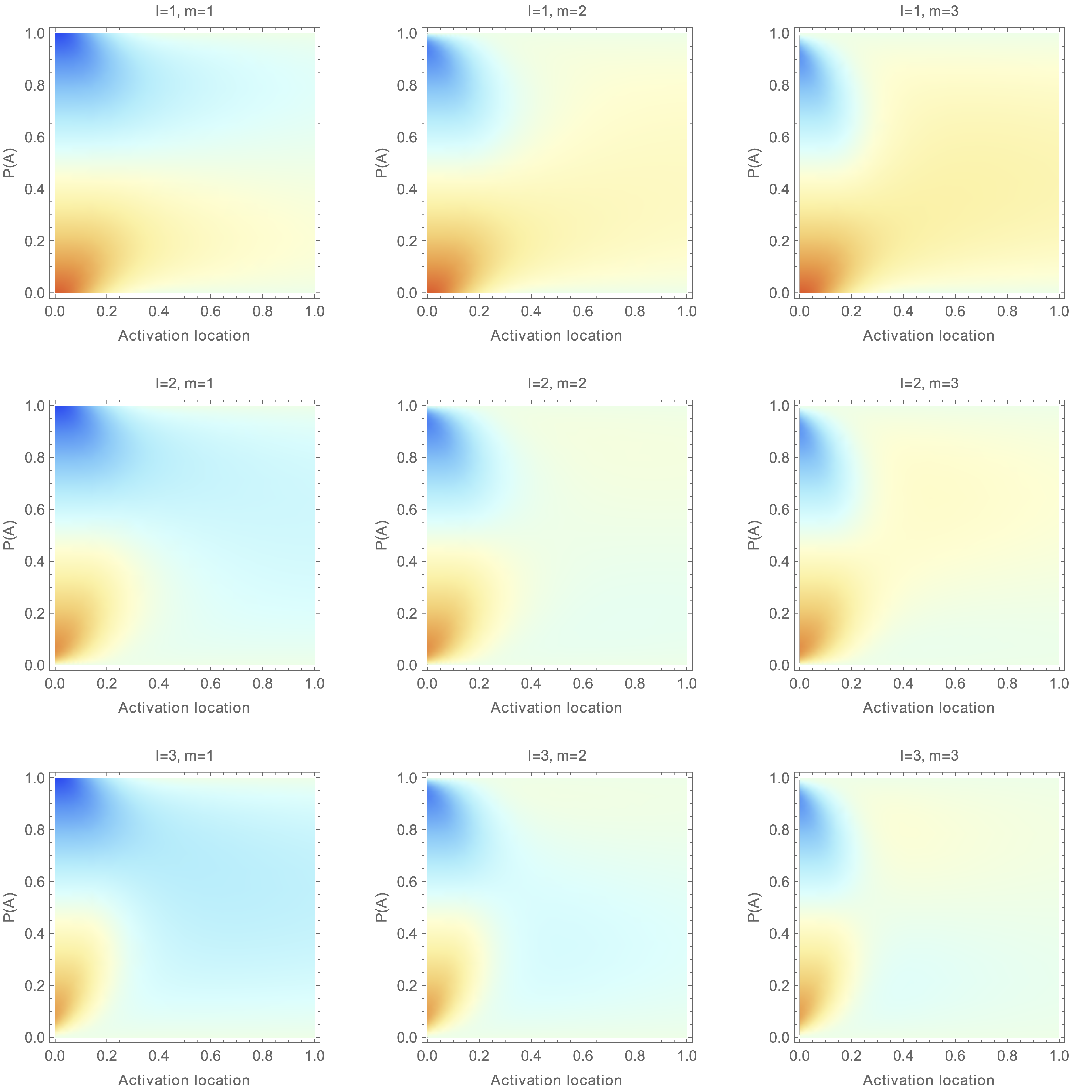}
\caption{Central transduction in activation $l$ and deactivation $m$ thresholds for number of points $K\sim\text{Poisson}(c=50)$, density kernel mass $p=1/10$, initializing activation threshold $k=8$, yielding $\xi_0\{A\}\simeq0.1438$}\label{fig:sifp}
\end{figure}
\FloatBarrier

\begin{figure}[h!]
\centering
\includegraphics[width=6.5in]{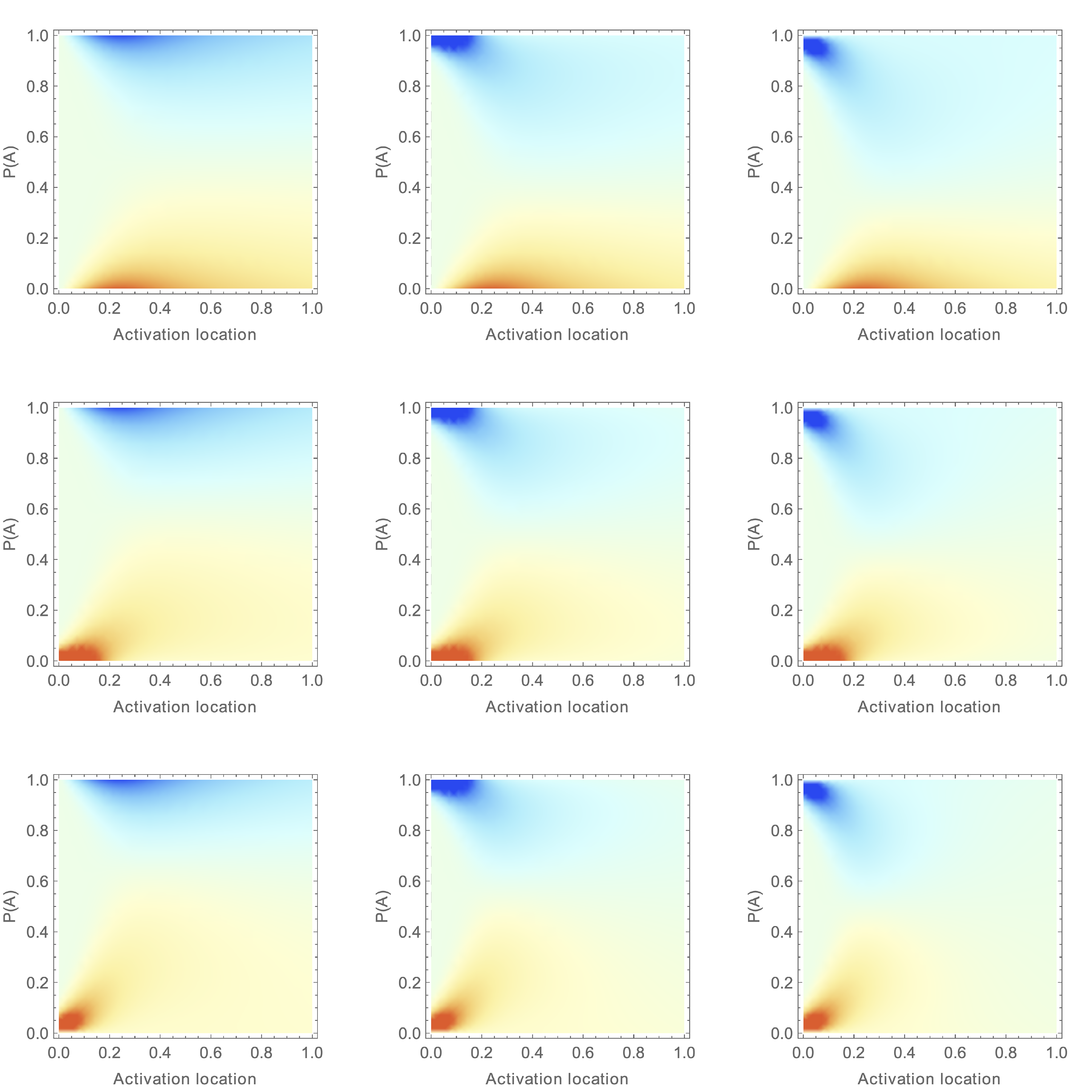}
\caption{Central transduction activation flux relative to activation in activation $l$ and deactivation $m$ thresholds for number of points $K\sim\text{Poisson}(c=50)$, density kernel mass $p=1/10$, initializing activation threshold $k=8$, yielding $\xi_0\{A\}\simeq0.1438$}\label{fig:sifp}
\end{figure}
\FloatBarrier

\begin{figure}[h!]
\centering
\includegraphics[width=6.5in]{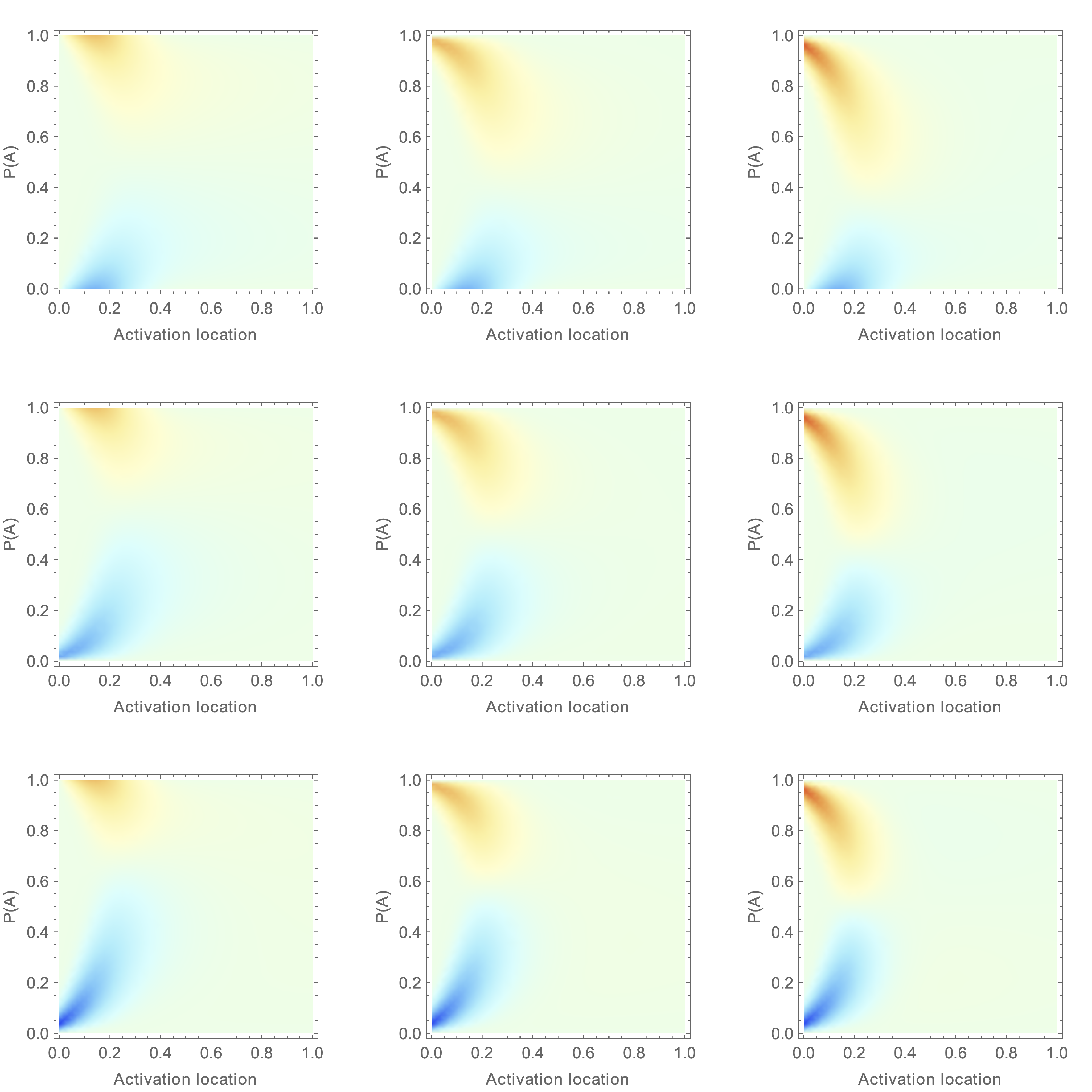}
\caption{Central transduction activation flux relative to location in activation $l$ and deactivation $m$ thresholds for number of points $K\sim\text{Poisson}(c=50)$, density kernel mass $p=1/10$, initializing activation threshold $k=8$, yielding $\xi_0\{A\}\simeq0.1438$}\label{fig:sifp}
\end{figure}
\FloatBarrier

\begin{figure}[h!]
\centering
\includegraphics[width=6.5in]{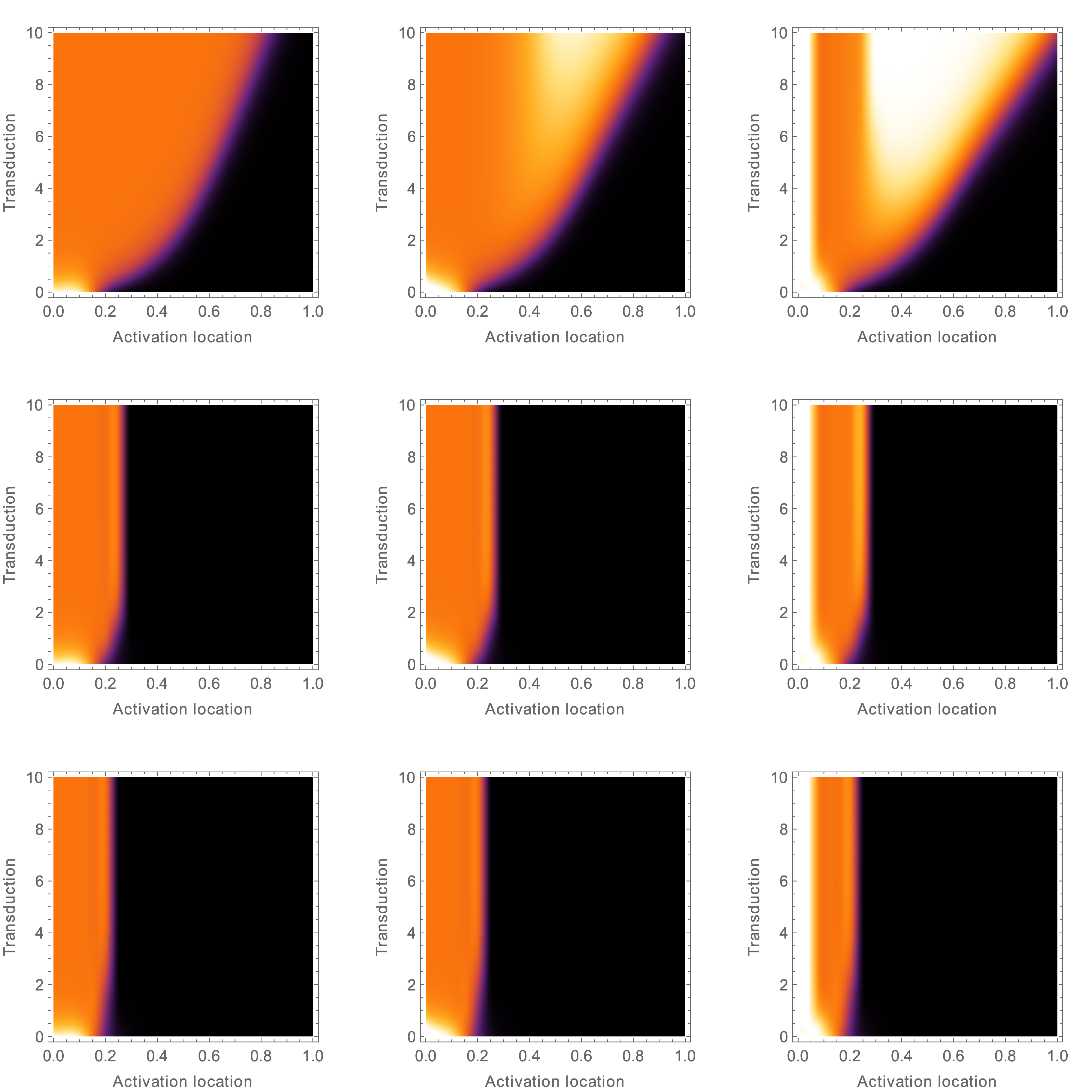}
\caption{Central transduction activation-deactivation in location and time  in activation $l$ and deactivation $m$ thresholds for number of points $K\sim\text{Poisson}(c=50)$, initializing activation threshold $k=13$}\label{fig:sifp}
\end{figure}
\FloatBarrier

\begin{figure}[h!]
\centering
\includegraphics[width=6.5in]{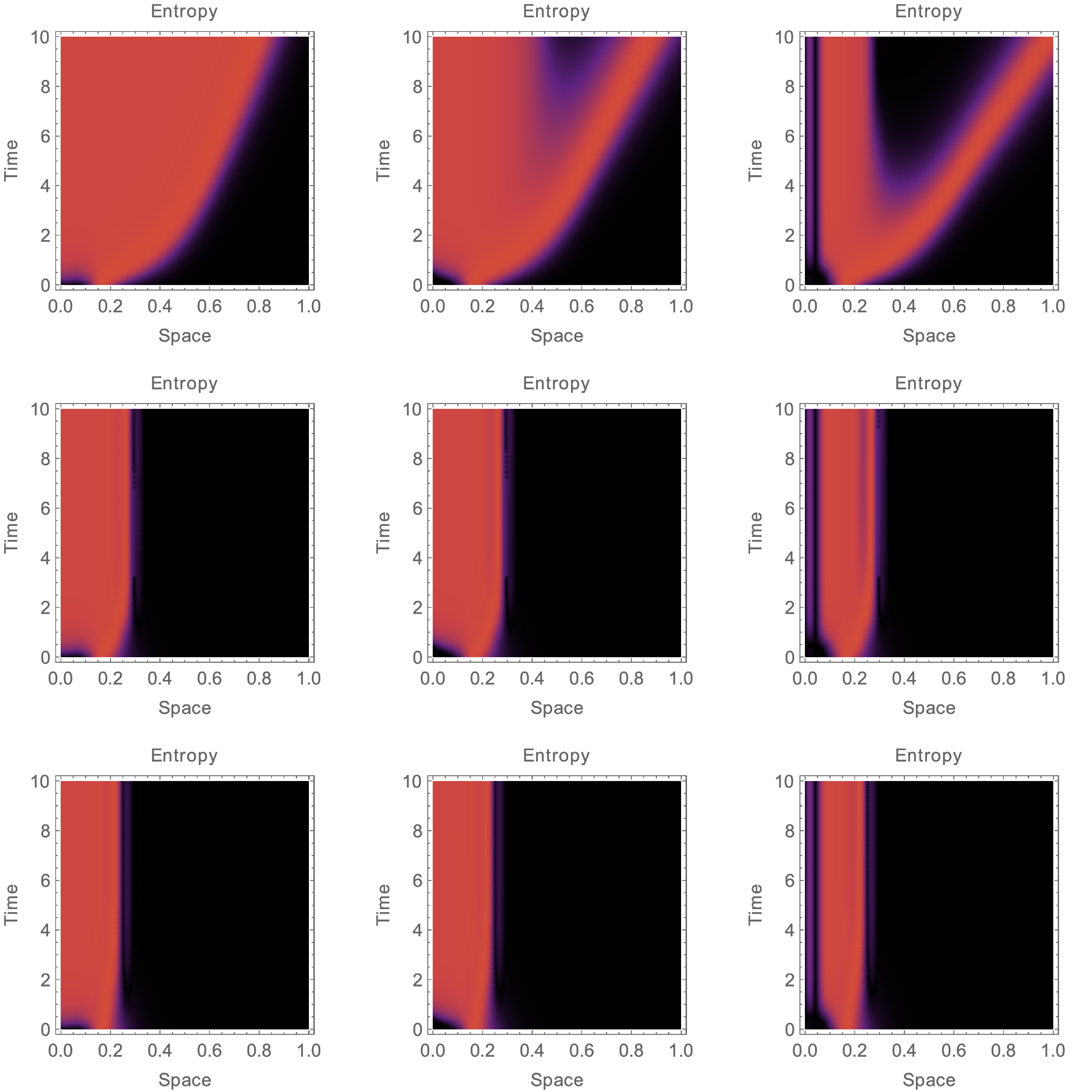}
\caption{Central transduction activation-deactivation entropy in location and time  in activation $l$ and deactivation $m$ thresholds for number of points $K\sim\text{Poisson}(c=50)$, initializing activation threshold $k=13$}\label{fig:sifp}
\end{figure}
\FloatBarrier

\begin{figure}[h!]
\centering
\includegraphics[width=6.5in]{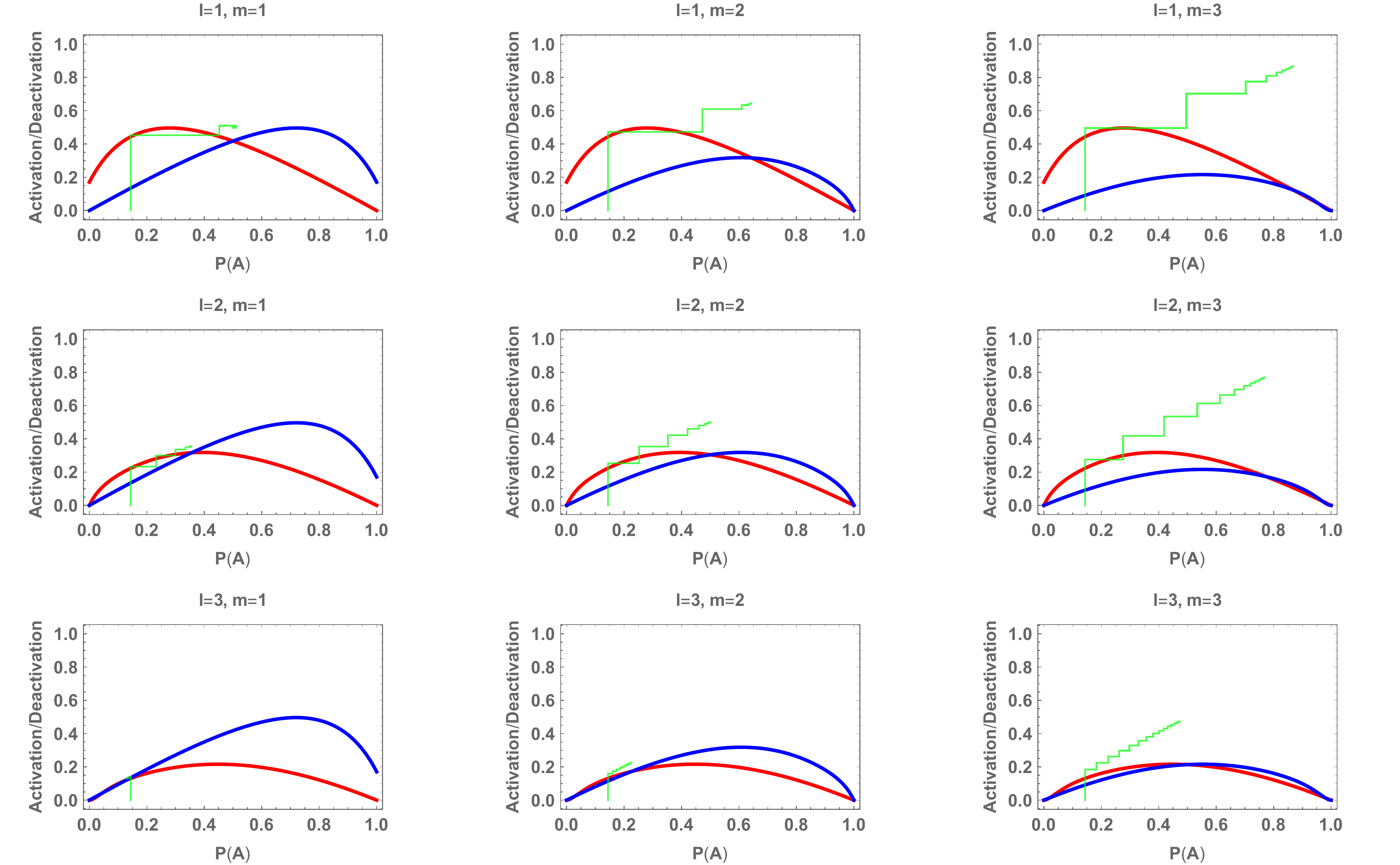}
\caption{Central transduction  in activation $l$ and deactivation $m$ thresholds for number of points $K\sim\text{Poisson}(c=50)$, density kernel mass $p=1/10$, initializing activation threshold $k=8$, yielding $\xi_0\{A\}\simeq0.1438$}\label{fig:sifp}
\end{figure}
\FloatBarrier

\begin{figure}[h!]
\centering
\includegraphics[width=6.5in]{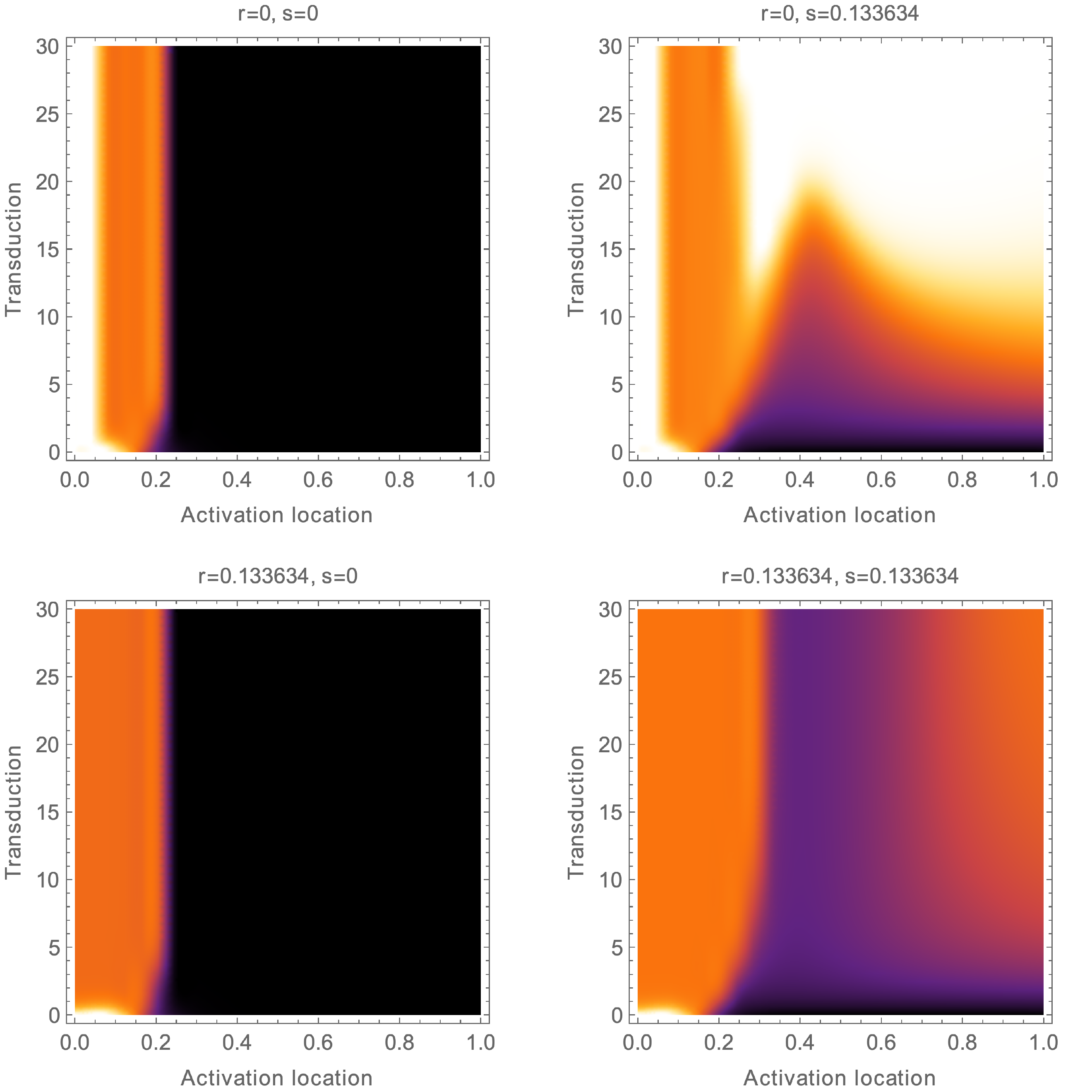}
\caption{Central transduction activation-deactivation in location for external activation $s$ and deactivation $r$ forcings $(r,s)\in\{0,0.1336\}^2$  for activation-deactivation thresholds $l=m=3$ for number of points $K\sim\text{Poisson}(c=50)$, initializing activation threshold $k=13$}\label{fig:sifp}
\end{figure}
\FloatBarrier

\clearpage
\section*{Central induction-transduction activation ($l=1, m=\infty$)} For these figures, we use the regular ITAD field equation ($C_1=0$) and  assume induced wave and energy equation dimensionalized constants to be unity. 

\begin{figure}[h!]
\centering
\includegraphics[width=5in]{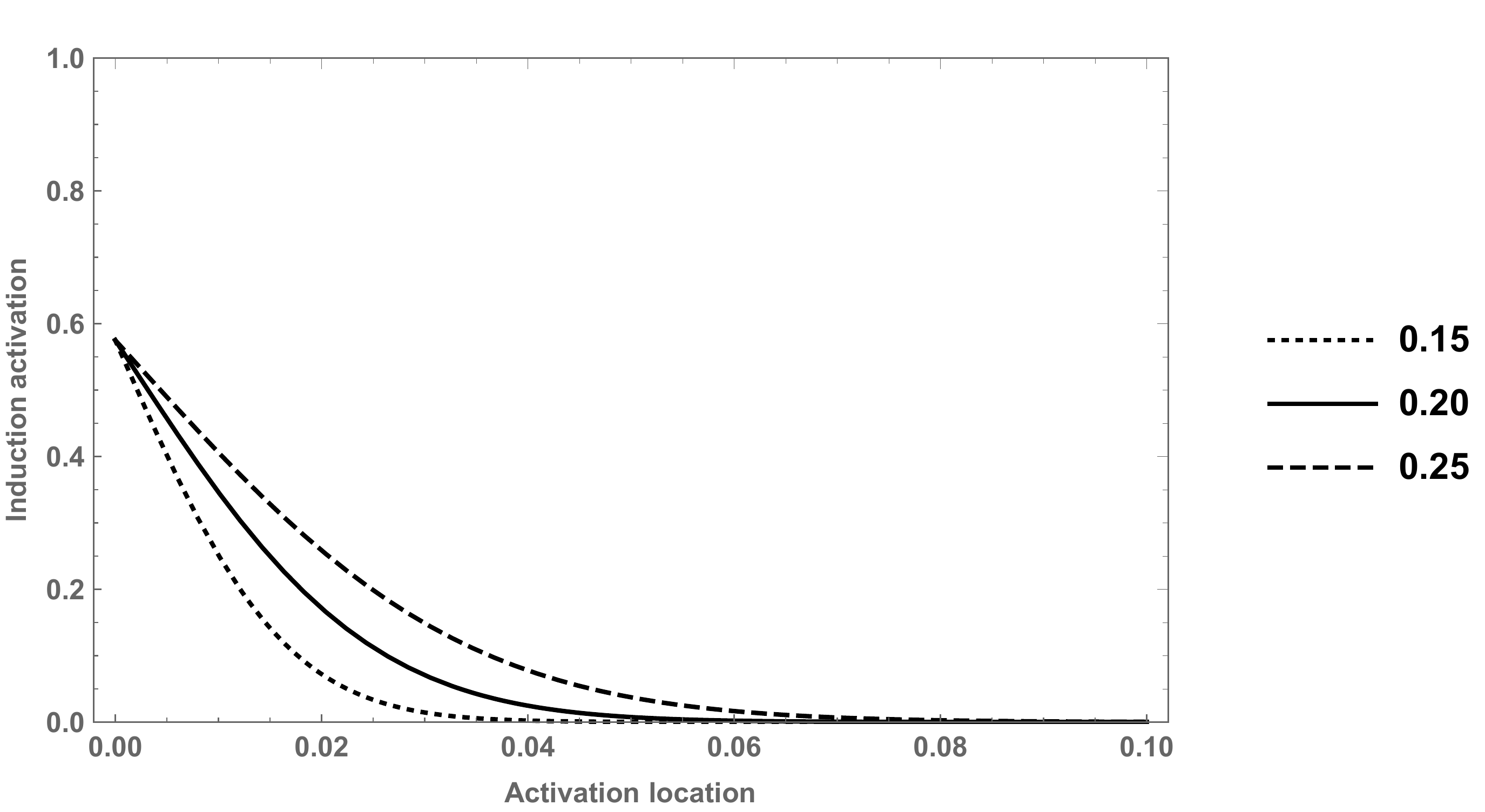}
\caption{Central induction activation for density kernel masses $0.15,0.2,0.25$ and induction activation threshold $k=50$ for $K\sim\text{Poisson}(c=50)$}\label{fig:nn}
\end{figure}

\begin{figure}[h!]
\centering
\includegraphics[width=7in]{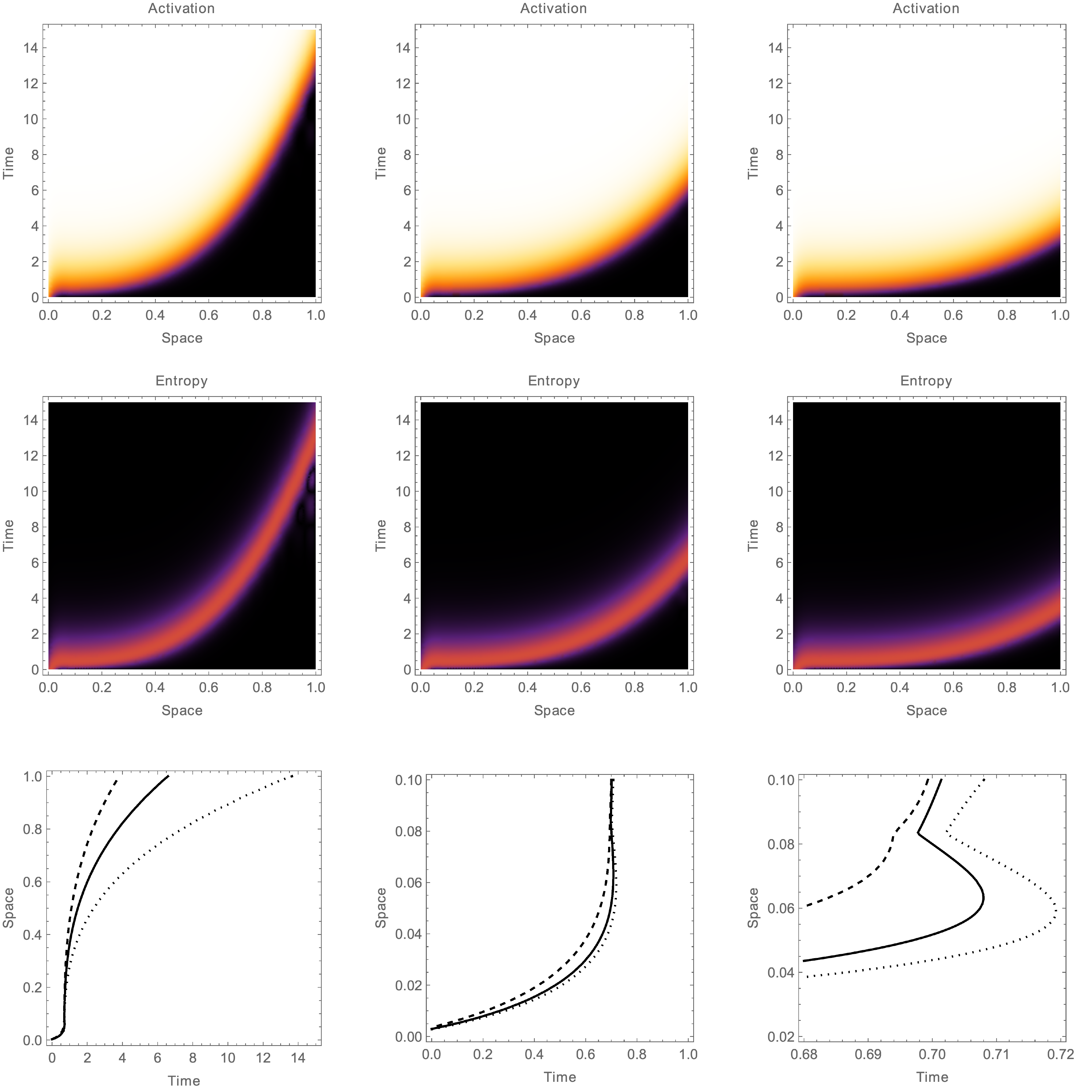}
\caption{Central induction-transduction activation for density kernel masses $0.15,0.2,0.25$ and induction activation threshold $k=50$ for $K\sim\text{Poisson}(c=50)$}\label{fig:nn}
\end{figure}

\begin{figure}[h!]
\centering
\includegraphics[width=6in]{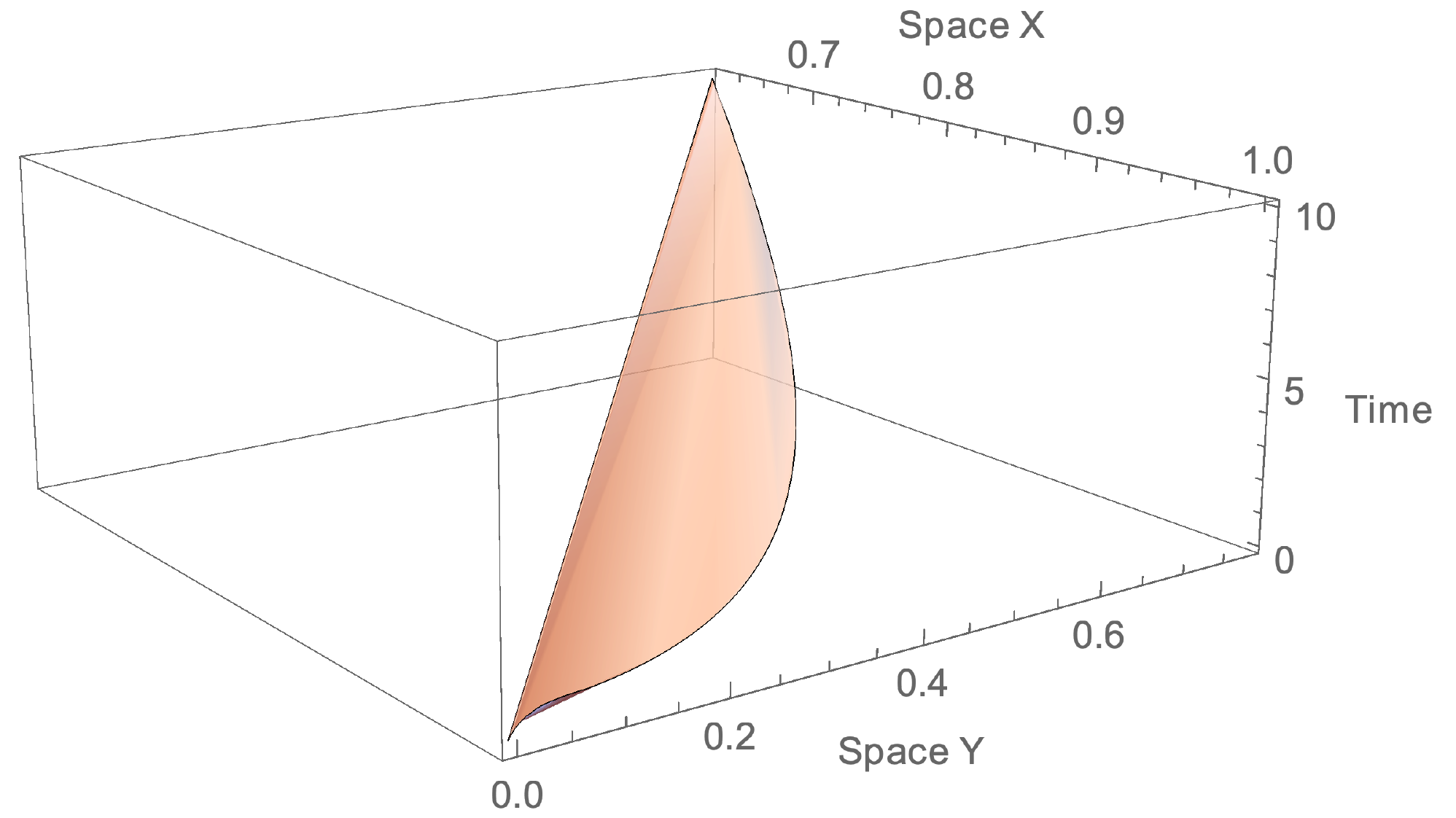}
\caption{Central induction-transduction activation stochastic space-time in coordinates $(\cos(x),\sin(x))$ for density kernel mass $0.15$ and induction activation threshold $k=50$ for $K\sim\text{Poisson}(c=50)$}\label{fig:nn}
\end{figure}

\begin{figure}[h!]
\centering
\includegraphics[width=6in]{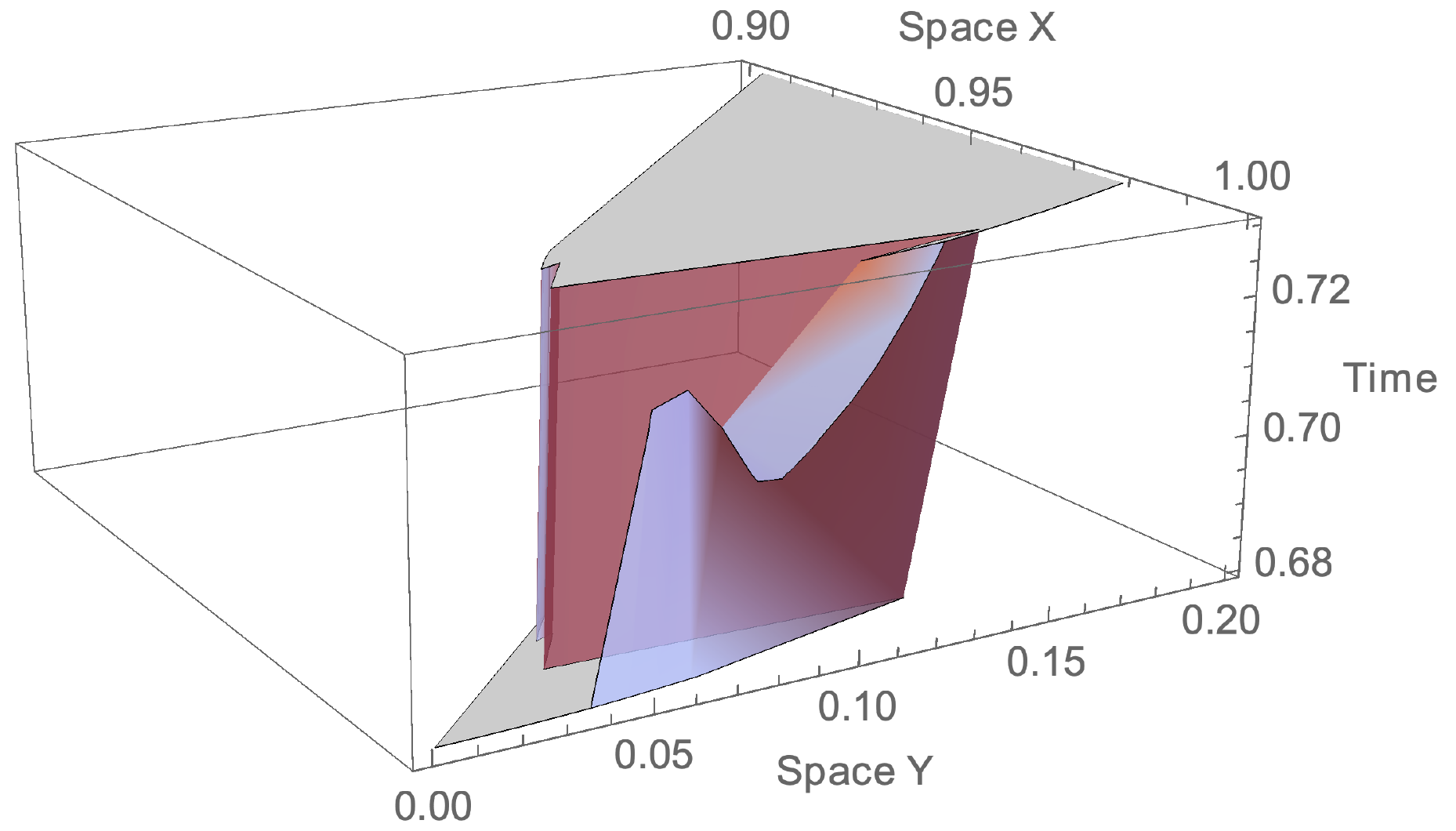}
\caption{Central induction-transduction activation stochastic space-time in coordinates $(\cos(x),\sin(x))$ for density kernel mass $0.15$ and induction activation threshold $k=50$ for $K\sim\text{Poisson}(c=50)$}\label{fig:nn}
\end{figure}

\begin{figure}[h!]
\centering
\includegraphics[width=6.5in]{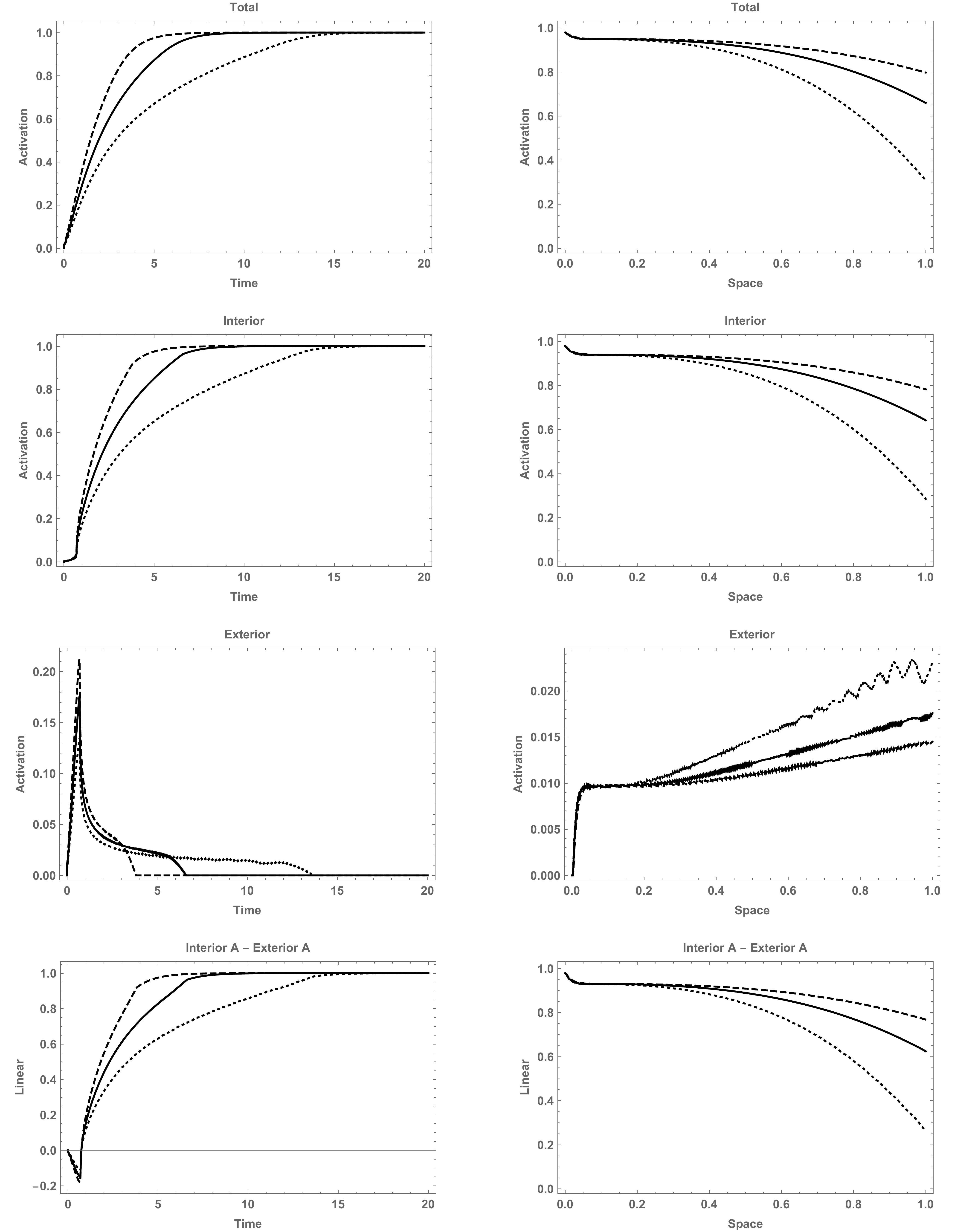}
\caption{Central induction-transduction activation relative to activation frontier for density kernel masses $0.15,0.2,0.25$ and induction activation threshold $k=50$ for $K\sim\text{Poisson}(c=50)$}\label{fig:nn}
\end{figure}

\begin{figure}[h!]
\centering
\includegraphics[width=6.5in]{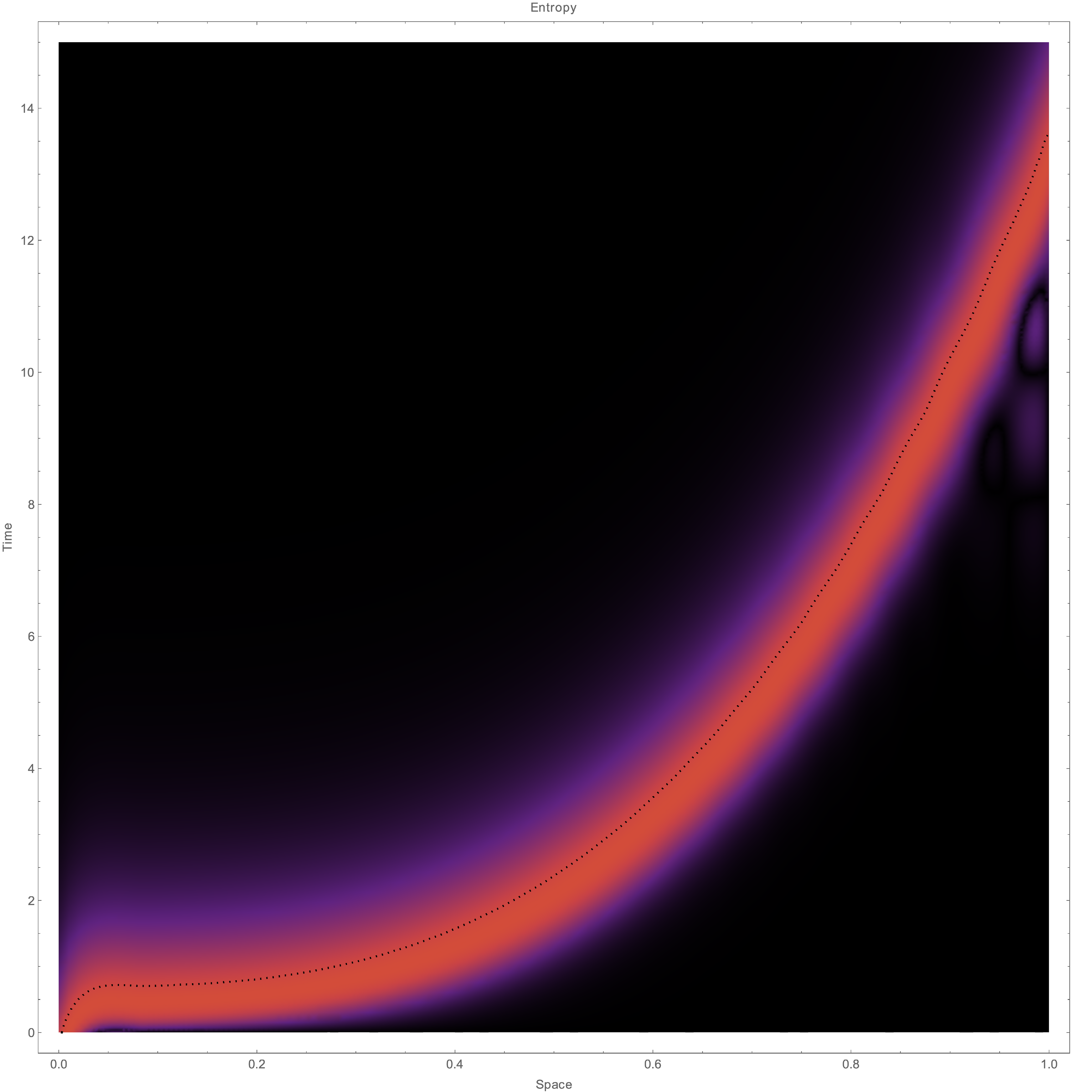}
\caption{Central induction-transduction activation entropy relative to activation frontier for density kernel mass $0.15$ and induction activation threshold $k=50$ for $K\sim\text{Poisson}(c=50)$}\label{fig:nn}
\end{figure}

\begin{figure}[h!]
\centering
\includegraphics[width=6.5in]{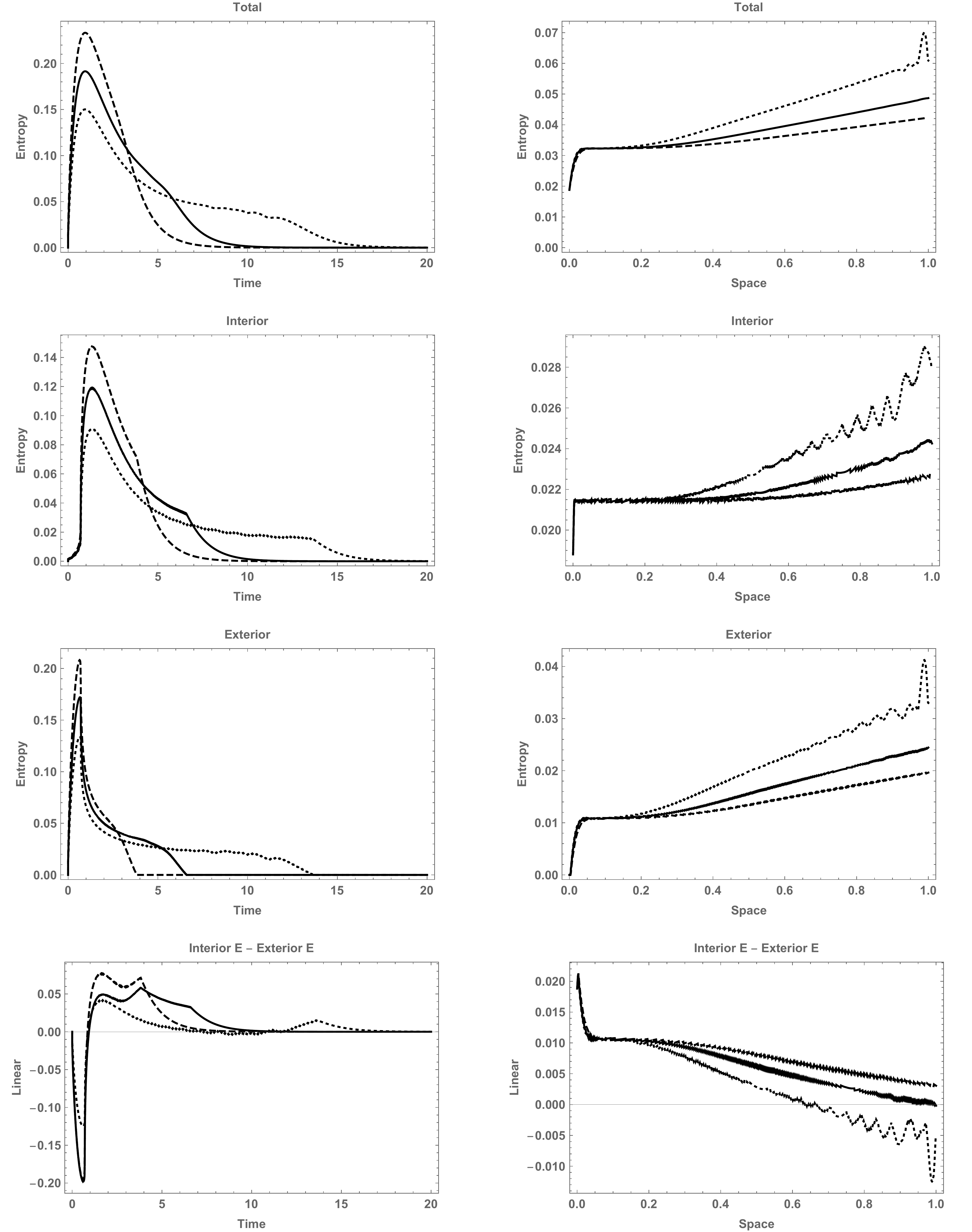}
\caption{Central induction-transduction activation entropies relative to activation frontier for density kernel masses $0.15,0.2,0.25$ and induction activation threshold $k=50$ for $K\sim\text{Poisson}(c=50)$}\label{fig:nn}
\end{figure}

\begin{figure}[h!]
\centering
\includegraphics[width=7in]{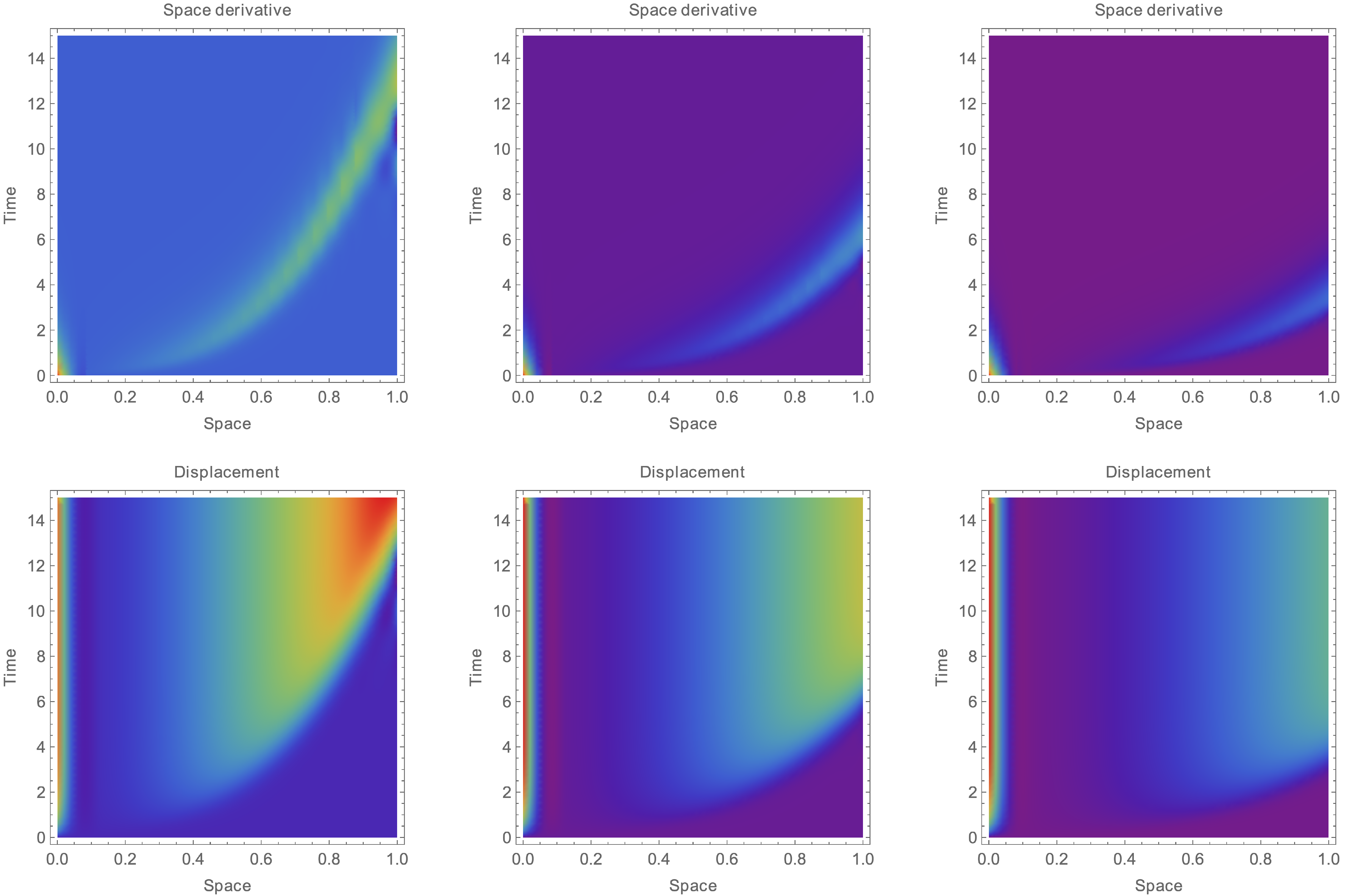}
\caption{Central induction-transduction activation wave propagation for density kernel masses $0.15,0.2,0.25$ and induction activation threshold $k=50$ for $K\sim\text{Poisson}(c=50)$}\label{fig:nn}
\end{figure}

\begin{figure}[h!]
\centering
\includegraphics[width=6.5in]{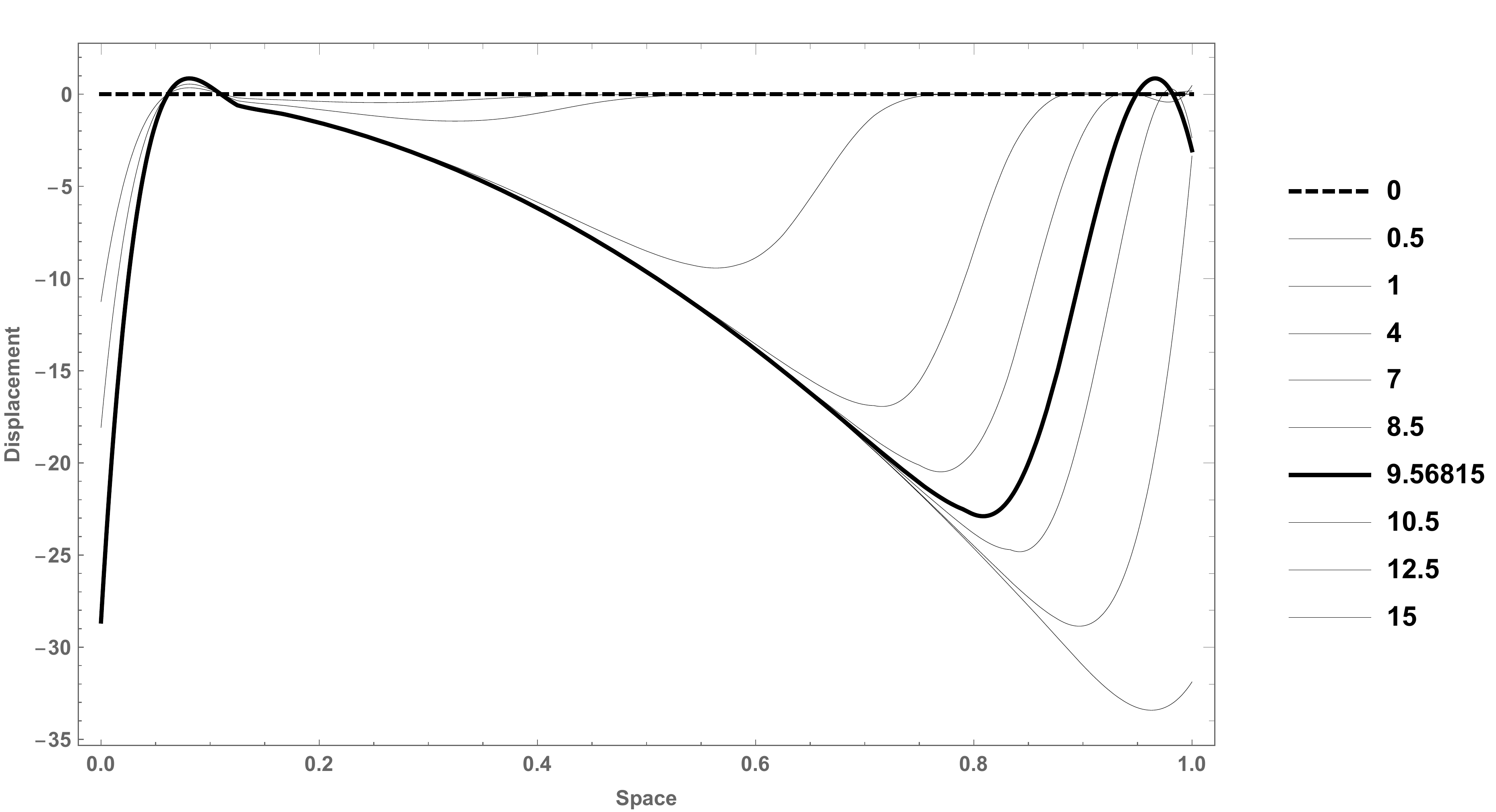}
\caption{Central induction-transduction activation wave propagation in time for density kernel mass $0.15$ and induction activation threshold $k=50$ for $K\sim\text{Poisson}(c=50)$}\label{fig:nn}
\end{figure}

\begin{figure}[h!]
\centering
\includegraphics[width=6.5in]{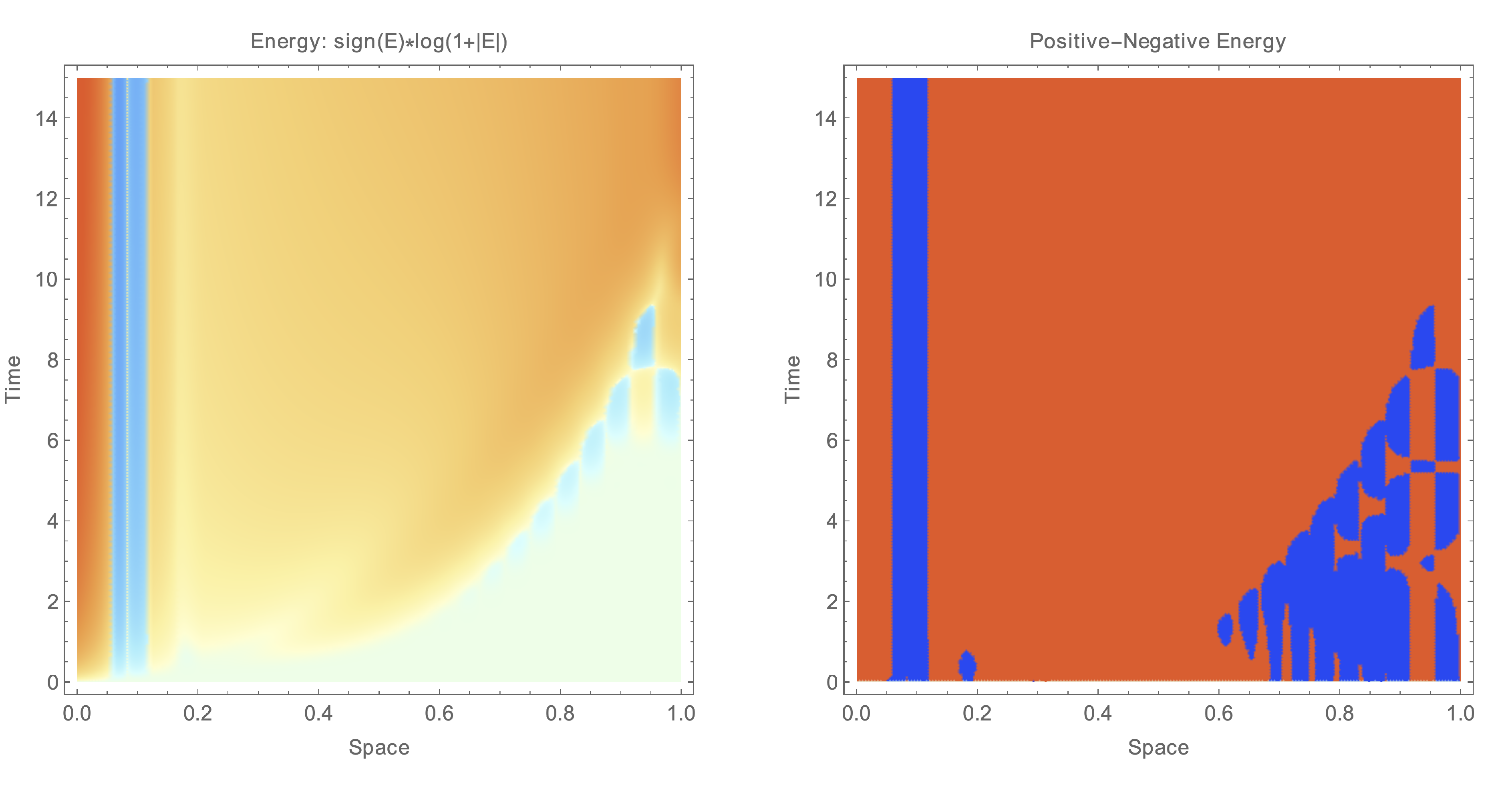}
\caption{Central induction-transduction activation energy in time for density kernel mass $0.15$ and induction activation threshold $k=50$ for $K\sim\text{Poisson}(c=50)$}\label{fig:nn}
\end{figure}

\clearpage
\newpage
\section{Subcentral induction-transduction activation-deactivation} For these figures, we take $E=[0,1]$, $\nu=\Leb$, and the density kernel as $f(x,y)=\exp_-a(x+y)$ on $E\times E$, where $a$ is chosen such that \[(\nu\times\nu)f = \int_{[0,1]^2}\D x\D y e^{-a(x+y)}= p\] and use the standard ITAD field equation, i.e. $C_1=0$. 

\begin{figure}[h!]
\centering
\begingroup
\captionsetup[subfigure]{width=5in,font=normalsize}
\subfloat[Activation with threshold $l\in\{1,2,3\}$\label{fig:r0}]{\includegraphics[width=3.5in]{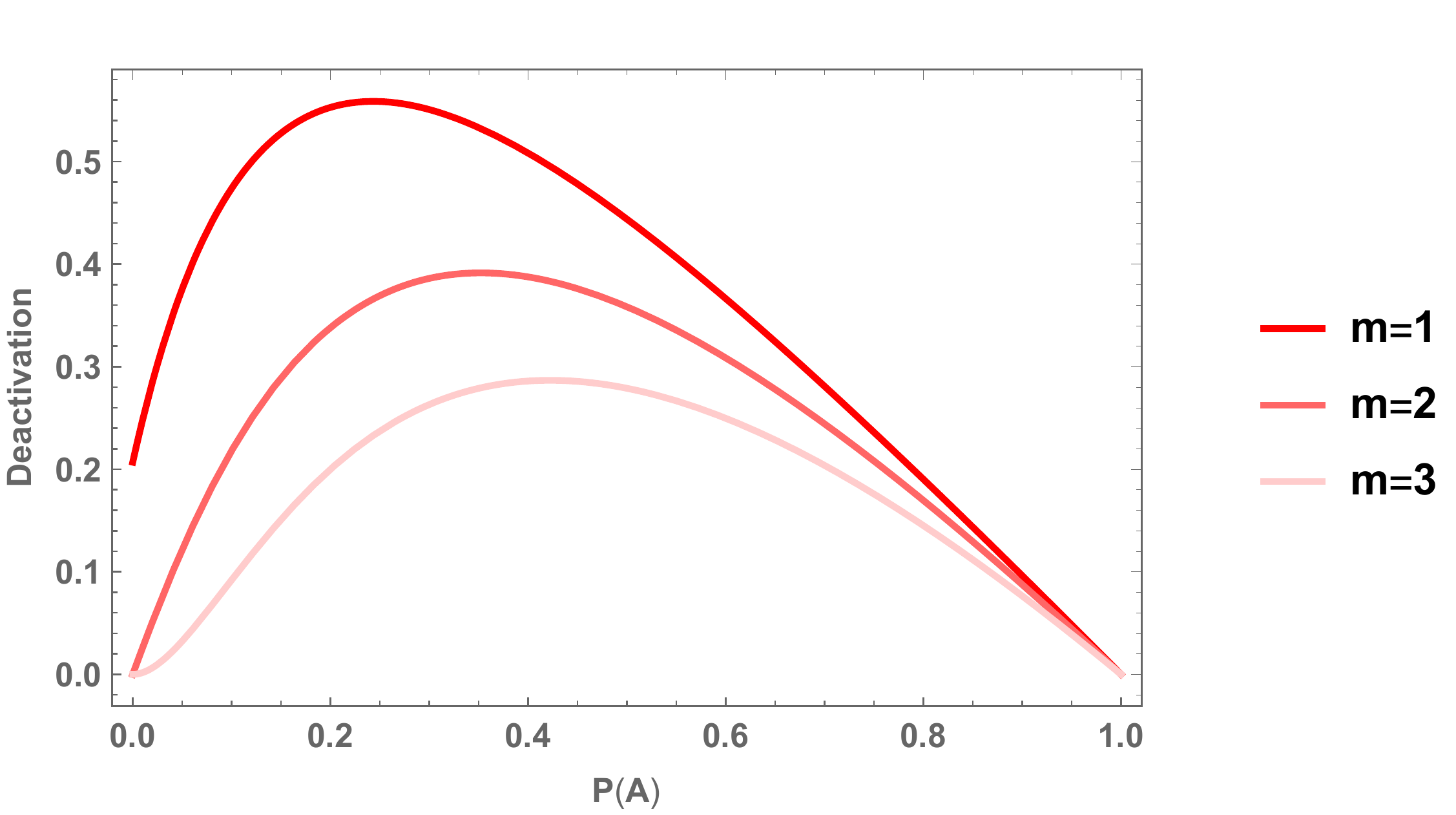}}
\subfloat[Deactivation with threshold $m\in\{1,2,3\}$\label{fig:r1}]{\includegraphics[width=3.5in]{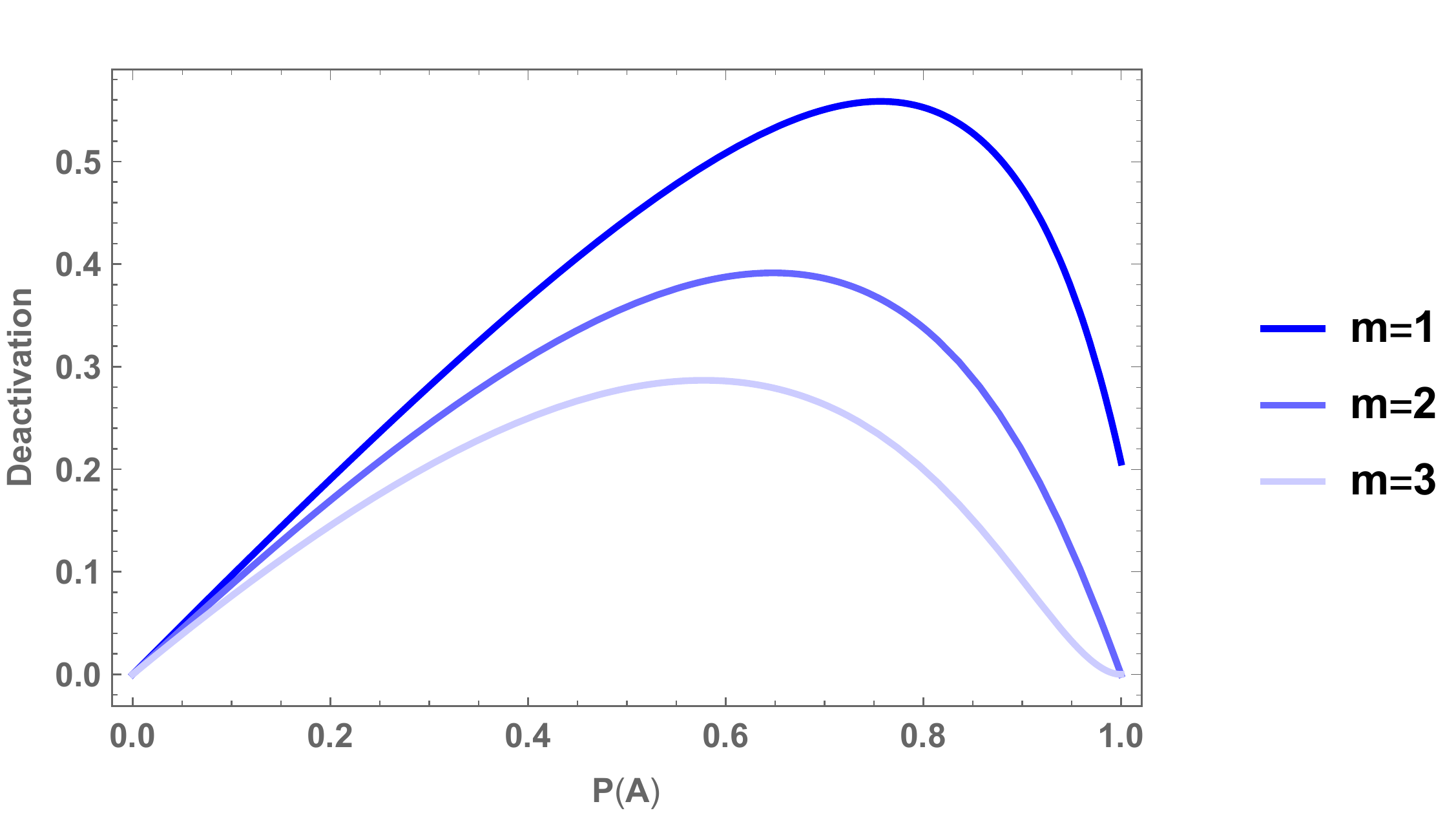}}\\
\subfloat[Net activation for $(l,m)\in\{1,2,3\}\times\{1,2,3\}$\label{fig:r2}]{\includegraphics[width=6.5in]{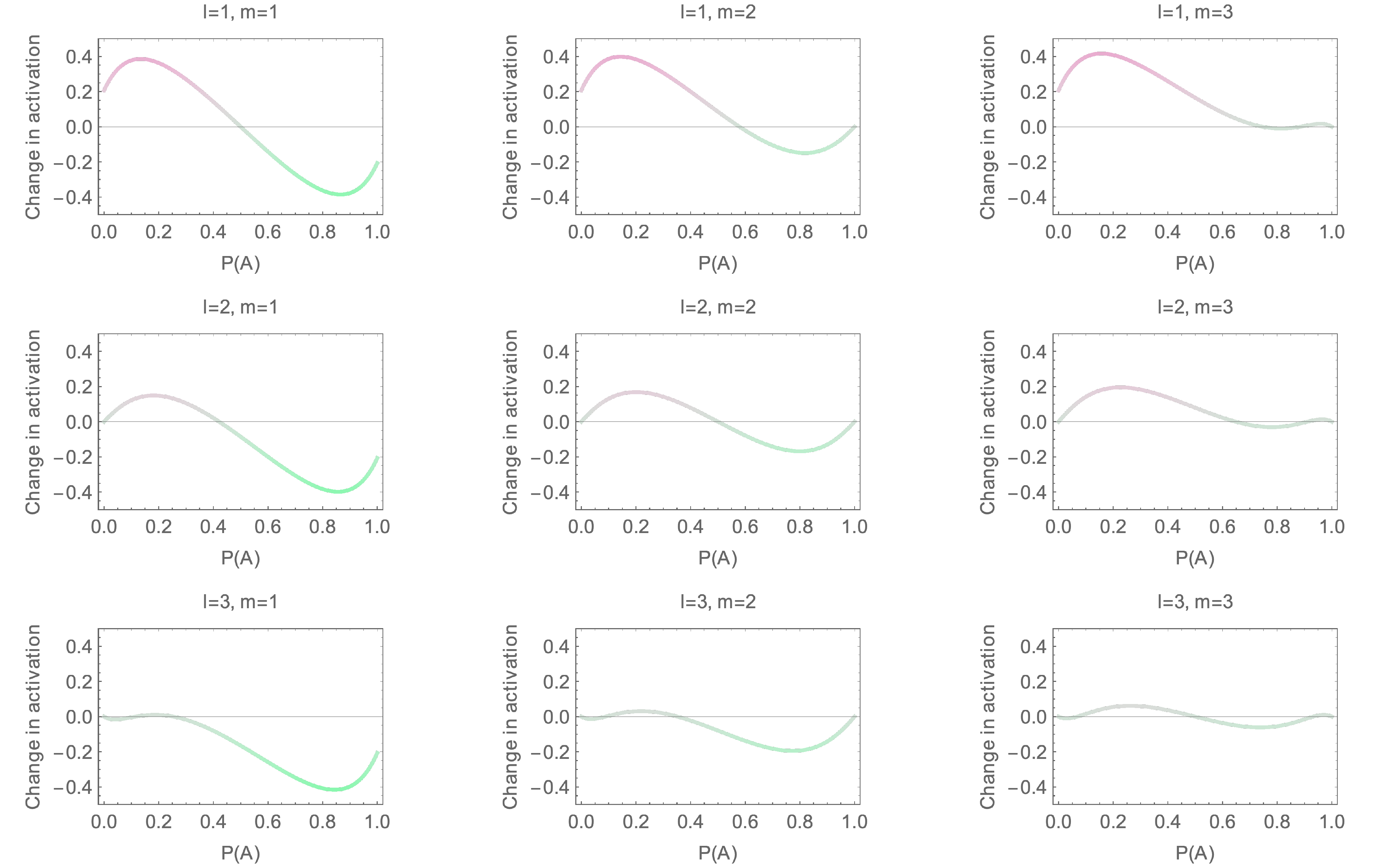}}
\endgroup
\caption{Subcentral transduction  in activation, deactivation, and net activation, in activation threshold $l$ and deactivation probability $r$ (blue is zero, red is one), for $K\sim\text{Poisson}(c=50)$ with density kernel mass $p=1/10$}\label{fig:orbitsdeact}
\end{figure}

\begin{figure}[h!]
\centering
\includegraphics[width=6.5in]{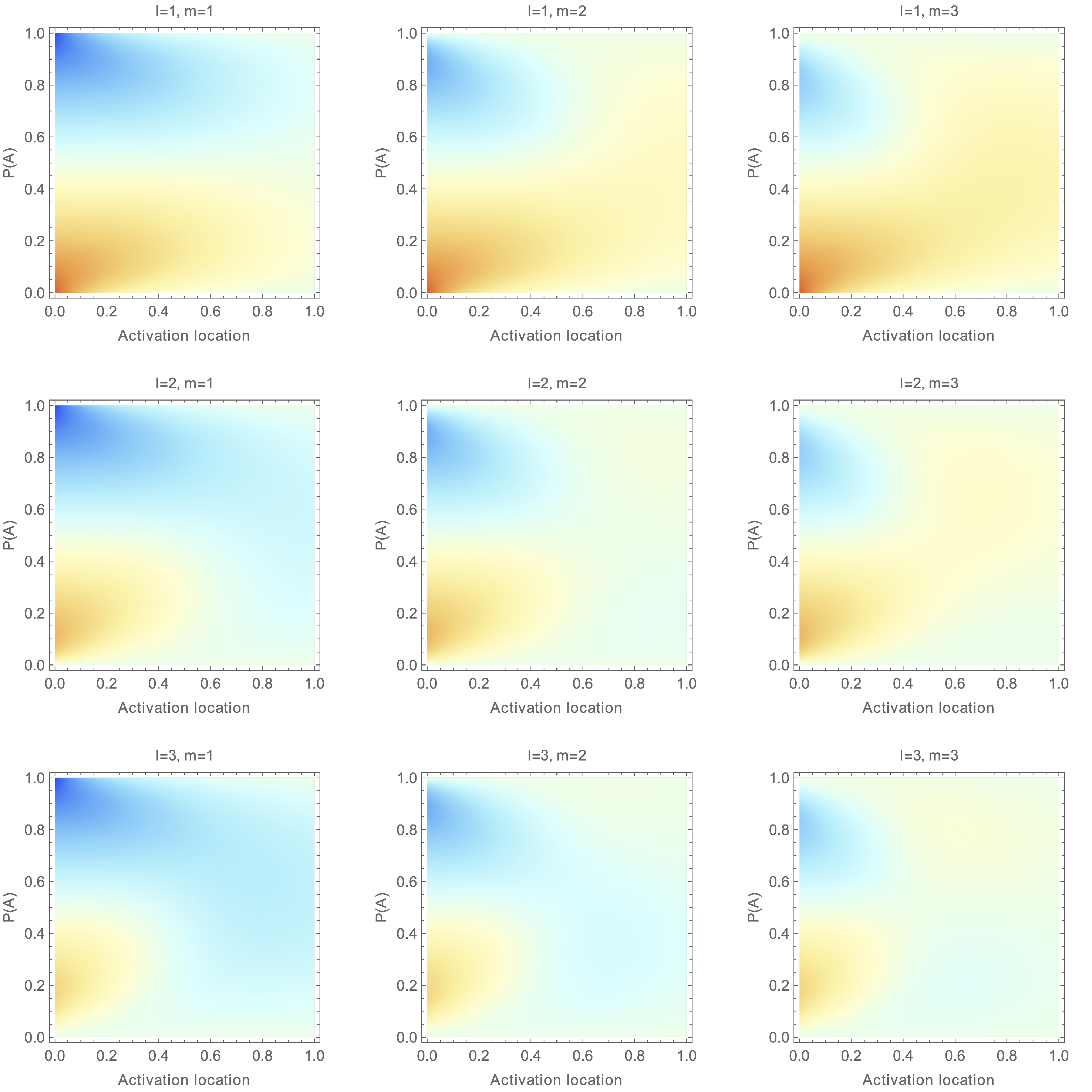}
\caption{Subcentral transduction  in activation $l$ and deactivation $m$ thresholds for number of points $K\sim\text{Poisson}(c=50)$, density kernel mass $p=1/10$, initializing activation threshold $k=8$, yielding $\xi_0\{A\}\simeq0.1438$}\label{fig:sifp}
\end{figure}
\FloatBarrier

\begin{figure}[h!]
\centering
\includegraphics[width=6.5in]{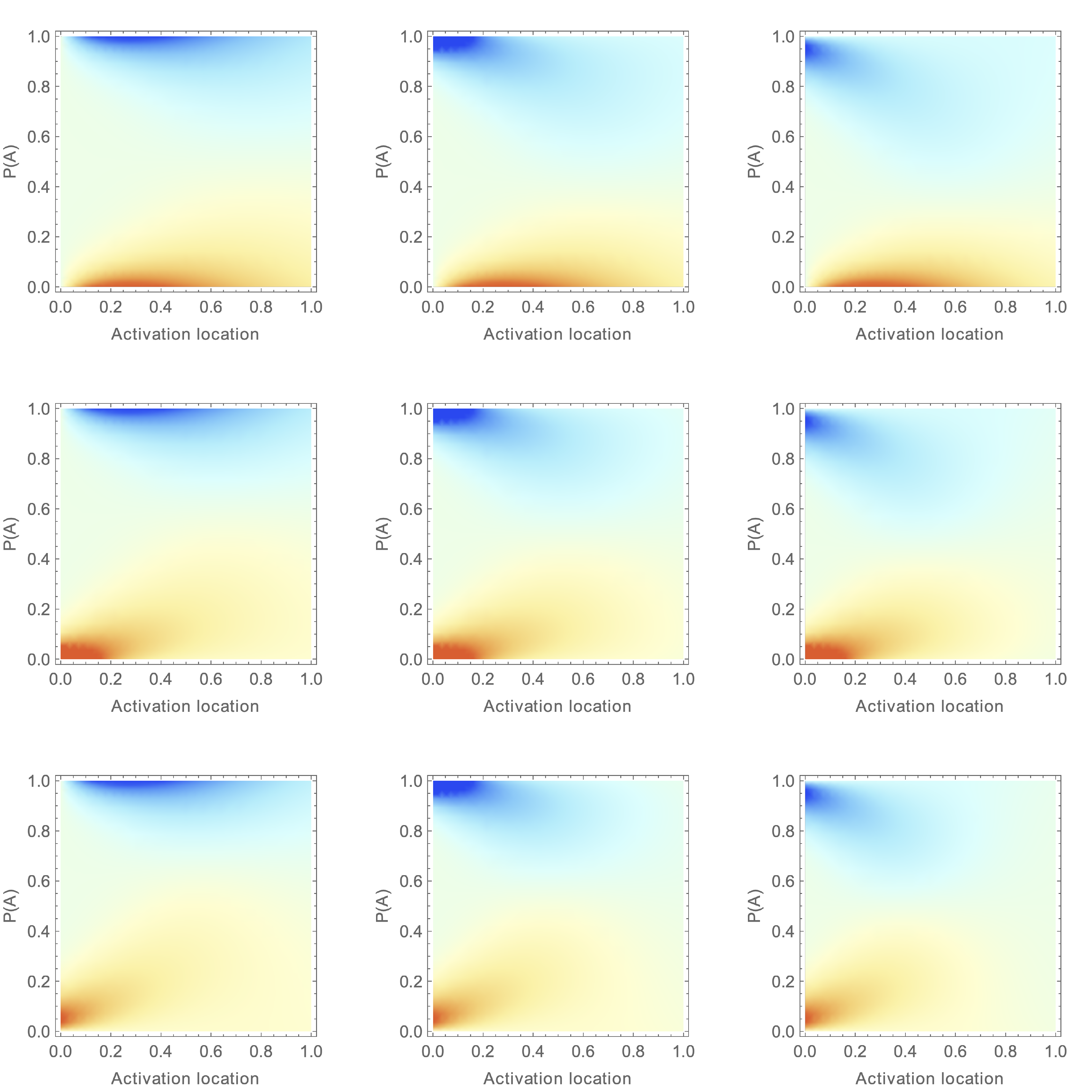}
\caption{Subcentral transduction activation flux relative to activation  in activation $l$ and deactivation $m$ thresholds for number of points $K\sim\text{Poisson}(c=50)$, density kernel mass $p=1/10$, initializing activation threshold $k=8$, yielding $\xi_0\{A\}\simeq0.1438$}\label{fig:sifp}
\end{figure}
\FloatBarrier

\begin{figure}[h!]
\centering
\includegraphics[width=6.5in]{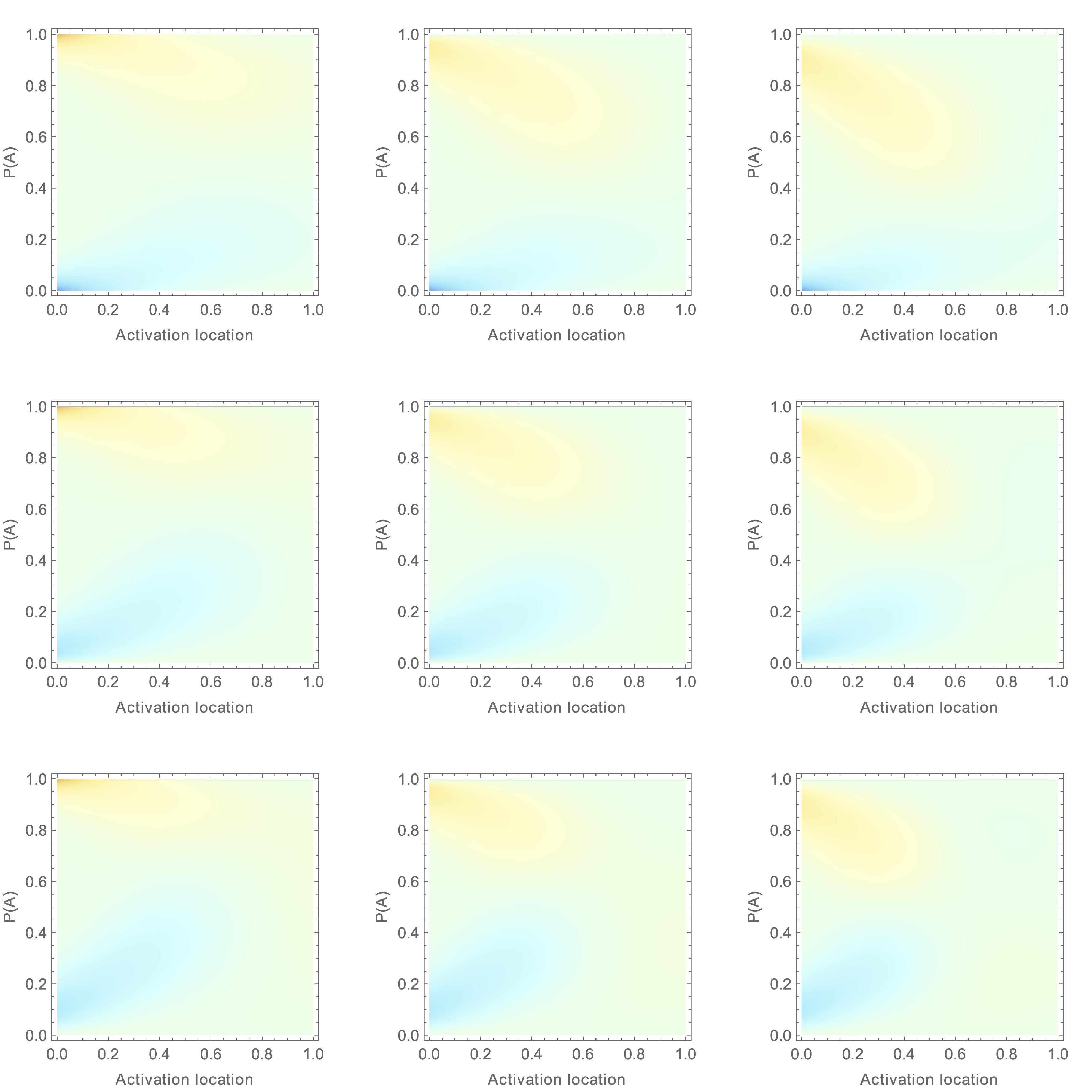}
\caption{Subcentral transduction activation flux relative to location  in activation $l$ and deactivation $m$ thresholds for number of points $K\sim\text{Poisson}(c=50)$, density kernel mass $p=1/10$, initializing activation threshold $k=8$, yielding $\xi_0\{A\}\simeq0.1438$}\label{fig:sifp}
\end{figure}
\FloatBarrier

\begin{figure}[h!]
\centering
\includegraphics[width=6.5in]{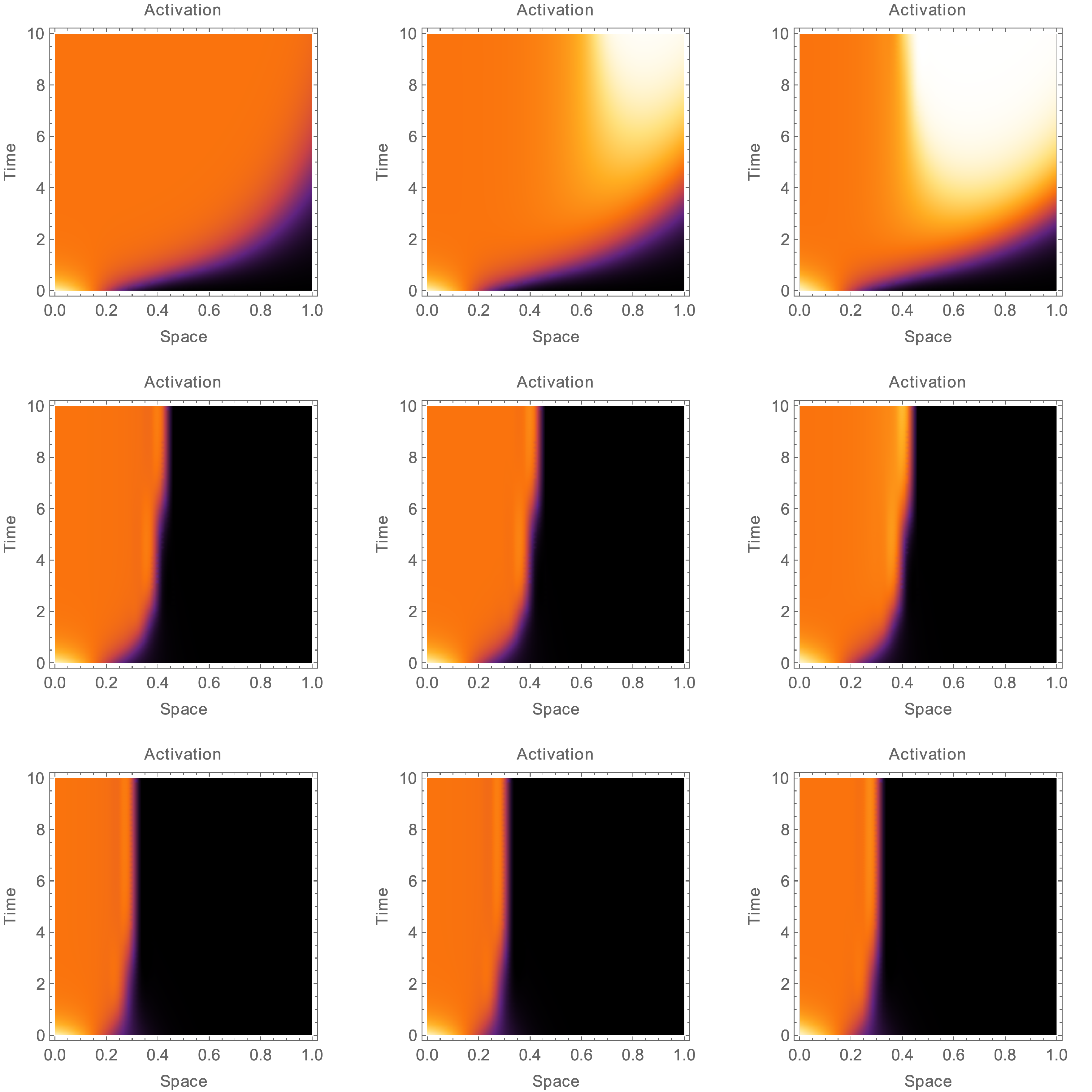}
\caption{Subcentral transduction activation-deactivation probability in location  in activation $l$ and deactivation $m$ thresholds for number of points $K\sim\text{Poisson}(c=50)$, initializing activation threshold $k=14$}\label{fig:sifp}
\end{figure}
\FloatBarrier

\begin{figure}[h!]
\centering
\includegraphics[width=6.5in]{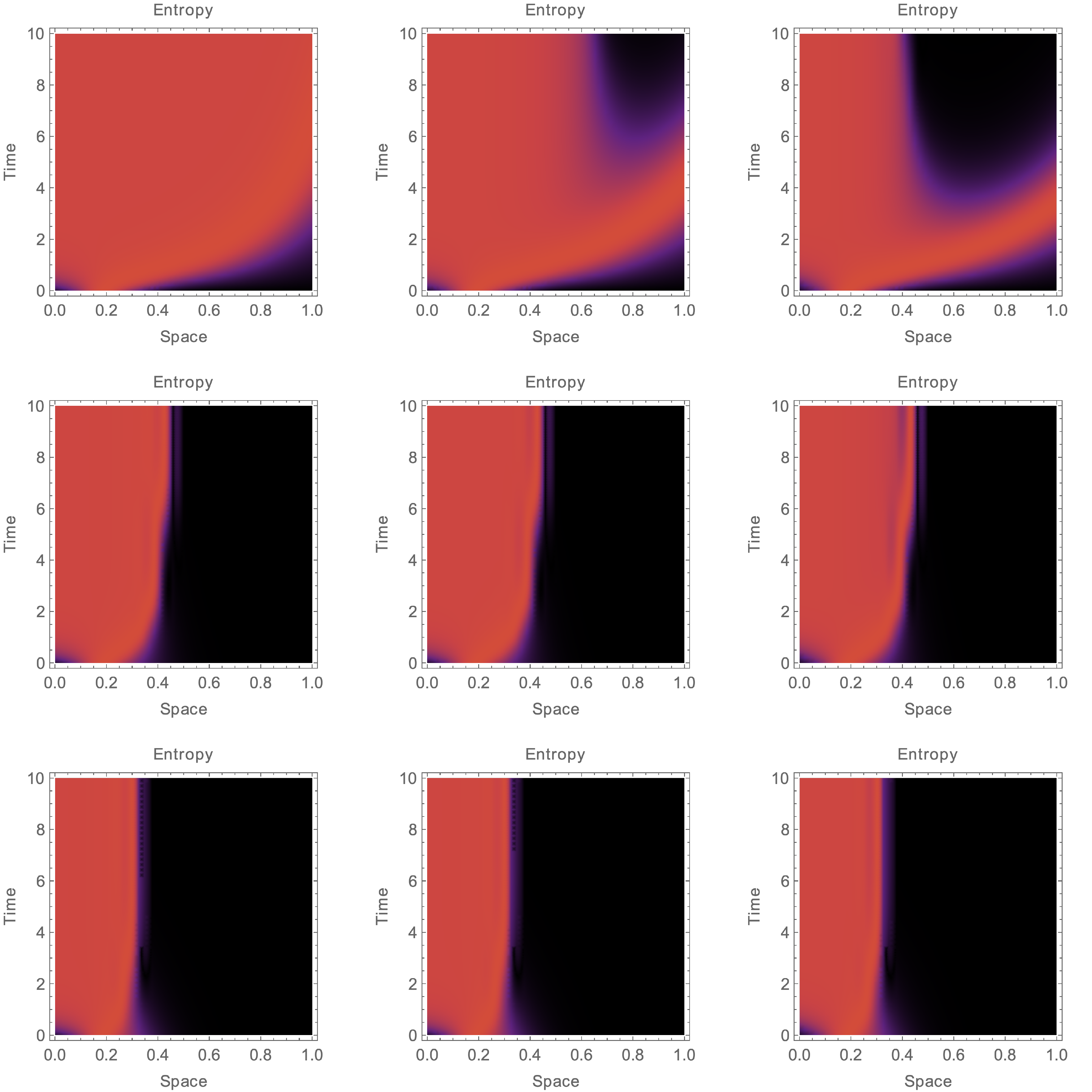}
\caption{Subcentral transduction activation-deactivation entropy in location  in activation $l$ and deactivation $m$ thresholds for number of points $K\sim\text{Poisson}(c=50)$, initializing activation threshold $k=14$}\label{fig:sifp}
\end{figure}
\FloatBarrier

\begin{figure}[h!]
\centering
\includegraphics[width=6.5in]{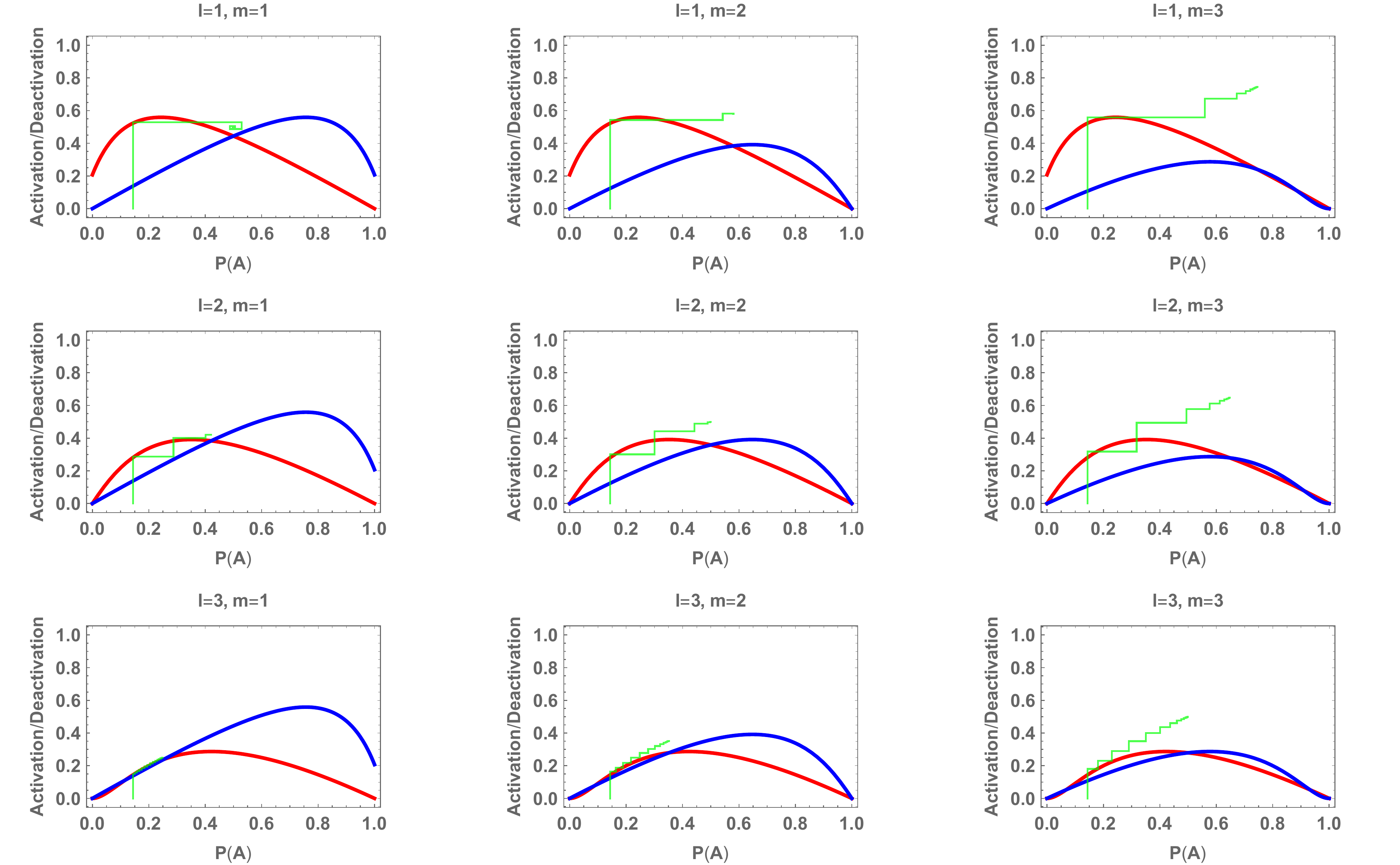}
\caption{Subcentral transduction  in activation $l$ and deactivation $m$ thresholds for number of points $K\sim\text{Poisson}(c=50)$, density kernel mass $p=1/10$,  initializing activation threshold $k=8$, yielding $\xi_0\{A\}\simeq0.1438$}\label{fig:sifp}
\end{figure}
\FloatBarrier

\clearpage
\newpage
\section*{Subcentral induction-transduction activation ($l=1,m=\infty$)} We use the standard ITAD field equation, i.e. $C_1=0$. 

\begin{figure}[h!]
\centering
\includegraphics[width=5in]{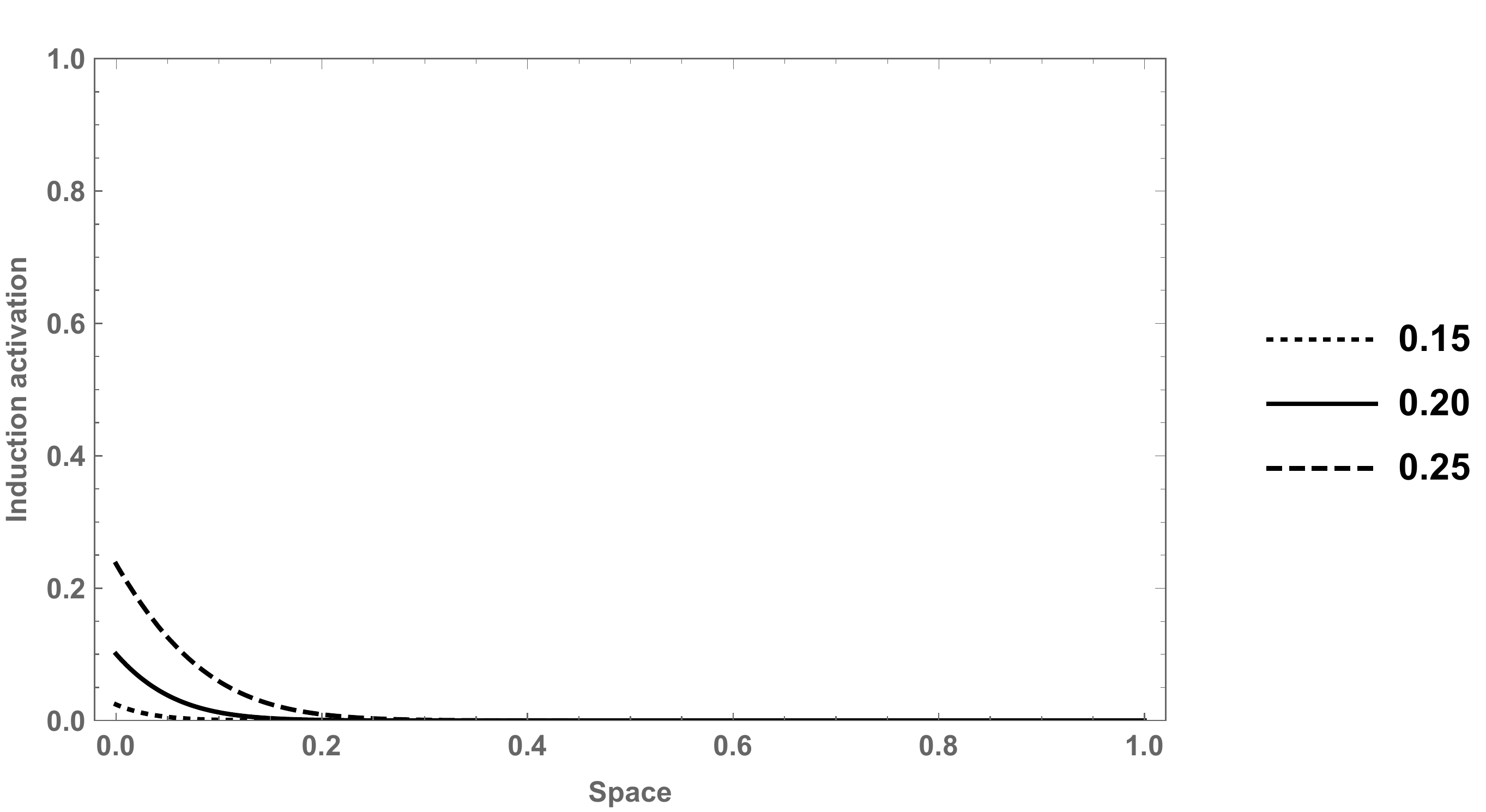}
\caption{Subcentral induction activation for density kernel masses $0.15,0.2,0.25$ and induction activation threshold $k=30$ for $K\sim\text{Poisson}(c=50)$}\label{fig:nn}
\end{figure}

\begin{figure}[h!]
\centering
\includegraphics[width=7in]{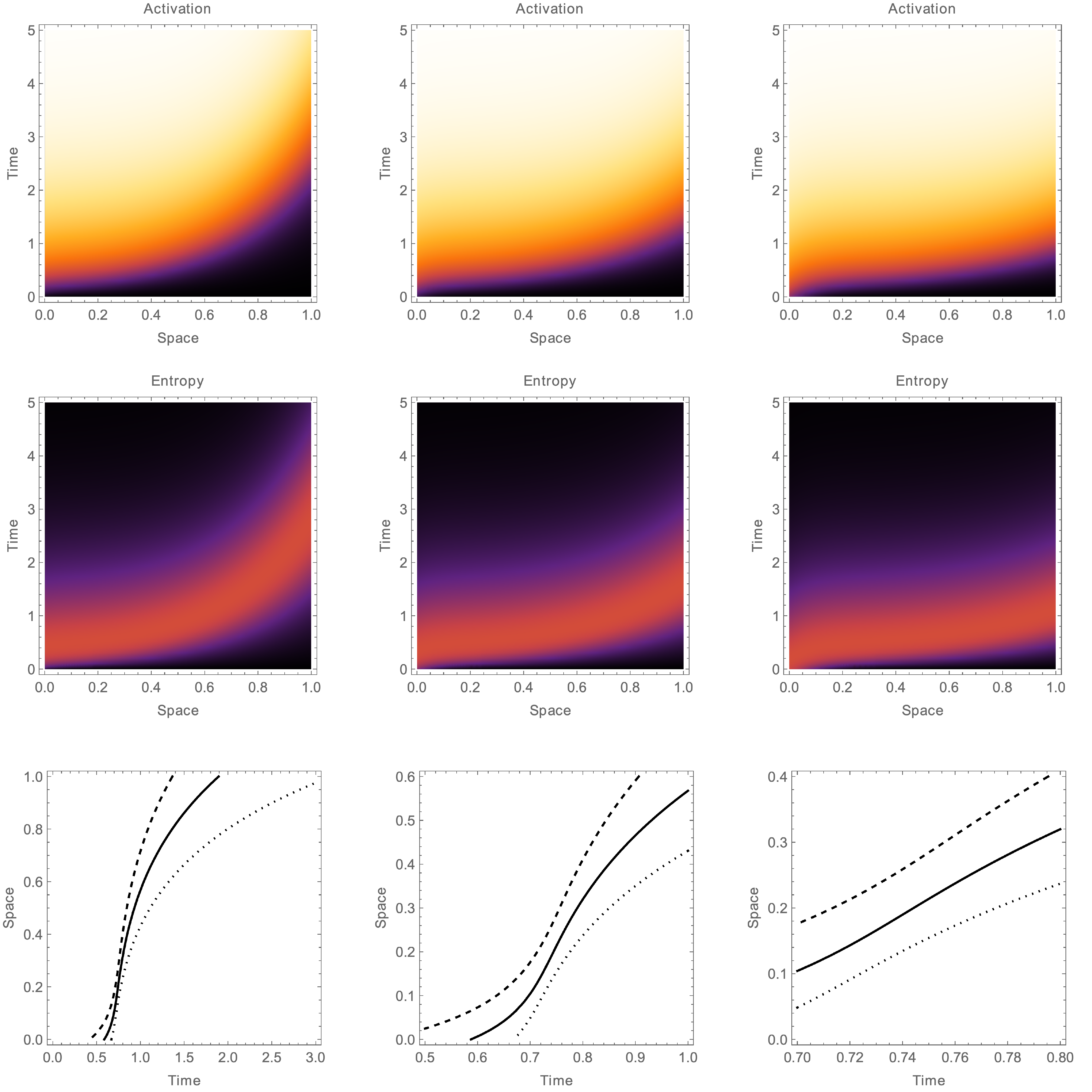}
\caption{Subcentral transduction activation for density kernel masses $0.15,0.2,0.25$ and induction activation threshold $k=30$ for $K\sim\text{Poisson}(c=50)$}\label{fig:nn}
\end{figure}

\clearpage
\newpage
\section{Decentral induction-transduction activation-deactivation} For this section, we take $E=[0,1]$, $\nu=\Leb$, and the (constant) density kernel as $f(x,y)=p$ on $E\times E$. Deactivation in this section is independent.

\subsection{Increasing transduction orbits of activation probability} 

We construct a discrete-time dynamical system on the label distribution. The initial condition is the fraction of points in activation as given by the activation model. The idea is that, at each iteration, the activation probability is incremented by the fraction of new activations.

Let $\xi_0\{A\}\equiv\xi\{D\}$. The recursive relation \[\xi_n\{D\} \equiv \xi_{n-1}\{D\} -\nu(p_{n-1}^B)\xi_{n-1}\{D\}\for n\ge 1\] where $p_{n}^B(x)\equiv\P_n(d_B(x,D)\ge l) + \P_n(d_B(x,D)=l-1)f(x,x)$ with $\P_n$ based on $\xi_n$ and $\P_0\equiv\P$ on $\xi$, describes the decay of deactivation probability by the normalized expected number of activations, a probability: $\xi_n\{D\}\rightarrow0$. We refer to the sequence $(\xi_n\{D\}: n\ge0)$ as the decaying orbit of deactivation probability. Similarly, we refer to the sequence $(\xi_n\{A\}: n\ge0)$  \[\xi_n\{A\} \equiv \xi_{n-1}\{A\} +\nu(p_{n-1}^B)(1-\xi_{n-1}\{A\})\for n\ge 1\] as the increasing orbit of activation probability, with $\xi_n\{A\}\rightarrow1$. Note that $p_{n}^B$ depends on $\xi_n$. The dynamics are driven by the change in the law of the deactivation degree function relative to $\xi_n$. This formulation describes the mean-field transduction of the activation. Below in Figure~\ref{fig:cm} we show the orbits of activation probability in transduction threshold $l\in\{1,2,3\}$ for Poisson $N$ with $c=50$, the constant density kernel with $p=1/10$, and initializing activation threshold $k=8$, giving $\xi_0\{A\}\simeq0.1438$. As $l$ increases, the dynamics slow. 

\begin{figure}[h!]
\centering
\includegraphics[width=4.5in]{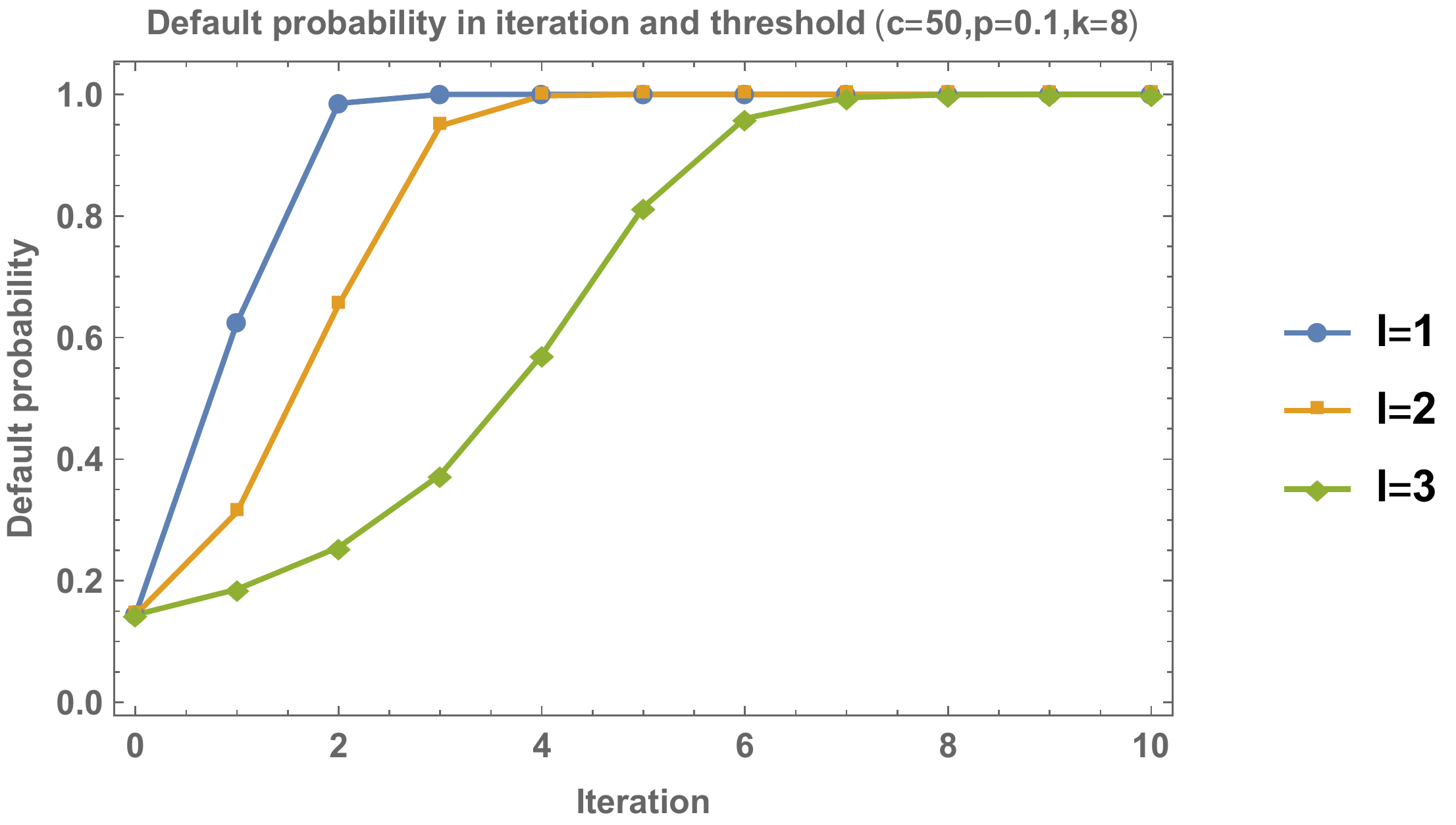}
\caption{Transduction $(\xi_n\{A\}:n=0,1,\dotsb,10)$ over 10 iterations for number of points $K\sim\text{Poisson}(c=50)$, constant density kernel with $p=1/10$, and $\xi_0\{A\}\simeq0.1438$, based on initializing activation threshold $k=8$}\label{fig:cm}
\end{figure}
\FloatBarrier

\subsection{Non-decreasing transduction orbits of activations}
A stochastic version is updating the labels of deactivations to activated based on the new activations at each iteration, i.e. each new activation increases/decreases the number of activations/deactivations by one. The corresponding degree sequences of the deactivations across iterations is non-decreasing, and the cumulative number of new activations across iterations is non-decreasing. Hence in the stochastic version, the total number of new activations across all iterations may be less than the total number of deactivations. 

\subsection{Point system example} We simulate the transduction process of the system relative to Poisson $N$ with $K\sim\text{Poisson}(c=50)$, the constant density kernel with $p=1/10$, initializing activation threshold $k=8$, and transduction activation threshold $l=2$. Below in Figure~\ref{fig:contag} we show the first five iterations of the transduction process, with outcome $K=50$, which gives the sequence of the number of activations as $(7,20,37,48,49,50)$ as determined by the graph. 

\begin{figure}[h!]
\centering
\begingroup
\captionsetup[subfigure]{width=5in,font=normalsize}
\subfloat[Iteration 0 (7 activations)\label{fig:i0}]{\includegraphics[width=2in]{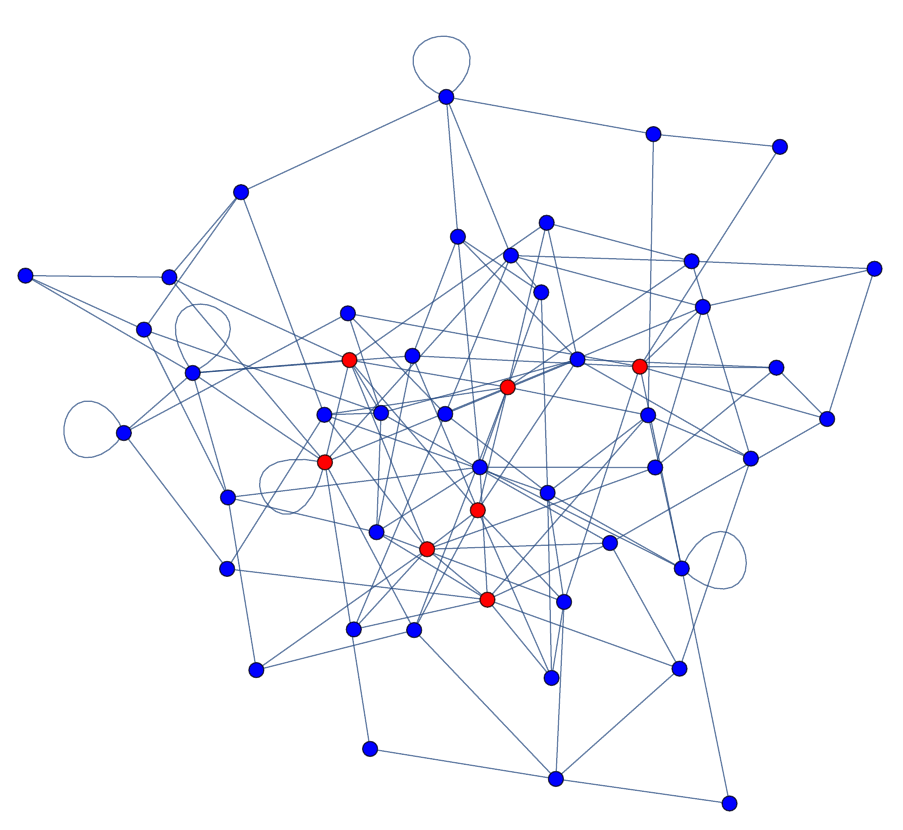}}
\subfloat[Iteration 1 (+13 activations)\label{fig:i0}]{\includegraphics[width=2in]{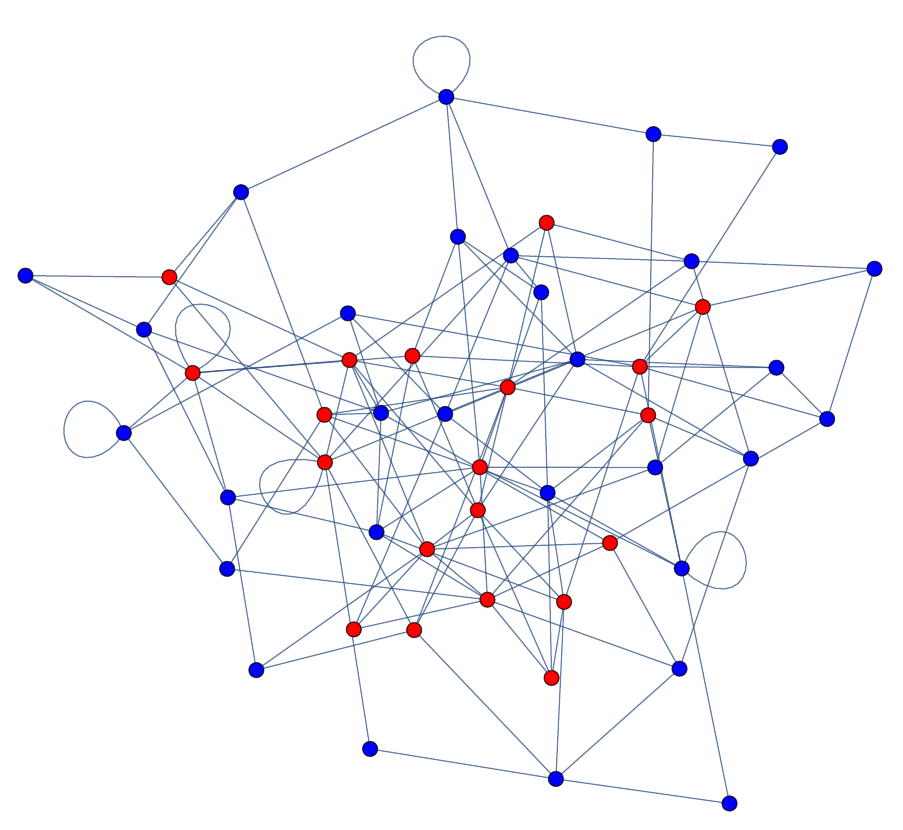}}
\subfloat[Iteration 2 (+17 activations)\label{fig:i0}]{\includegraphics[width=2in]{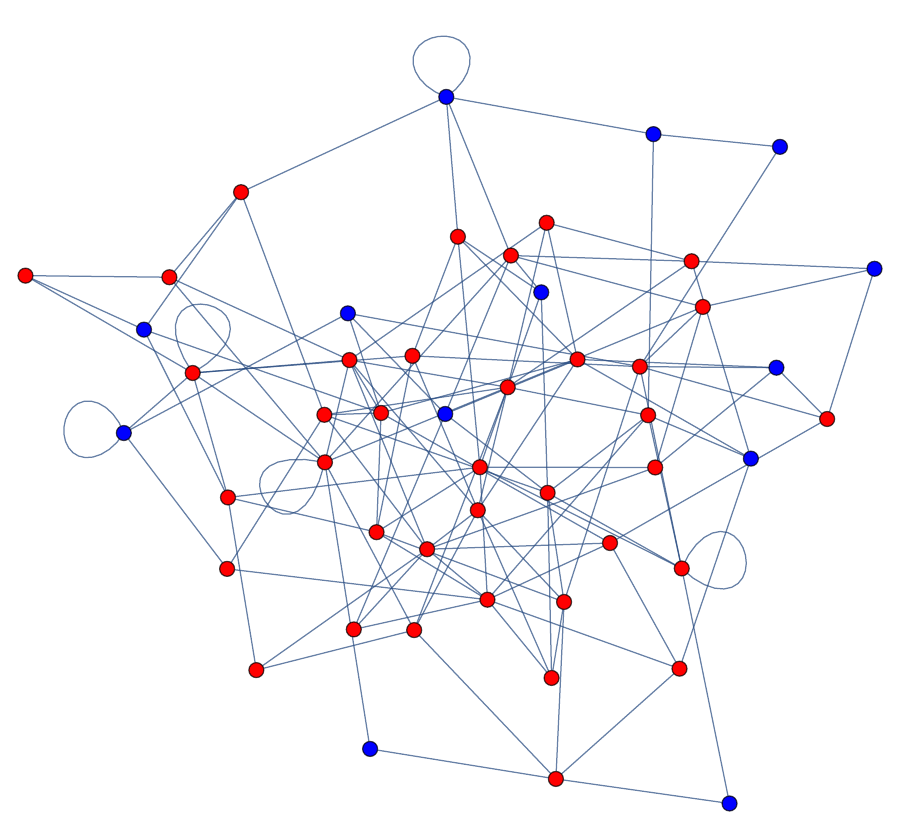}}\\
\subfloat[Iteration 3 (+11 activations)\label{fig:i0}]{\includegraphics[width=2in]{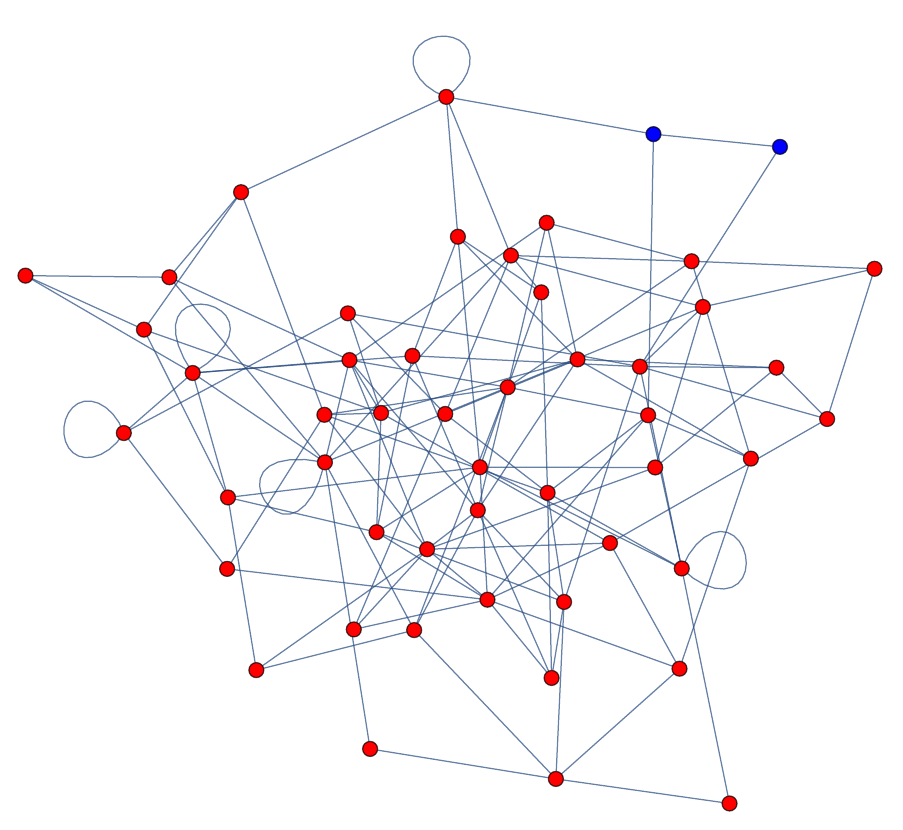}}
\subfloat[Iteration 4 (+1 activation)\label{fig:i0}]{\includegraphics[width=2in]{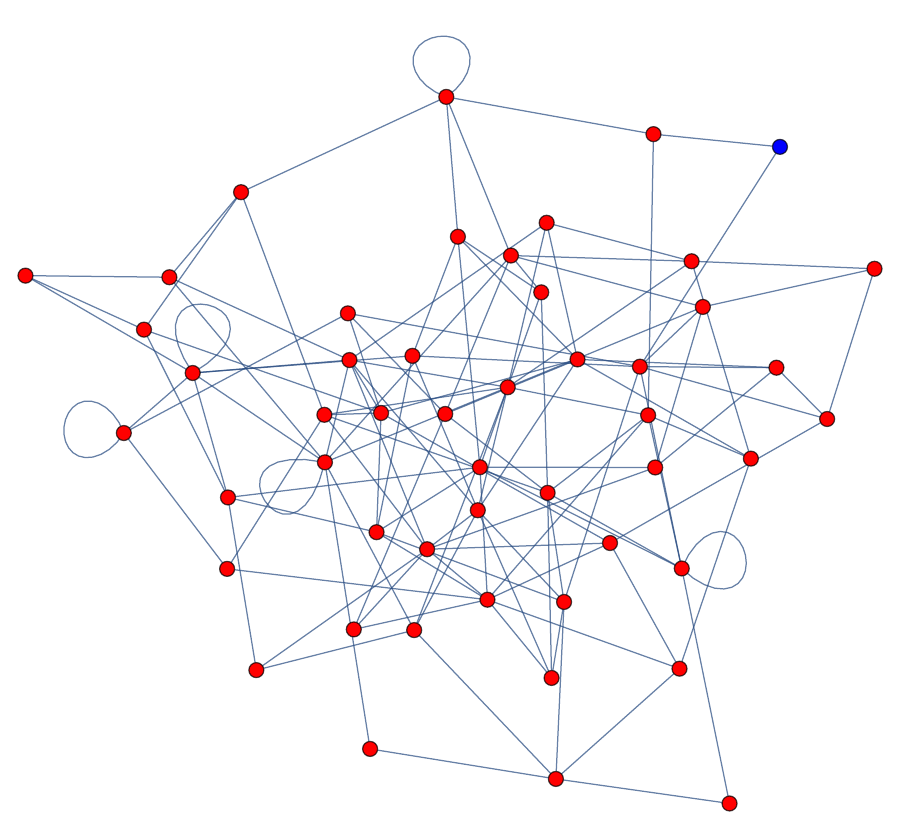}}
\subfloat[Iteration 5 (+1 activation)\label{fig:i0}]{\includegraphics[width=2in]{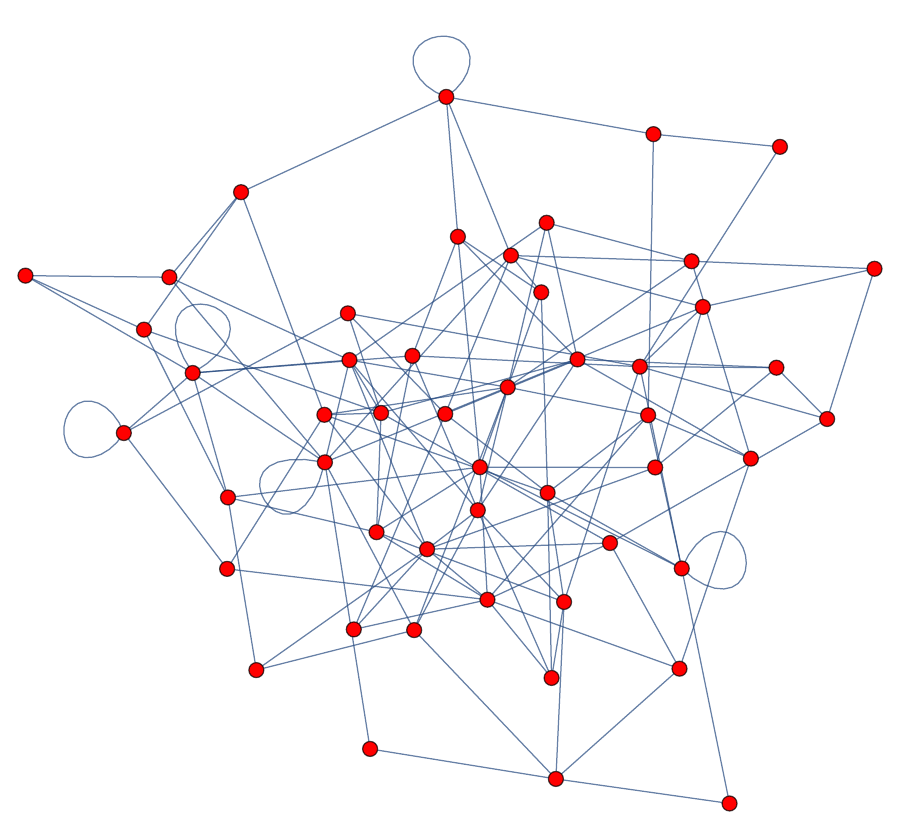}}\\
\endgroup
\caption{Simulation of activation transduction over random graph (blue is deactivation, red is activation) for Poisson $N$ with realization $K=50$, constant density kernel with $p=1/10$, initializing activation threshold $k=8$, and transduction activation threshold $l=2$, yielding simulation activation count sequence $(7,20,37,48,49,50)$}\label{fig:contag}
\end{figure}

\FloatBarrier

\subsection{Number of deactivations ($A\rightarrow D$ deactivating vertices)} We define the number of deactivations as \[R = \sum_i^K\ind{E}(X_i)\Lambda(Y_i)\] where $\Lambda(y)\sim\text{Bernoulli}(r\ind{}(y=A))$ with deactivation probability $r$. The mean is given by \[\E R = cr\int_F\xi(\D y)\ind{}(y=A)=cr\xi\{A\}\]

\subsection{Decentral transduction under deactivation and hysteresis bifurcations} 

\begin{figure}[h!]
\centering
\includegraphics[width=5in]{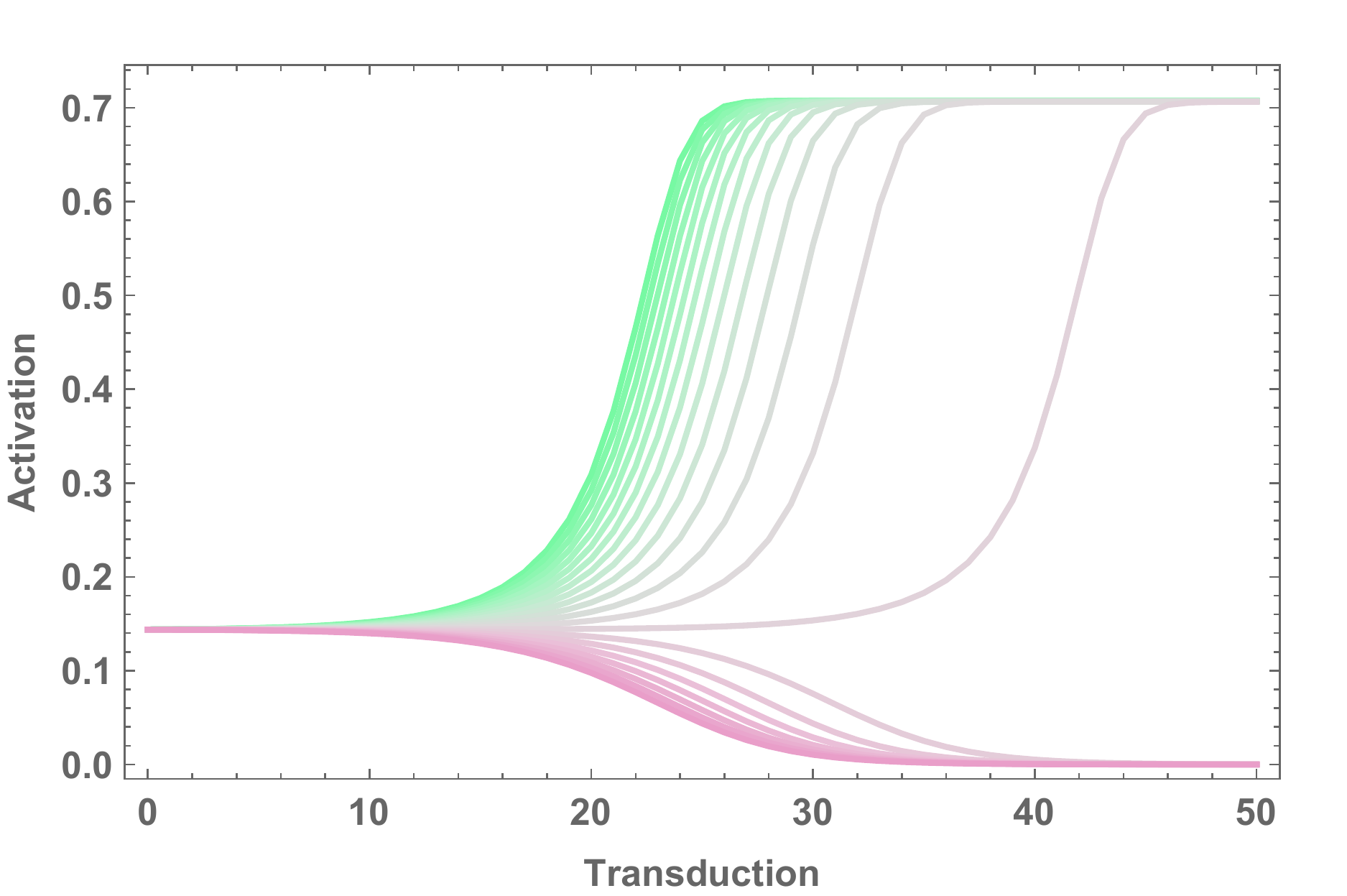}
\caption{Bifurcation and latency of transduction $(\xi_n\{A\}:n=0,1,\dotsb,50)$ over 50 iterations for number of points $K\sim\text{Poisson}(c=50)$, constant density kernel with $p=1/10$, and $\xi_0\{A\}\simeq0.1438$, based on initializing activation threshold $k=8$, transduction threshold $l=3$, and varying deactivation probability $r$ from $0.291$ (blue) to $0.293$ (red)}\label{fig:bifur}
\end{figure}

\begin{figure}[h!]
\centering
\includegraphics[width=5in]{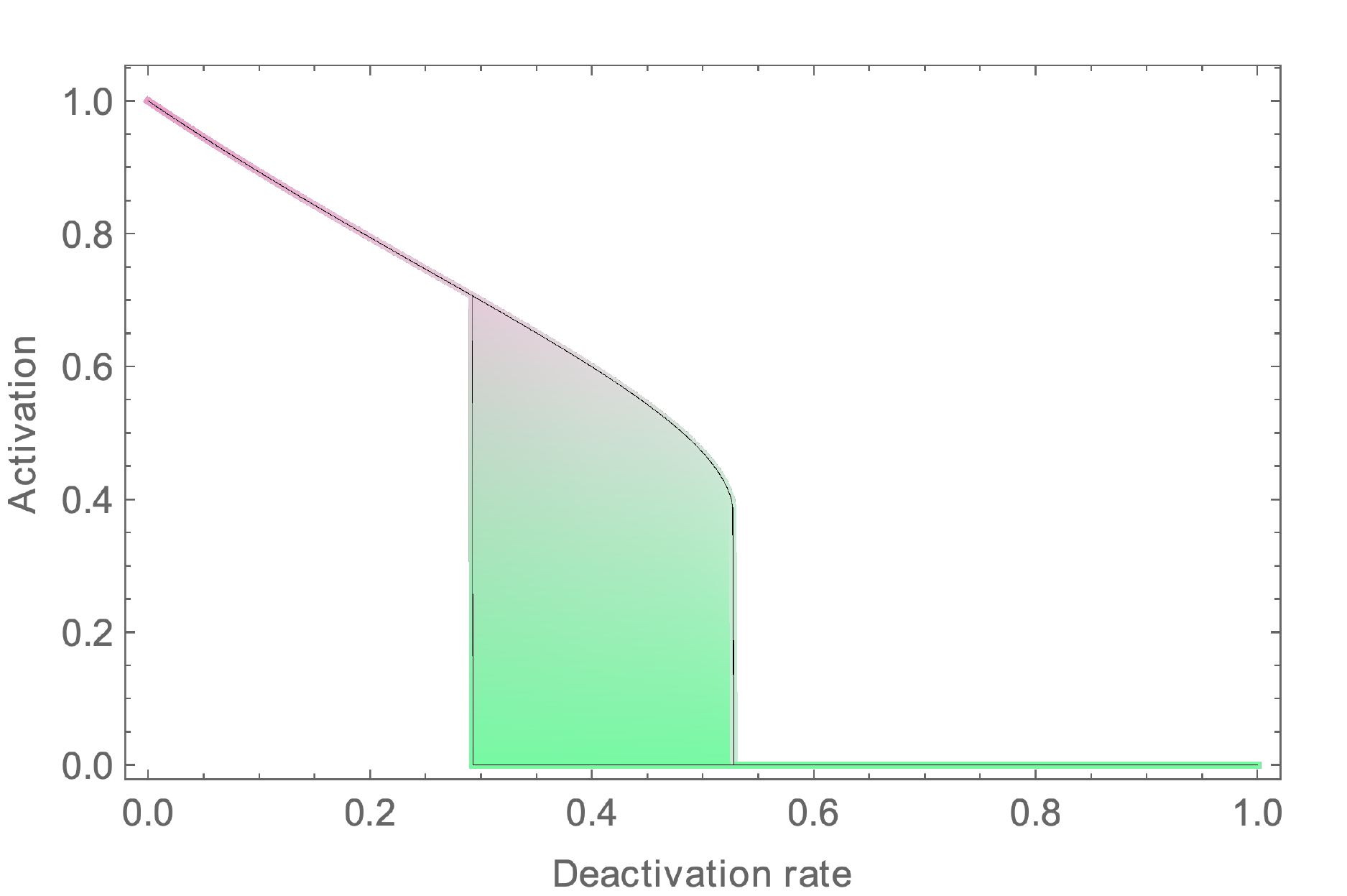}
\caption{Hysteresis of fixed-points of transduction in deactivation probability $r$ for number of points $K\sim\text{Poisson}(c=50)$, constant density kernel with $p=1/10$, and $\xi_0\{A\}\simeq0.1438$, based on initializing activation threshold $k=8$ and transduction threshold $l=3$}\label{fig:fp}
\end{figure}

We add deactivation to the recursion \[\xi_n\{D\} \equiv \xi_{n-1}\{D\} -\nu(p_{n-1}^B)\xi_{n-1}\{D\} + r(1-\xi_{n-1}\{D\})\for n\ge 1\] with $p_{n}^B(x)$ as before, with activation recursion \[\xi_n\{A\}\equiv \xi_{n-1}\{A\} + \nu(p_{n-1}^B)(1-\xi_{n-1}\{A\}) - r\xi_{n-1}\{A\}\for n\ge 1\] Recall $p_n^B$ is non-linear in $\xi_n$. We plot the forward-orbits of activation probability $(\xi_n\{A\}: n)$ with deactivation in Figure~\ref{fig:orbitsdeact} for varying $l$. The dependence on final state in $(l,r)$ is marked. Firstly, even with complete deactivation $r=1$ each iteration, positive fixed-points of activation exist for $l=1,2$. This illustrates that the probability of point activation is dominated by the transduction activation threshold $l$. For the highest activation threshold $l=3$, the activation probability converges to zero. Moreover, the activation probability for $l=3$ exhibits phase transition, in particular a bifurcation, of a jump-discontinuity in deactivation probability $r$, depicted below in the forward orbits (Figure~\ref{fig:bifur}). Moreover, the bifurcation is of the hysteresis type, shown in the fixed-points (Figure~\ref{fig:fp}).

\begin{figure}[h!]
\centering
\begingroup
\captionsetup[subfigure]{width=5in,font=normalsize}
\subfloat[Activation for activation threshold $l\in\{1,2,3\}$\label{fig:r0}]{\includegraphics[width=3.5in]{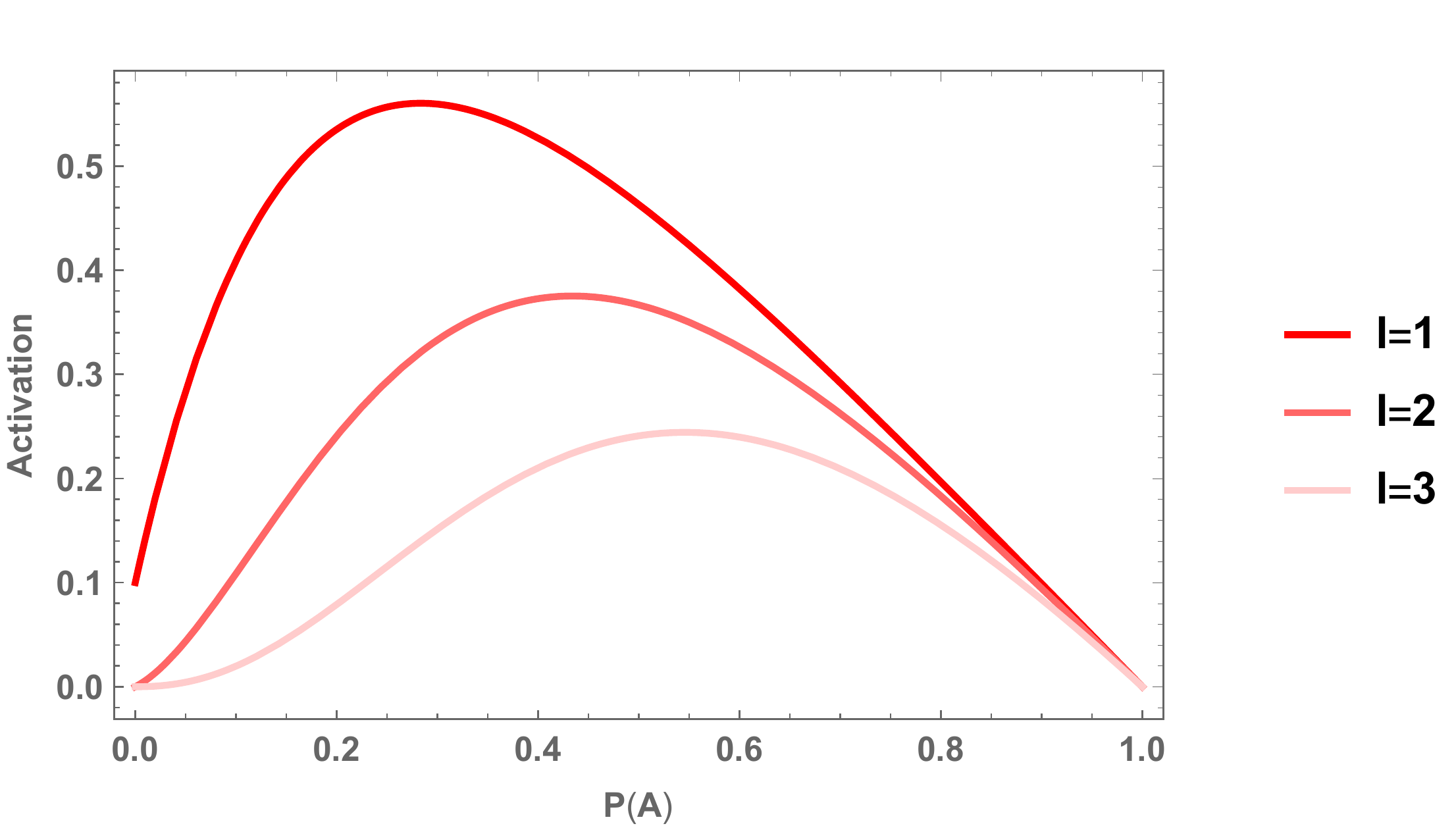}}
\subfloat[Deactivation \label{fig:r1}]{\includegraphics[width=3.5in]{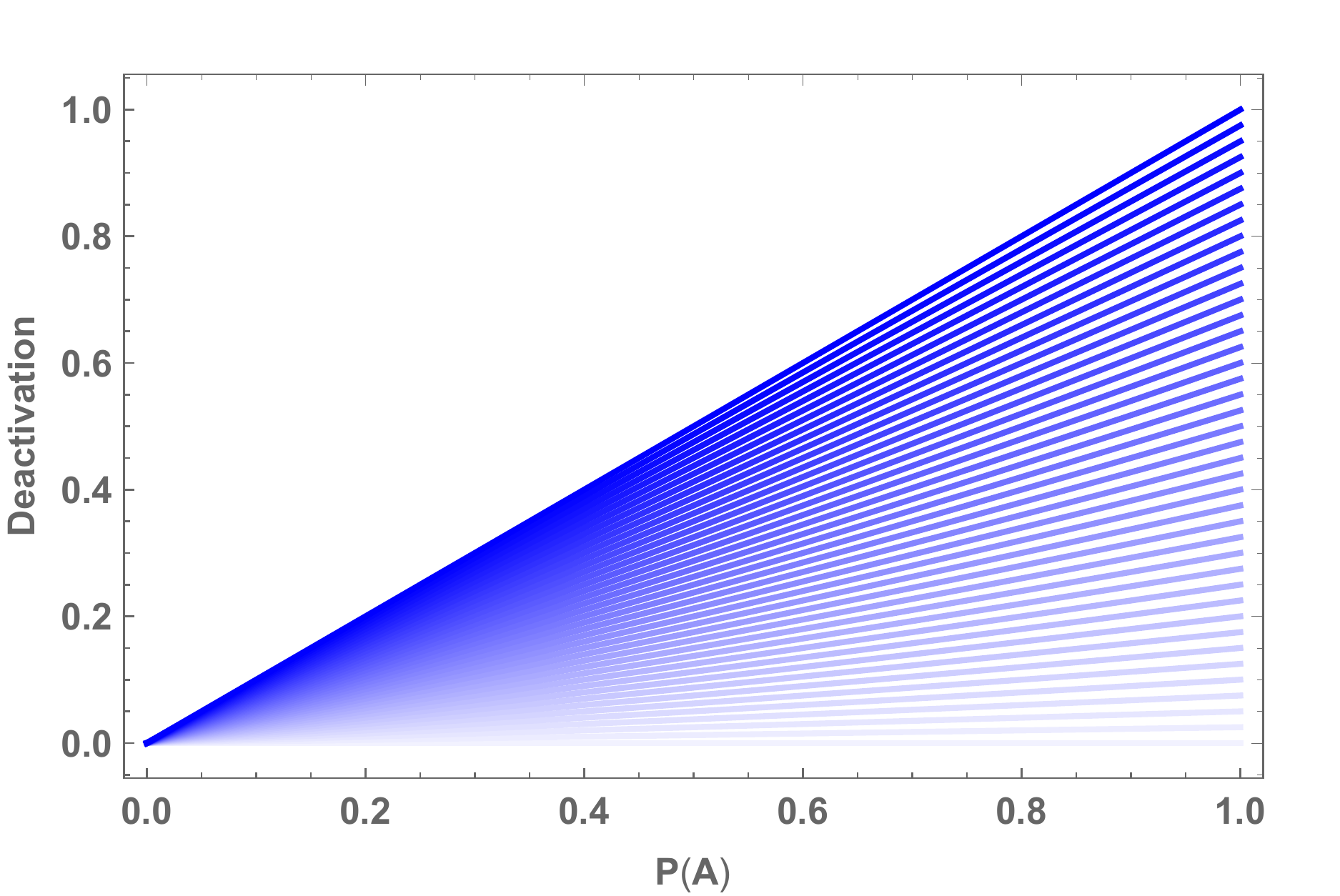}}\\
\subfloat[Net activation for activation threshold $l=3$\label{fig:r2}]{\includegraphics[width=3.5in]{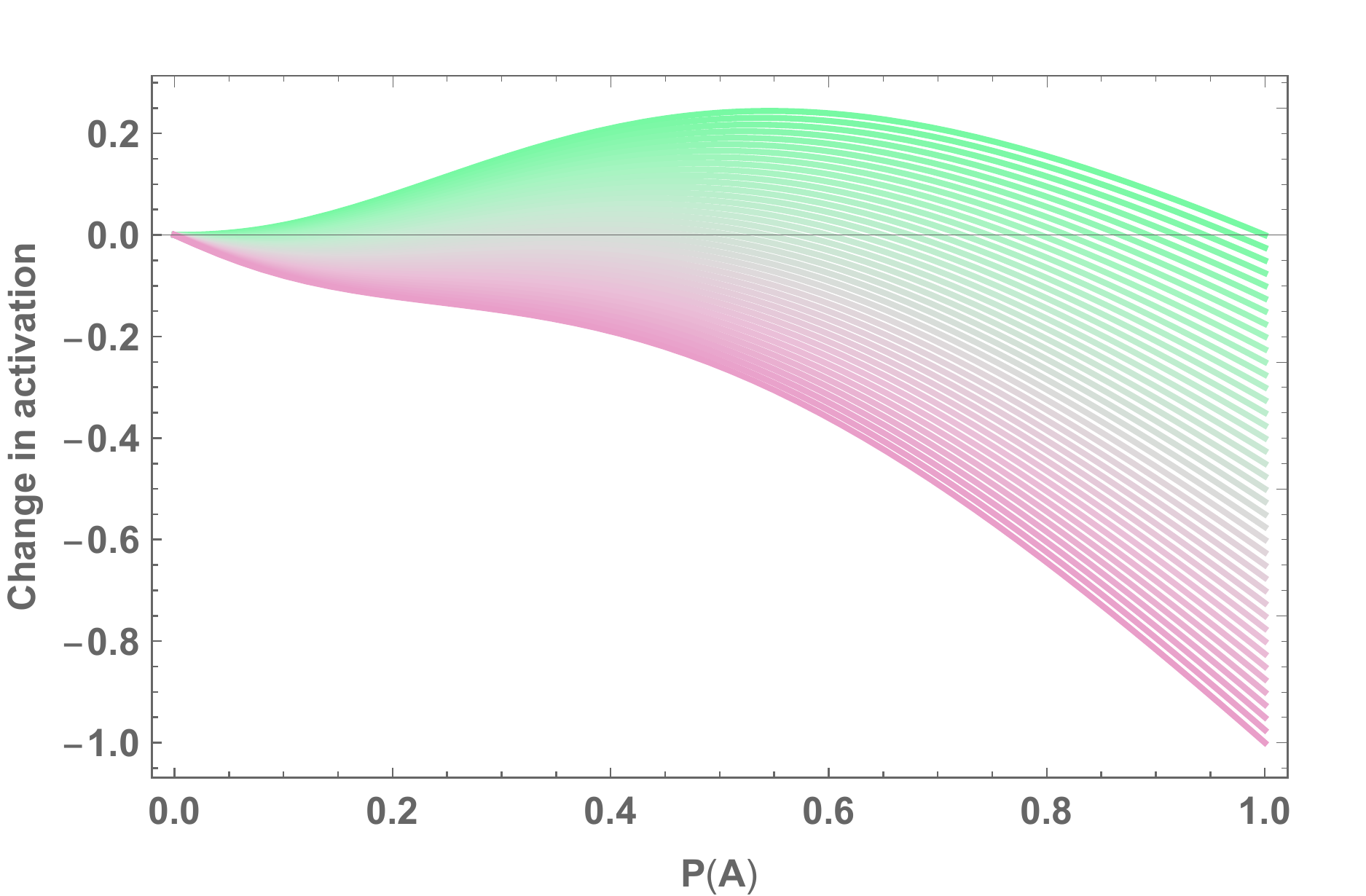}}
\subfloat[Net activation for activation threshold $l=3$ (zoomed)\label{fig:r2}]{\includegraphics[width=3.5in]{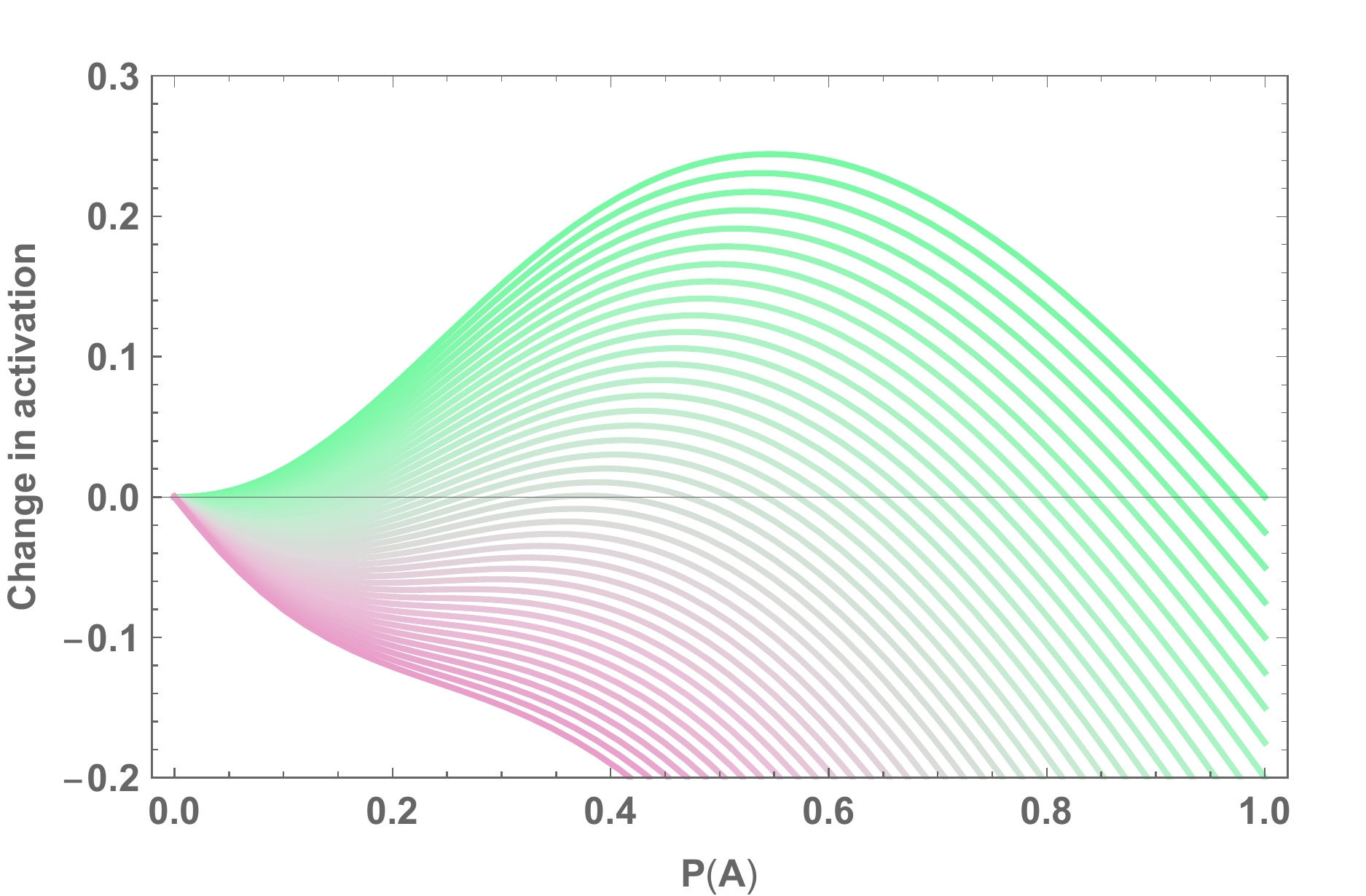}}\\
\subfloat[Net activation for activation threshold $l=3$ (zoomed)\label{fig:r2}]{\includegraphics[width=3.5in]{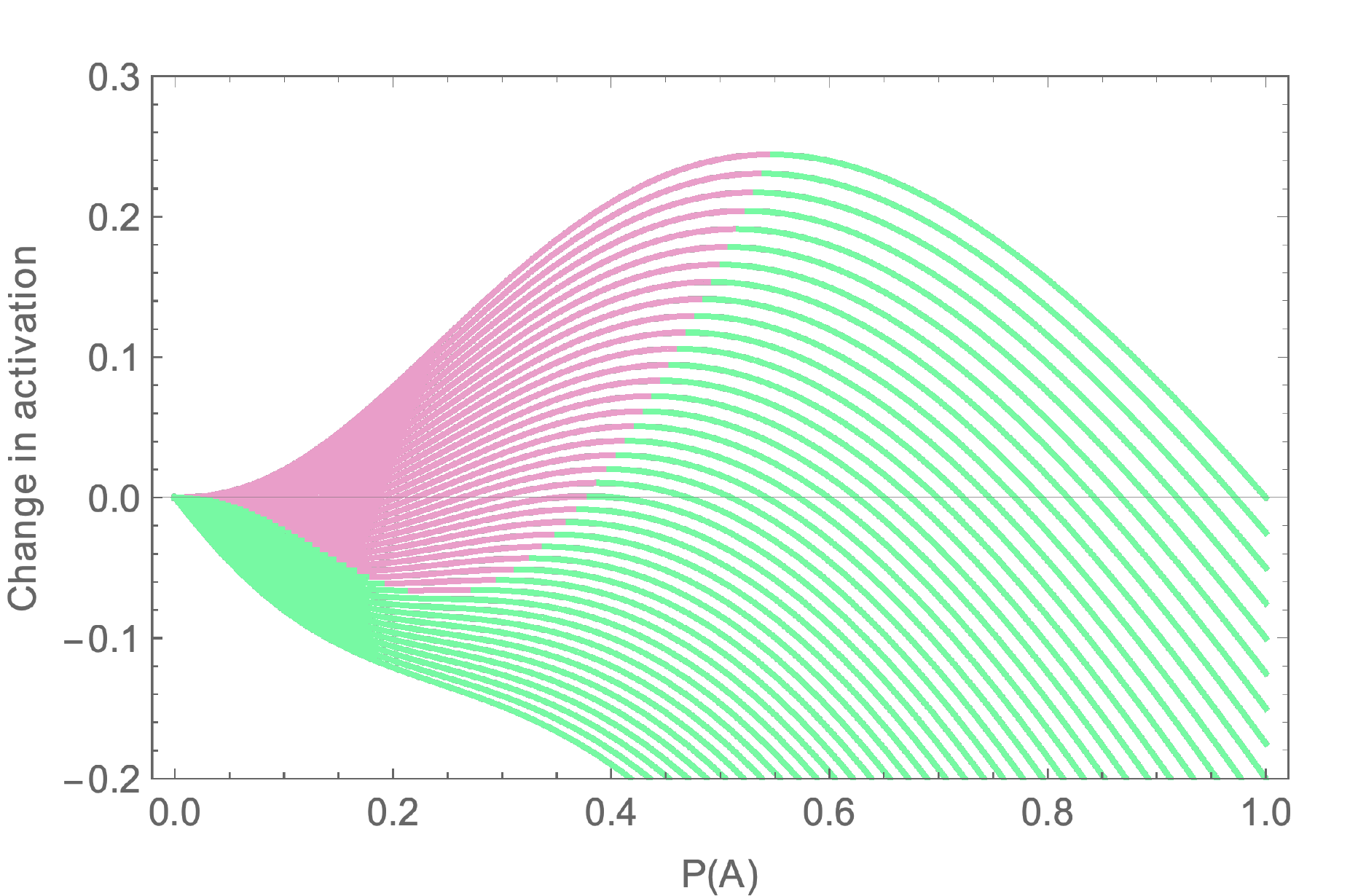}}\\
\endgroup
\caption{Activation, deactivation, and net activation in activation threshold $l$ and deactivation probability $r$ (blue is zero, red is one) for $K\sim\text{Poisson}(c=50)$ and the constant density kernel with $p=1/10$}\label{fig:orbitsdeact}
\end{figure}

\begin{figure}[h!]
\centering
\includegraphics[width=6.5in]{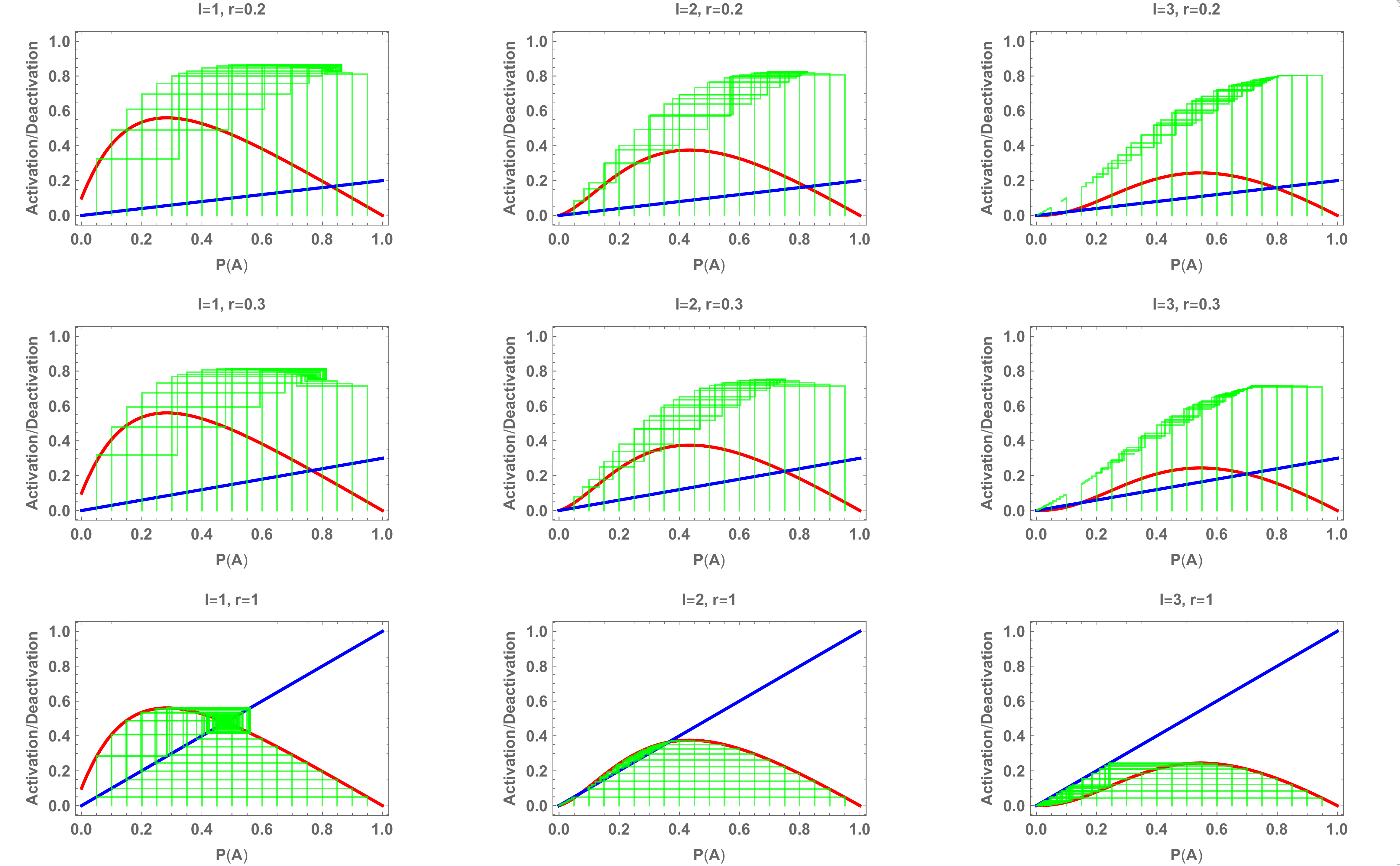}
\caption{Activation, deactivation in activation threshold $l\in\{1,2,3\}$ and deactivation rate $r\in\{0.2,0.3,1\}$ for number of points $K\sim\text{Poisson}(c=50)$, constant density kernel with $p=1/10$,  initializing activation threshold $k=8$, yielding $\xi_0\{A\}\simeq0.1438$}\label{fig:sifp}
\end{figure}
\FloatBarrier

\begin{figure}[h!]
\centering
\includegraphics[width=6.5in]{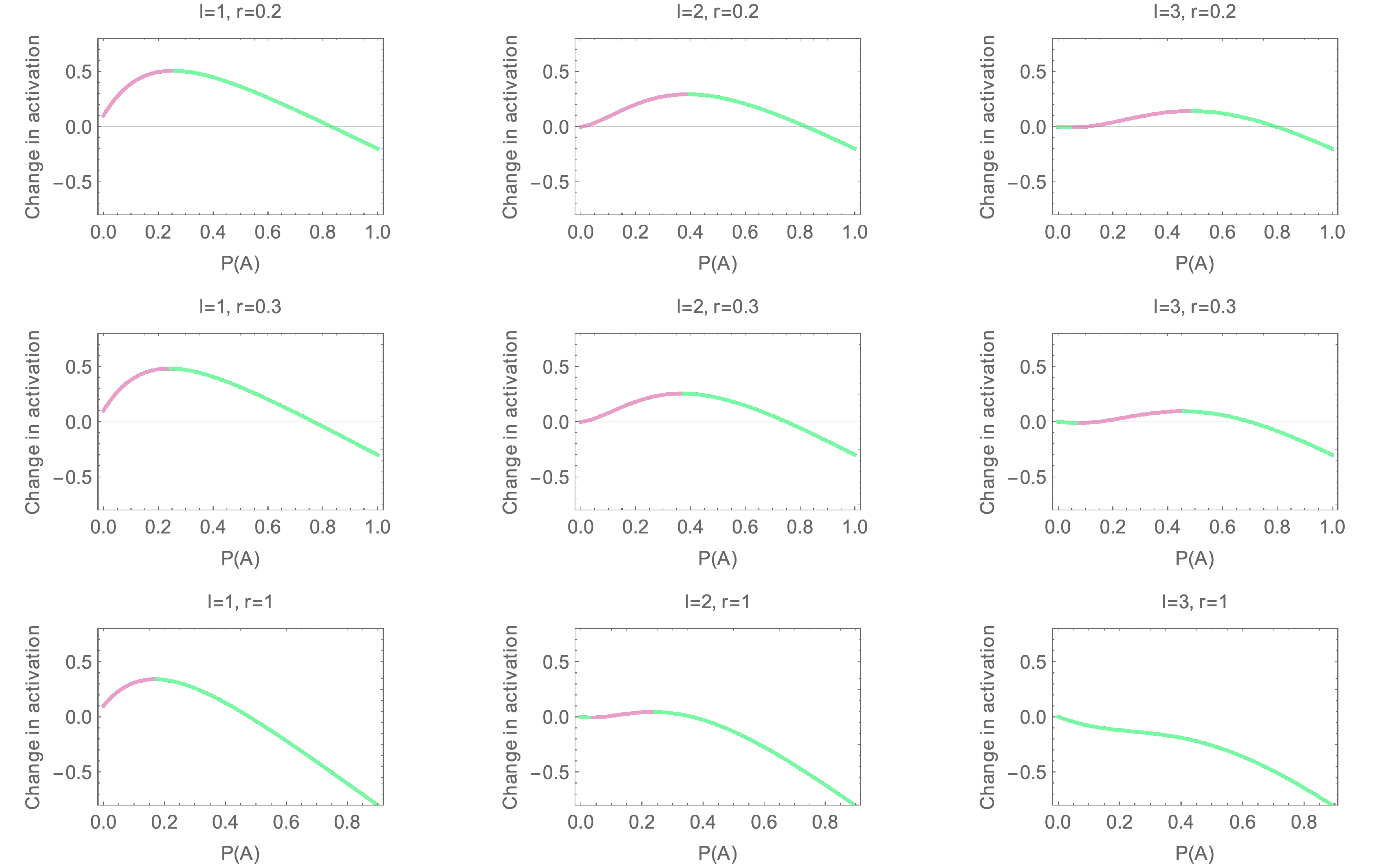}
\caption{Net activation in activation threshold $l\in\{1,2,3\}$ threshold and deactivation rate $r\in\{0.2,0.3,1\}$ for number of points $K\sim\text{Poisson}(c=50)$, constant density kernel with $p=1/10$,  initializing activation threshold $k=8$, yielding $\xi_0\{A\}\simeq0.1438$}\label{fig:sifp}
\end{figure}
\FloatBarrier

\begin{figure}[h!]
\centering
\begingroup
\captionsetup[subfigure]{width=5in,font=normalsize}
\subfloat[$r=0.2$\label{fig:r0}]{\includegraphics[width=3.5in]{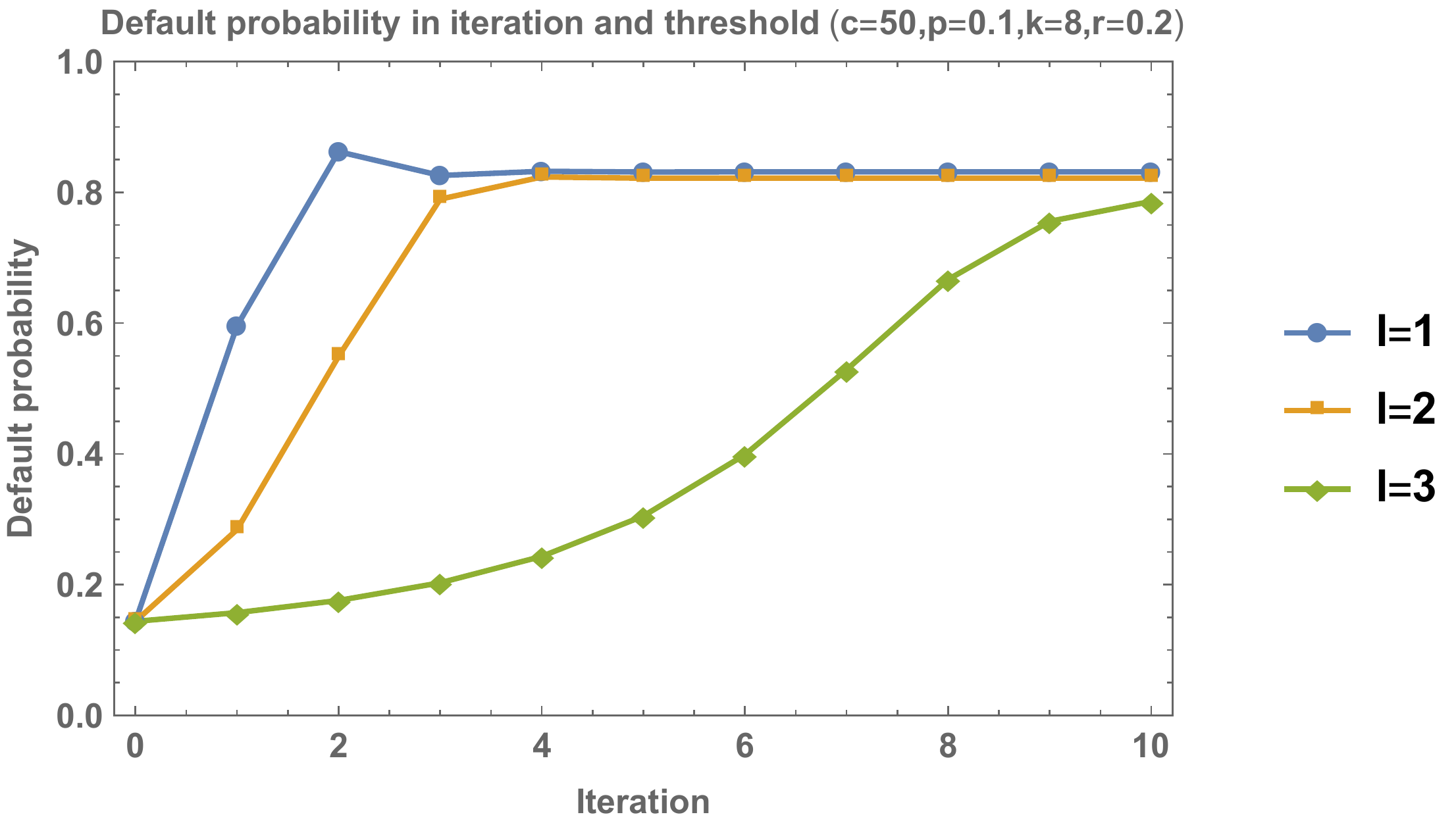}}
\subfloat[$r=0.3$\label{fig:r1}]{\includegraphics[width=3.5in]{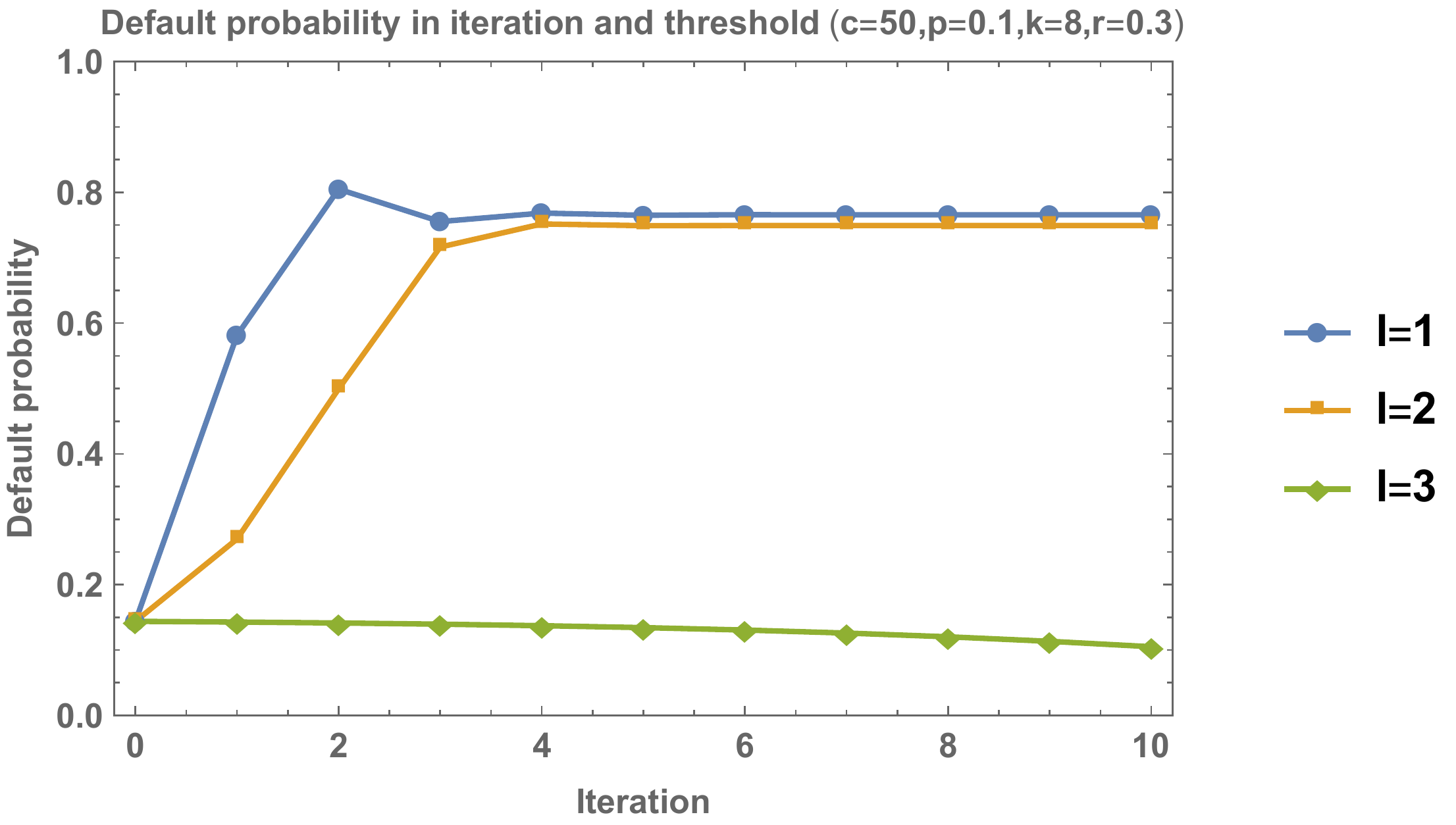}}\\
\subfloat[$r=1$\label{fig:r2}]{\includegraphics[width=3.5in]{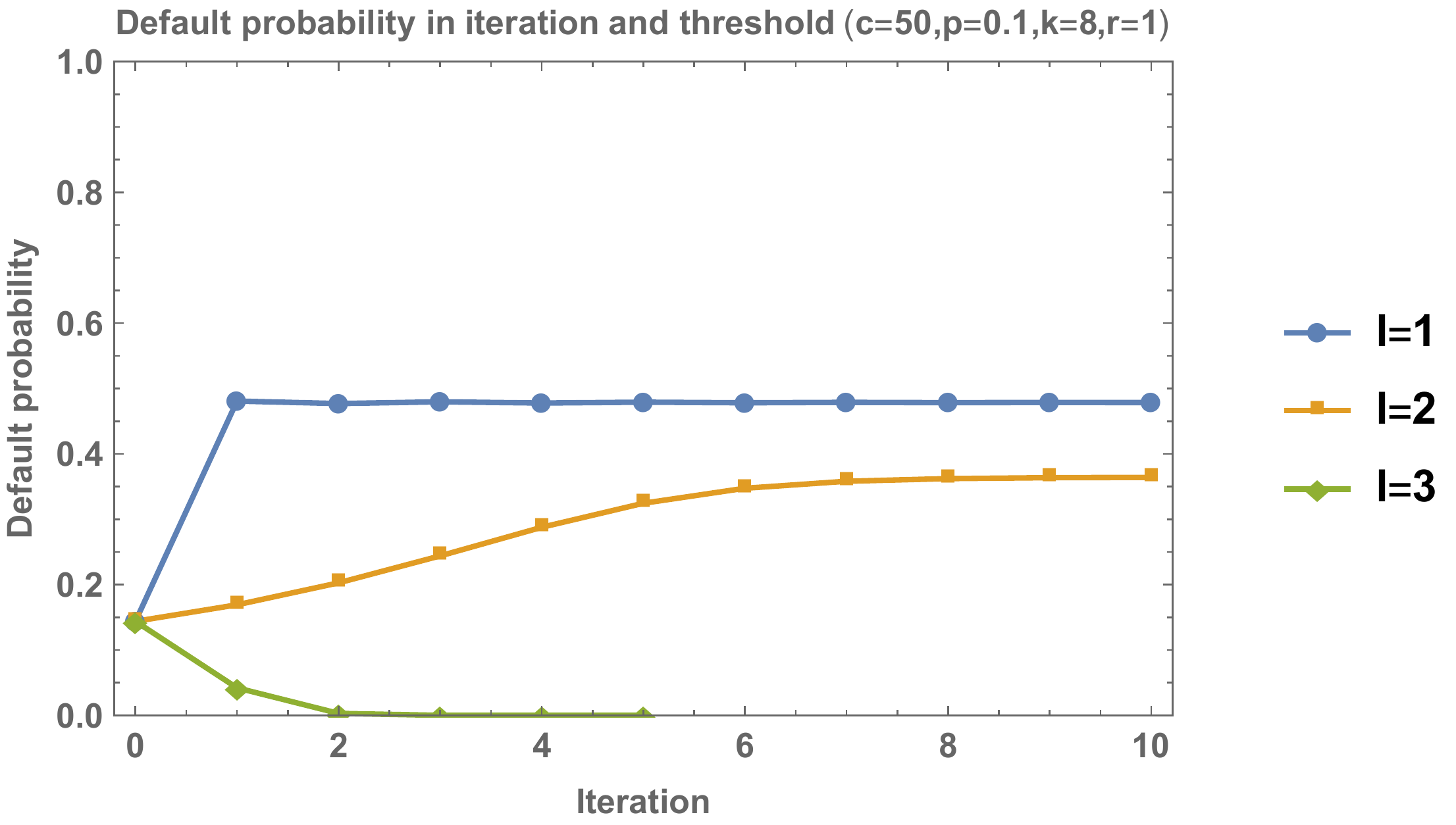}}\\
\endgroup
\caption{Transduction $(\xi_n\{A\}:n=0,1,\dotsb,10)$ over 10 iterations for number of points $K\sim\text{Poisson}(c=50)$, constant density kernel with $p=1/10$, and $\xi_0\{A\}\simeq0.1438$, based on initializing activation threshold $k=8$, with deactivation probability $r\in\{0.2,0.3,1\}$}\label{fig:orbitsdeact}
\end{figure}

\FloatBarrier

\FloatBarrier

\section{Decentral transduction activation-deactivation} 

Here both activation and deactivation are density based (`dual' model). The first example is intrinsic polarization, whereas the second contains external forcings on activation and deactivation. 

\subsection{Pure intrinsic} The activation recursion is given by \[\xi_n\{A\}\equiv \xi_{n-1}\{A\} + \nu(p_{n-1}^B)(1-\xi_{n-1}\{A\}) - \nu(q_{n-1}^B)\xi_{n-1}\{A\}\for n\ge 1\] where $q_{n}^B(x)=\P_n(d_B(x,A)\ge m)+\P_n(d_B(x,A)= m-1)f(x,x)$. 

\begin{figure}[h!]
\centering
\begingroup
\captionsetup[subfigure]{width=5in,font=normalsize}
\subfloat[Activation with threshold $l\in\{1,2,3\}$\label{fig:r0}]{\includegraphics[width=3.5in]{activationcurves.pdf}}
\subfloat[Deactivation with threshold $m\in\{1,2,3\}$\label{fig:r1}]{\includegraphics[width=3.5in]{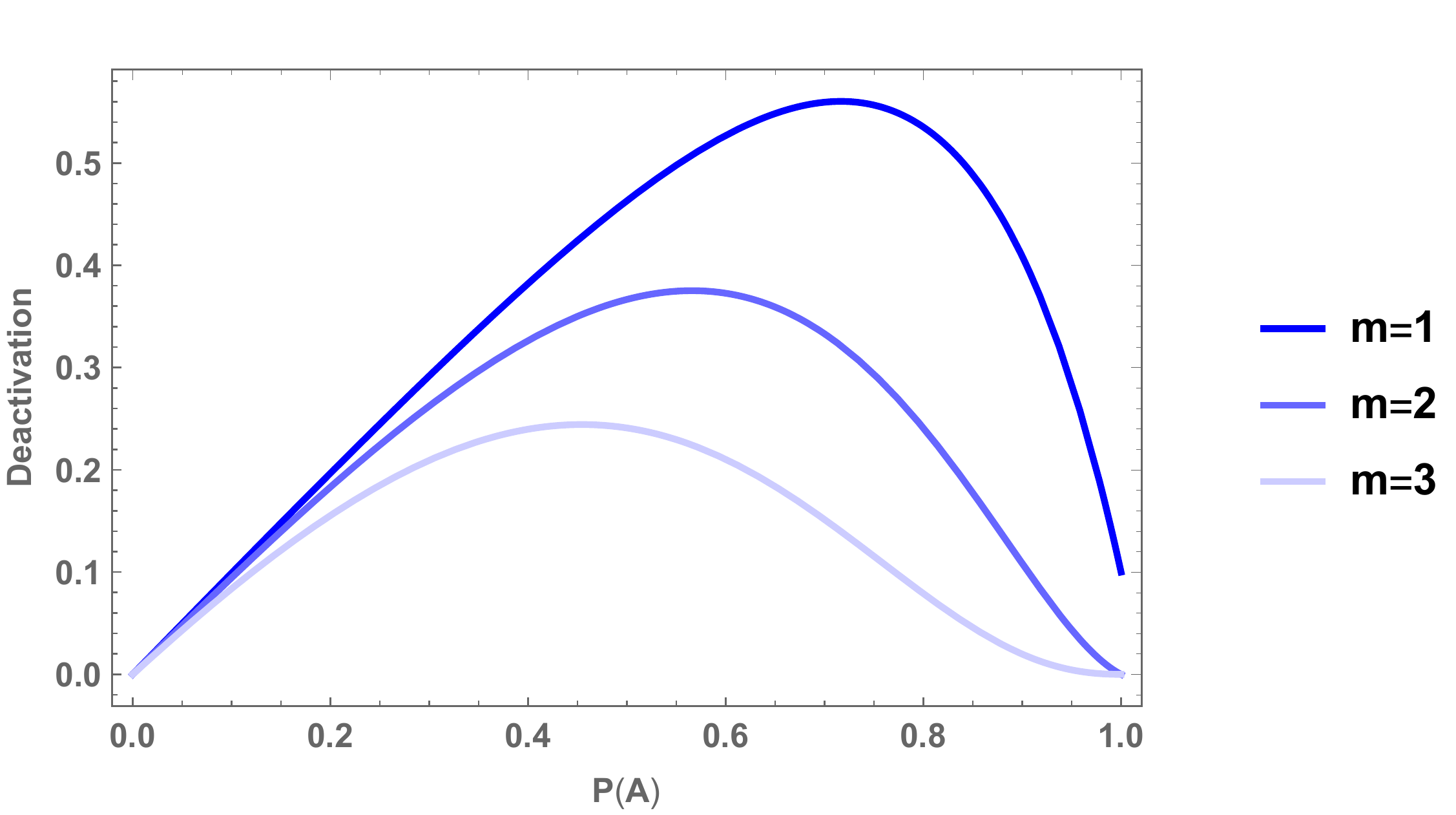}}\\
\subfloat[Net activation for $(l,m)\in\{1,2,3\}\times\{1,2,3\}$\label{fig:r2}]{\includegraphics[width=6.5in]{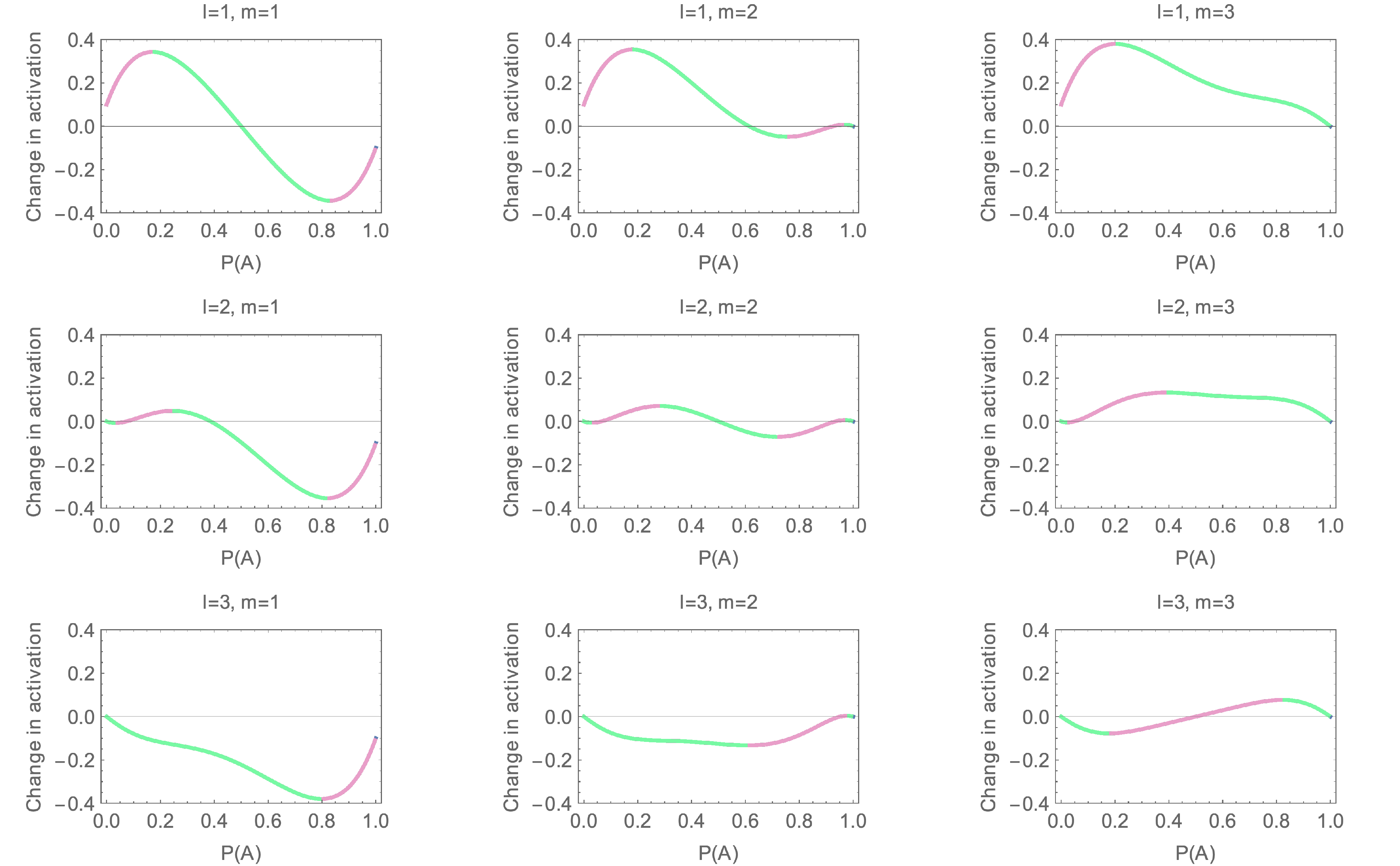}}
\endgroup
\caption{Polarization in activation, deactivation, and net activation, in activation threshold $l$ and deactivation probability $r$ (blue is zero, red is one), for $K\sim\text{Poisson}(c=50)$ and the constant density kernel with $p=1/10$}\label{fig:orbitsdeact}
\end{figure}

\begin{figure}[h!]
\centering
\includegraphics[width=6.5in]{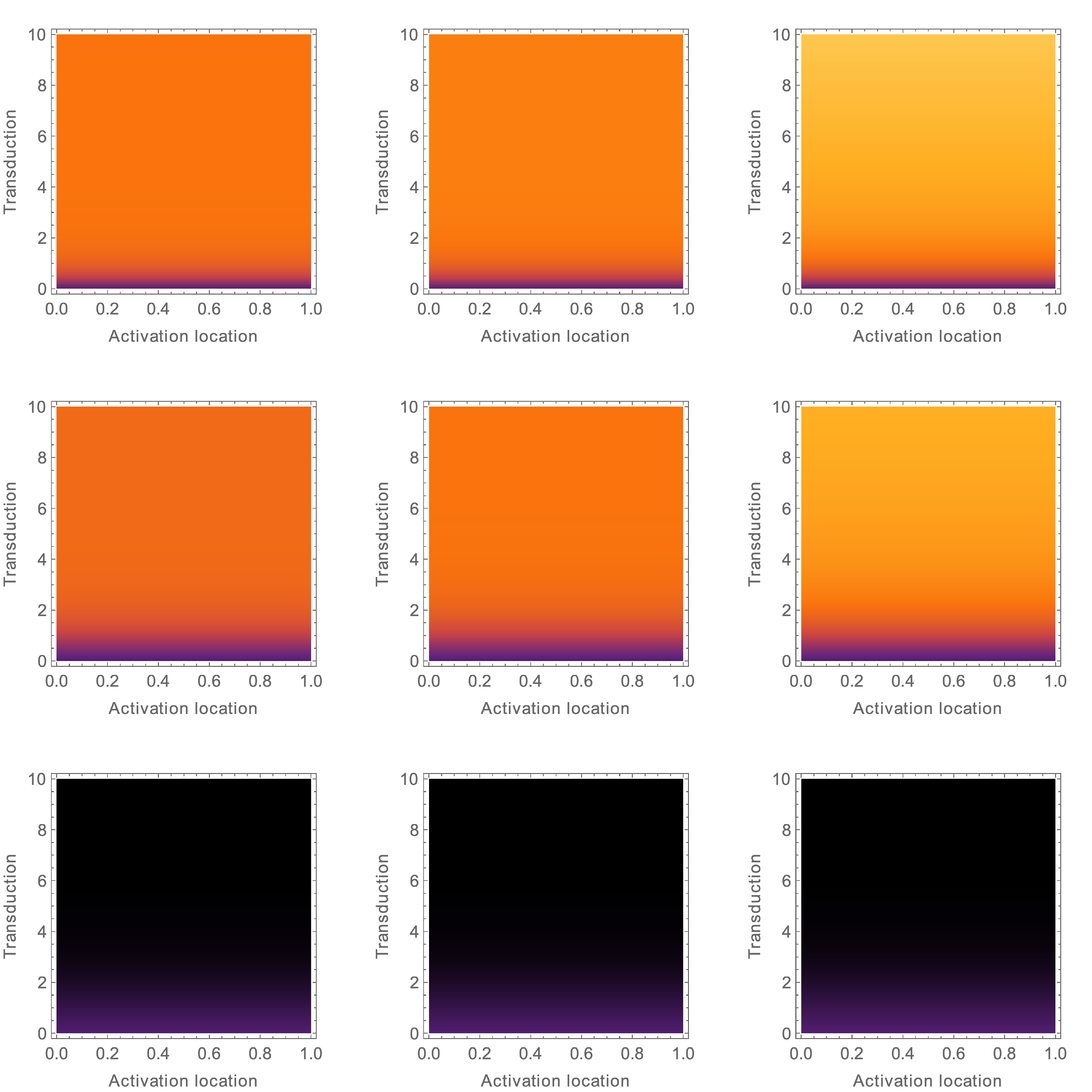}
\caption{Decentral transduction activation-deactivation in location  in activation $l$ and deactivation $m$ thresholds for number of points $K\sim\text{Poisson}(c=50)$, initializing activation threshold $k=8$}\label{fig:sifp}
\end{figure}
\FloatBarrier

\begin{figure}[h!]
\centering
\includegraphics[width=6.5in]{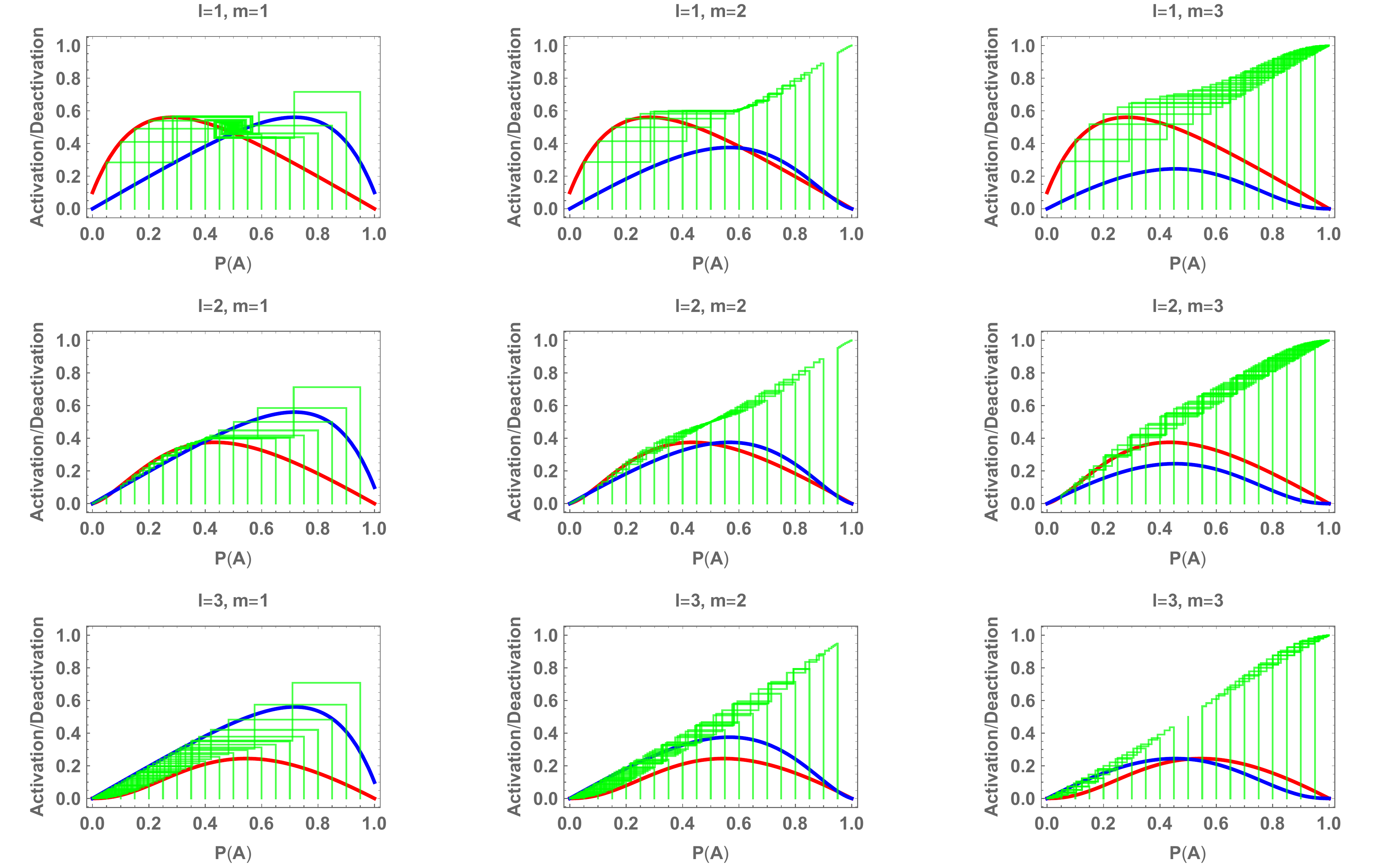}
\caption{Polarization in activation $l$ and deactivation $m$ thresholds for number of points $K\sim\text{Poisson}(c=50)$, constant density kernel with $p=1/10$,  initializing activation threshold $k=8$, yielding $\xi_0\{A\}\simeq0.1438$}\label{fig:sifp}
\end{figure}
\FloatBarrier

\begin{figure}[h!]
\centering
\begingroup
\captionsetup[subfigure]{width=5in,font=normalsize}
\subfloat[\label{fig:r0}]{\includegraphics[width=3.5in]{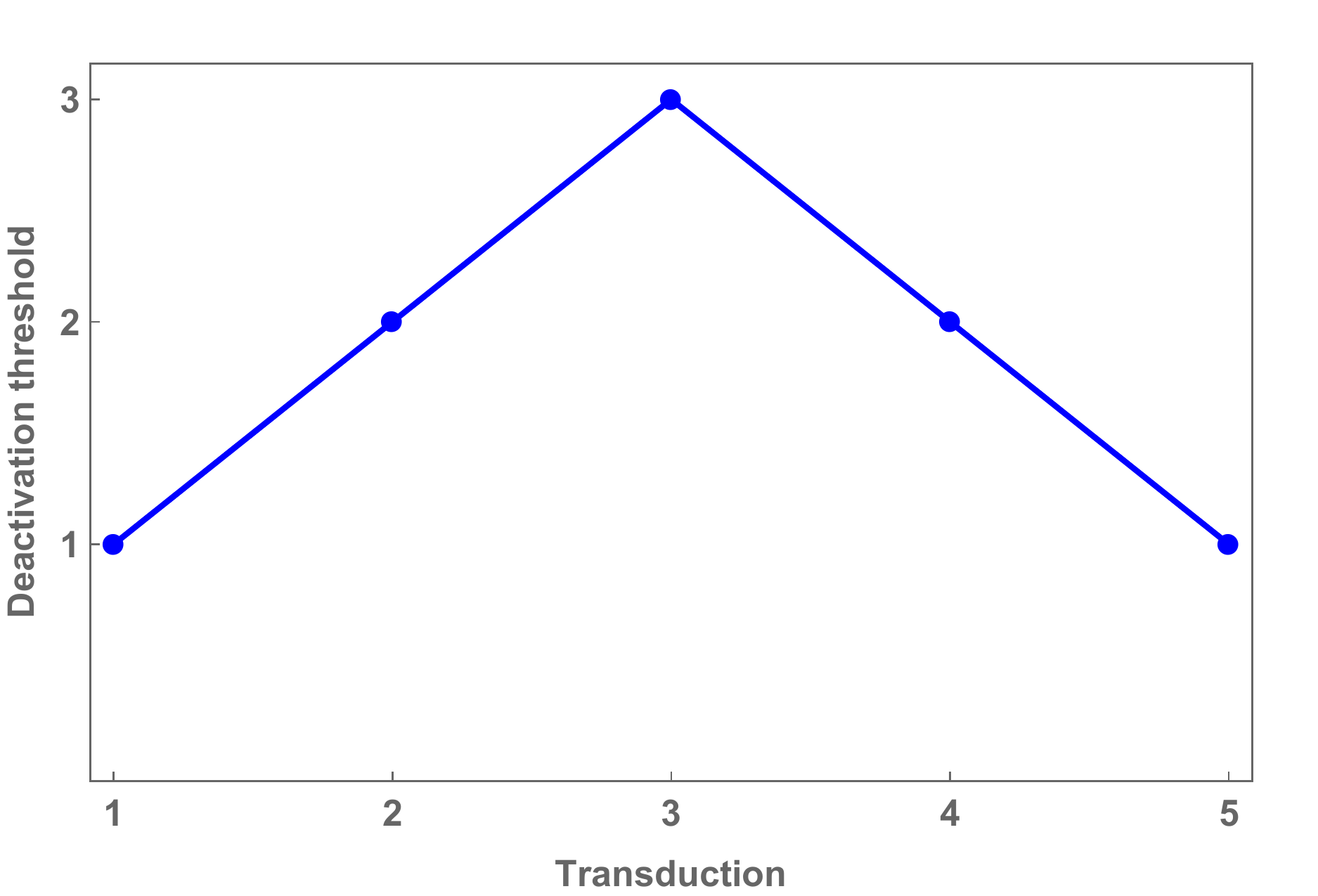}}
\subfloat[\label{fig:r1}]{\includegraphics[width=3.5in]{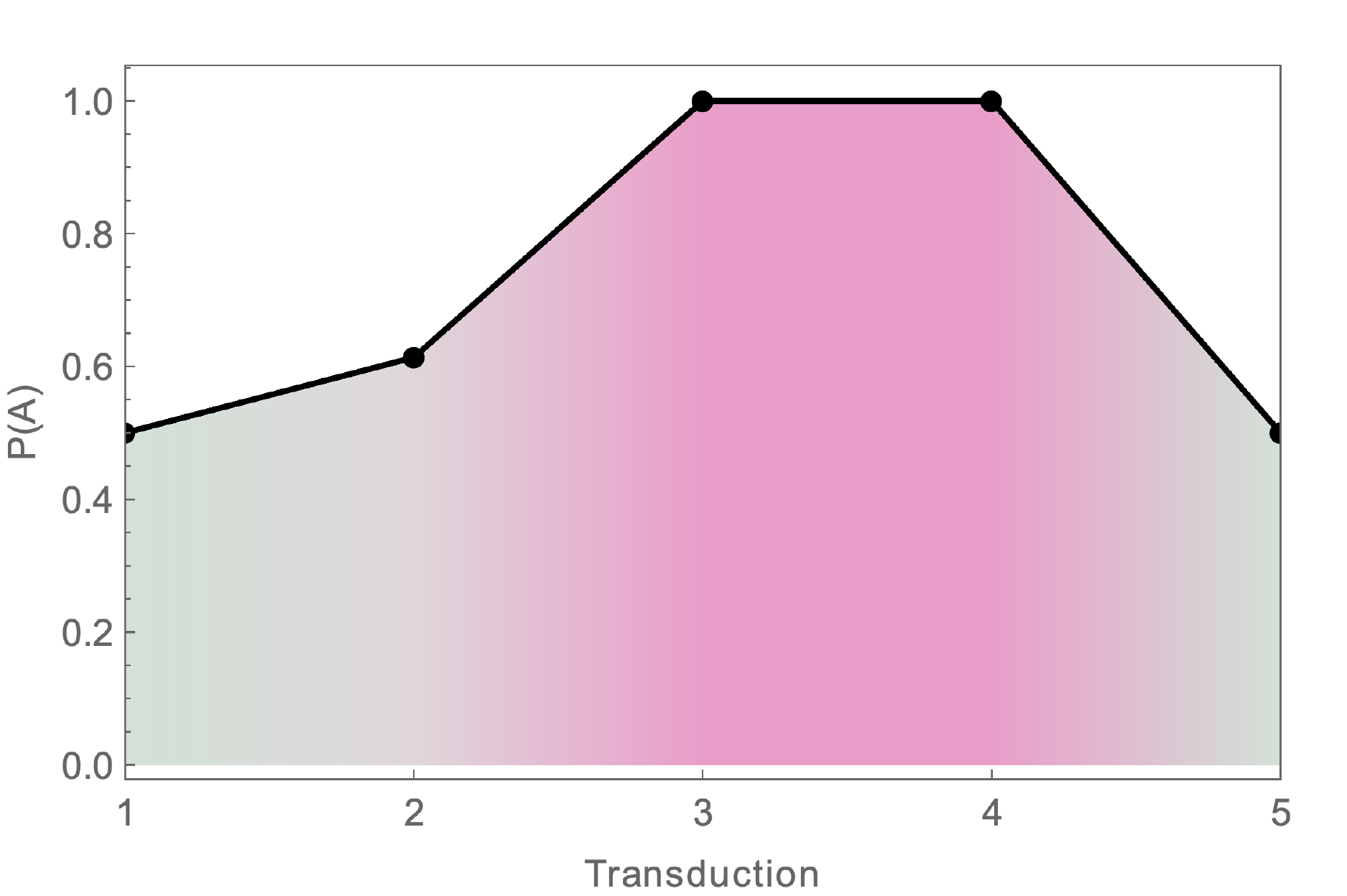}}\\
\subfloat[\label{fig:r2}]{\includegraphics[width=3.5in]{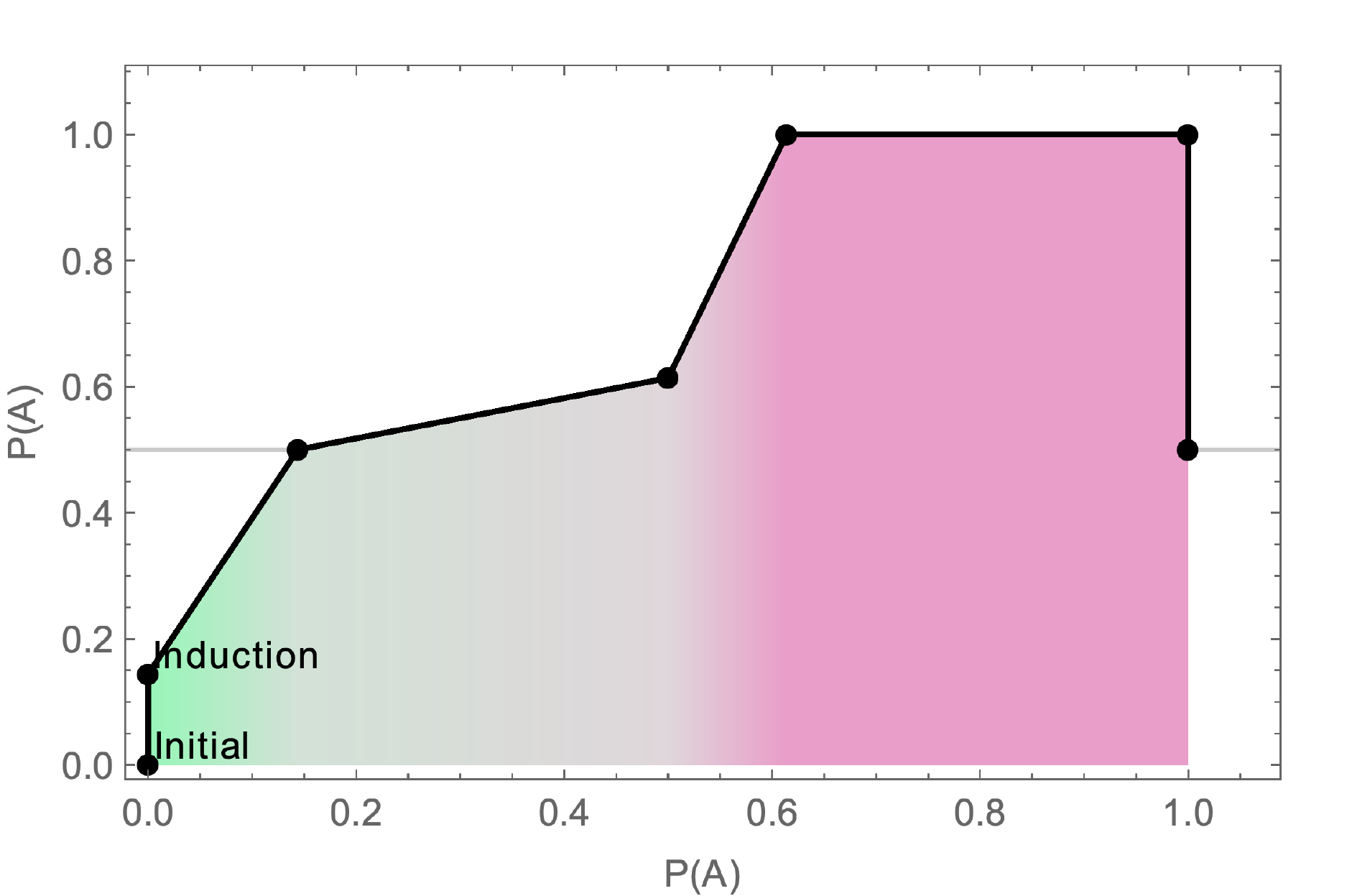}}
\subfloat[\label{fig:r2}]{\includegraphics[width=3.5in]{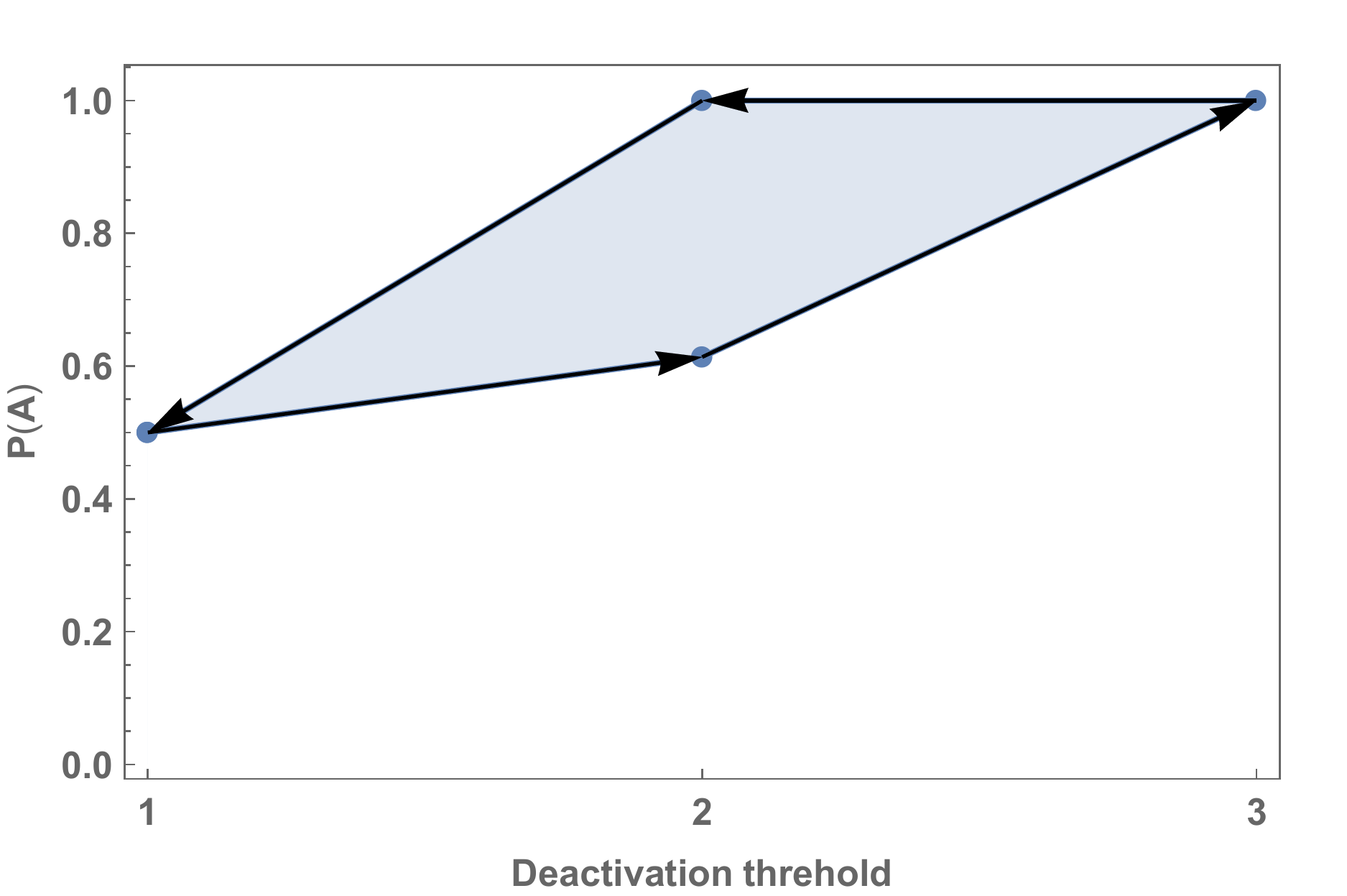}}\\
\endgroup
\caption{Hysteresis of polarization activation for fixed activation threshold $l=1$ and varying deactivation thresholds with sequence $(1,2,3,2,1)$; with $K\sim\text{Poisson}(c=50)$, the constant density kernel with $p=1/10$, and induction by initializing activation threshold $k=8$, yielding $\xi_0\{A\}\simeq0.1438$}\label{fig:orbitsdeact}
\end{figure}

\begin{figure}[h!]
\centering
\begingroup
\captionsetup[subfigure]{width=5in,font=normalsize}
\subfloat[\label{fig:r0}]{\includegraphics[width=3.5in]{hysteresispolarizationthree.pdf}}
\subfloat[\label{fig:r1}]{\includegraphics[width=3.5in]{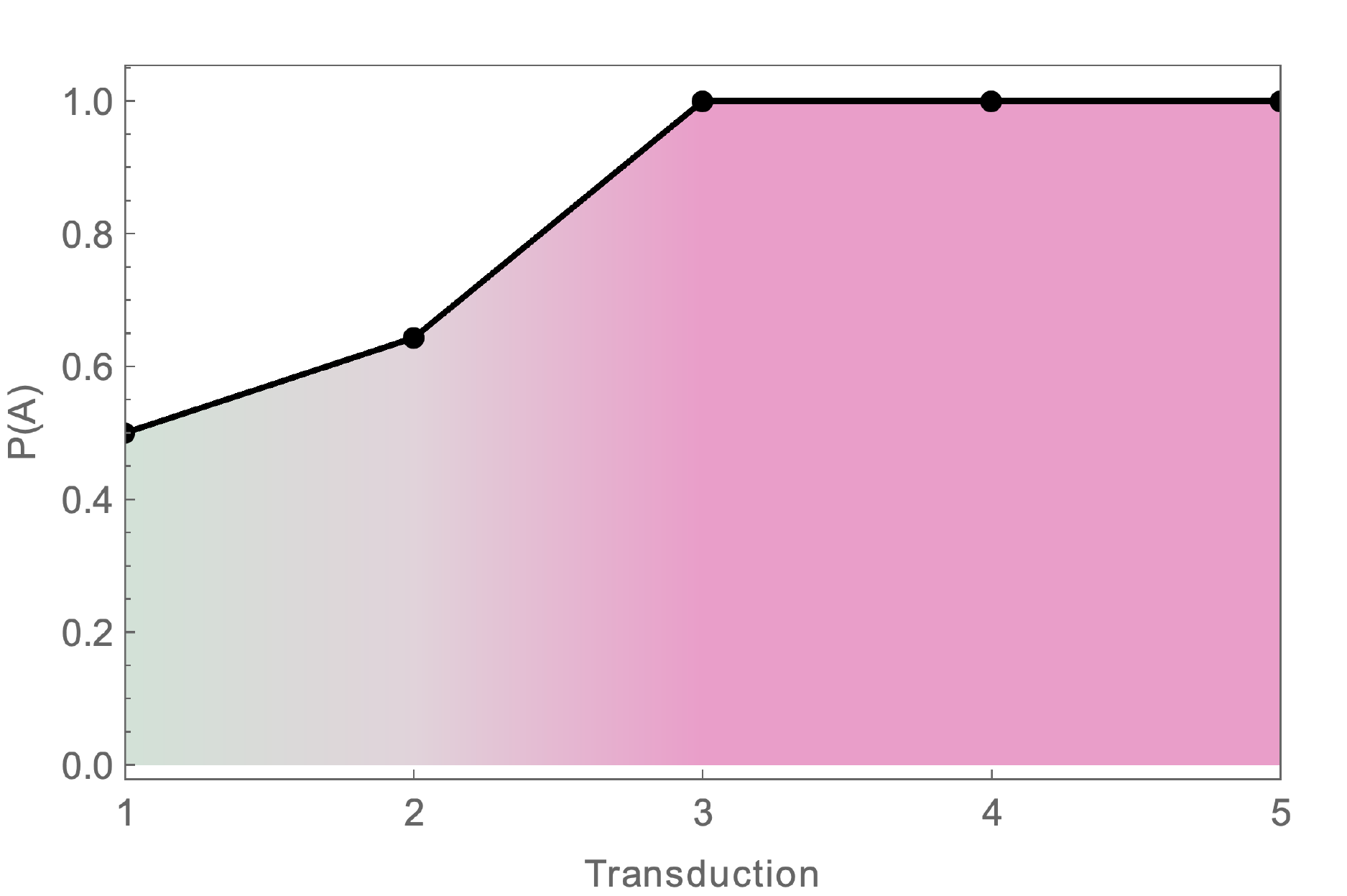}}\\
\subfloat[\label{fig:r2}]{\includegraphics[width=3.5in]{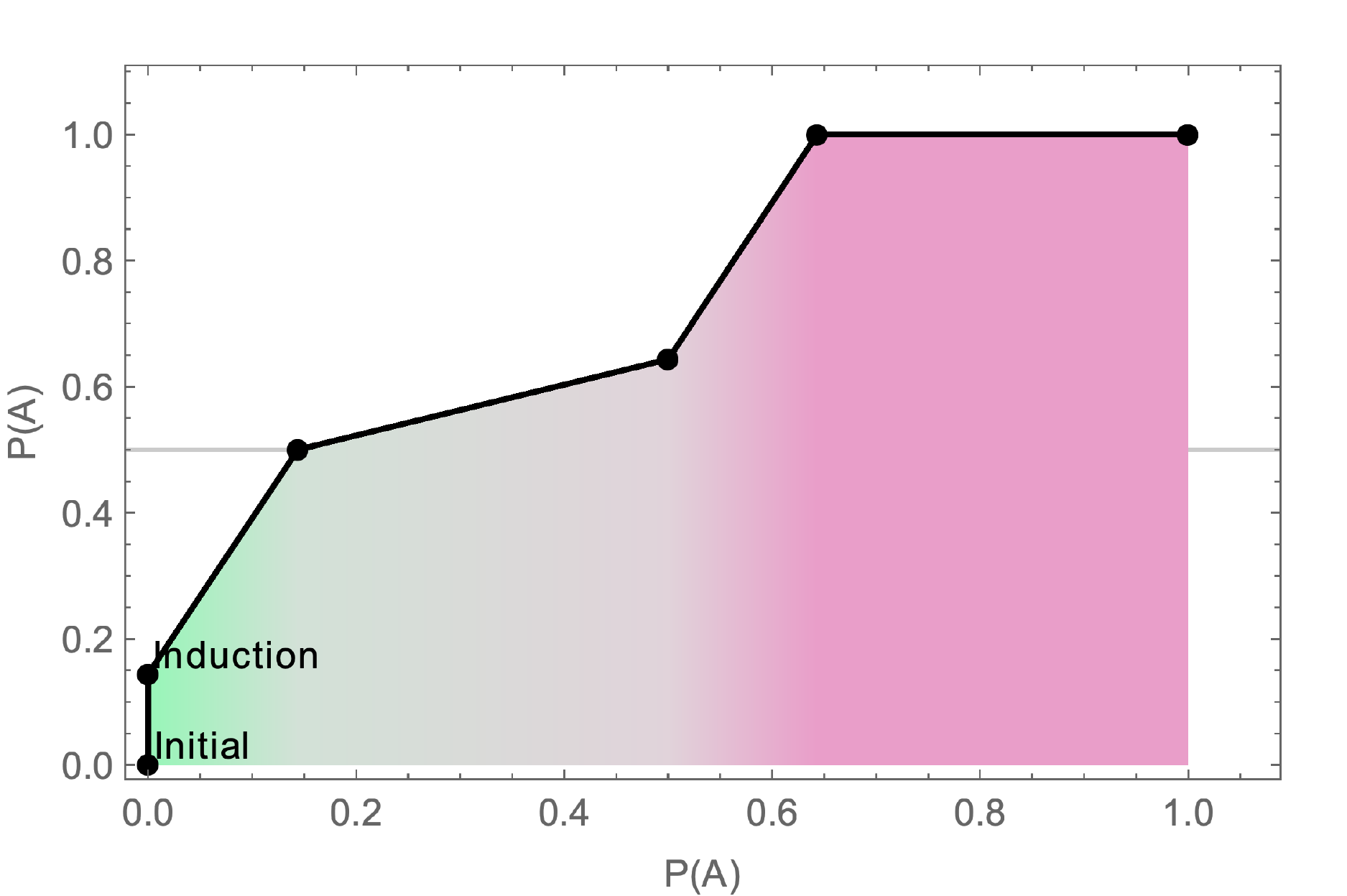}}
\subfloat[\label{fig:r2}]{\includegraphics[width=3.5in]{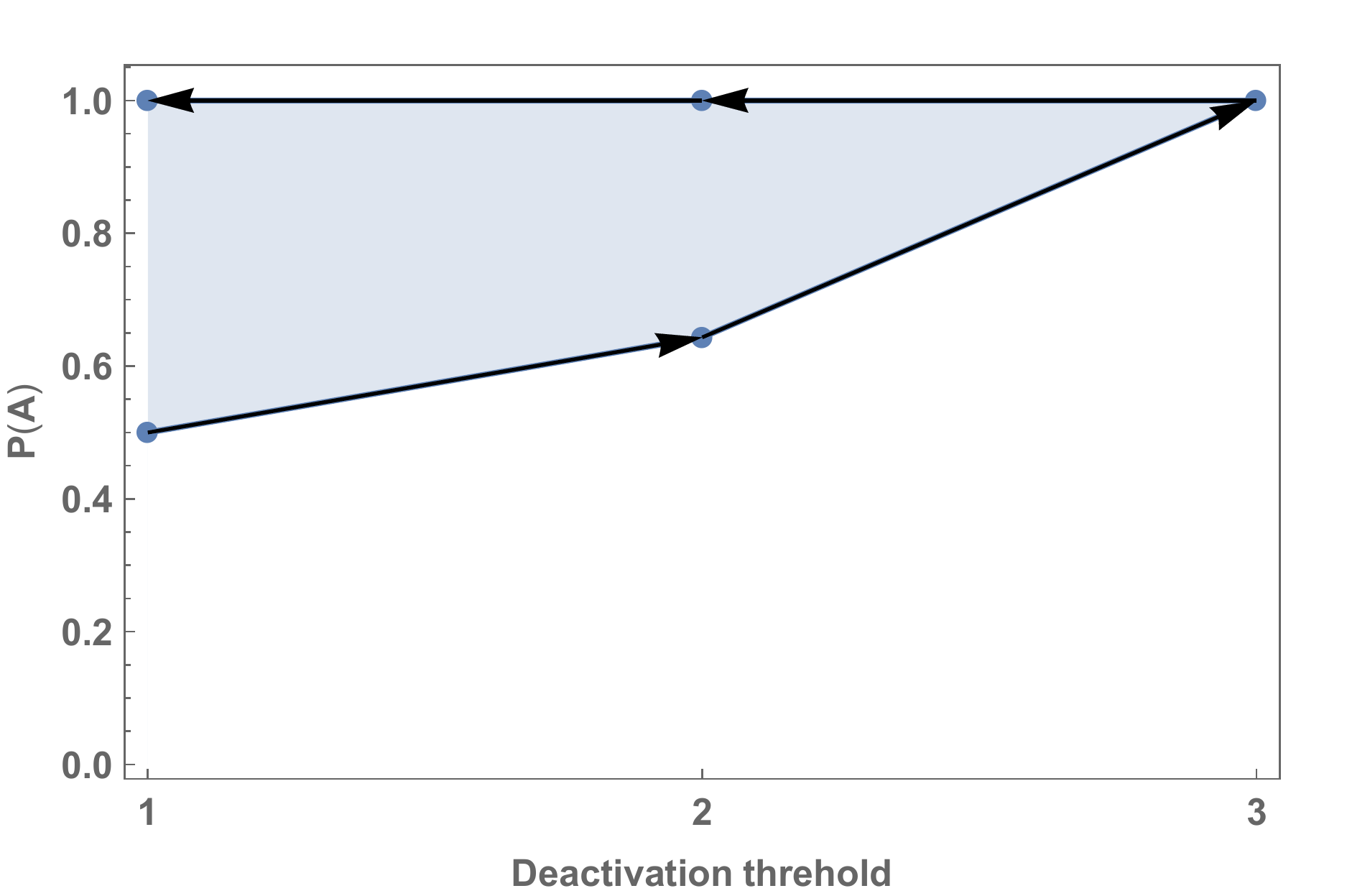}}\\
\endgroup
\caption{Hysteresis of polarization activation for no self-activation-deactivation and fixed activation threshold $l=1$ and varying deactivation thresholds with sequence $(1,2,3,2,1)$; with $K\sim\text{Poisson}(c=50)$, the constant density kernel with $p=1/10$, and induction by initializing activation threshold $k=8$, yielding $\xi_0\{A\}\simeq0.1438$}\label{fig:orbitsdeact}
\end{figure}

\begin{figure}[h!]
\centering
\begingroup
\captionsetup[subfigure]{width=5in,font=normalsize}
\subfloat[\label{fig:r0}]{\includegraphics[width=3.5in]{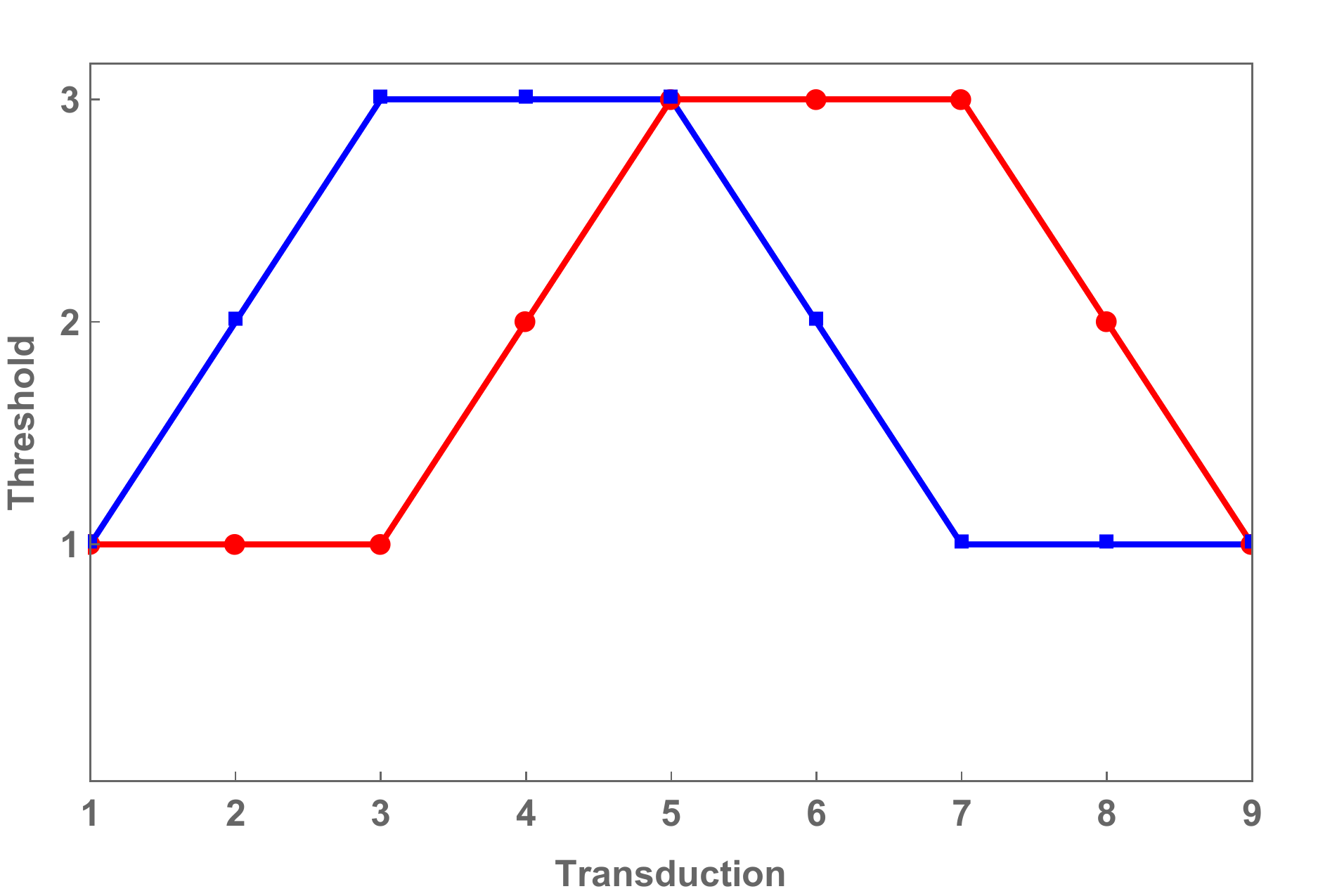}}
\subfloat[\label{fig:r0}]{\includegraphics[width=3.5in]{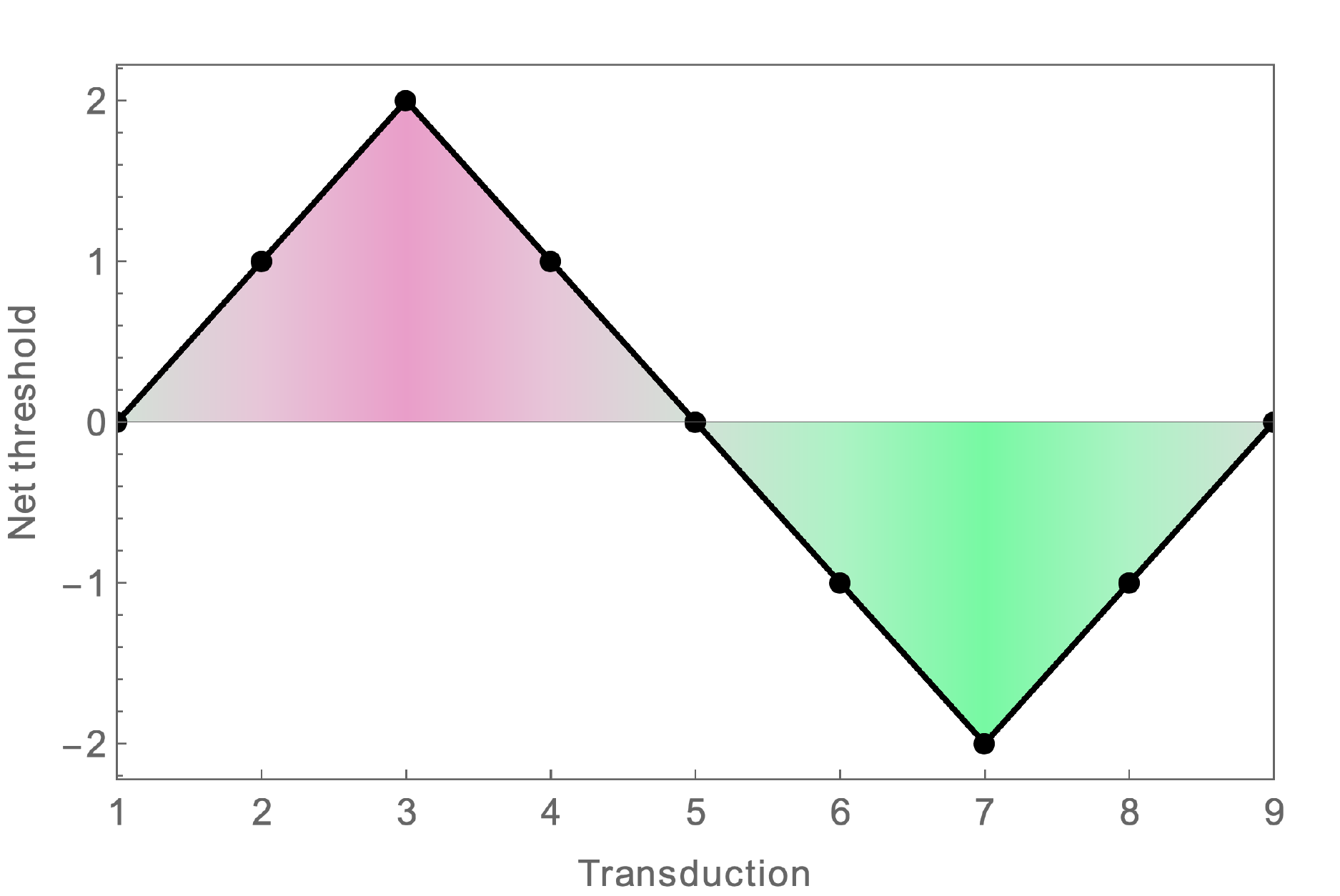}}\\
\subfloat[\label{fig:r0}]{\includegraphics[width=3.5in]{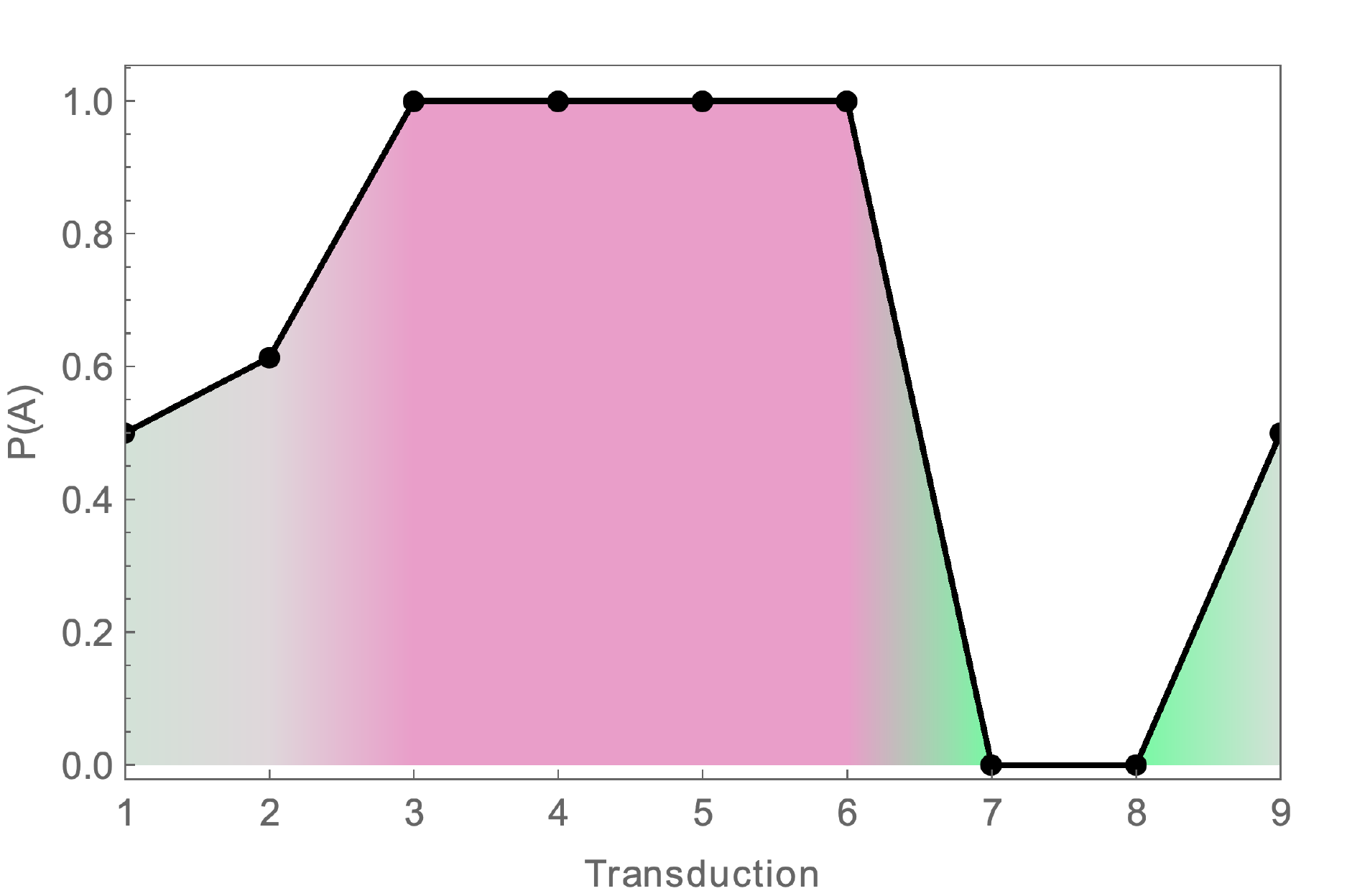}}
\subfloat[\label{fig:r0}]{\includegraphics[width=3.5in]{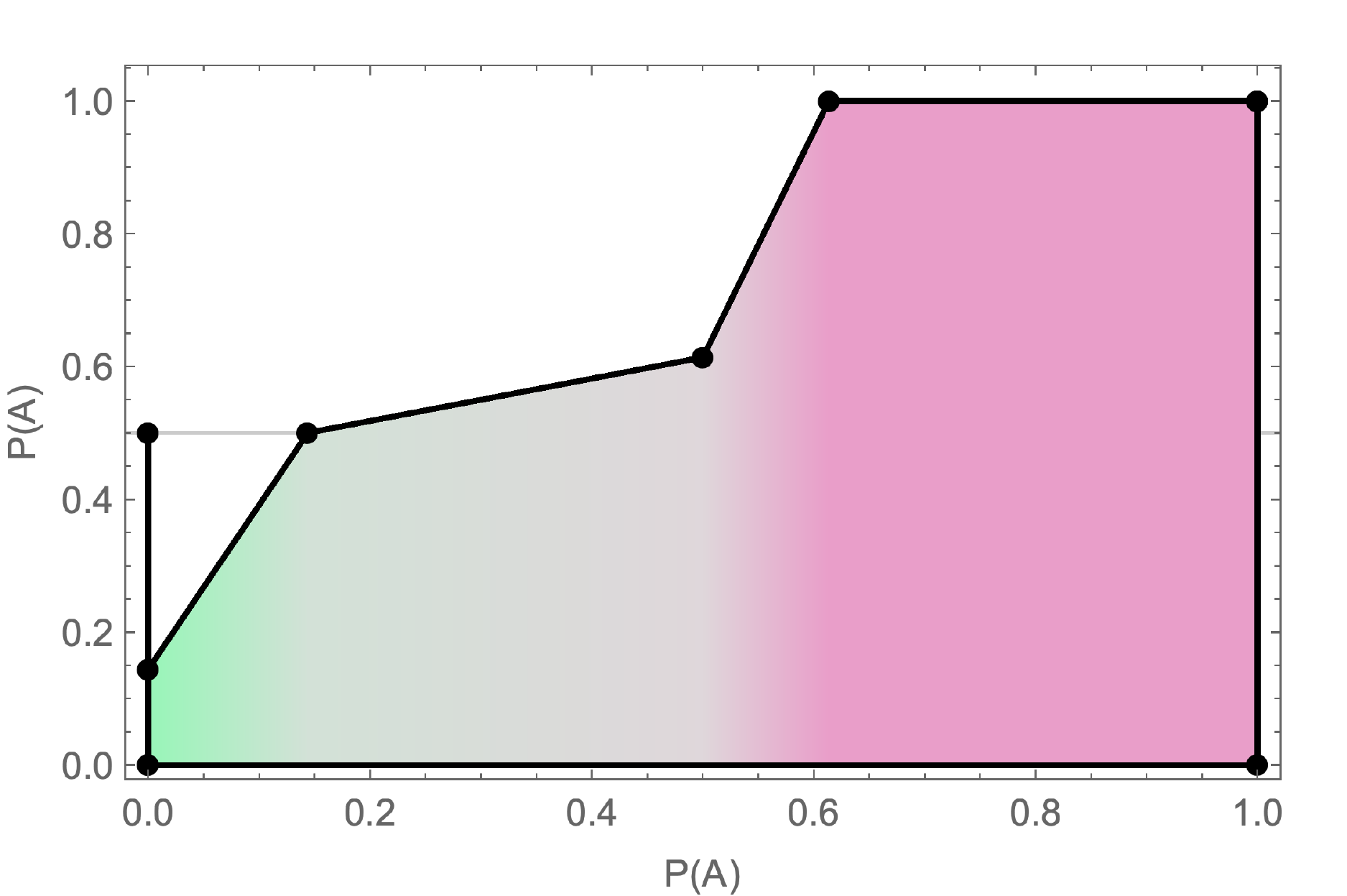}}\\
\subfloat[\label{fig:r0}]{\includegraphics[width=4.5in]{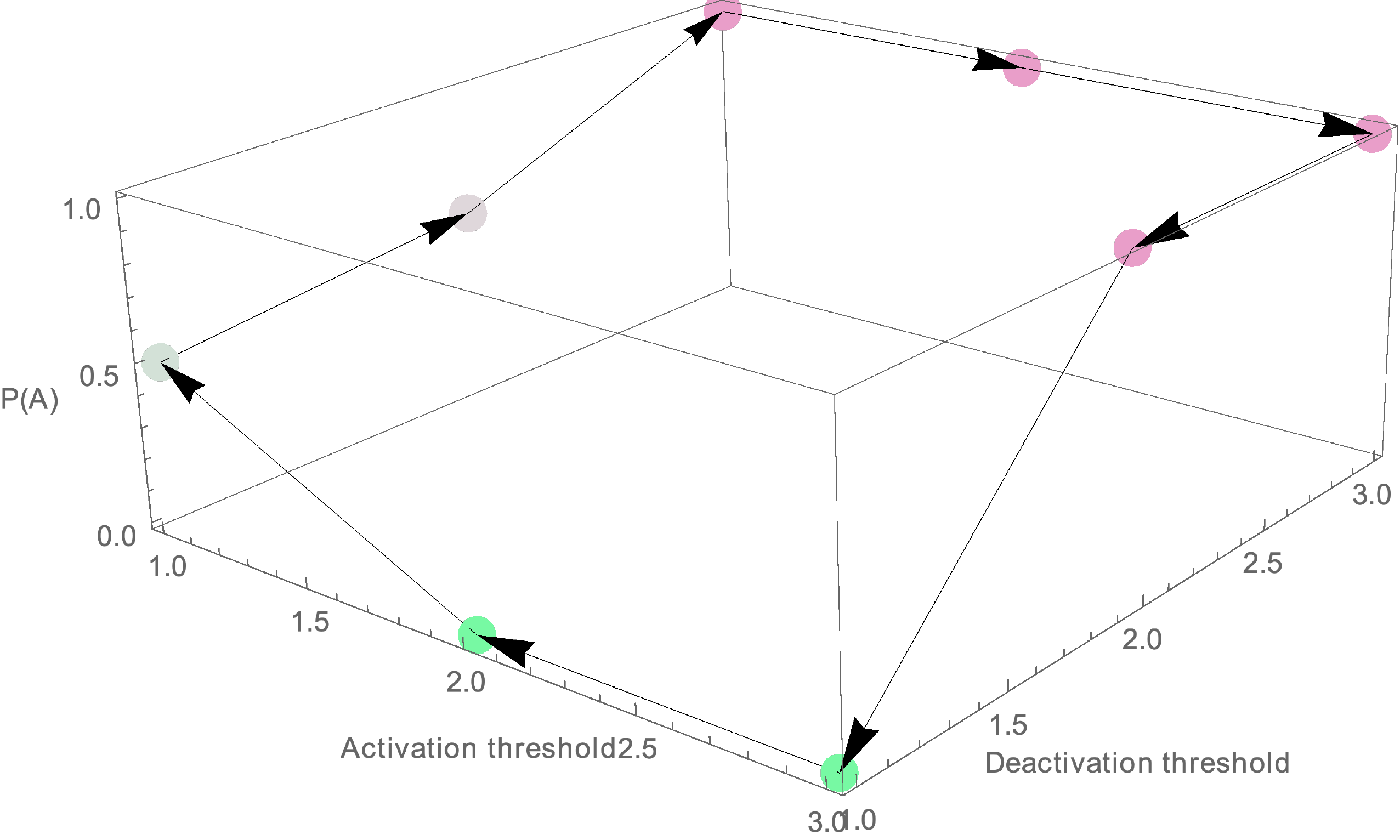}}
\endgroup
\caption{Hysteresis of polarization activation for activation-deactivation threshold sequence $((1,1),(1,2),(1,3),(2,3),(3,3),(3,2),(3,1),(2,1),(1,1))$ and self-activation-deactivation; with $K\sim\text{Poisson}(c=50)$, the constant density kernel with $p=1/10$, and induction by initializing activation threshold $k=8$, yielding $\xi_0\{A\}\simeq0.1438$}\label{fig:orbitsdeact}
\end{figure}

\begin{figure}[h!]
\centering
\begingroup
\captionsetup[subfigure]{width=5in,font=normalsize}
\subfloat[\label{fig:r0}]{\includegraphics[width=3.5in]{hysteresispolarizationthreedimensionsseq.pdf}}
\subfloat[\label{fig:r0}]{\includegraphics[width=3.5in]{hysteresispolarizationthreedimensionschange.pdf}}\\
\subfloat[\label{fig:r0}]{\includegraphics[width=3.5in]{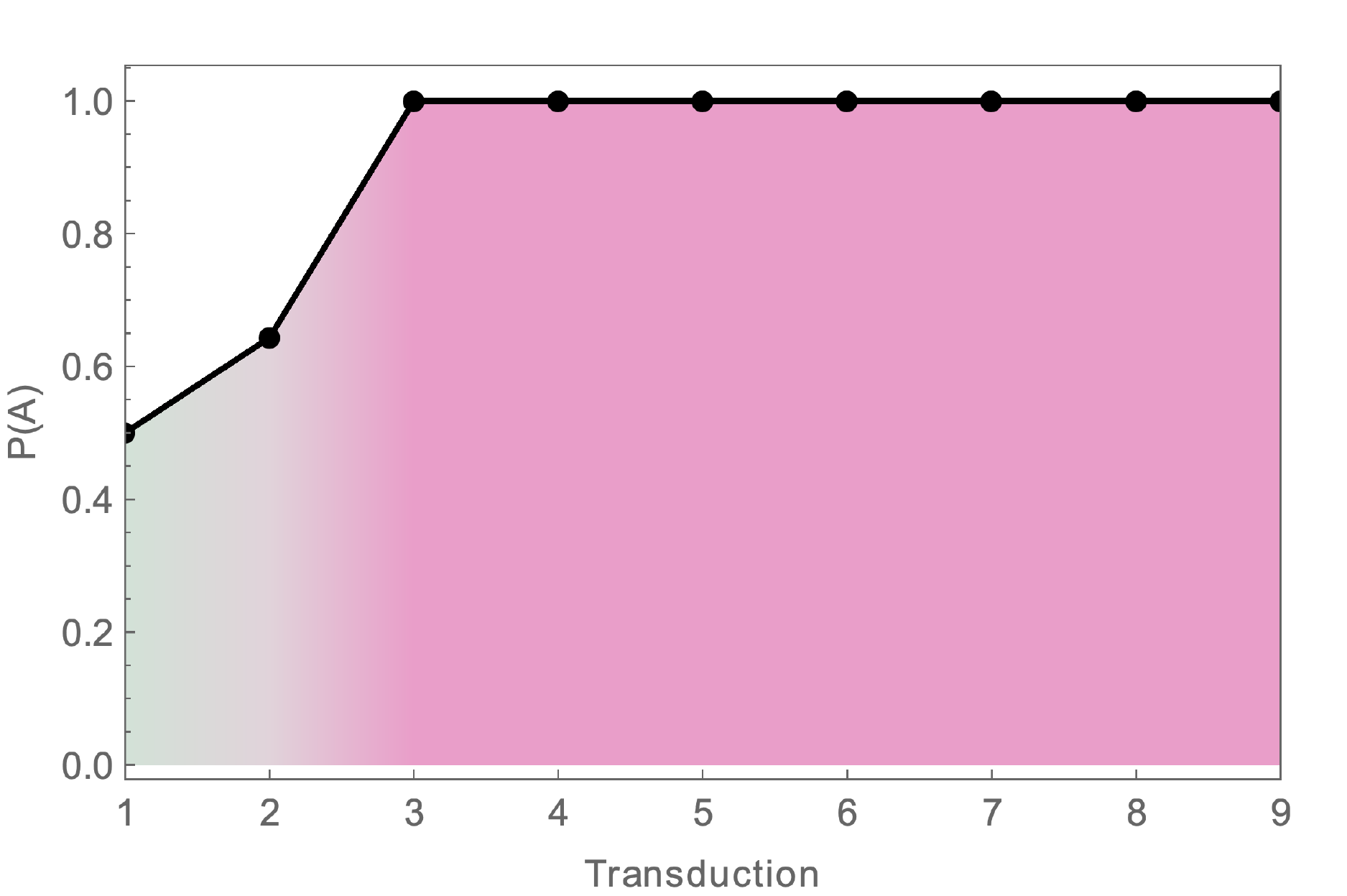}}
\subfloat[\label{fig:r0}]{\includegraphics[width=3.5in]{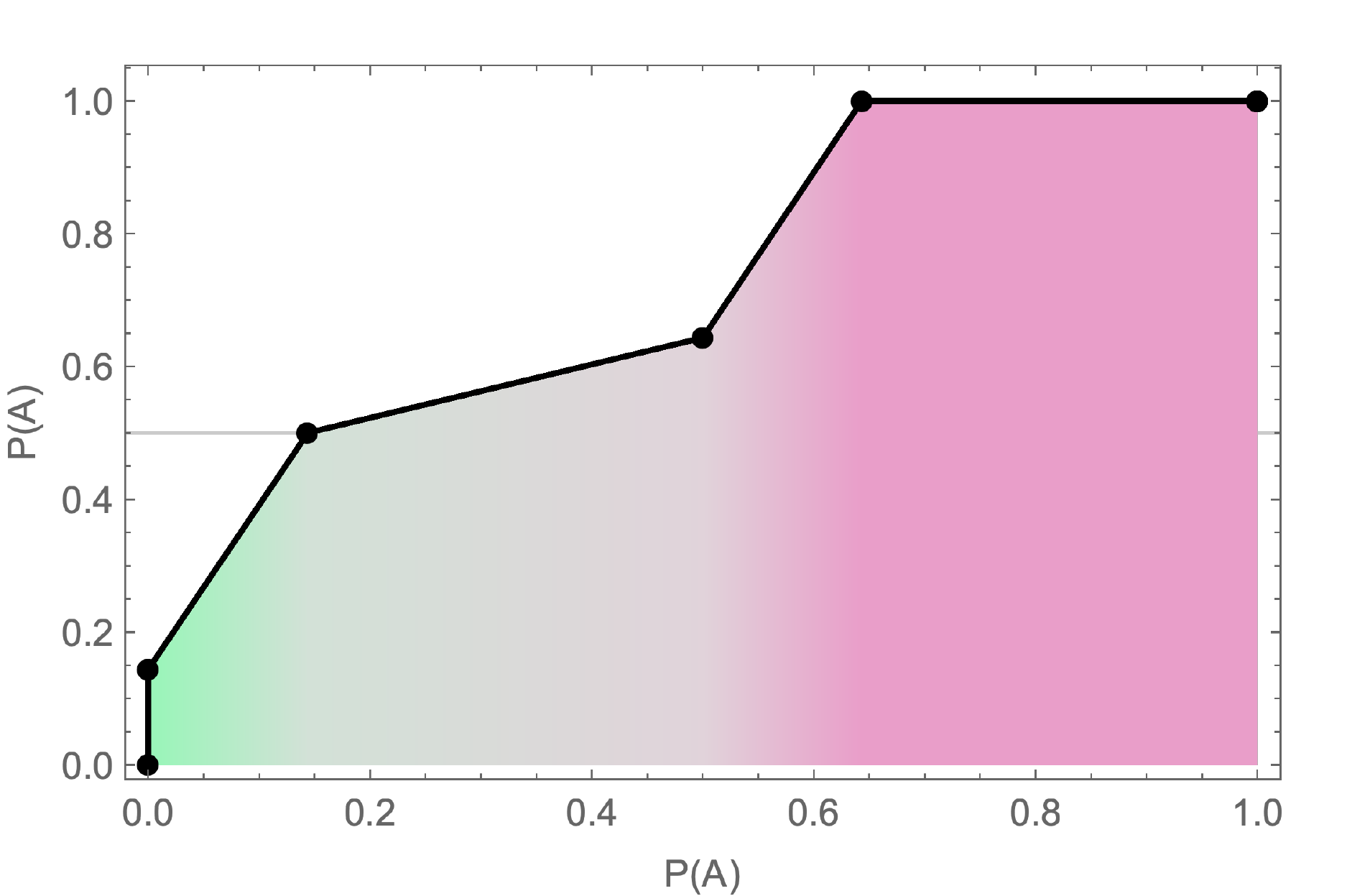}}\\
\subfloat[\label{fig:r0}]{\includegraphics[width=4.5in]{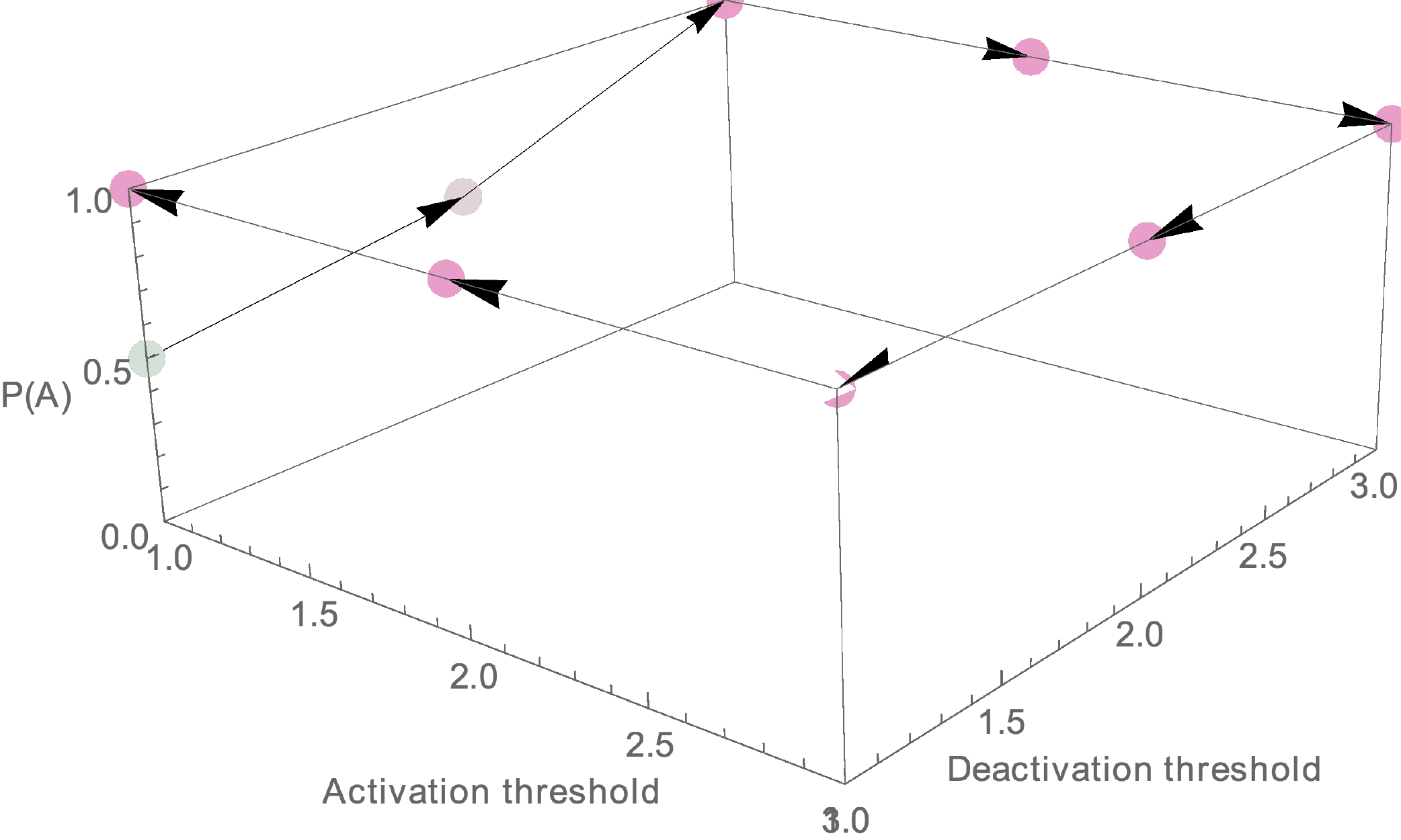}}
\endgroup
\caption{Hysteresis of polarization activation for activation-deactivation threshold sequence $((1,1),(1,2),(1,3),(2,3),(3,3),(3,2),(3,1),(2,1),(1,1))$ and no self-activation-deactivation; with $K\sim\text{Poisson}(c=50)$, the constant density kernel with $p=1/10$, and induction by initializing activation threshold $k=8$, yielding $\xi_0\{A\}\simeq0.1438$}\label{fig:orbitsdeact}
\end{figure}

\begin{figure}[h!]
\centering
\begingroup
\captionsetup[subfigure]{width=5in,font=normalsize}
\subfloat[\label{fig:r0}]{\includegraphics[width=3.5in]{hysteresispolarizationthreedimensionsseq.pdf}}
\subfloat[\label{fig:r0}]{\includegraphics[width=3.5in]{hysteresispolarizationthreedimensionschange.pdf}}\\
\subfloat[\label{fig:r0}]{\includegraphics[width=3.5in]{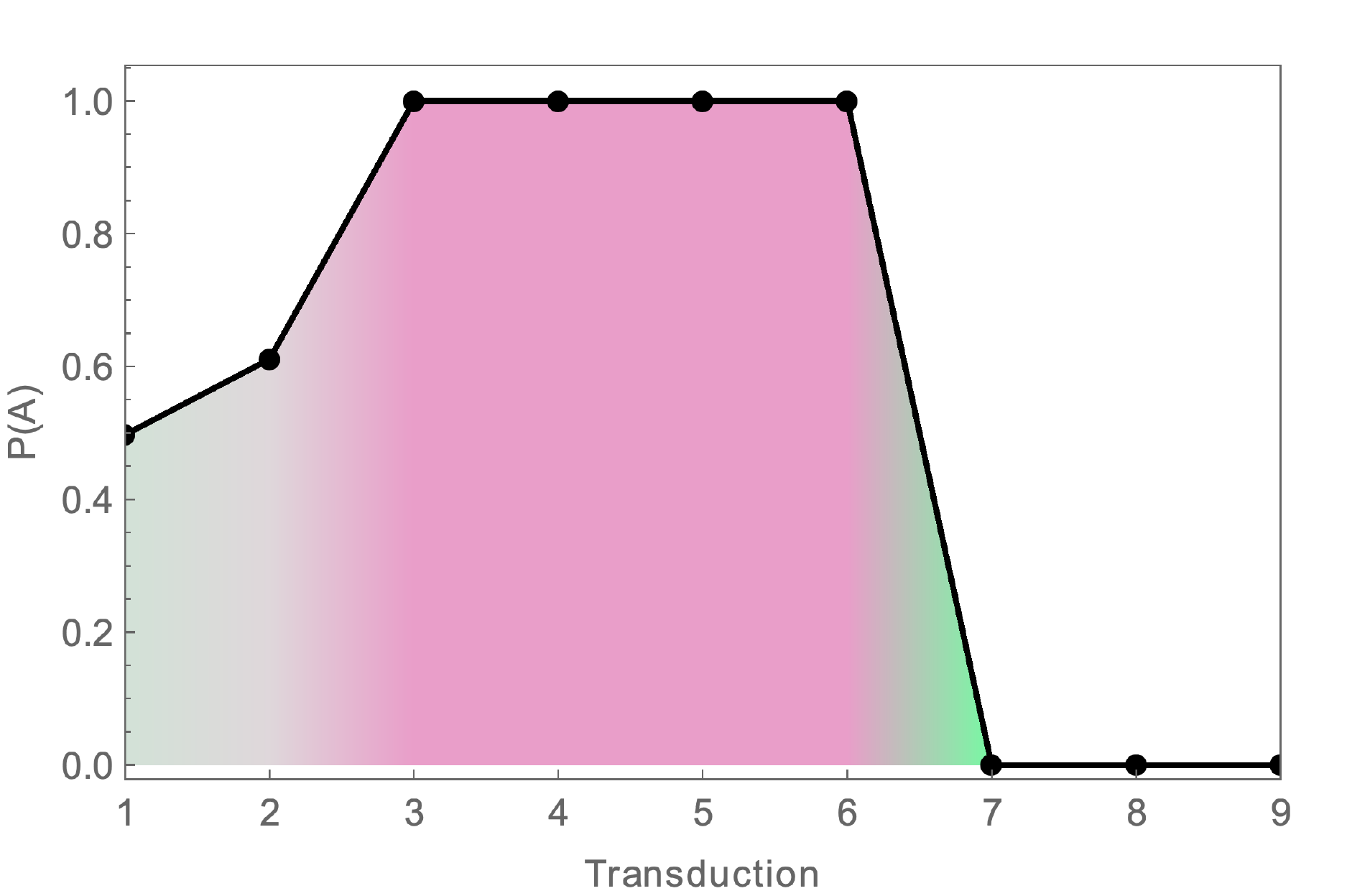}}
\subfloat[\label{fig:r0}]{\includegraphics[width=3.5in]{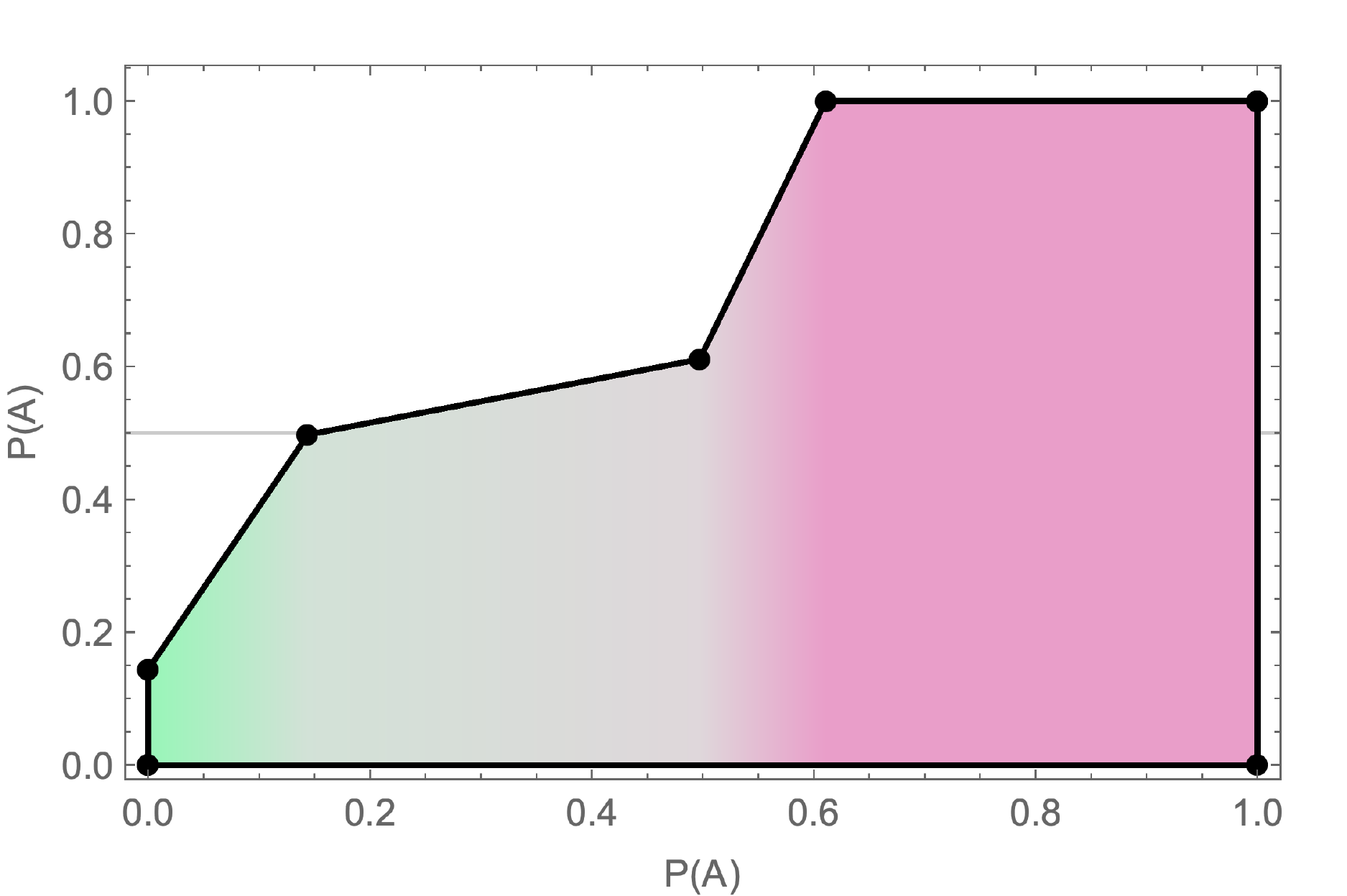}}\\
\subfloat[\label{fig:r0}]{\includegraphics[width=4.5in]{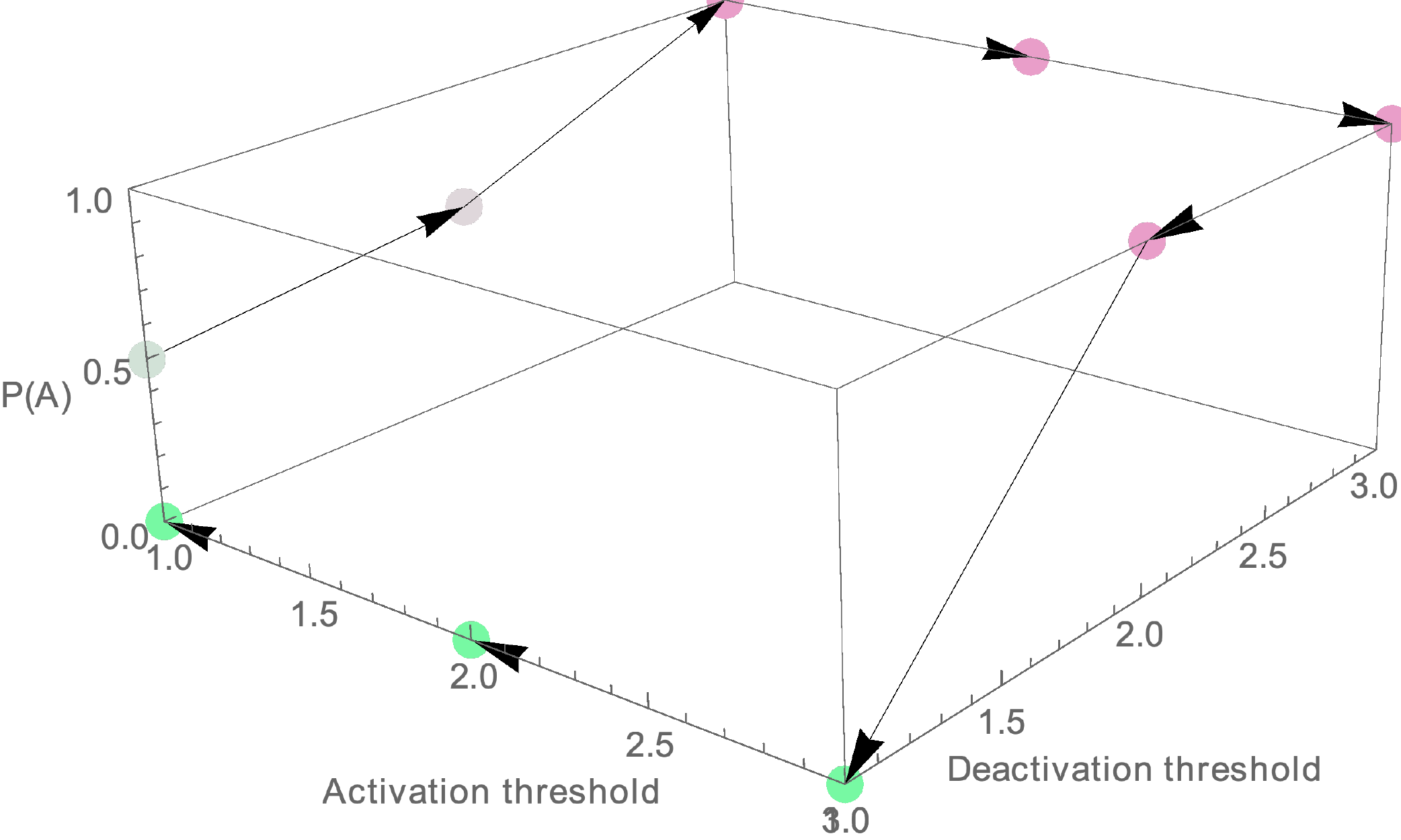}}
\endgroup
\caption{Hysteresis of polarization activation for activation-deactivation threshold sequence $((1,1),(1,2),(1,3),(2,3),(3,3),(3,2),(3,1),(2,1),(1,1))$ and self-deactivation; with $K\sim\text{Poisson}(c=50)$, the constant density kernel with $p=1/10$, and induction by initializing activation threshold $k=8$, yielding $\xi_0\{A\}\simeq0.1438$}\label{fig:orbitsdeact}
\end{figure}

\begin{figure}[h!]
\centering
\includegraphics[width=6.5in]{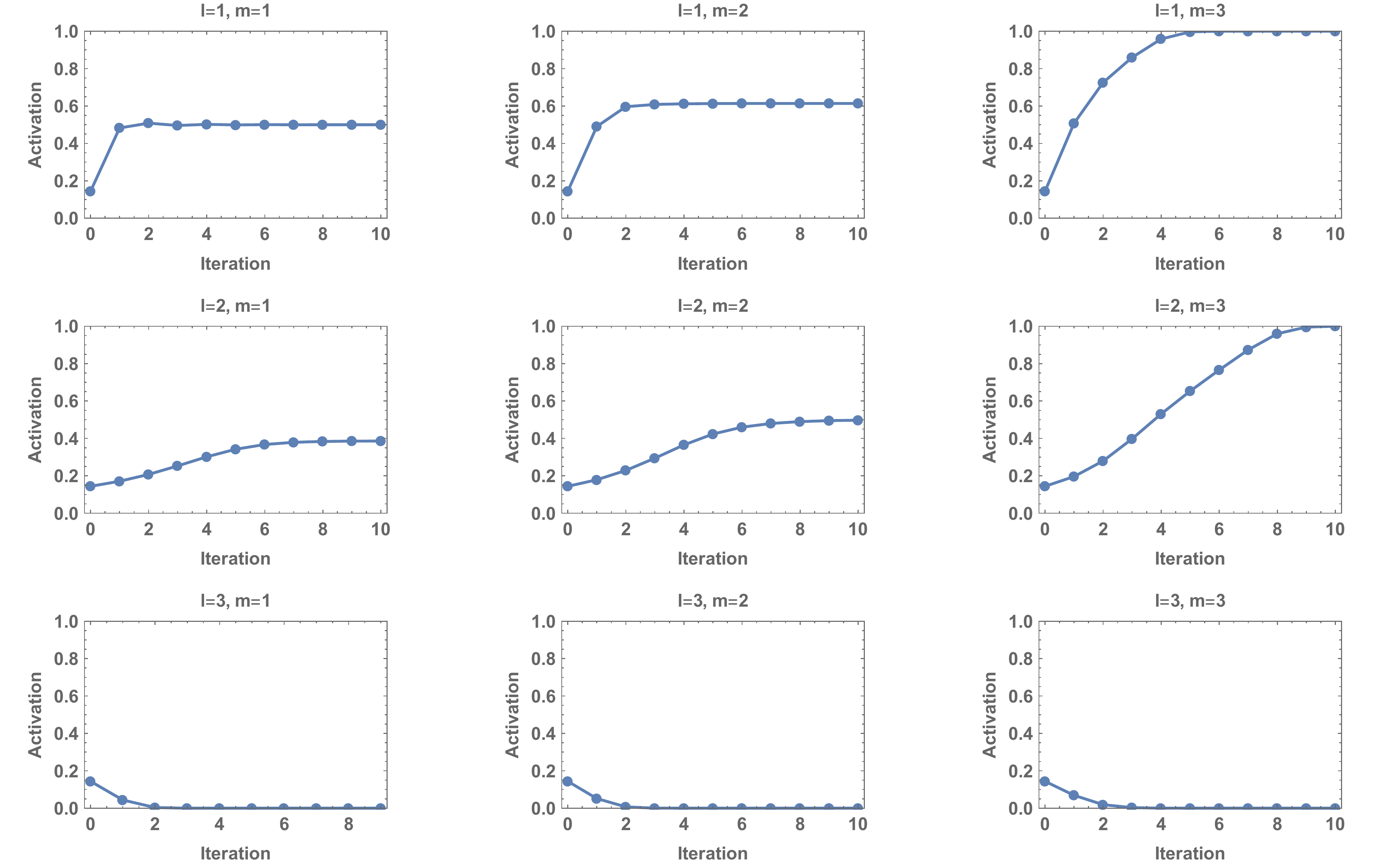}
\caption{Transduction polarization in activation $l$ and deactivation $m$ thresholds for number of points $K\sim\text{Poisson}(c=50)$, constant density kernel with $p=1/10$,  initializing activation threshold $k=8$, yielding $\xi_0\{A\}\simeq0.1438$}\label{fig:sifp}
\end{figure}
\FloatBarrier

\newpage
\subsection{Mixed intrinsic and external} The activation recursion is given by \[\xi_n\{A\}\equiv \xi_{n-1}\{A\} + (\nu(p_{n-1}^B)+s)(1-\xi_{n-1}\{A\}) - (\nu(q_{n-1}^B)+r)\xi_{n-1}\{A\}\for n\ge 1\] where $r,s\in[0,1]$ chosen such that $\nu(q_{n-1}^B)+r\le 1$ and  $\nu(p_{n-1}^B)+s\le 1$ hold for all $\xi\{A\}$. This regime can admit periodic orbits. 

\begin{figure}[h!]
\centering
\includegraphics[width=6.5in]{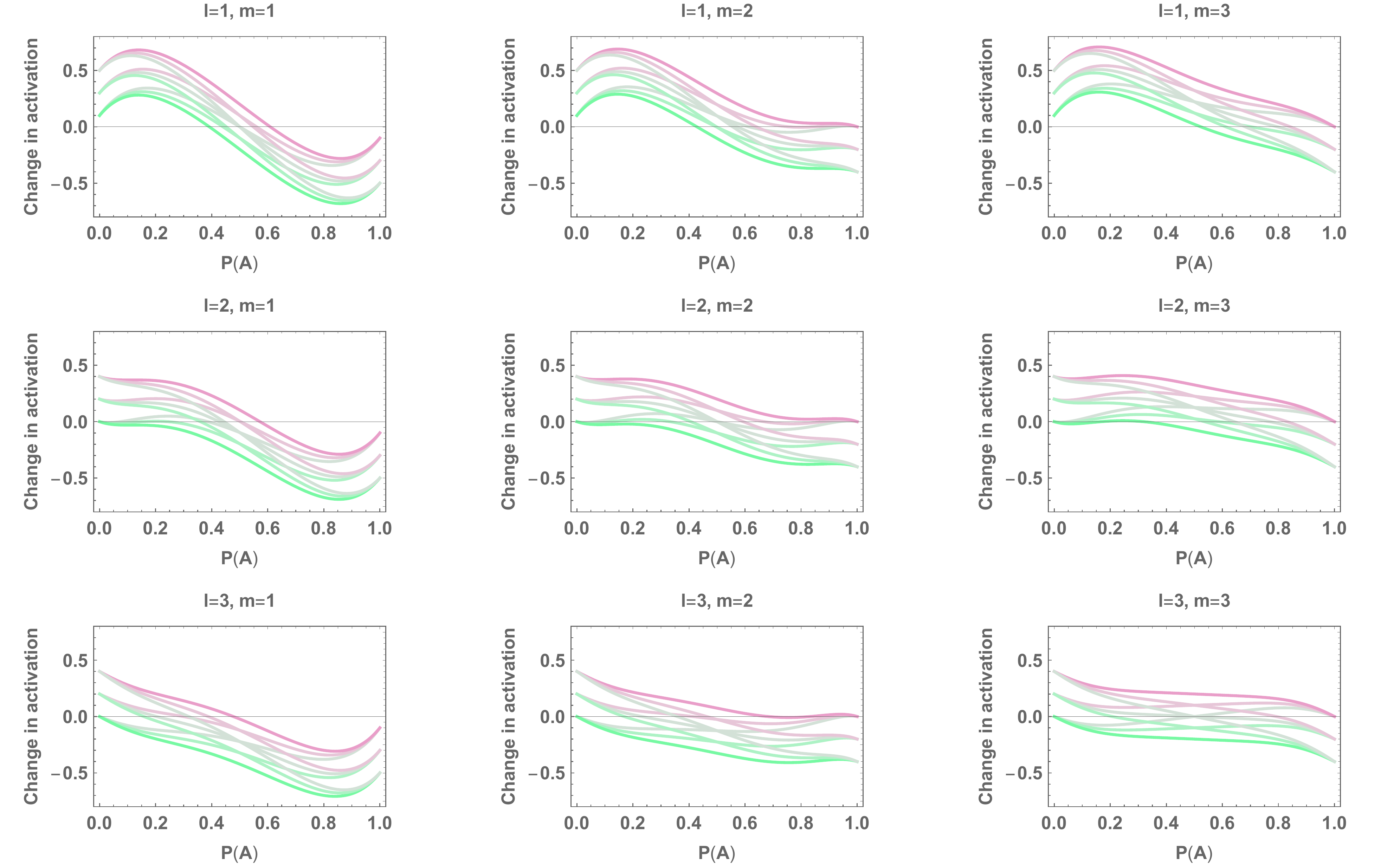}
\caption{Polarization net activation, in activation-deactivation thresholds $(l,m)\in\{1,2,3\}\times\{1,2,3\}$, deactivation probability $r$, and activation probability $s$, for $K\sim\text{Poisson}(c=50)$ and the constant density kernel with $p=1/10$}\label{fig:sifp}
\end{figure}
\FloatBarrier

\begin{figure}[h!]
\centering
\includegraphics[width=6.5in]{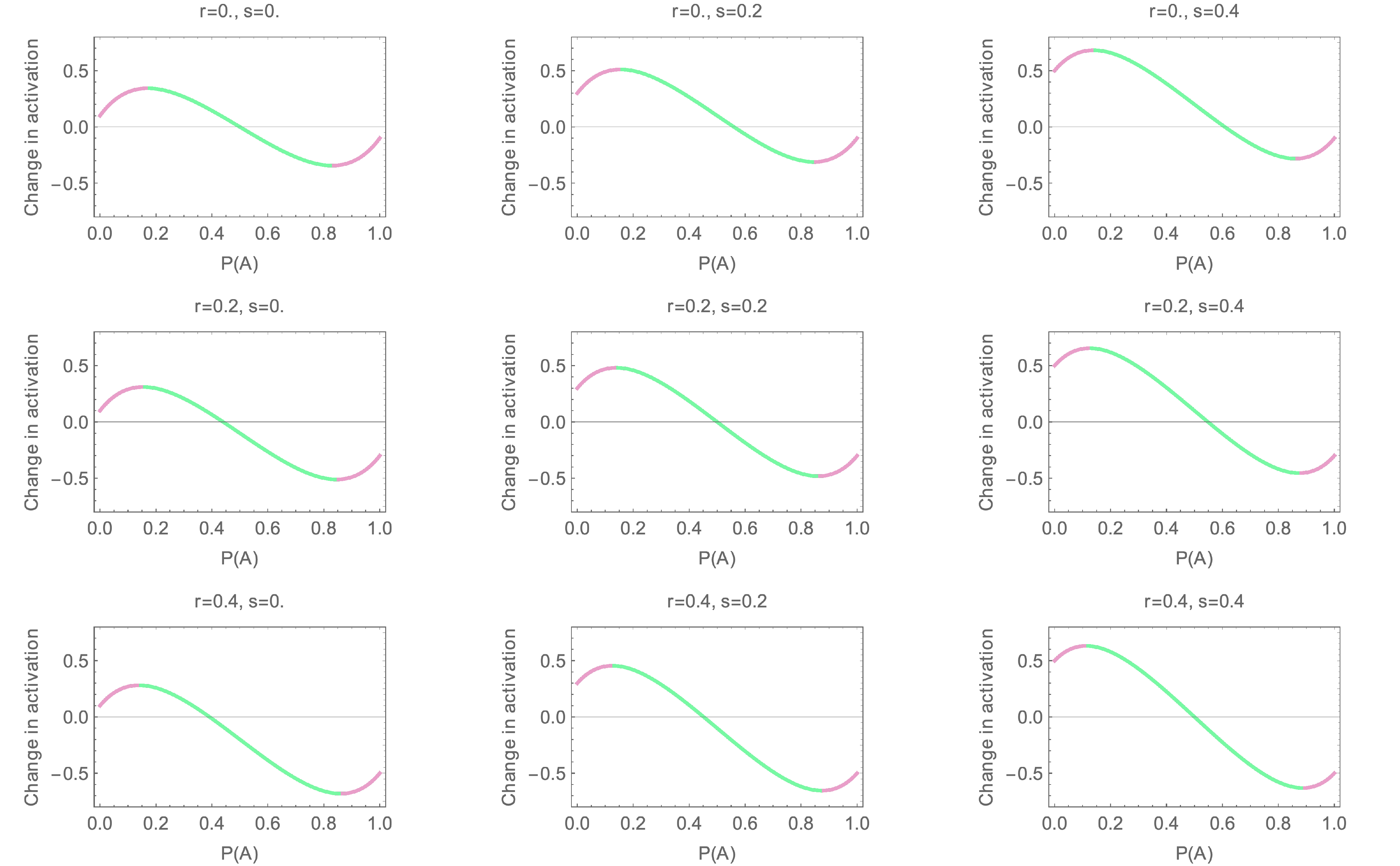}
\caption{Polarization net activation in deactivation $r$ and activation $s$ probabilities for activation-deactivation thresholds $(l,m)=(1,1)$, $K\sim\text{Poisson}(c=50)$ and the constant density kernel with $p=1/10$}\label{fig:sifp}
\end{figure}
\FloatBarrier

\begin{figure}[h!]
\centering
\includegraphics[width=6.5in]{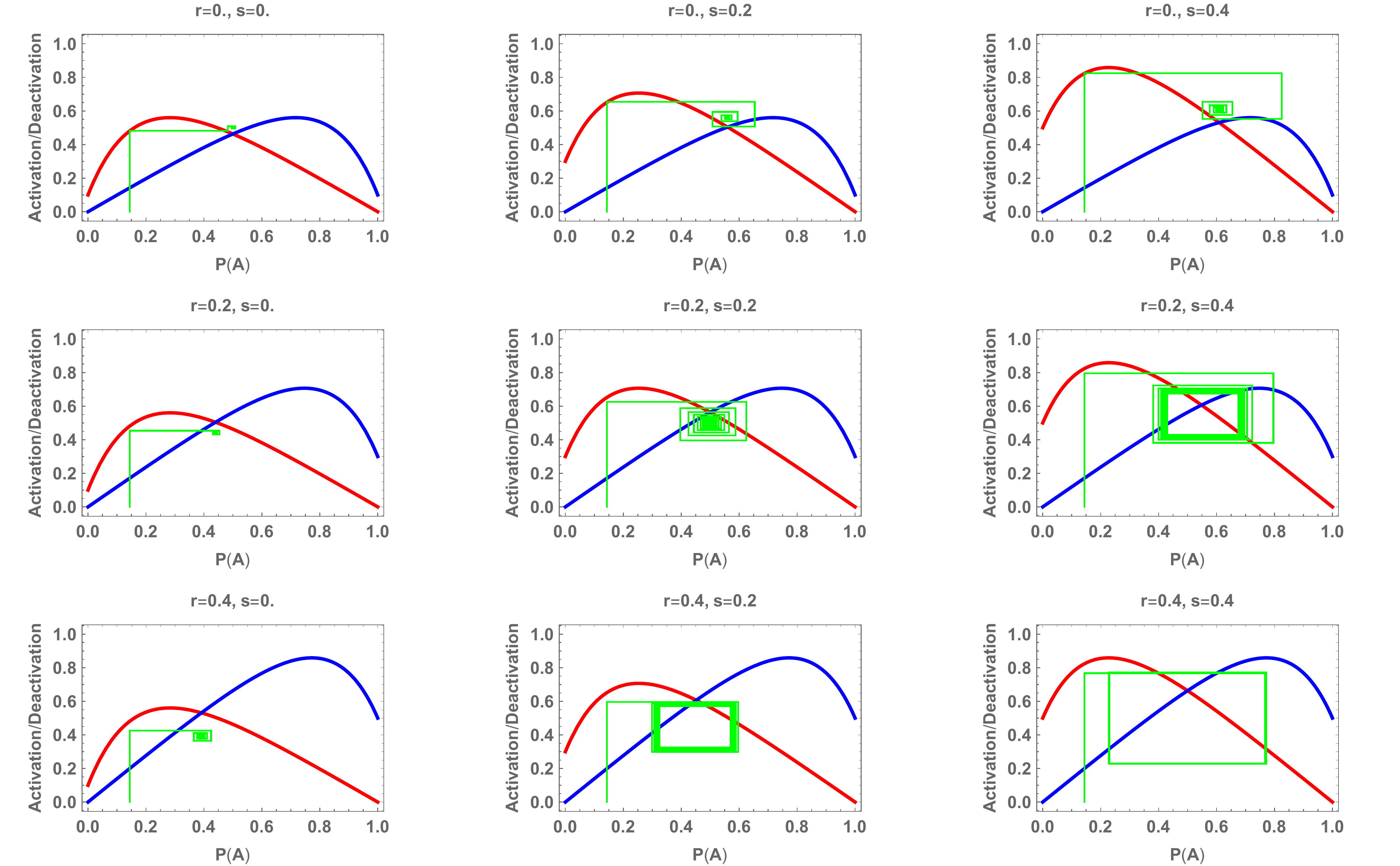}
\caption{Polarization net activation in deactivation $r$ and activation $s$ probabilities for activation-deactivation thresholds $(l,m)=(1,1)$, $K\sim\text{Poisson}(c=50)$ and the constant density kernel with $p=1/10$}\label{fig:sifp}
\end{figure}
\FloatBarrier

\begin{figure}[h!]
\centering
\includegraphics[width=6.5in]{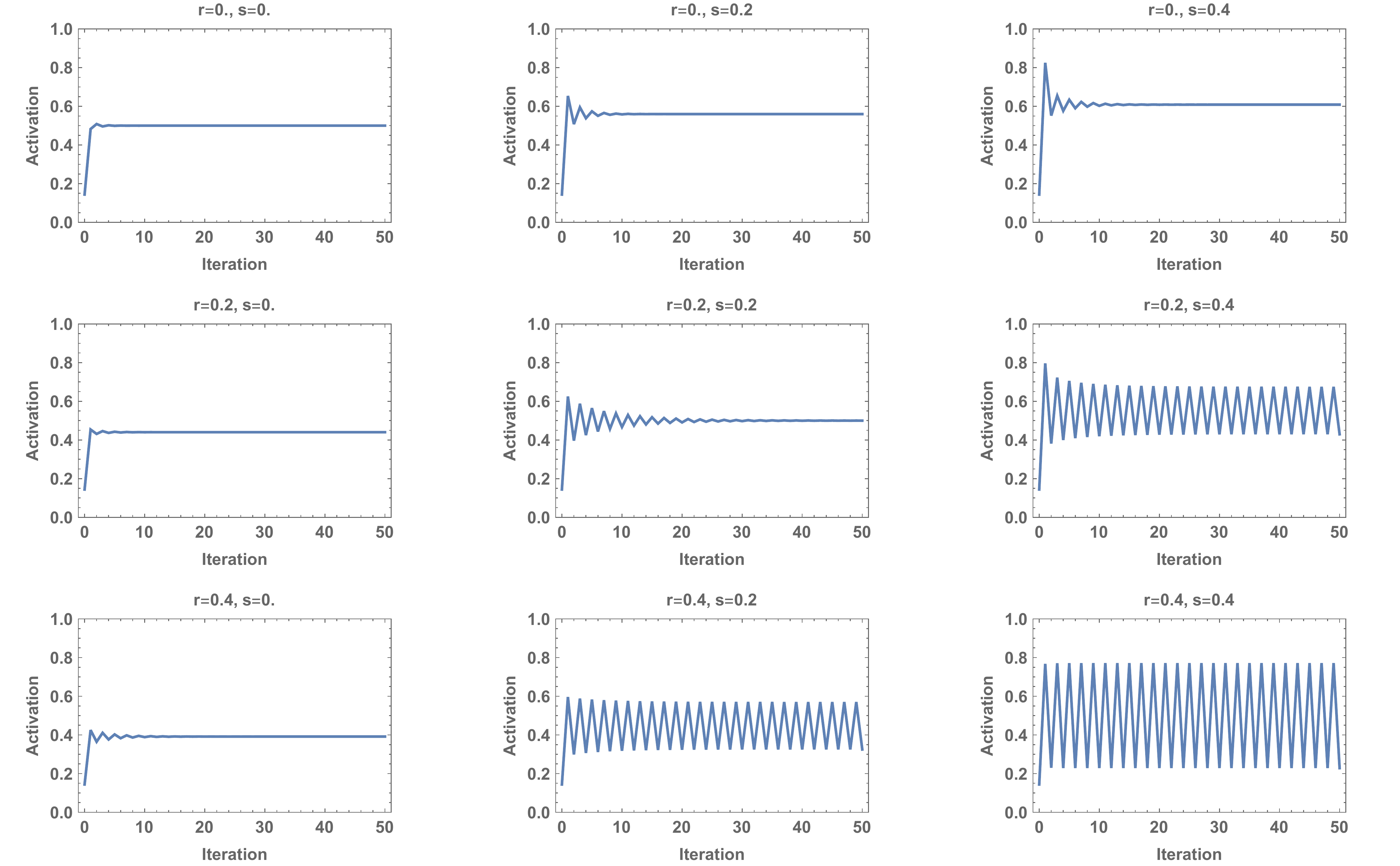}
\caption{Polarization and periodicity in free activation $s$ and deactivation $r$ for activation and deactivation thresholds $l=m=1$, number of points $K\sim\text{Poisson}(c=50)$, constant density kernel with $p=1/10$, initializing activation threshold $k=8$, yielding $\xi_0\{A\}\simeq0.1438$}\label{fig:sifp}
\end{figure}
\FloatBarrier

\newpage
\section{Decentral transduction of chaos by activation annihilation and superfluity} Chaos is transduced with complete deactivation (annihilation of activation; $r=1$) and superfluous activation (renormalization).

\subsection{Transduction} Chaotic orbits of transduction are defined by \[\xi_n\{A\} \equiv \wp\,\nu(p_{n-1}^A)(1-\xi_{n-1}\{A\})\for n\ge 1\] where $\wp>1$ is the chaotic renormalizing (`chaos') factor \[\wp^{-1} \equiv \sup_{\xi\{A\}>0}\nu(p_{0}^B)(1-\xi\{A\})<1\] Recall that $\nu(p_{0}^B)$ depends on $\xi\{A\}$. Indicating $\varepsilon \equiv \wp-1$ the recursion \[\xi_n\{A\} \equiv (1+\varepsilon)\nu(p_{n-1}^B)(1-\xi_{n-1}\{A\})\for n\ge 1\] shows that the origination of chaos is the additional (superfluous) source of activation probability caused by renormalization. That is, for $\varepsilon=0$, there is no chaos.  

For Poisson $K$, constant density kernel $p$, and transduction threshold $l=1$, the chaos factor is given by 
\[\wp = \frac{c p W\left(\frac{e^{cp+1}}{1-p}\right)}{\left(W\left(\frac{e^{cp+1}}{1-p}\right)-1\right)^2}\] where $W$ is the product-logarithm function. 

\begin{figure}[h!]
\centering
\includegraphics[width=4.5in]{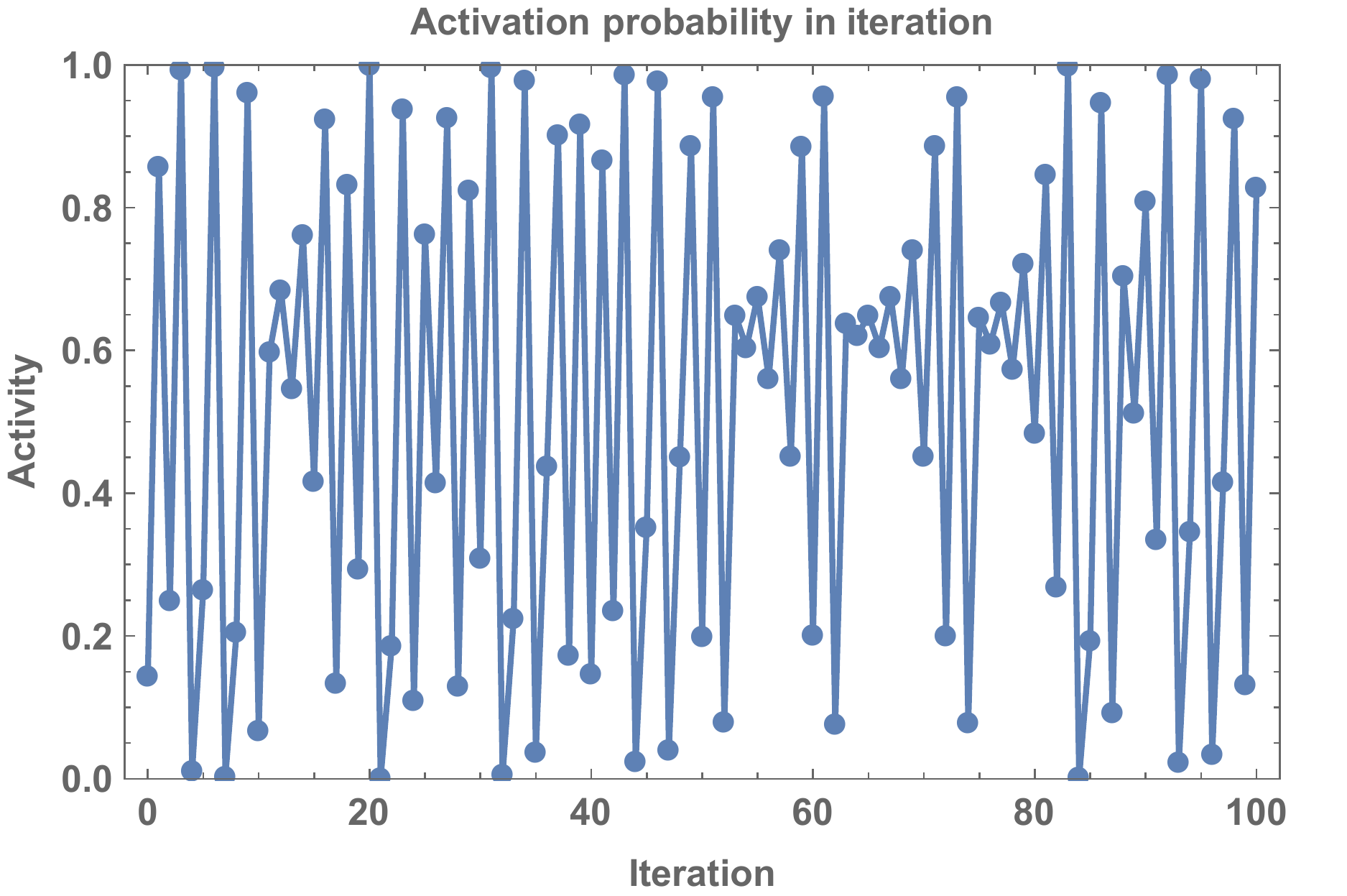}
\caption{Chaos in the forward orbits of transduction for number of points $K\sim\text{Poisson}(c=50)$, constant density kernel with $p=1/10$, $\xi_0\{A\}\simeq0.1438$, initializing activation threshold $k=8$, transduction threshold $l=1$, and renormalization factor $\wp\simeq 1.7849$}\label{fig:chaosport}
\end{figure}

In Figure~\ref{fig:chaosport}. We plot the forward orbits of transduction for the running example, with transduction threshold $l=1$ and renormalization factor $\wp\simeq 1.7849$. They are chaotic.

\begin{figure}[h!]
\centering
\begingroup
\captionsetup[subfigure]{width=5in,font=normalsize}
\subfloat[Superfluous activation for activation threshold $l\in\{1,2,3\}$\label{fig:r0}]{\includegraphics[width=4in]{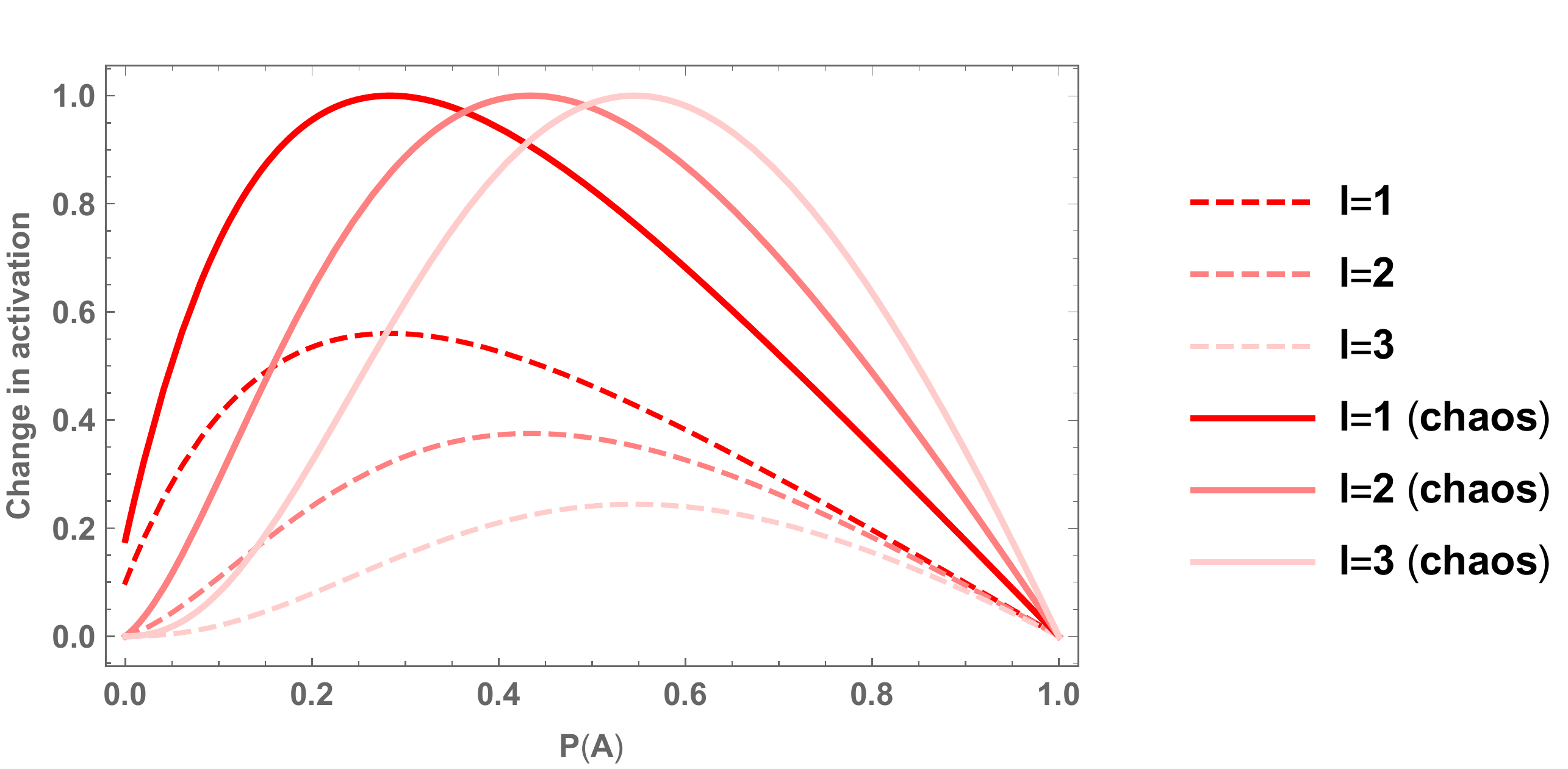}}
\subfloat[Activation annihilation \label{fig:r1}]{\includegraphics[width=3in]{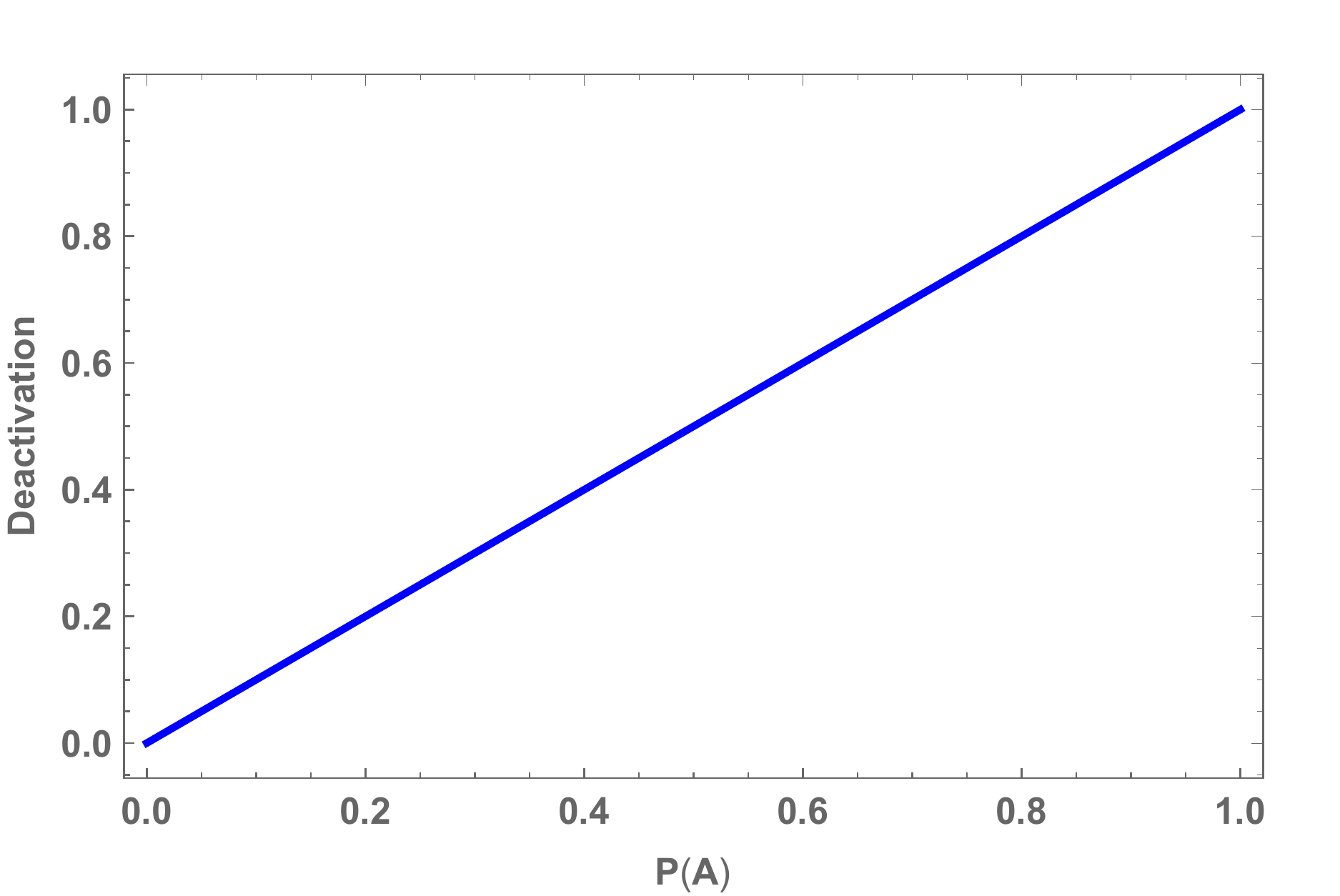}}\\
\subfloat[Net activation for activation threshold $l\in\{1,2,3\}$\label{fig:r2}]{\includegraphics[width=3in]{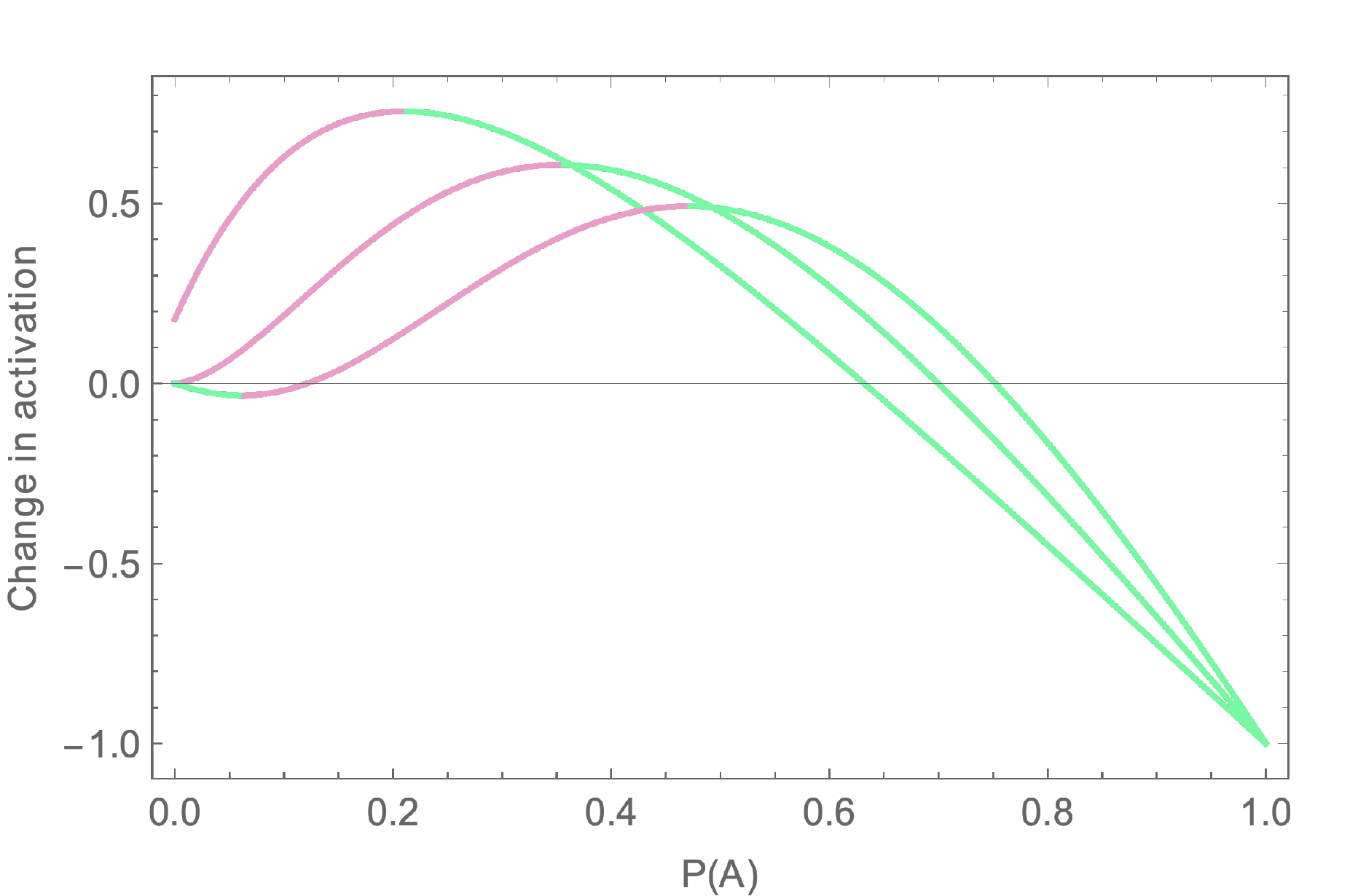}}\\
\endgroup
\caption{Chaos of superfluous activation, activation annihilation, and net activation in activation threshold $l$ for $K\sim\text{Poisson}(c=50)$ and the constant density kernel with $p=1/10$}\label{fig:orbitsdeact}
\end{figure}

\begin{figure}[h!]
\centering
\begingroup
\captionsetup[subfigure]{width=5in,font=normalsize}
\subfloat[$l=1$\label{fig:r0}]{\includegraphics[width=3in]{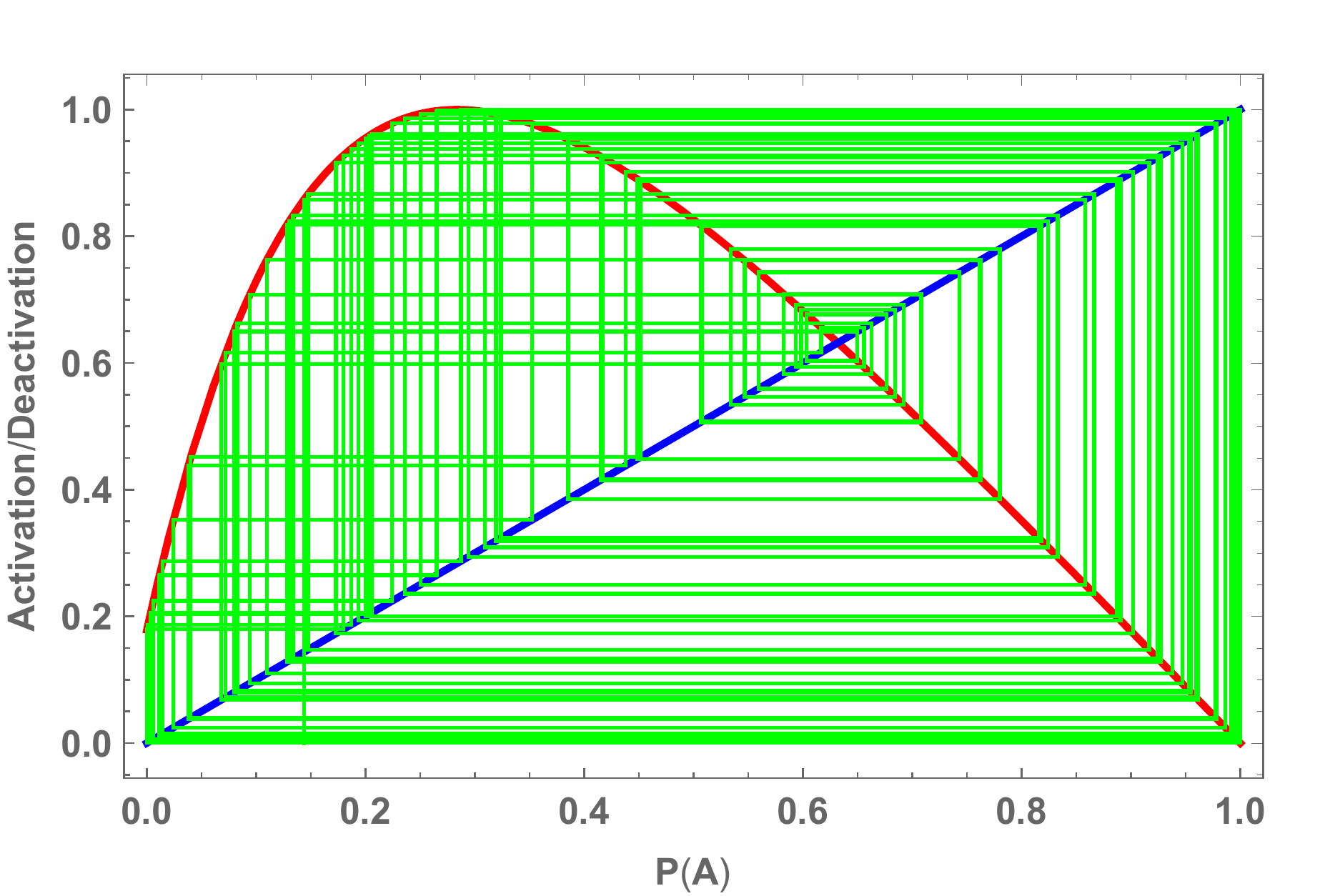}}
\subfloat[$l=2$\label{fig:r1}]{\includegraphics[width=3in]{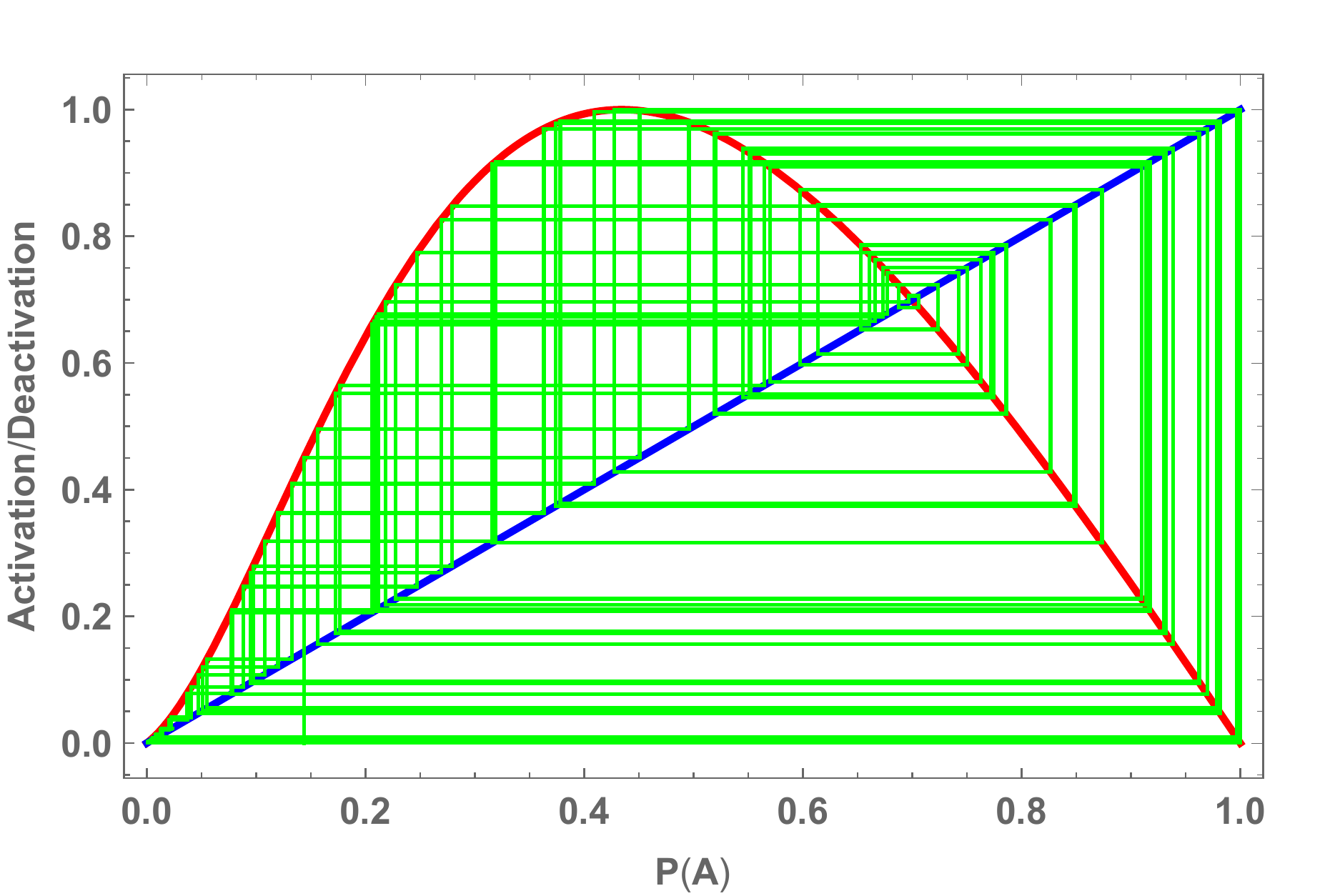}}\\
\subfloat[$l=3$\label{fig:r2}]{\includegraphics[width=3in]{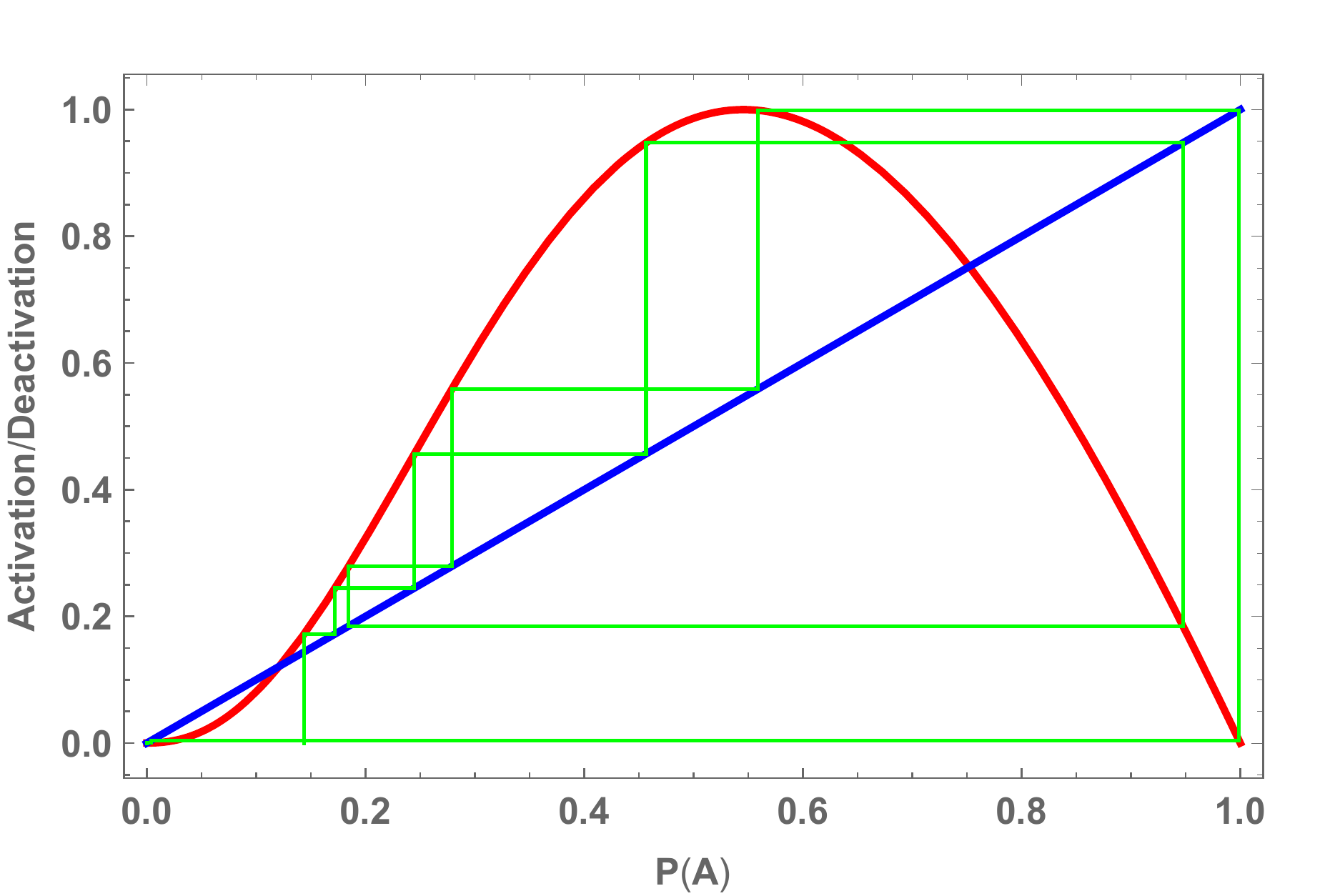}}
\subfloat[Logistic map $\beta=4$\label{fig:r2}]{\includegraphics[width=3in]{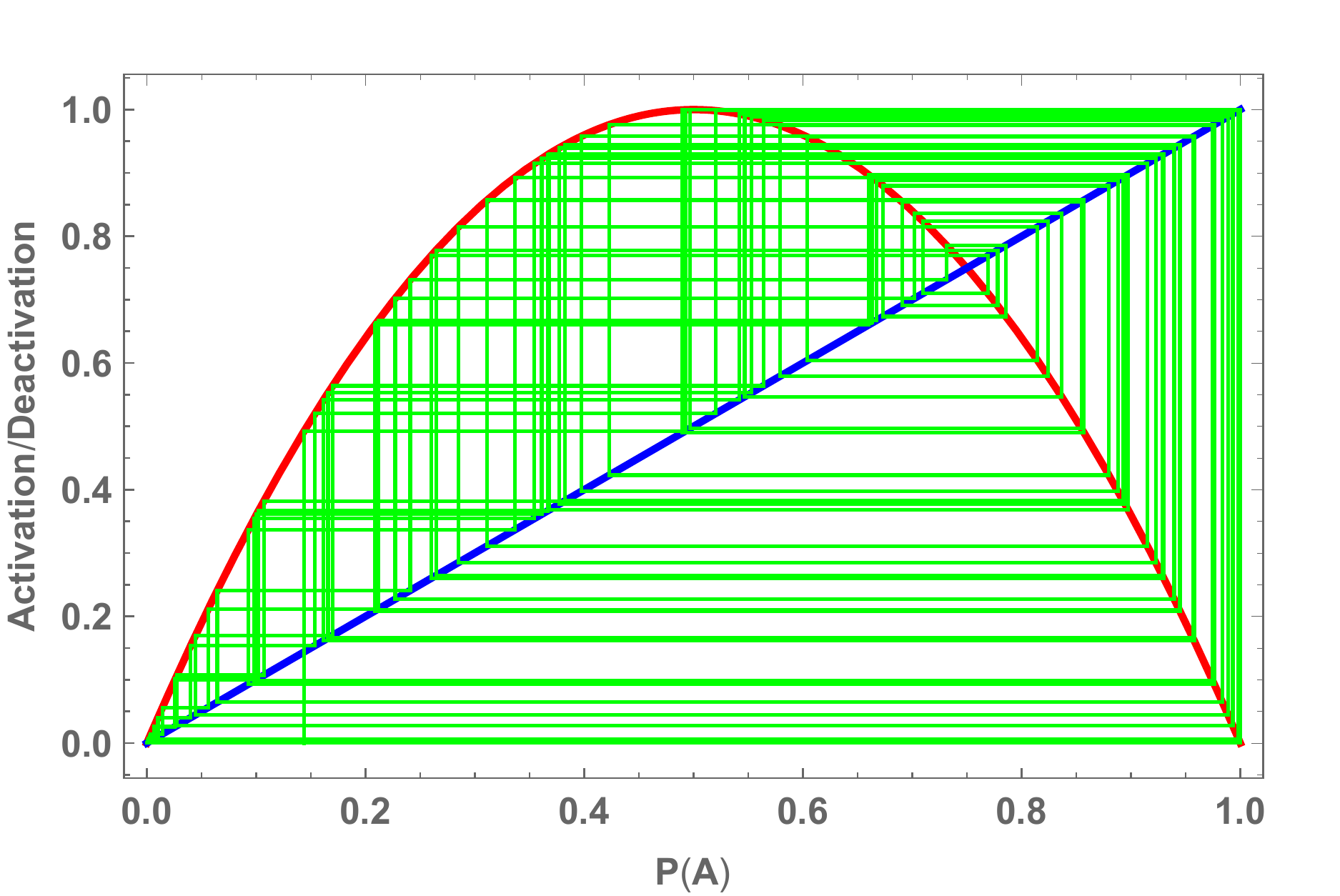}}
\endgroup
\caption{Cobweb diagrams in activation threshold $l$ for 100 iterations, initial condition $\xi_0\{A\}\simeq0.1438$, $K\sim\text{Poisson}(c=50)$ and the constant density kernel with $p=1/10$}\label{fig:orbitsdeact}
\end{figure}

\FloatBarrier

\subsection{Special cases} The susceptible-infected-susceptible (SIS) compartmental model for spread of infectious disease (activation) with recovery (deactivation) is given by recursion relation \[\xi_n\{A\} \equiv \xi_{n-1}\{A\} + \beta\xi_{n-1}\{A\}(1-\xi_{n-1}\{A\}) - r\xi_{n-1}\{A\}\for n\ge 1\] with parameters $\beta$ and $r$ in $[0,1]$. The ITAD mean-field reduces to the SIS compartmental model for a single point $K\sim\text{Dirac}(1)$, the constant density kernel with $\beta\in[0,1]$, $\nu(p_{n-1}^B)=p_{n-1}^B=\P(d_B(D)\ge 1)=1-\psi_D(0)=1-\psi(1-\beta\xi\{A\})=\beta\xi\{A\}$ where $d_B$ has pgf $\psi_D(t)=\psi(1-\beta\xi\{A\}+\beta\xi\{A\}t)=1-\beta\xi\{A\}+\beta\xi\{A\}t$. 

 The SIS compartmental model is the logistic map for $r=1$ and $\beta\in[0,4]$ \[\xi_n\{A\} \equiv \beta\xi_{n-1}\{A\}(1-\xi_{n-1}\{A\})\for n\ge 1\] The logistic map is chaotic for $\beta=4$. 
 
 \clearpage
\newpage
\section{Local induction-transduction activation-deactivation} For these figures, we take $E=[0,1]$, $\nu=\Leb$, and the density kernel as $f(x,y)=q\ind{}(|x-y|\le r)$ on $E\times E$, where we set $q=1/2$ and choose $r$ such that \[(\nu\times\nu)f = \frac{1}{2}\int_{[0,1]^2}\D x\D y \ind{}(|x-y|\le r)= p\] and use the standard ITAD field equation, i.e. $C_1=0$. 

\begin{figure}[h!]
\centering
\begingroup
\captionsetup[subfigure]{width=5in,font=normalsize}
\subfloat[Activation with threshold $l\in\{1,2,3\}$\label{fig:r0}]{\includegraphics[width=3.5in]{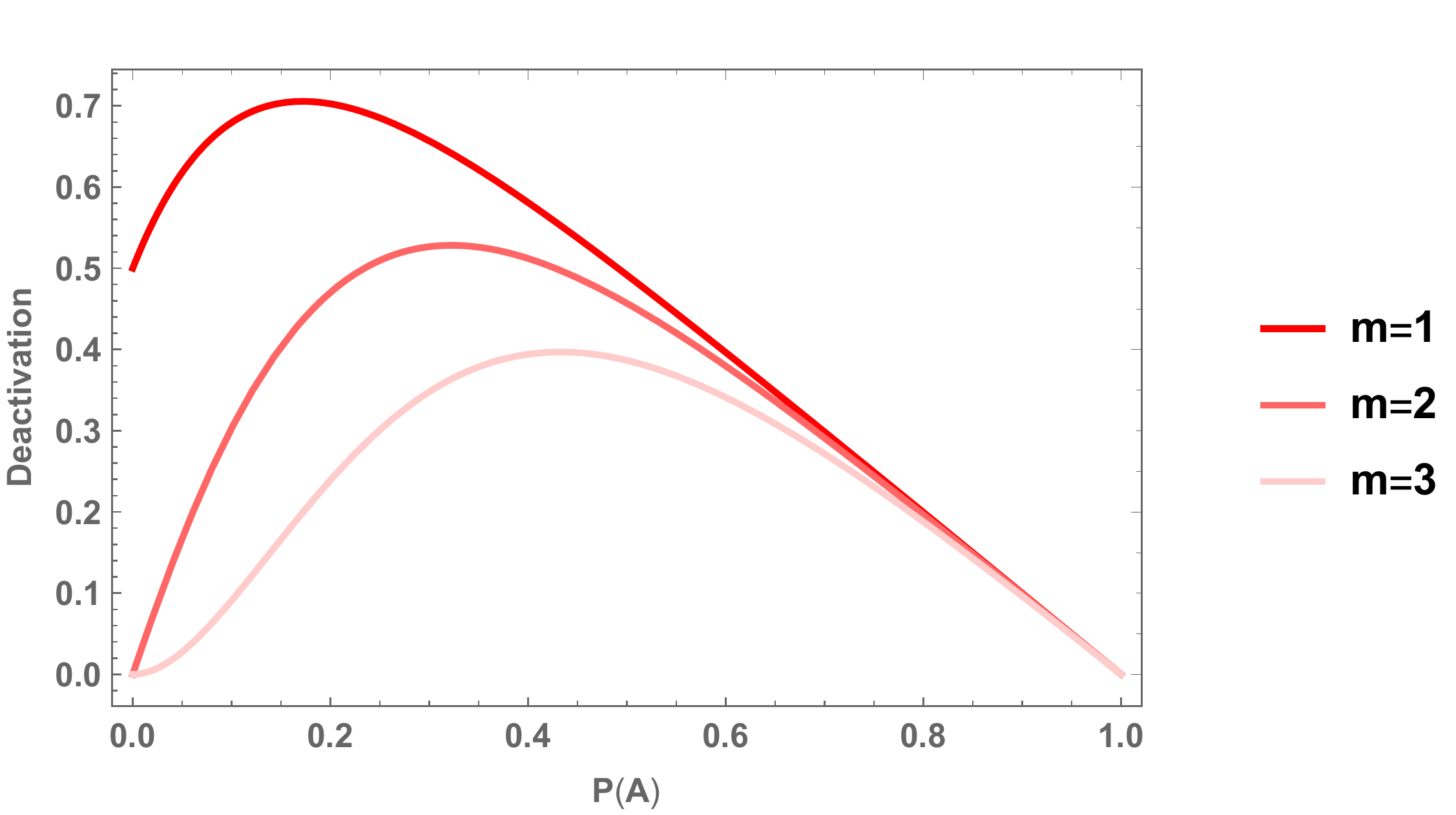}}
\subfloat[Deactivation with threshold $m\in\{1,2,3\}$\label{fig:r1}]{\includegraphics[width=3.5in]{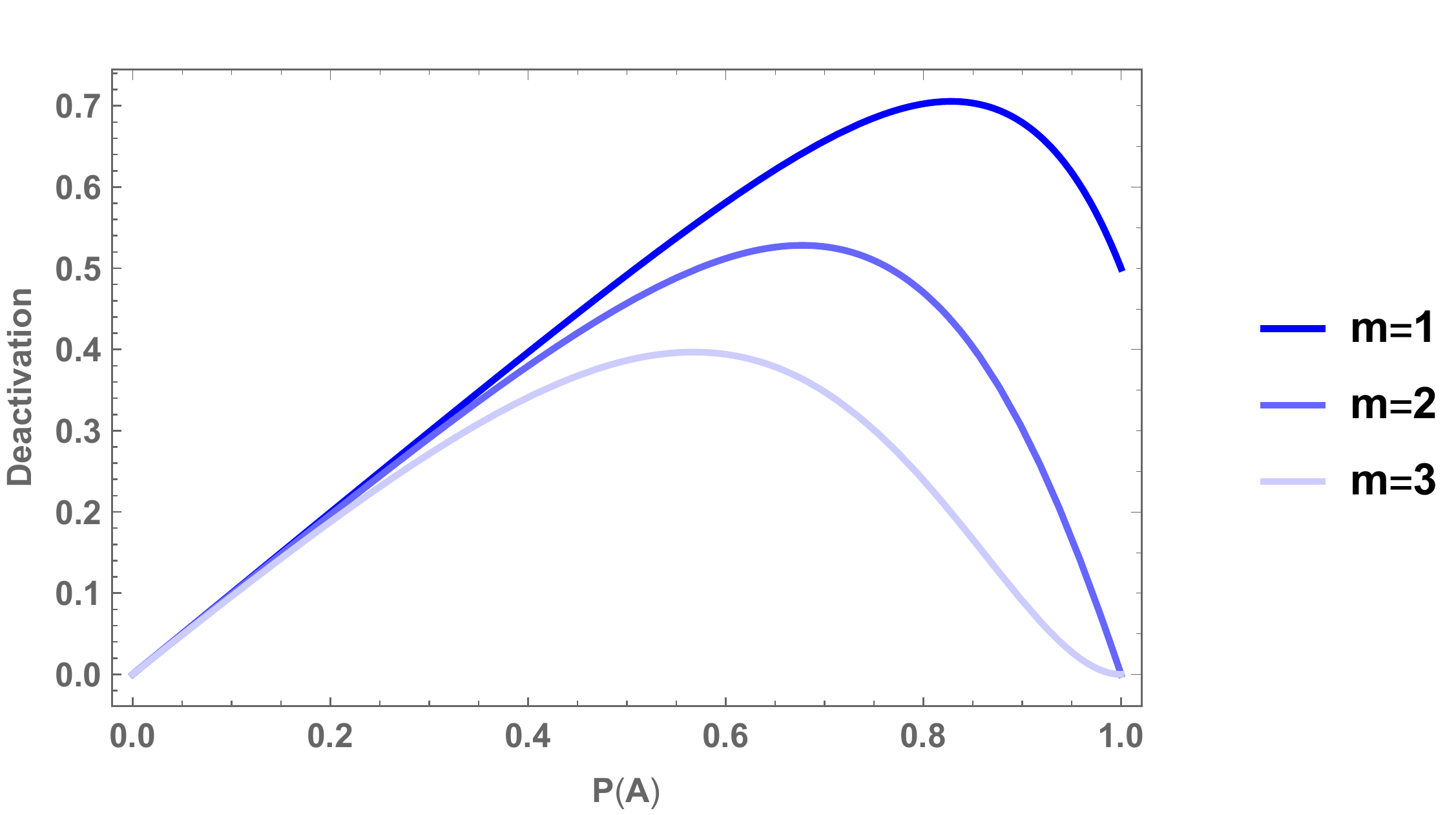}}\\
\subfloat[Net activation for $(l,m)\in\{1,2,3\}\times\{1,2,3\}$\label{fig:r2}]{\includegraphics[width=6.5in]{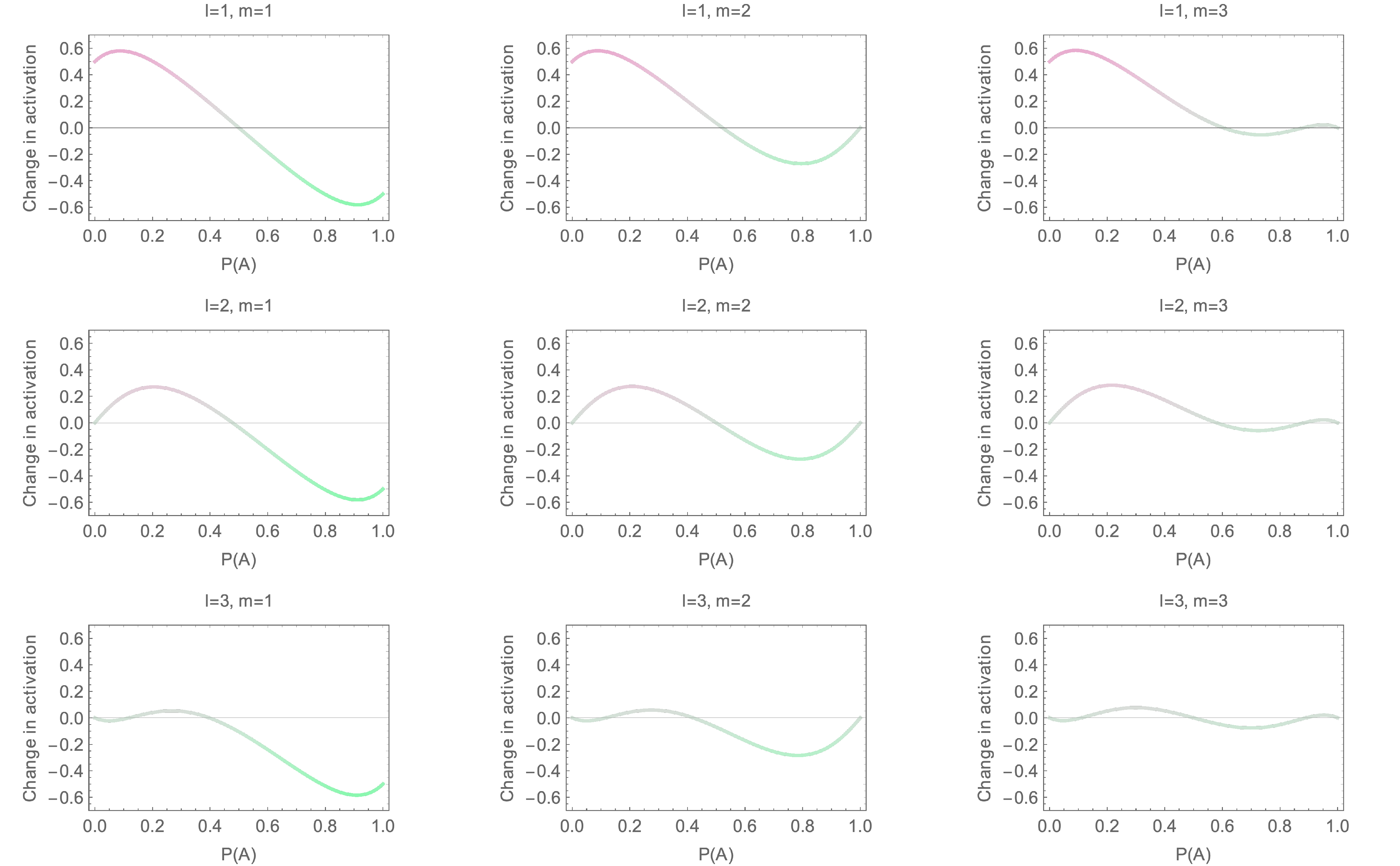}}
\endgroup
\caption{Local transduction  in activation, deactivation, and net activation, in activation threshold $l$ and deactivation probability $r$ (blue is zero, red is one), for $K\sim\text{Poisson}(c=50)$ and kernel $q\ind{}(|x-y|\le r)$ for $q=0.5$ and $r$ such that mass is $p=1/10$}\label{fig:orbitsdeact}
\end{figure}

\begin{figure}[h!]
\centering
\includegraphics[width=6.5in]{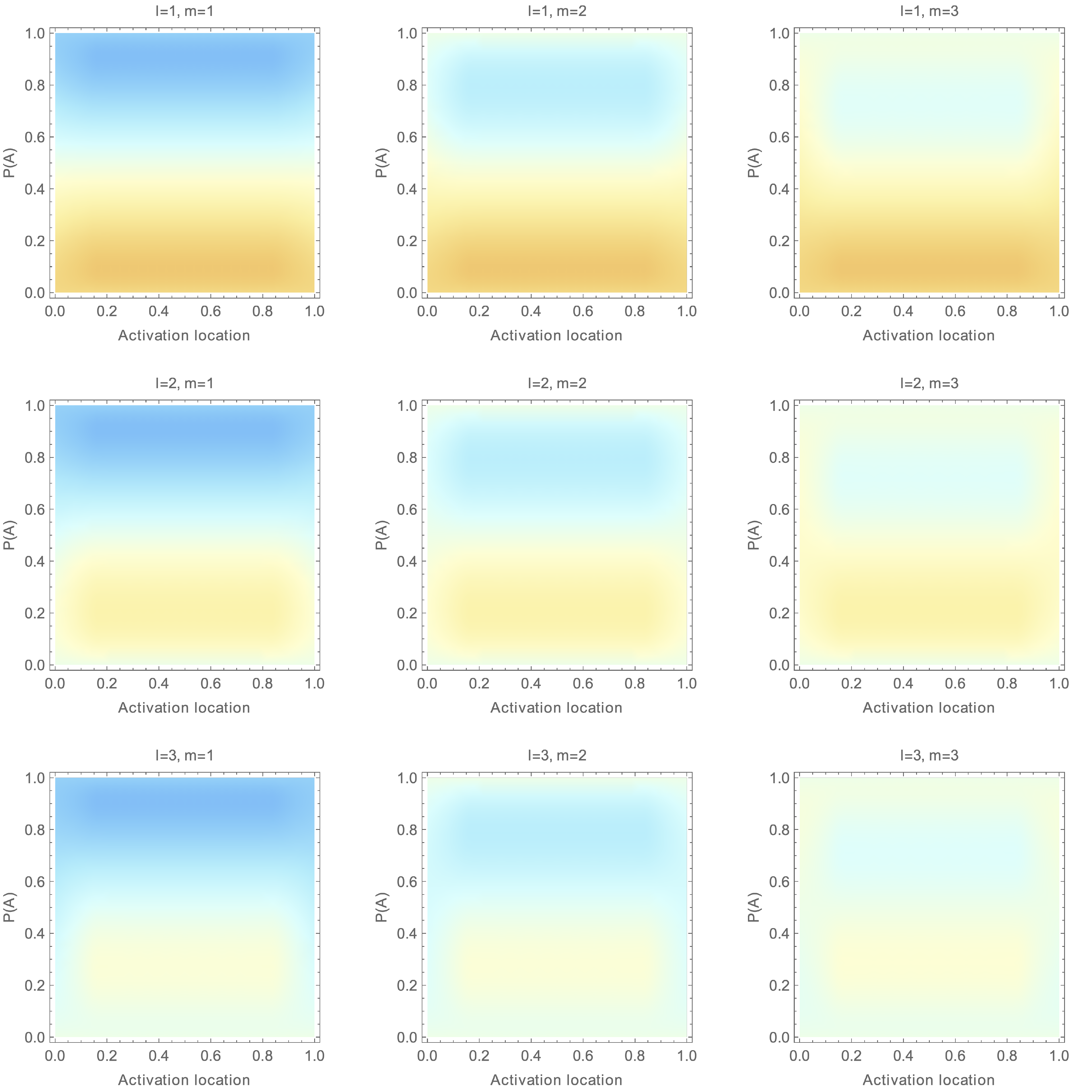}
\caption{Local transduction  in activation $l$ and deactivation $m$ thresholds for number of points $K\sim\text{Poisson}(c=50)$, kernel $q\ind{}(|x-y|\le r)$ for $q=0.5$ and $r$ such that mass is $p=1/10$}\label{fig:sifp}
\end{figure}
\FloatBarrier

\begin{figure}[h!]
\centering
\includegraphics[width=6.5in]{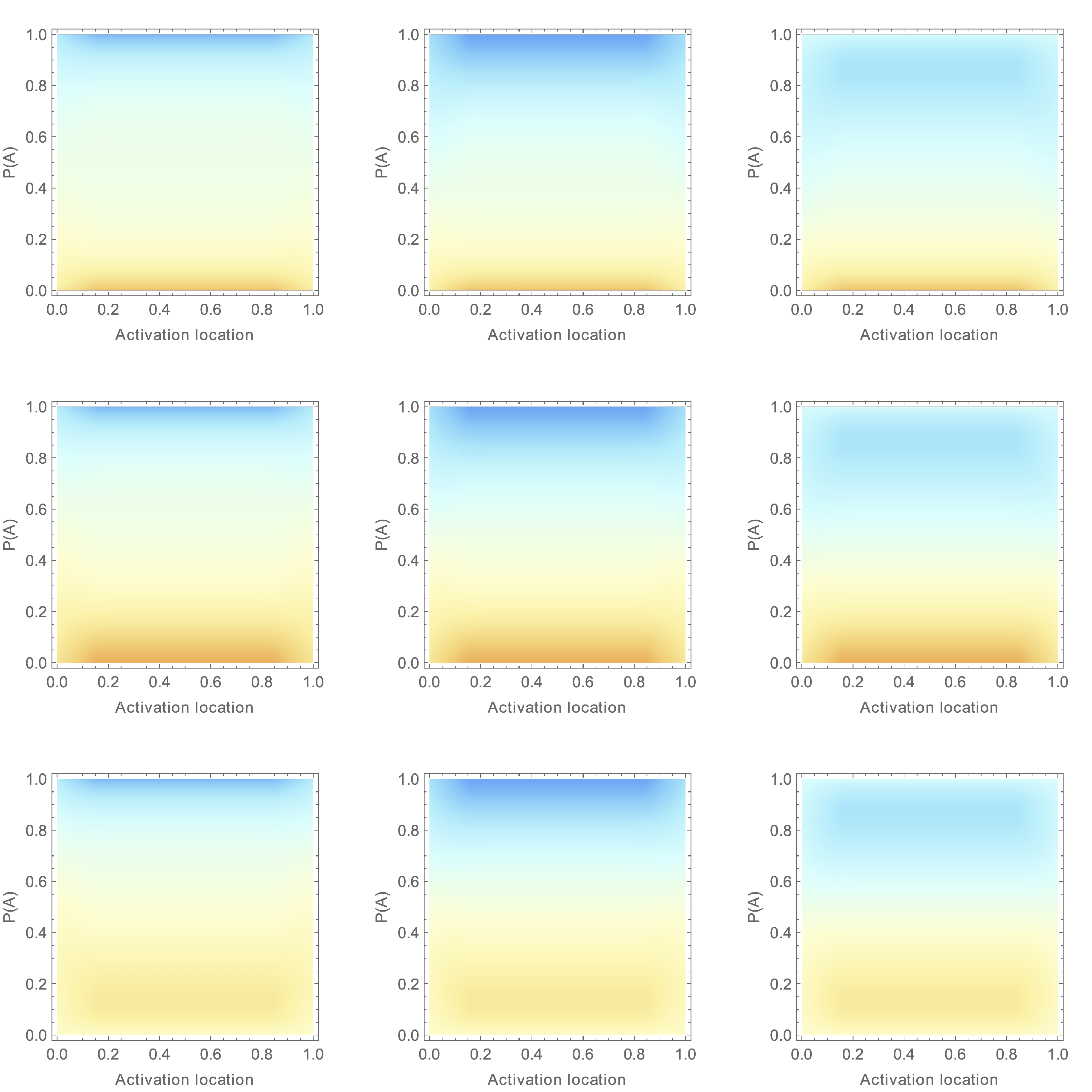}
\caption{Local transduction activation flux relative to activation  in activation $l$ and deactivation $m$ thresholds for number of points $K\sim\text{Poisson}(c=50)$, density kernel mass $p=1/10$, initializing activation threshold $k=8$, yielding $\xi_0\{A\}\simeq0.1438$}\label{fig:sifp}
\end{figure}
\FloatBarrier

\begin{figure}[h!]
\centering
\includegraphics[width=6.5in]{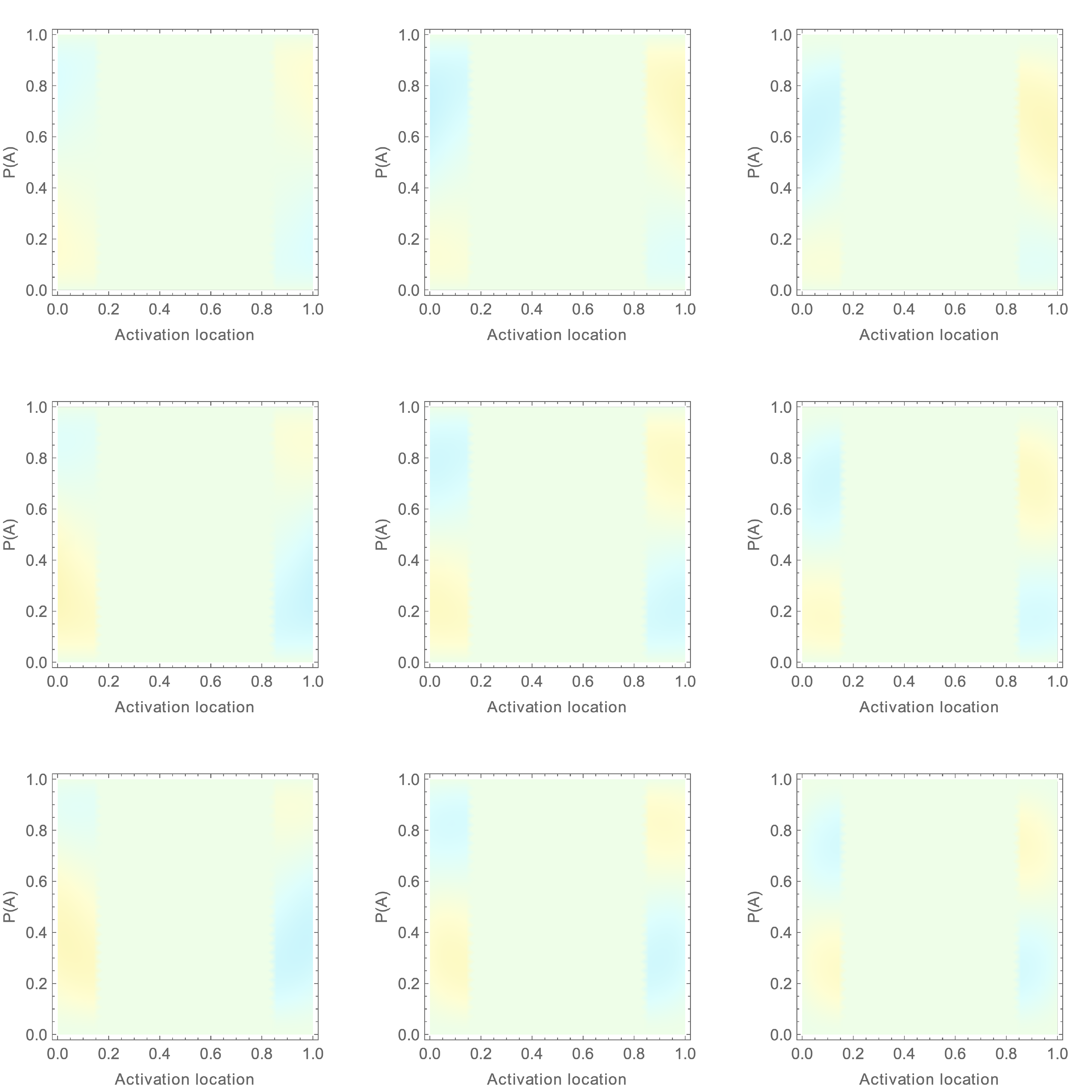}
\caption{Local transduction activation flux relative to location  in activation $l$ and deactivation $m$ thresholds for number of points $K\sim\text{Poisson}(c=50)$, density kernel mass $p=1/10$, initializing activation threshold $k=8$, yielding $\xi_0\{A\}\simeq0.1438$}\label{fig:sifp}
\end{figure}
\FloatBarrier

\begin{figure}[h!]
\centering
\includegraphics[width=6.5in]{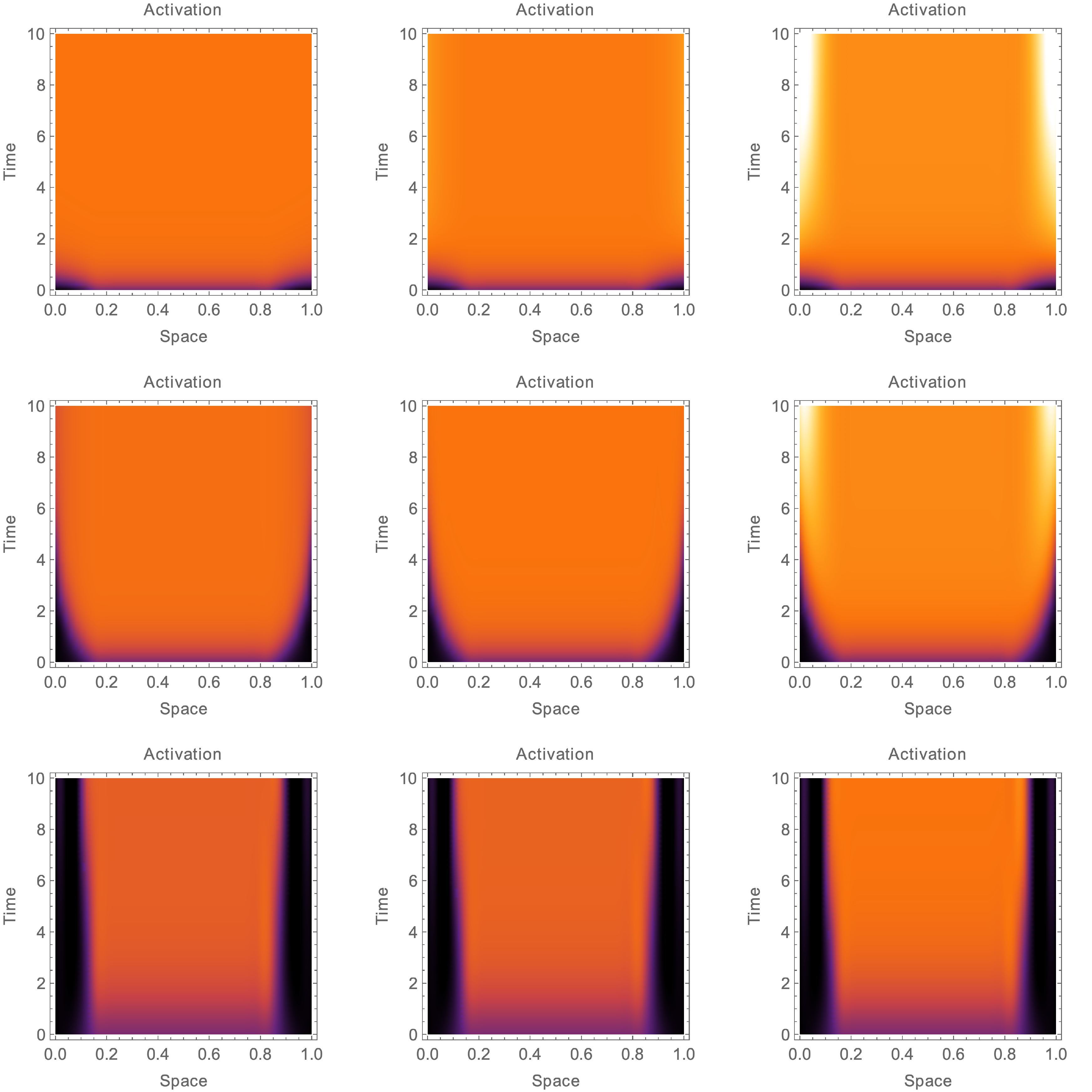}
\caption{Local transduction activation-deactivation probability in location  in activation $l$ and deactivation $m$ thresholds for number of points $K\sim\text{Poisson}(c=50)$, initializing activation threshold $k=14$}\label{fig:sifp}
\end{figure}
\FloatBarrier

\begin{figure}[h!]
\centering
\includegraphics[width=6.5in]{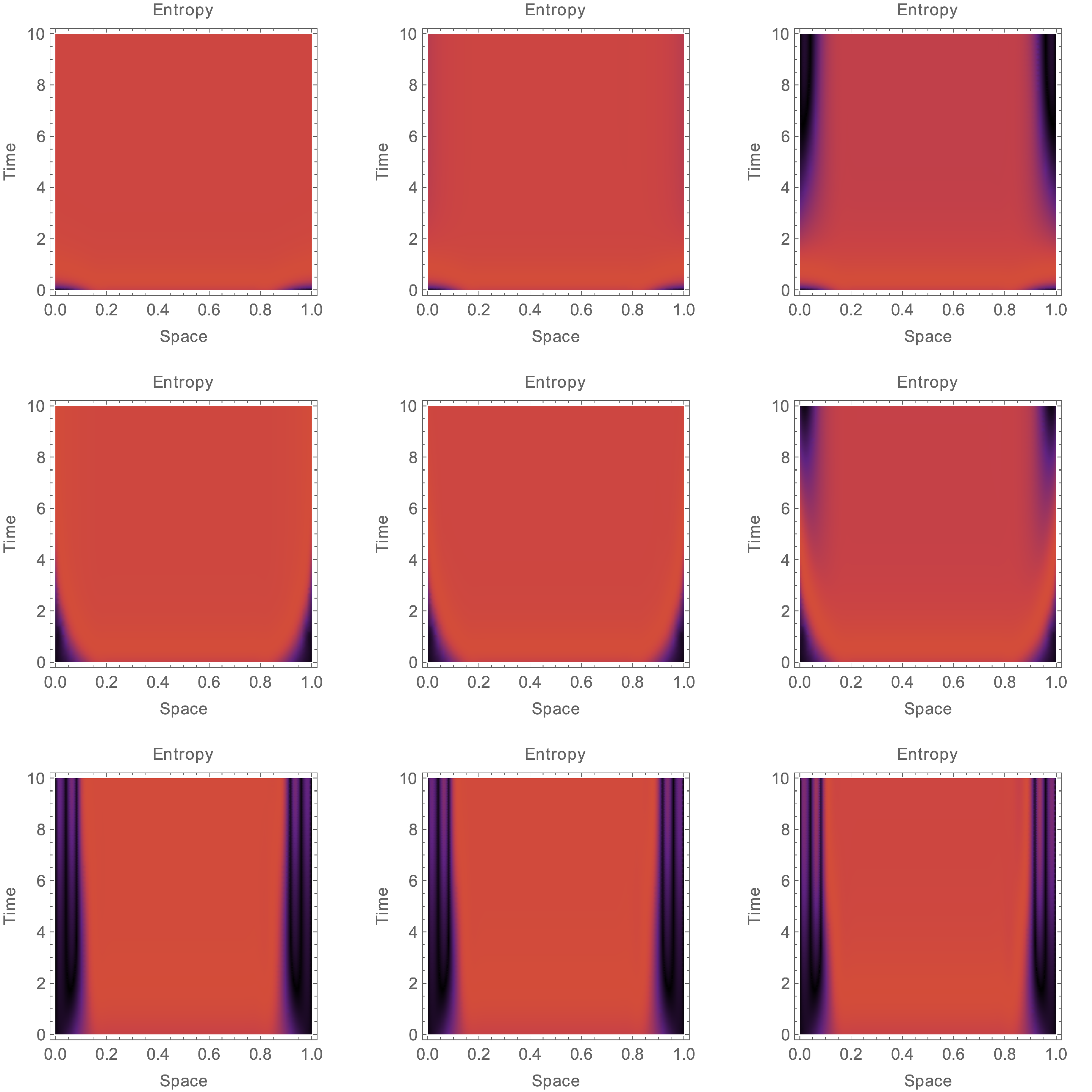}
\caption{Local transduction activation-deactivation entropy in location  in activation $l$ and deactivation $m$ thresholds for number of points $K\sim\text{Poisson}(c=50)$, initializing activation threshold $k=14$}\label{fig:sifp}
\end{figure}
\FloatBarrier

\begin{figure}[h!]
\centering
\includegraphics[width=6.5in]{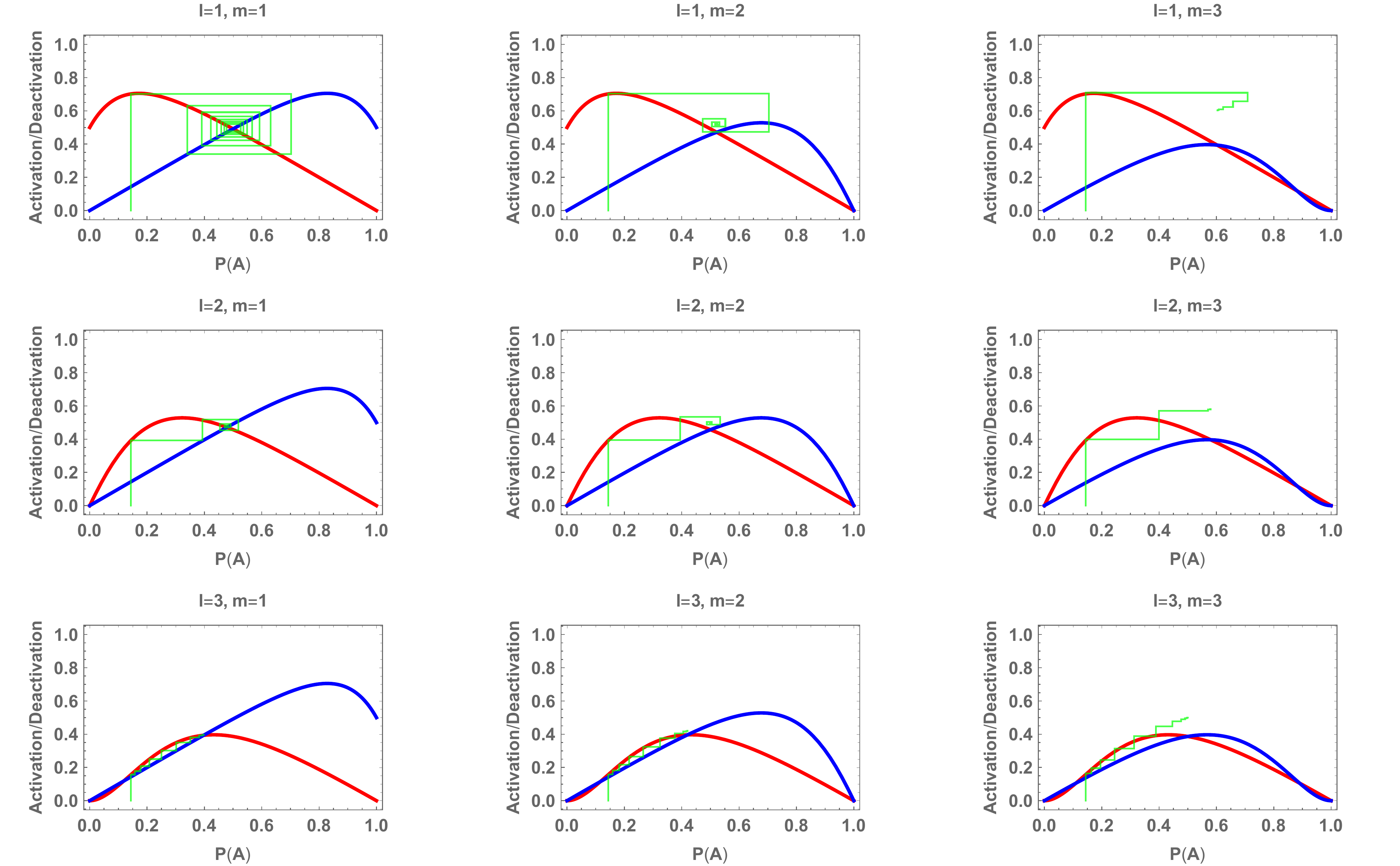}
\caption{Local transduction  in activation $l$ and deactivation $m$ thresholds for number of points $K\sim\text{Poisson}(c=50)$, kernel $q\ind{}(|x-y|\le r)$ for $q=0.5$ and $r$ such that mass is $p=1/10$, initial $\xi_0\{A\}\simeq0.1438$}\label{fig:sifp}
\end{figure}
\FloatBarrier

\clearpage
\section*{Local induction-transduction activation ($l=1,m=\infty$)} For these figures, we assume the standard ITAD field equation ($C_1=0$) and the induced wave and energy equation dimensionalized constants to be unity.

\begin{figure}[h!]
\centering
\includegraphics[width=5in]{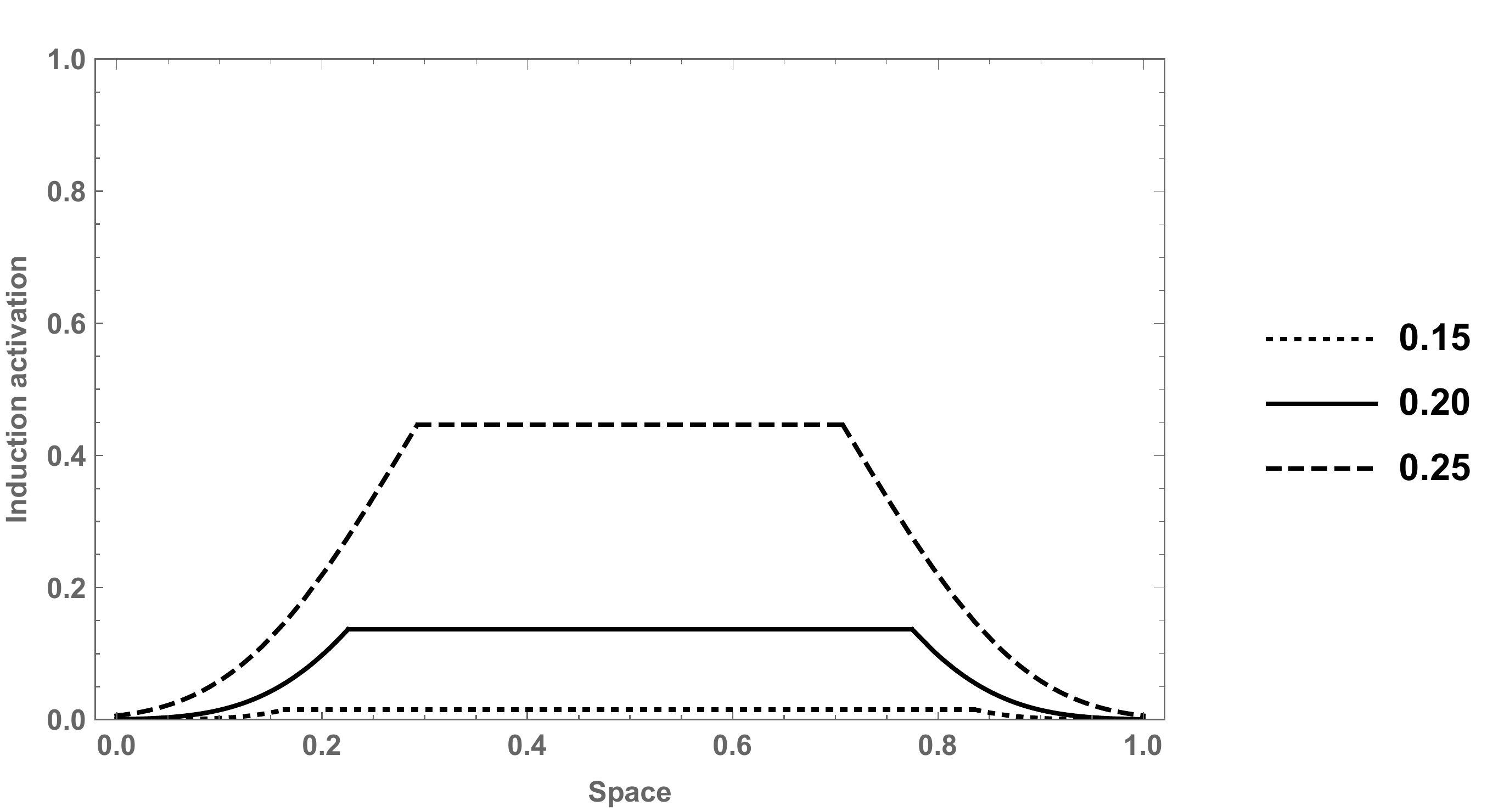}
\caption{Local induction activation for density kernel masses $0.15,0.2,0.25$ and induction activation threshold $k=15$ for $K\sim\text{Poisson}(c=50)$}\label{fig:nn}
\end{figure}

\begin{figure}[h!]
\centering
\includegraphics[width=7in]{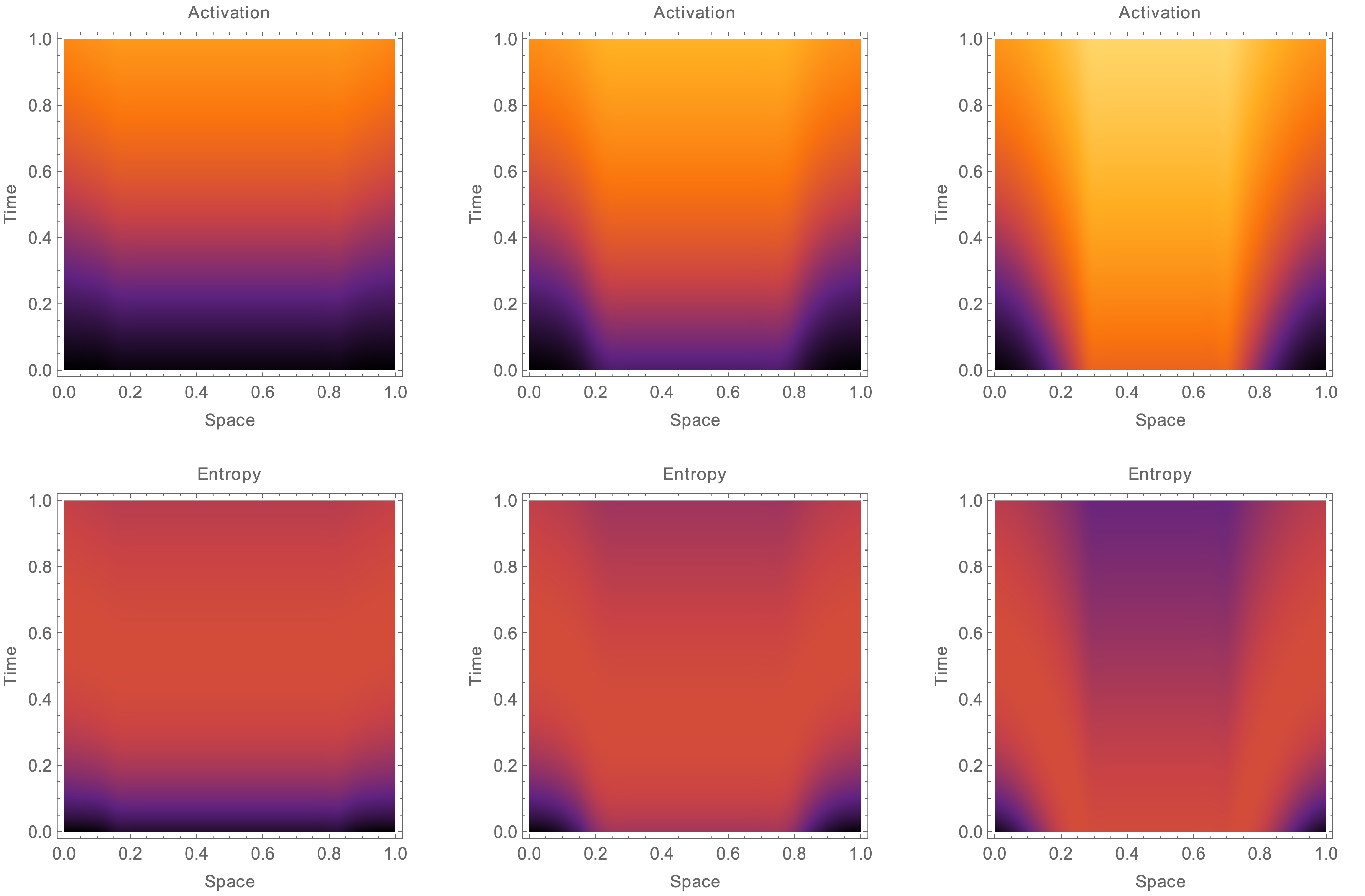}
\caption{Local transduction activation for density kernel masses $0.15,0.2,0.25$ and induction activation threshold $k=15$ for $K\sim\text{Poisson}(c=50)$}\label{fig:nn}
\end{figure}

\begin{figure}[h!]
\centering
\includegraphics[width=7in]{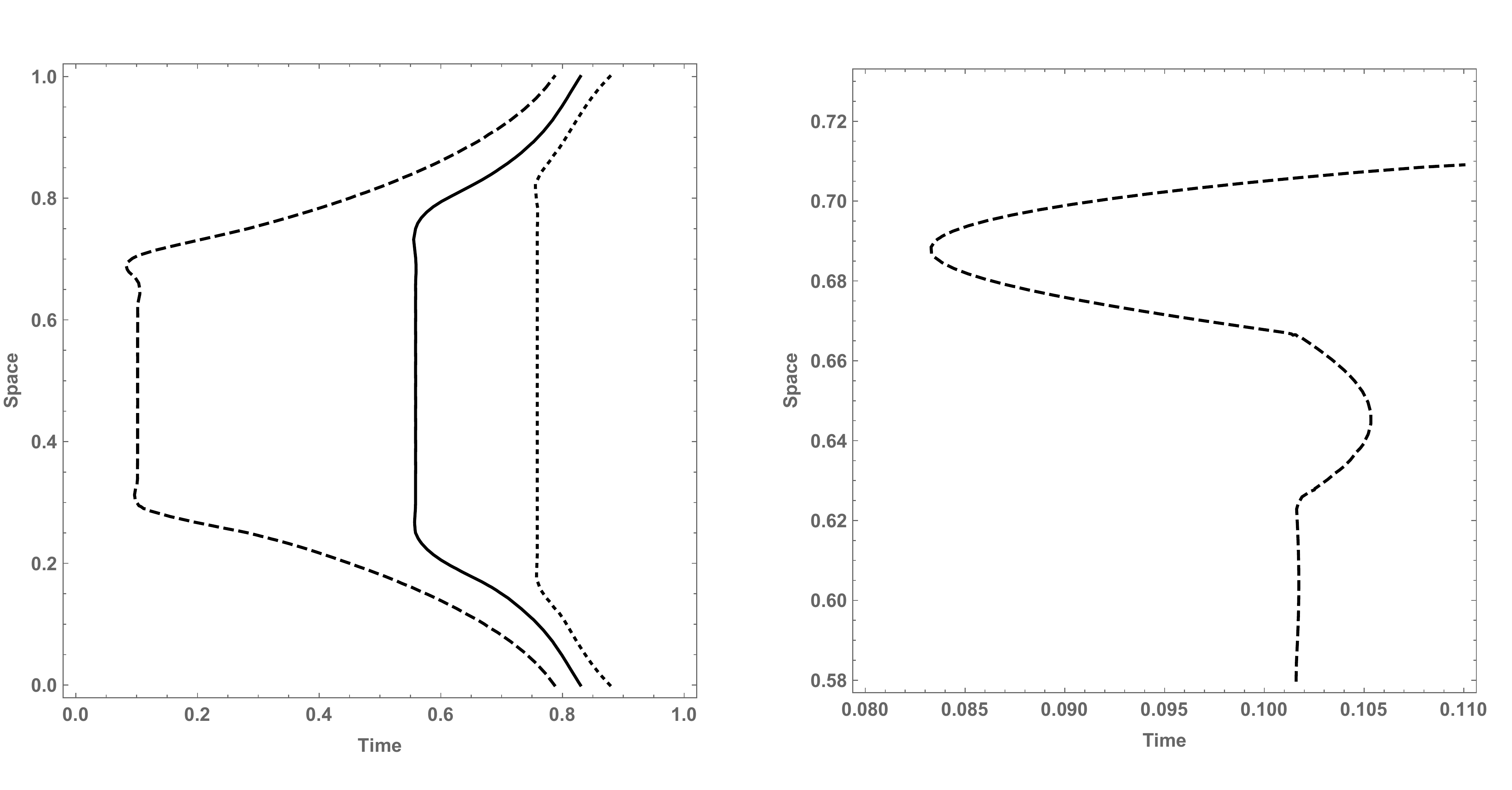}
\caption{Local transduction activation for density kernel masses $0.15,0.2,0.25$ and induction activation threshold $k=15$ for $K\sim\text{Poisson}(c=50)$}\label{fig:nn}
\end{figure}

\begin{figure}[h!]
\centering
\includegraphics[width=6.5in]{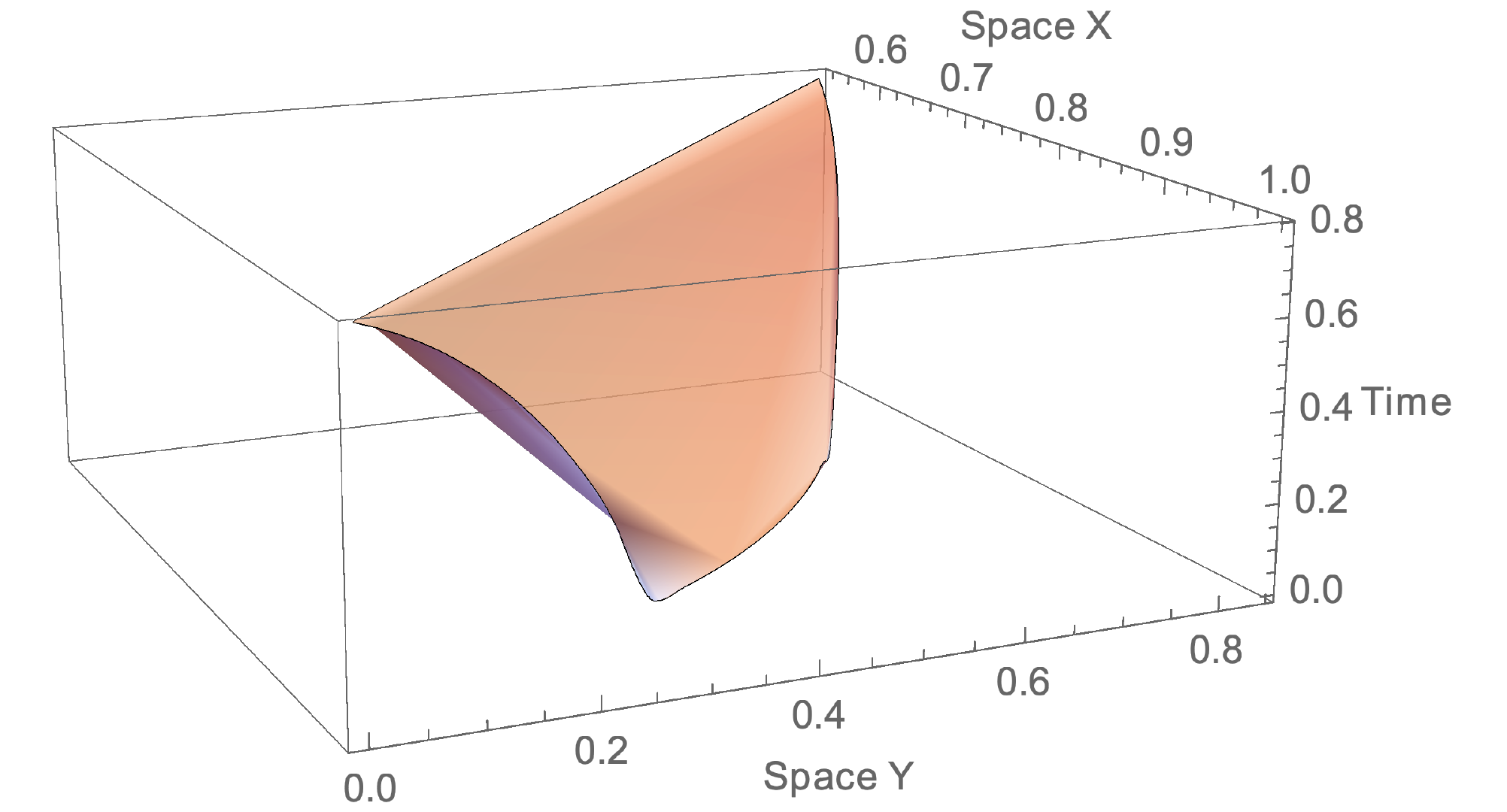}
\caption{Local transduction activation stochastic space-time in coordinates $(\cos(x),\sin(x))$ for density kernel mass $0.25$ and induction activation threshold $k=15$ for $K\sim\text{Poisson}(c=50)$}\label{fig:nn}
\end{figure}

\begin{figure}[h!]
\centering
\includegraphics[width=6.5in]{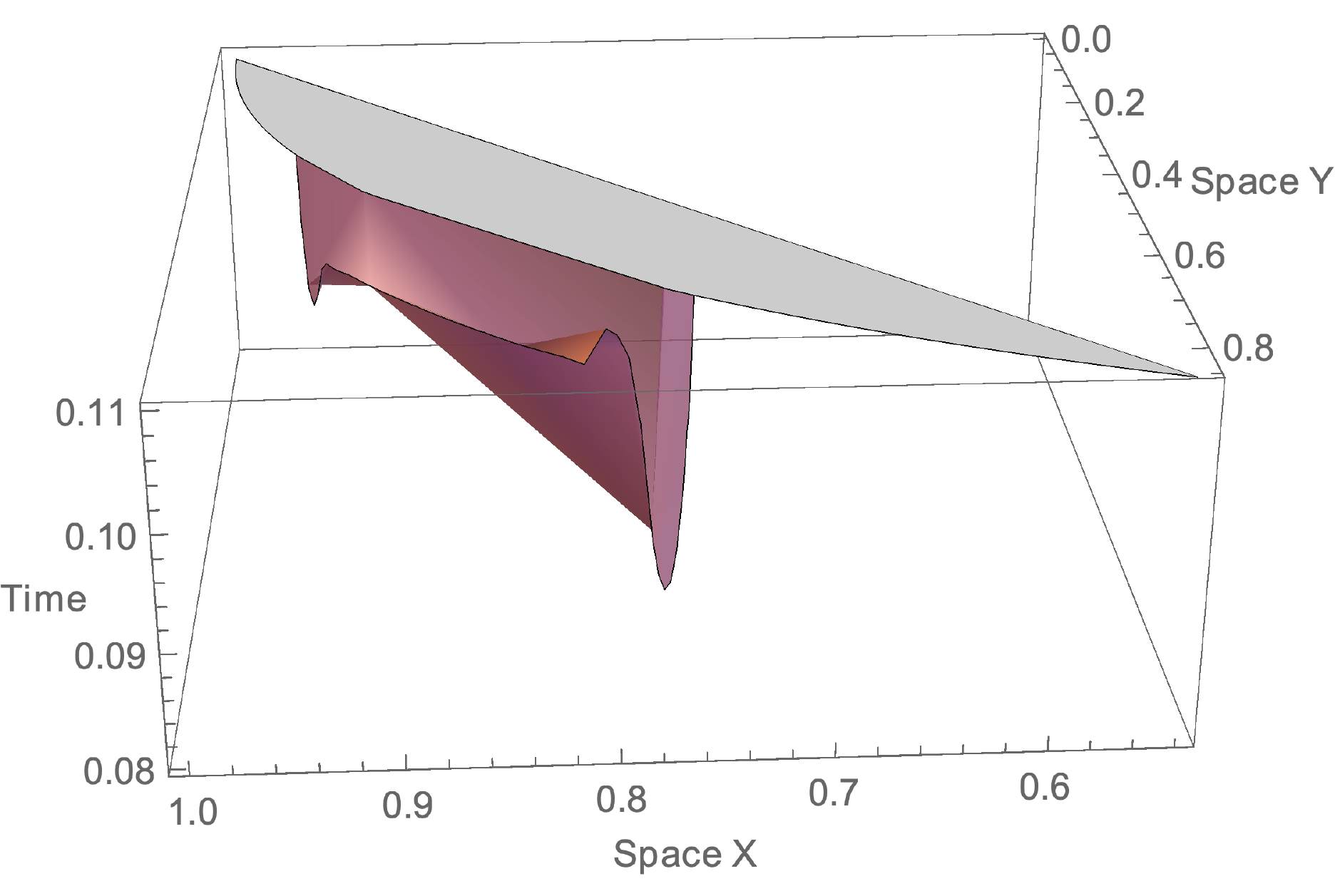}
\caption{Local transduction activation stochastic space-time in coordinates $(\cos(x),\sin(x))$ for density kernel mass $0.25$ and induction activation threshold $k=15$ for $K\sim\text{Poisson}(c=50)$}\label{fig:nn}
\end{figure}

\end{document}